\pdfoutput=1
\documentclass{article} %

\usepackage{iclr2027_conference,times}

\iclrfinalcopy

\usepackage{amsmath,amsfonts,amssymb,amsthm}

\usepackage{prettyref}
\usepackage{mathrsfs}
\usepackage{graphicx}
\usepackage{wrapfig}
\usepackage{subfig}
\usepackage{bbold}
\usepackage{tabu}
\usepackage{comment}
\usepackage{MnSymbol}
\usepackage{multirow}
\usepackage{booktabs}
\usepackage{enumitem}
\usepackage{algorithm}
\usepackage[table]{xcolor}
\usepackage{tikz}
\usepackage[noend]{algpseudocode}
\usepackage[colorlinks,allcolors=gray,hypertexnames=true]{hyperref}

\newtheorem{assume}{\TE{Assumption}}[section]
\newtheorem{define}{\TE{Definition}}[section]
\newtheorem{theorem}{\TE{Theorem}}[section]
\newtheorem{lemma}[theorem]{\TE{Lemma}}
\newtheorem{remark}[theorem]{\TE{Remark}}
\newtheorem{proposition}[theorem]{\TE{Proposition}}
\newtheorem{corollary}[theorem]{\TE{Corollary}}

\newrefformat{fig}{Figure~\ref{#1}}
\newrefformat{par}{Section~\ref{#1}}
\newrefformat{appen}{Appendix~\ref{#1}}
\newrefformat{sec}{Section~\ref{#1}}
\newrefformat{sub}{Section~\ref{#1}}
\newrefformat{table}{Table~\ref{#1}}
\newrefformat{ass}{Assumption~\ref{#1}}
\newrefformat{alg}{Algorithm~\ref{#1}}
\newrefformat{def}{Definition~\ref{#1}}
\newrefformat{thm}{Theorem~\ref{#1}}
\newrefformat{cor}{Corollary~\ref{#1}}
\newrefformat{lem}{Lemma~\ref{#1}}
\newrefformat{step}{Step~\ref{#1}}
\newrefformat{ln}{Line~\ref{#1}}
\newrefformat{rem}{Remark~\ref{#1}}
\newrefformat{prop}{Proposition~\ref{#1}}
\newrefformat{eq}{Equation~\ref{#1}}
\newrefformat{pb}{Problem~\ref{#1}}
\newrefformat{it}{Item~\ref{#1}}
\newrefformat{te}{Term~\ref{#1}}
\def\Eqref Eq:#1:{\eqref{eq:#1}}
\newrefformat{Eq}{Equation~\Eqref#1:}

\newcommand{\defeq}{\mathrel{\mathop:}=}
\newcommand{\R}{\mathbb{R}}
\newcommand{\TE}[1]{\textbf{#1}}
\newcommand{\prox}{\text{prox}}
\newcommand{\argmin}[1]{\underset{#1}{\text{argmin}}\;}
\newcommand{\ST}{\text{s.t.}}
\newcommand{\STQ}{\ST\quad}
\newenvironment{eqaligned}{\begin{equation}\begin{aligned}}{\end{aligned}\end{equation}}

\newcommand{\IAll}{\mathcal{I}}
\newcommand{\IO}{\mathcal{I}_\text{o}}
\newcommand{\IC}{\mathcal{I}_\text{c}}
\newcommand{\Crit}
{\text{Crit}}
\newcommand{\CritWKKT}
{\Crit_\text{w-KKT}}
\newcommand{\CritKKT}
{\Crit_\text{KKT}}
\newcommand{\CritbFJ}
{\Crit_\text{b-FJ}}

\usepackage{xcolor}
\definecolor{Black}{rgb}{0.0, 0.0, 0.0}
\definecolor{Blue} {rgb}{0.2, 0.2, 0.8}
\definecolor{Red}  {rgb}{0.8, 0.2, 0.2}
\definecolor{Green}{rgb}{0.2, 0.8, 0.2}
\newcommand{\rev}[1]{#1}
\newenvironment{revision}{}{}

\title{Parallelizable Interior-Point ADMM for Block-Structured Nonconvex Programs}

\author{Zherong Pan \\
Meta Reality Labs \\
\texttt{zrpan@meta.com}}

\begin{document}

\maketitle
\lhead{Preprint}

\begin{abstract}
Physics simulation is an increasingly important computational primitive for AI, and local--global/ADMM splitting exposes its parallel structure but lacks matching finite-time convergence guarantees. We close this gap with two ADMM variants, Interior-Point ADMM (IP-ADMM), for block-structured nonconvex composite programs with smooth nonlinear inequality and affine constraints. A one-sided finite-activation wall keeps each constraint split state strictly feasible while the original-space infeasibility vanishes, giving closed-form constraint updates, one proximal call per objective block, and one equality-constrained quadratic consensus solve per iteration. \rev{To our knowledge, this makes IP-ADMM the first single-loop method that updates nonconvex inequality constraints locally and in parallel with finite-time convergence guarantees.} Ordinary IP-ADMM attains a CQ-free best-iterate weak-KKT rate of $O(K^{-1/3})$ or a banded Fritz--John rate of $O(K^{-1/4})$, which converts to KKT at $O(K^{-1/4})$ under a uniform constraint qualification; these are nearby-point certificates that reduce to conventional KKT for smooth blocks. A checkpoint variant, IP-ADMM-Ckpt, improves the asymptotic KKT exponent to $O(K^{-1/3})$. We also show how the bounds scale with the objective- and constraint-block counts and the coupling norm.
\end{abstract}

\section{Introduction}
Physics is becoming a computational foundation for AI.  Differentiable simulators place mechanics inside learning and control loops~\cite{belbuteperes2018endtoend,qiao2020scalable,hu2020difftaichi}; physics-informed networks and learned mesh simulators inject governing equations or physical structure into predictive models~\cite{raissi2019physics,pfaff2021learning}; and world models and embodied agents learn or exploit interaction dynamics for planning and control~\cite{hafner2020dream,savva2019habitat}.  \rev{Learned models also define the physics itself: neural signed distance fields (SDFs) serve as contact and collision geometry in differentiable rigid-body simulation and robot control~\cite{park2019deepsdf,strecke2021diffsdfsim,koptev2023njsdf,bukhari2025composite}.}  In all of these settings, we run simulation many times.  Each run must be accurate enough to respect physical structure, and fast enough to supply the experience or gradients that the surrounding AI system needs.

Parallelism is a main way to reach that efficiency.  Physical models naturally break into local constitutive, subsystem, and contact computations that are coupled through shared state, and high-throughput robot learning already benefits from running many physics environments at the same time~\cite{makoviychuk2021isaacgym}. Local--global and ADMM solvers are a promising direction along this line of research. These methods recast the simulation as a nonconvex constrained optimization problem, and expose parallelism by splitting the global model into local ones. This strategy has been very successful in practice for modeling deformable bodies, multibody systems, and contact-rich mechanics~\cite{narain2016admm,overby2017admm,lee2023subsystem,zhu2024frictional,liu2026distributed}. These methods establish the practical value of decomposition. Unfortunately, existing finite-time convergence analyses still cannot explain their empirical success: standard nonconvex ADMM mainly treats objective blocks under affine coupling, whereas methods for nonlinear constraints update all constraints jointly, through a penalty, a linearization, or an inner solve~\cite{jiagrimmer2025,yang2025singleloop,elbourkhissi2025}.

We ask whether a parallelizable splitting method can provide finite-time convergence guarantees for this nonlinear constrained model. To answer this question, we analyze two variants of the ADMM algorithm. Both introduce one split state $y_i=A_ix$ per local model, where each local model is either an objective term or a nonlinear inequality constraint. Our key idea is to apply a one-sided finite-activation quadratic wall to each inequality block, which keeps each split state strictly inside its feasible region. The original variable need not satisfy $f_i(A_ix)\geq0$ at every iteration; instead, the consensus error and the original-space infeasibility vanish as the iterations proceed. Because the split states stay interior, our method is a variant of the interior-point method, which we call IP-ADMM.
\rev{Our paper makes two contributions.
\begin{itemize}[leftmargin=*,nosep]
\item \textbf{A practical explanation of the local constraint-update approach.} Local--global and ADMM physics solvers already split the model into local blocks, update each block independently, and solve one consensus problem~\cite{narain2016admm,overby2017admm,lee2023subsystem,zhu2024frictional,liu2026distributed}; our problem setting and assumptions follow this structure. IP-ADMM keeps one proximal call per objective block, a closed-form update per constraint block, and one quadratic consensus solve per iteration, and provides the finite-time original-problem guarantees these methods lack. To our knowledge, it is the first single-loop method that updates nonconvex inequality constraints locally and in parallel with such guarantees (\prettyref{table:prior}). Because a simulator defines a family of problems as the discretization is refined (\prettyref{fig:block-scaling-cloth}), our bounds keep the objective-block count $N$, the constraint-block count $M$, and the coupling norm $\|A\|_2$ explicit.
\item \textbf{A theoretical extension of nonconvex ADMM to nonlinear inequality constraints.} Existing nonconvex ADMM analyses allow only affine coupling~\cite{yuan2025admm,maia2026adaptive}. Following the continuation idea of IPDS-ADMM~\cite{yuan2025admm}, ordinary IP-ADMM attains CQ-free best-iterate w-KKT and b-FJ rates of $O(K^{-1/3})$ and $O(K^{-1/4})$, the latter with the same exponent as the double-loop method of~\cite{jiagrimmer2025}. Under Uniform Band-MFCQ, b-FJ converts to KKT at $O(K^{-1/4})$, and the checkpoint variant IP-ADMM-Ckpt improves this to $O(K^{-1/3})$ after its first checkpoint at iteration $K_0$, which can be large. \prettyref{table:rates} summarizes this guarantee ladder, and \prettyref{table:prior} compares it with prior methods (\prettyref{appen:prior}).
\end{itemize}}

\section{Related Work}
\begin{wraptable}{r}{0.55\textwidth}
\vspace{-40px}
\centering
\scriptsize
\setlength{\tabcolsep}{2pt}
\renewcommand{\arraystretch}{1.03}
\caption{Physics regimes for our certificate hierarchy.}
\label{table:certificate-physics}
\begin{tabular}{@{}p{0.13\textwidth}p{0.37\textwidth}@{}}
\toprule
Certificate & Physics regime and precedent \\
\midrule
w-KKT & Force balance; forces at any distance \\
b-FJ, fixed $\epsilon$ & Zero force beyond a fixed band~\cite{li2020ipc} \\
b-FJ, $\epsilon\to0$ & Shrinking penalty layers~\cite{harmon2009acm} \\
KKT & MFCQ-qualified contact~\cite{anitescuhart2004} \\
\bottomrule
\end{tabular}
\end{wraptable}
\paragraph{Physics+AI and parallel physics simulation.}
Physics enters AI through several complementary ways.  Differentiable engines expose simulation to end-to-end learning and control~\cite{belbuteperes2018endtoend,qiao2020scalable,hu2020difftaichi}.  Physics-informed networks encode governing-equation residuals, and graph-based learned simulators use physical discretizations with local message passing~\cite{raissi2019physics,pfaff2021learning}.  Latent world models and embodied platforms support prediction and interaction at scale~\cite{hafner2020dream,savva2019habitat,makoviychuk2021isaacgym}. \rev{Learned geometry can also define the constraints of a simulation: neural SDFs~\cite{park2019deepsdf} supply the contact geometry of nonconvex bodies in differentiable rigid-body simulation~\cite{strecke2021diffsdfsim}, neural density fields supply contact in differentiable simulation of learned objects~\cite{lecleach2023dano}, and neural distance fields supply collision avoidance in robot control and navigation~\cite{koptev2023njsdf,bukhari2025composite}. \prettyref{appen:neural} uses such a field as the contact constraint of \prettyref{pb:prob}.} These developments motivate scalable physics simulators that make good use of parallel processors. Local--global and ADMM formulations are a promising direction along this line. These algorithms separate constitutive, subsystem, or contact handling from global consensus updates in projective dynamics, multibody systems, soft-body contact, and distributed affine-body dynamics~\cite{narain2016admm,overby2017admm,lee2023subsystem,zhu2024frictional,liu2026distributed}. This literature motivates our architecture. Our main goal is to establish finite-time original-problem certification for the nonlinear program in \prettyref{pb:prob} under its local-consensus split.

\paragraph{Contact and admissible trajectories.}
Collision avoidance is an important source of nonlinear and often nonconvex constraints in mechanics: feasible motion depends on pose-dependent distances, changing closest features, and self-contact geometry.  Robust cloth contact and continuous collision detection seek noninterpenetrating motion over a time interval~\cite{bridson2002cloth,harmon2009acm,wang2021ccd}.  IPC uses a barrier-augmented incremental potential and a CCD-aware line search from a collision-free state; C-IPC and rigid IPC extend this admissible-trajectory philosophy to codimensional and reduced-coordinate rigid contact~\cite{li2020ipc,li2021cipc,ferguson2021rigidipc}. BC-ADMM provides a parallelizable algorithm for a class of biconvex functions, but its rate certifies a fixed-barrier relaxation, its local subproblems are solved iteratively, and it does not establish finite-time guarantees for certificates of the original problem~\cite{pan2026bcadmm}.

\paragraph{Nonconvex splitting and functional constraints.}
Nonconvex ADMM commonly allows nonconvex or nonsmooth objective blocks while keeping affine coupling. IPDS-ADMM provides a closely aligned continuation precedent: it combines a polynomially increasing penalty with inverse-penalty smoothing and a square-root-penalty-rescaled multiplier~\cite{yuan2025admm}. Adaptive proximal ADMM instead uses progress tests to gate multiplier updates~\cite{maia2026adaptive}. For nonlinear functional constraints,~\cite{jiagrimmer2025} obtain normalized Fritz--John certificates without a CQ, and KKT certificates under a uniform strong-MFCQ condition for the proximal subproblems, while Hinder and Ye analyze a nonconvex log-barrier method from a strictly feasible start~\cite{hinder2024worst}. ABIP and ACVI combine ADMM and interior-point ideas for linear programs and constrained variational inequalities, respectively~\cite{lin2021abip,yang2023acvi}. \rev{However, these algorithms update their nonlinear constraints jointly rather than locally. Among the methods compared in \prettyref{table:prior} (\prettyref{appen:prior}), IP-ADMM is the only one that combines a single loop, local updates of nonconvex constraints, and a finite-time rate.}

\paragraph{Certificates and qualification.}
\prettyref{table:certificate-physics} maps our increasingly informative certificate hierarchy to contact regimes.  \rev{Weak-KKT (w-KKT) gives consensus and force balance without complementarity. Normalized banded Fritz--John (b-FJ) at accuracy $\epsilon$ additionally supplies complementarity, nontrivial normalized multipliers, and support of every nonzero force on constraints with $f_i\leq\epsilon$, without a CQ. A fixed band is IPC's accuracy notion: its repulsion vanishes beyond the target distance $\hat d$, so a frictionless IPC solution is a b-FJ point at accuracy $\hat d$ (\prettyref{prop:finite-band}), although a fixed band does not converge to zero-band stationarity~\cite{li2020ipc}. IP-ADMM instead shrinks the band, paralleling ACM's nested penalty layers, whose thickness decreases while stiffness increases~\cite{harmon2009acm}.} Under Uniform Band-MFCQ, a sufficiently accurate b-FJ certificate converts to KKT, which removes the Fritz--John scaling ambiguity. Anitescu and Hart provide a narrow mechanics precedent: a pointed-friction-cone condition implies MFCQ for their rigid-contact quadratic program~\cite{anitescuhart2004}. Thus the cited physics papers instantiate these regimes, while our analysis defines and connects the exact certificates; the CQ-free FJ versus CQ-dependent KKT boundary follows the broader optimization literature~\cite{jiagrimmer2025,hinder2024worst}.

\section{Problem Setting and Simulation Instances}
We study a block-separable composite program. It has smooth (possibly nonconvex) inequality-constraint functions $f_i$, smooth (possibly nonconvex) objective terms $h_i$, convex nonsmooth proximable terms $\psi_i$, affine equality constraints through $C$, and affine coupling through the maps $A_i$.  Here $N,M\geq1$, $\IAll=\{1,\ldots,N+M\}$, $\IO=\{1,\ldots,N\}$, $\IC=\IAll\setminus\IO$, $x\in\R^{n_x}$, $A_i\in\R^{n_i\times n_x}$ for $i\in\IAll$, $A$ denotes the vertical stack of $A_1,\ldots,A_{N+M}$, and $\varphi_i\defeq h_i+\psi_i$ for $i\in\IO$.  The problem is
\begin{eqaligned}
\label{pb:prob}
\argmin{x}\ \sum_{i\in\IO}\varphi_i(A_ix)\qquad\STQ f_i(A_ix)\geq0\ (i\in\IC),\quad Cx=d.
\end{eqaligned}

\begin{assume}[Constraint regularity]
\label{ass:f}
For every $i\in\IC$, $f_i$ has an $L_i$-Lipschitz gradient.
\end{assume}
\begin{assume}[Bounded constraint gradients]
\label{ass:bg}
For every $i\in\IC$, $\|\nabla f_i(y)\|\leq M_i$ for all $y$.
\end{assume}
\begin{assume}[Feasibility]
\label{ass:feasible}
The initial point $x^1$ is strictly feasible: $Cx^1=d$ and $m\defeq\min_{i\in\IC}f_i(A_ix^1)>0$.  A point is feasible when $Cx=d$ and $f_i(A_ix)\geq0$ for every $i\in\IC$.
\end{assume}

\begin{assume}[Composite objective blocks]
\label{ass:obj}
For each $i\in\IO$, $h_i$ is continuously differentiable, possibly nonconvex, globally $D_i$-Lipschitz, and has an $H_i$-Lipschitz gradient.  The function $\psi_i$ is proper, closed, convex, finite-valued on $\R^{n_i}$, globally $C_i^\psi$-Lipschitz, and $\prox_{t\psi_i}$ is exactly evaluable for every $t>0$.  We also assume the lower bound $\sum_{i\in\IO}\varphi_i(A_ix)\geq0$ for every $x$.
\end{assume}

These are the only four standing assumptions. For later use, define $\widehat C_i\defeq D_i+C_i^\psi$ for $i\in\IO$. \prettyref{appen:algebraic} records the global estimates implied by these assumptions. \prettyref{pb:prob} overlaps two established problem classes but matches neither exactly. Affine-coupled nonconvex ADMM allows rich composite objective blocks, but it does not directly represent the nonlinear inequalities we use here~\cite{yuan2025admm,maia2026adaptive}. Functional-constraint methods handle nonlinear inequalities and give FJ/KKT guarantees under their own assumptions~\cite{jiagrimmer2025,hinder2024worst}, but they do not exploit this local-consensus decomposition. \rev{IP-ADMM keeps the local-consensus decomposition of the former while handling the nonlinear inequalities of the latter, and certifies Fritz--John points without a CQ (\prettyref{table:rates}).}

\subsection*{Motivating Applications}
\begin{figure*}[t]
\vspace{-10px}
\centering
\setlength{\tabcolsep}{0pt}
\begin{tabular}{cccc}
\includegraphics[width=0.242\linewidth,trim=34 53 34 53,clip]{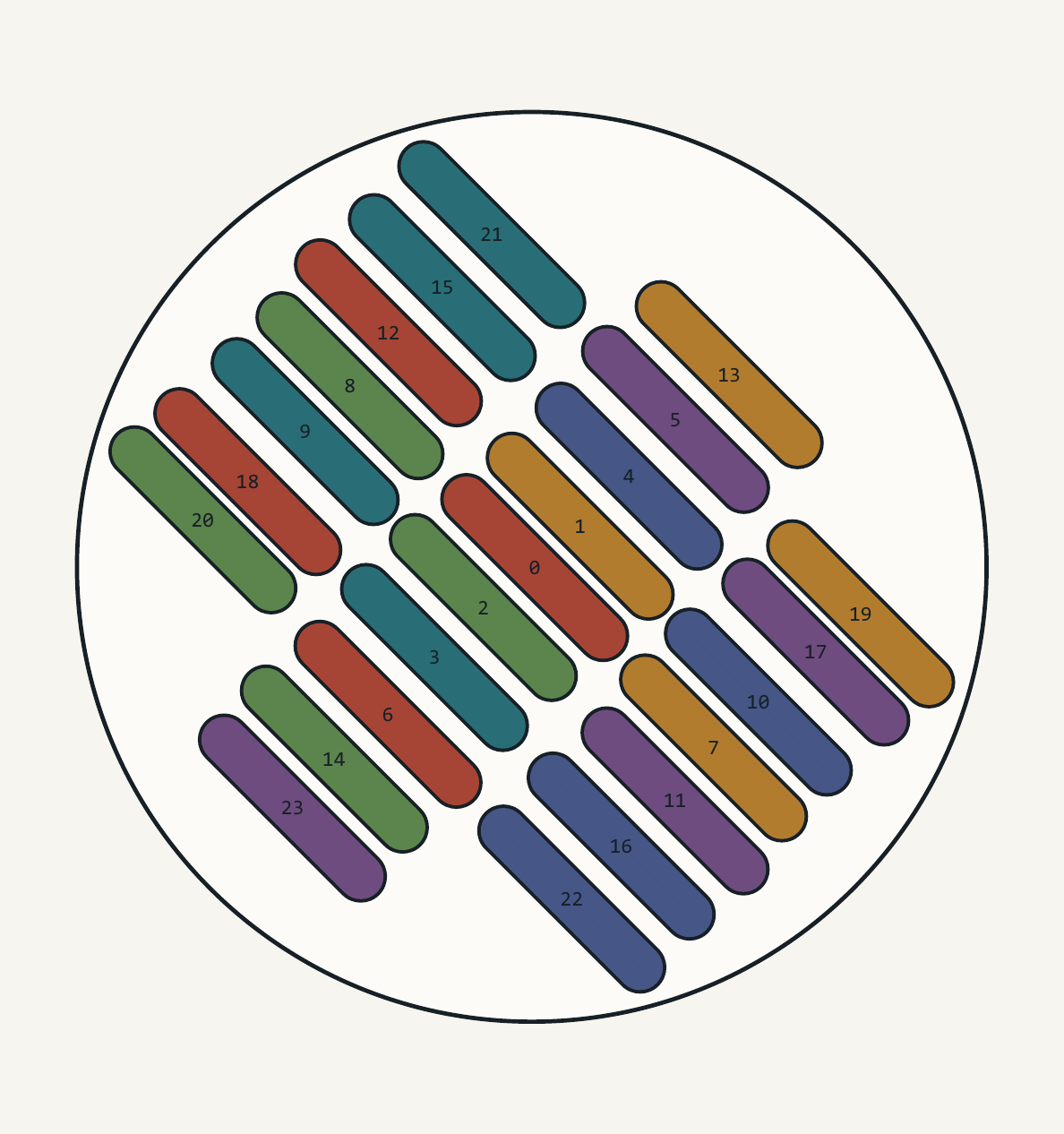}&
\includegraphics[width=0.242\linewidth,trim=34 53 34 53,clip]{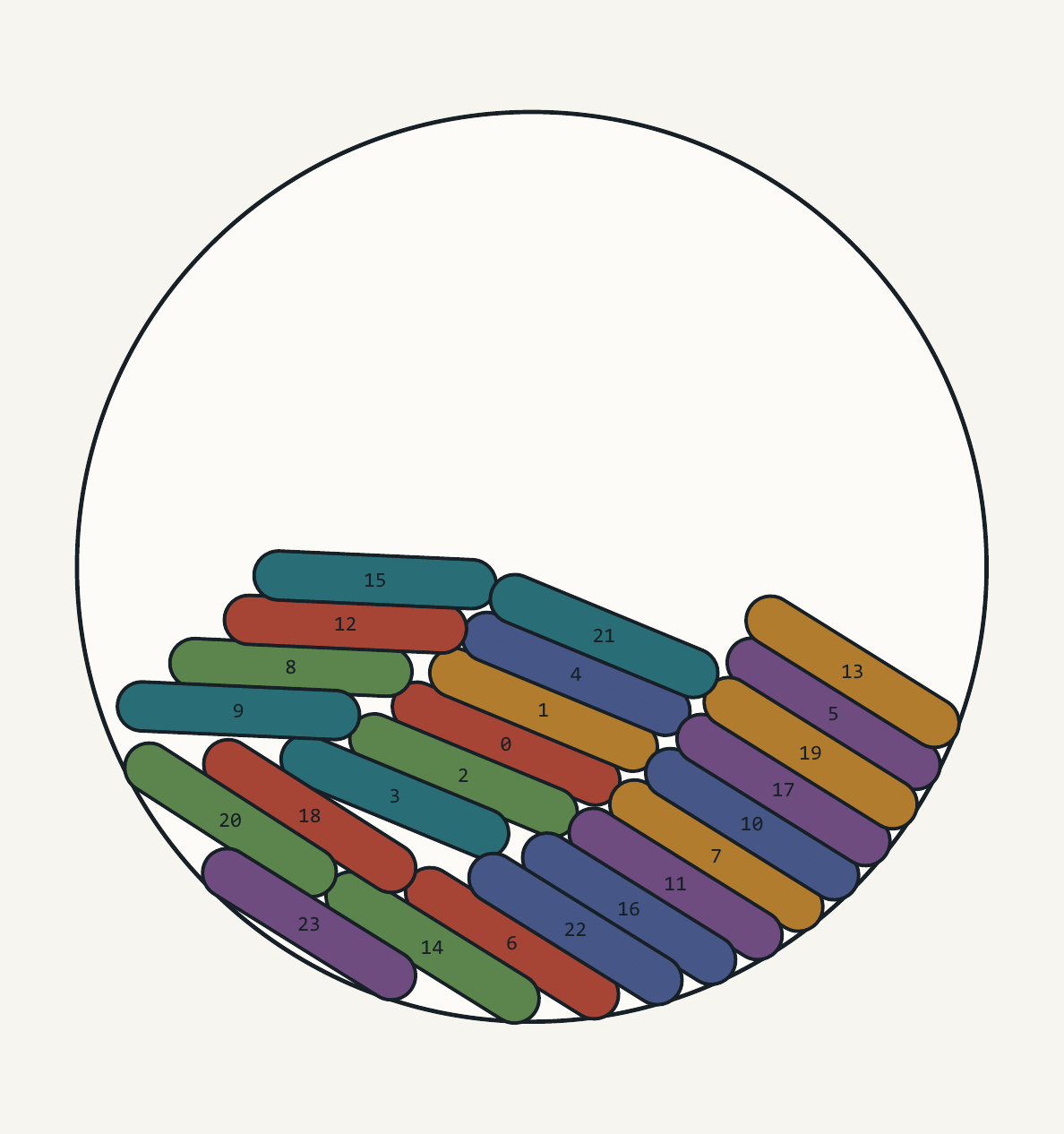}&
\includegraphics[width=0.242\linewidth,trim=34 53 34 53,clip]{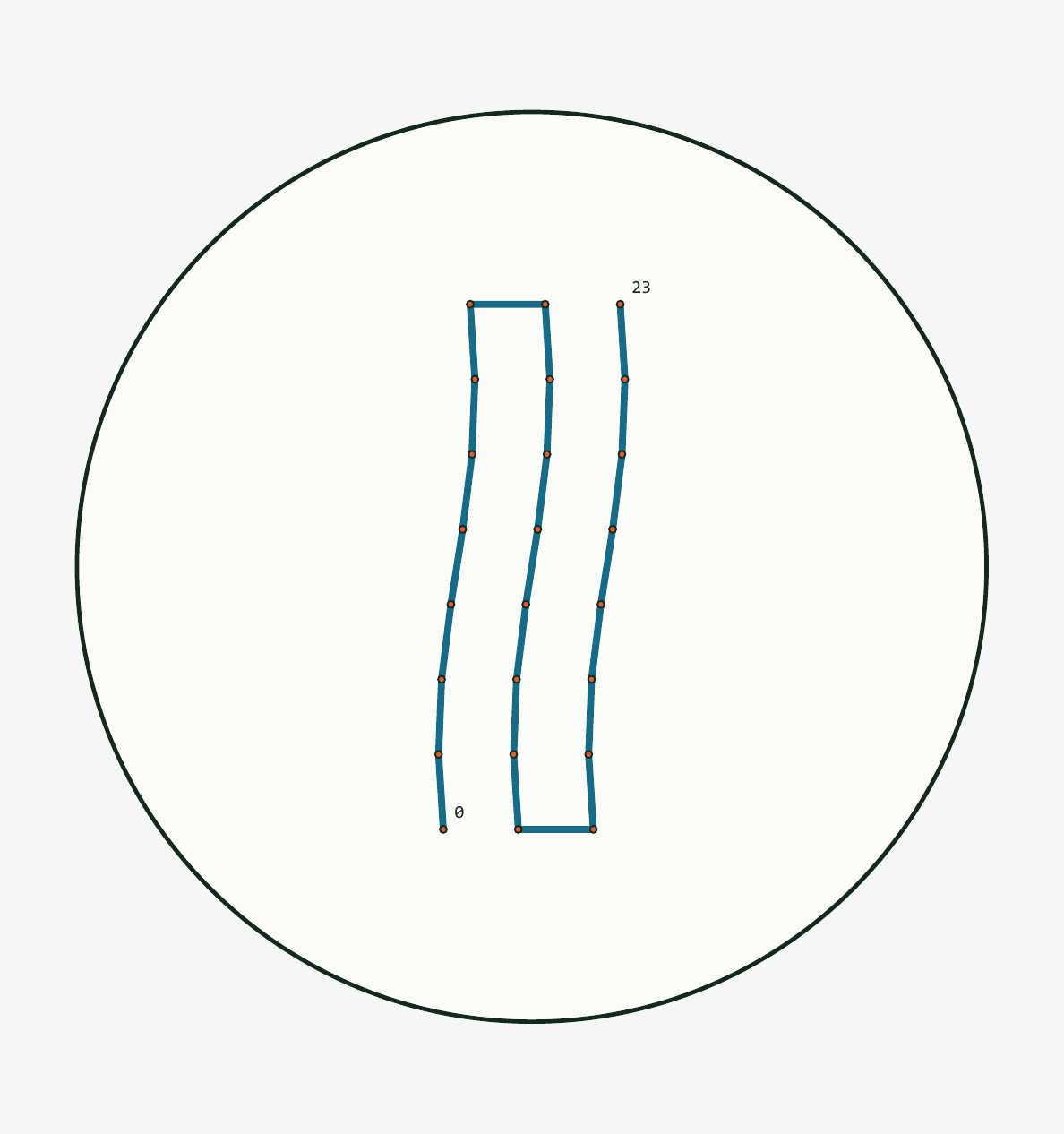}&
\includegraphics[width=0.242\linewidth,trim=34 53 34 53,clip]{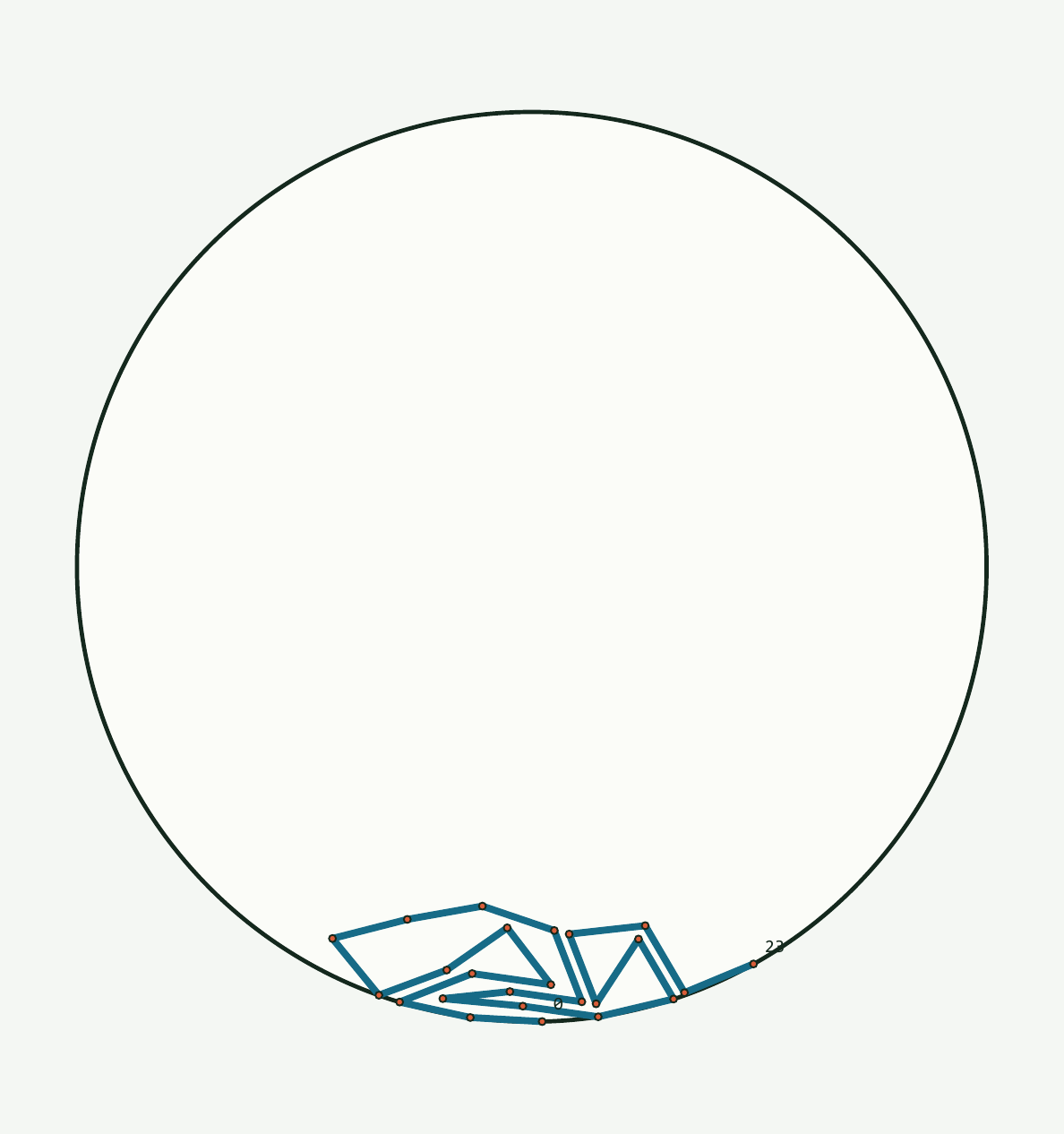}
\end{tabular}
\caption{Initial and final ordinary IP-ADMM configurations for capsule stacking (left pair) and cloth resting (right pair) at $n=24$.}
\label{fig:applications}
\vspace{-5px}
\end{figure*}

It is well known that physics simulation problems can be recast as nonconvex constrained optimization~\cite{narain2016admm,overby2017admm}, and we show that our problem setting applies directly to such problems and solves them in parallel. As~\prettyref{fig:applications} illustrates, we compute the resting configurations of rigid and deformable objects, where (self-)collision constraints are the main source of nonconvexity. Detailed formulations and checks against the four standing assumptions appear in \prettyref{appen:applications}, where we apply a series of constraint mollifications so the problem conforms to our assumptions. These mollifications are standard techniques; see, e.g.,~\cite{pan2026bcadmm,li2020ipc}.
\TE{Capsule stacking:} We minimize gravitational potential over $n$ planar capsule poses inside a circular container. The $3n$ pose variables form $N=n$ objective blocks; directed capsule--capsule collision-clearance constraints and capsule--container boundary collision constraints form $M=2n^2$ inequality blocks. All blocks use local selectors $A_i$, there are no affine equalities, and initialization is strictly feasible.
\TE{Cloth resting:} We minimize gravity plus stretch and two-hop bending energies for $n$ planar cloth nodes, modeling an elastic 2D piece of cloth as described in~\cite{10.1145/566654.566624}. The $2n$ coordinates give $N=3n-3$ objective blocks; cloth self-collision constraints over nonincident node--segment pairs and cloth--container boundary collision constraints give $M=(n-1)(n-2)+n$ inequality blocks.

\section{IP-ADMM(-Ckpt) and Convergence Analysis}
We develop two variants of the ADMM algorithm for \prettyref{pb:prob}.  Both smooth the convex nonsmooth objective components with Moreau envelopes and represent each nonlinear inequality by a one-sided finite-activation wall, yielding independent local updates followed by one equality-constrained quadratic consensus update.  Ordinary IP-ADMM uses polynomial continuation, similar to IPDS-ADMM~\cite{yuan2025admm}, and provides CQ-free best-iterate w-KKT and b-FJ guarantees, whereas IP-ADMM-Ckpt promotes its continuation state only at successful checkpoints and, under Uniform Band-MFCQ, returns a latest-checkpoint KKT certificate with an improved asymptotic rate once its first checkpoint completes. We first define their common model and algorithm and then state the mode-specific convergence results; complete constants and proofs are deferred to the appendices. \rev{\prettyref{table:rates} in \prettyref{appen:prior} summarizes these guarantees, and \prettyref{table:prior} compares them with prior methods.} IP-ADMM rows select the best certified iterate in the stated window, whereas IP-ADMM-Ckpt returns the latest checkpoint.

\paragraph{Algorithm Design.} We begin by defining our proximal operators. For the convex component in objective, we use the standard Moreau envelopes below
\begin{eqaligned}
(\psi_i)_\mu(y)\defeq\min_z\left\{\psi_i(z)+\frac{1}{2\mu}\|z-y\|^2\right\},\qquad P_{i,\mu}^\psi(y)\defeq\prox_{\mu\psi_i}(y).
\end{eqaligned}
For each mode $\mathsf m\in\{\mathsf{IP},\mathsf{Ckpt}\}$, bind its proximal and wall parameters by
\begin{eqaligned}
(\theta_{1,\mathsf m},\theta_{2,\mathsf m},G_{\mathsf m})
\defeq
\begin{cases}
(\theta_1,\theta_2,G),&\mathsf m=\mathsf{IP},\\
(\theta_1^{\rm cp},\theta_2^{\rm cp},G^{\rm cp}),&\mathsf m=\mathsf{Ckpt}.
\end{cases}
\end{eqaligned}
The mode-specific values are fixed in~\prettyref{appen:parameters}.  For a constraint and $L>0$, define the active one-sided wall by
\begin{eqaligned}
g_{\mathsf m}(s,L)=
\begin{cases}
G_{\mathsf m}-\sqrt{2LG_{\mathsf m}}\,s,&s<0,\\
\frac{L}{2}\max(\sqrt{2G_{\mathsf m}/L}-s,0)^2,&s\geq0.
\end{cases}
\end{eqaligned}
Write $\beta(L)\defeq\sqrt{2G_{\mathsf m}/L}$ for the active wall width.  In the rest of the algorithm and in each mode-specific proof, $\mathsf m$ is fixed by the selected mode. We therefore unify objective and constraint blocks as
\begin{eqaligned}
\bar f_i(y,L,\mu)\defeq
\begin{cases}
h_i(y)+(\psi_i)_\mu(y),&i\in\IO,\\
g_{\mathsf m}(f_i(y),L),&i\in\IC.
\end{cases}
\end{eqaligned}
The calculus of Moreau envelopes is well understood, and the activation wall follows from standard convex analysis.  We summarize the useful properties in~\prettyref{appen:algebraic}.  Introducing split variables $y_i=A_ix$, stacked as $y\defeq(y_i)_{i\in\IAll}$ in the block order of $A$, with consensus multipliers $\lambda\defeq(\lambda_i)_{i\in\IAll}$ stacked the same way, gives the wall-penalized split formulation used by the algorithm
\begin{eqaligned}
\argmin{x,y}\ \sum_{i\in\IO}\varphi_i(y_i)+\sum_{i\in\IC}g_{\mathsf m}(f_i(y_i),L)\qquad\STQ y=Ax,\quad Cx=d.
\end{eqaligned}
For the algorithm, we smooth the convex objective components and use the smoothed augmented Lagrangian
\begin{eqaligned}
\label{eq:Lagrangian}
\mathcal L(x,y,\lambda,L,\mu,\rho)
=\sum_{i\in\IAll}\bar f_i(y_i,L,\mu)
+\lambda^T(y-Ax)+\frac\rho2\|y-Ax\|^2.
\end{eqaligned}

At iteration $k$, let $d_i^k\defeq y_i^{k+1}-y_i^k$ denote the block step and $e_i^k\defeq\lambda_i^k+\rho^k(y_i^k-A_ix^k)$ the current consensus force; the stacked step is $d^k\defeq(d_i^k)_{i\in\IAll}$, obtained by concatenating the blocks in the order $i=1,\ldots,N+M$, so that $\|d^k\|^2=\sum_{i\in\IAll}\|d_i^k\|^2$. For a constraint block, $\ell_i^k(u)\defeq f_i(y_i^k)+\langle\nabla f_i(y_i^k),u-y_i^k\rangle$ is its affine linearization.  On an objective block, $w_i^{k+1}\defeq\nabla(\psi_i)_{\mu^k}(y_i^{k+1})$ is the smoothed convex force; on a constraint block, $\iota_i^{k+1}\defeq g_{\mathsf m}'(\ell_i^k(y_i^{k+1}),L^k)$ is the scalar wall force. The model-error curvature $\tau_i^k$, model smoothness $\zeta_i^k$, and resulting block force $q_i^k$ are consolidated as
\begin{eqaligned}
(\tau_i^k,\zeta_i^k,q_i^k)\defeq
\begin{cases}
(H_i,1/\mu^k,\nabla h_i(y_i^k)+w_i^{k+1}),&i\in\IO,\\
(\sqrt{2L^kG_{\mathsf m}}\,L_i,L^kM_i^2,\iota_i^{k+1}\nabla f_i(y_i^k)),&i\in\IC.
\end{cases}
\end{eqaligned}
Set $\mathsf B_i^k\defeq\rho^k+\tau_i^k$, the curvature that majorizes the linearized block Lagrangian, and $m_i^k\defeq\theta_{1,\mathsf m}\mathsf B_i^k$, the strongly convex proximal curvature.  The block models are
\begin{eqaligned}
\label{eq:block-models}
\mathcal M_i^k(u)\defeq
\begin{cases}
(\psi_i)_{\mu^k}(u)+\langle\nabla h_i(y_i^k)+e_i^k,u-y_i^k\rangle+\dfrac{m_i^k}{2}\|u-y_i^k\|^2,&i\in\IO,\\[1mm]
g_{\mathsf m}(\ell_i^k(u),L^k)+\langle e_i^k,u-y_i^k\rangle+\dfrac{m_i^k}{2}\|u-y_i^k\|^2,&i\in\IC.
\end{cases}
\end{eqaligned}
Both are strongly convex. Linearizing the smooth block terms makes the first model evaluable with one convex proximal call and gives a closed form for the second model.

For an objective block, set $t_i^k\defeq1/m_i^k$ and $z_{o,i}^k\defeq y_i^k-t_i^k(\nabla h_i(y_i^k)+e_i^k)$.  The one-call interface in~\prettyref{lem:one-call-envelope-prox} returns $y_i^{k+1}$ together with the proximal point $p_i^{k+1}$ and force $w_i^{k+1}$; its exact identities are recorded in \prettyref{eq:prox-oracle}.  The constraint model has the three-case closed form of~\prettyref{lem:added_algebraic_3}; hence neither branch requires an iterative inner solve. In Ckpt mode a failed test freezes the stored dual and all continuation parameters but retains the new primal variables; the virtual multiplier is discarded.  A successful step records its checkpoint certificate tuple with pre-promotion parameters, then commits the multiplier and promotes the continuation state.  Thus the method has one global primal trajectory, and no certified latest-checkpoint output exists before the first success.

\begin{figure*}[t]
\vspace{-10px}
\centering
\begin{tabular}{cc}
\includegraphics[width=0.4\textwidth]{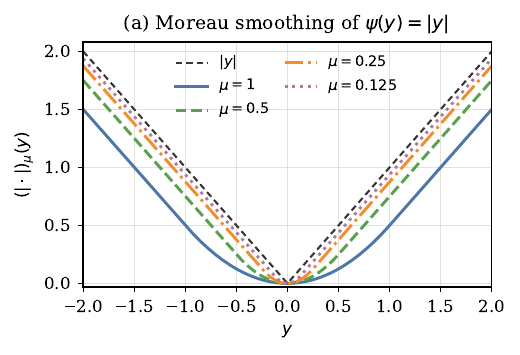}&
\includegraphics[width=0.4\textwidth]{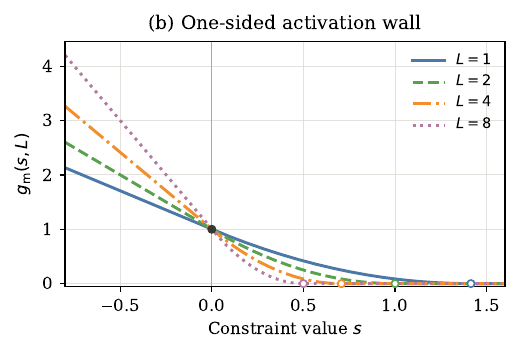}
\end{tabular}
\vspace{-10px}
\caption{Geometry induced by the continuation schedule in~\prettyref{eq:IP-schedules}, shown under the illustrative normalization $\rho^0=L^0=\eta=G=1$.  Panel~(a) uses $\mu\in\{1,1/2,1/4,1/8\}$, and panel~(b) uses $L\in\{1,2,4,8\}$; these values represent scheduled models rather than iterations for a particular exponent, while $p$ determines how quickly the iterates traverse them.  In panel~(b), open markers indicate the zero endpoints $\beta=\sqrt{2/L}$.}
\vspace{-10px}
\label{fig:IP-schedules}
\end{figure*}

Our two modes below use the same primal update with positive continuation parameters $(\rho^k,L^k,\mu^k)$ with the only difference being the scheduling of these parameters. In the ordinary IP-ADMM mode, choose $\sigma\in(0,2)$, $\rho^0,L^0,\eta>0$, $\xi\in(0,1]$, and $p\in[0,1)$, and set
\begin{eqaligned}
\label{eq:IP-schedules}
\rho^k=\rho^0(1+\xi k^p),\qquad L^k=L^0(1+\xi k^p),\qquad \mu^k=\frac{\eta}{\rho^k}.
\end{eqaligned}
The common factor in $(\rho^k,L^k)$ keeps $\rho^k/L^k=\rho^0/L^0$ fixed, whereas $\mu^k=\eta/[\rho^0(1+\xi k^p)]$.  Thus, for the model convex term $\psi(y)=|y|$ in the left panel of~\prettyref{fig:IP-schedules}, the quadratic transition region of its Moreau envelope contracts like $1/(1+\xi k^p)$ and the envelope approaches $|y|$.  Meanwhile the wall in the right panel has feasible-side activation width
\begin{eqaligned}
\beta^k\defeq\beta(L^k)=\sqrt{2G_{\mathsf m}/L^k}.
\end{eqaligned}
In the ordinary mode $G_{\mathsf m}=G$, so $\beta^k=\sqrt{2G/(L^0(1+\xi k^p))}$.  Its boundary height remains $g_{\mathsf m}(0,L^k)=G$, but for each fixed violation $s<0$, $g_{\mathsf m}(s,L^k)=G-\sqrt{2L^0G(1+\xi k^p)}\,s$ grows with the schedule.  Hence $p>0$ moves through increasingly sharp objective and constraint models, while $p=0$ uses one fixed member of each family.

The increasing-penalty/decreasing-smoothing continuation is inspired by IPDS-ADMM~\cite{yuan2025admm}. We depart from their method by additionally increasing the nonlinear-constraint wall curvature, so its activation band shrinks while consensus and reconstruction errors vanish.  The IP-ADMM-Ckpt mode instead freezes $(\rho^k,L^k,\mu^k)$ between successful checkpoints and uses the geometric update shown in \prettyref{alg:IP-ADMM}. Its gated dual update is loosely motivated by the progress-tested multiplier updates of~\cite{maia2026adaptive}; the retained-primal and full continuation-state semantics below are specific to our nonlinear-inequality checkpoint analysis.

\begin{algorithm}
\caption{\label{alg:IP-ADMM}IP-ADMM(-Ckpt)}
\begin{algorithmic}
\Require{$x^1$ satisfying~\prettyref{ass:feasible}; $\mathsf m\in\{\mathsf{IP},\mathsf{Ckpt}\}$, $\chi>1$, $\ell_{\rm cp}>0$, and $\tau_{\rm cp}>0$}
\State $y_i^1\gets A_ix^1$, $\lambda^1\gets0$; initialize $(\rho^1,L^1,\mu^1)$ for the selected mode
\For{$k=1,2,\ldots$}
\For{$i=1,\ldots,N+M$}
\State $y_i^{k+1}\gets\argmin{u}\mathcal M_i^k(u)$
\EndFor
\State $\displaystyle x^{k+1}\gets\argmin{Cx=d}\left\{
\mathcal L(x,y^{k+1},\lambda^k,L^k,\mu^k,\rho^k)
+(\theta_{2,\mathsf m}\rho^k/2)\|A(x-x^k)\|^2\right\}$
\State $r^{k+1}\gets y^{k+1}-Ax^{k+1}$
\If{$\mathsf m=\mathsf{Ckpt}$}
\State $\widehat\lambda^{k+1}\gets\lambda^k+\rho^kr^{k+1}$
\State $E^{k+1}\gets\rho^k\sqrt{\|y^{k+1}-y^k\|^2+\|A(x^{k+1}-x^k)\|^2}$
\If{$E^{k+1}>\tau_{\rm cp}/\rho^k$}
\State $(\lambda^{k+1},\rho^{k+1},L^{k+1},\mu^{k+1})\gets(\lambda^k,\rho^k,L^k,\mu^k)$
\Else
\State record the step-$k$ tuple with $(\rho^k,L^k,\mu^k)$ as the latest-checkpoint output
\State $(\lambda^{k+1},\rho^{k+1},L^{k+1},\mu^{k+1})\gets(\widehat\lambda^{k+1},\chi\rho^k,\chi^2L^k,\mu^k/\chi)$
\EndIf
\Else
\State $\lambda^{k+1}\gets\lambda^k+\sigma\rho^kr^{k+1}$ and set $(\rho^{k+1},L^{k+1},\mu^{k+1})$ by~\prettyref{eq:IP-schedules}
\EndIf
\EndFor
\end{algorithmic}
\end{algorithm}

\begin{remark}
The $x$-update in~\prettyref{alg:IP-ADMM} minimizes an equality-constrained quadratic that depends on $x$ only through $Ax$ subject to $Cx=d$, so the minimizer need not be unique when $\ker A\cap\ker C\neq\{0\}$; in that case any minimizing $x^{k+1}$ may be selected, since $Ax^{k+1}$ is determined and every quantity used by the algorithm and its analysis is invariant to the unobserved component of $x$ in $\ker A\cap\ker C$.
\end{remark}

\paragraph{Family-level block scaling.} Finite-time analysis of a single problem usually hides the block-count dependence in its constants. In physics simulation, however, users change the problem dimension by refining the discretization, which gives a family of problems. We therefore make that dependence explicit through the following convention.
\begin{assume}[Family-level block-scaling convention]
\label{ass:block-scaling}
For the optional family-level rate interpretation, consider a sequence of instances with $N\to\infty$ and $M\geq1$, and fix $p\in(0,1)$.  All constants below are uniform in the family member and block index.  Every element of $\{\widehat C_i,C_i^\psi,H_i,\varphi_i(A_ix^1):i\in\IO\}
\cup\{M_i,L_i:i\in\IC\}\cup\{m\}$ is $\mathcal O(1)$, $\dim(x)=\mathcal O(N)$, and $M=\mathcal O(N^2)$. Moreover, there are constants $s_A\in\mathbb N$ and $0<c_A\leq C_A<\infty$ such that $\|\operatorname{vec}(A_i)\|_0\leq s_A$, and $c_A\leq\|\operatorname{vec}(A_i)\|_1\leq C_A$ for all $i\in\IAll$. Finally, $\varphi_i(A_ix^1)$ and $\widehat C_i$ ($i\in\IO$), $M_i$ ($i\in\IC$), and $m$ are bounded below by positive constants uniform in the family.
\end{assume}
\begin{figure}[H]
\vspace{-5px}
\centering
\begin{tikzpicture}[
  x=1cm,
  y=1cm,
  line cap=round,
  line join=round,
  strand/.style={draw=black!76, line width=0.8pt},
  ground/.style={draw=black!48, line width=0.65pt},
  clothnode/.style={circle, draw=black!75, fill=white, line width=0.55pt, inner sep=1.25pt}
]
\begin{scope}[xshift=-0.20cm]
\draw[strand]
  plot coordinates {
    (0.32,2.77) (0.76,2.91) (1.28,2.88) (1.88,2.78)
    (2.48,2.44) (2.64,1.98) (2.14,1.70) (1.39,1.64) (0.59,1.82)
  };
\foreach \p in {
  (0.32,2.77),(0.76,2.91),(1.28,2.88),(1.88,2.78),(2.48,2.44),
  (2.64,1.98),(2.14,1.70),(1.39,1.64),(0.59,1.82)
}
  \node[clothnode] at \p {};
\end{scope}

\begin{scope}[xshift=0.70cm]
\draw[strand]
  plot coordinates {
    (3.70,2.78) (4.30,2.95) (4.98,2.95) (5.66,2.87) (6.42,2.58)
    (7.05,2.18) (6.72,1.82) (5.92,1.63) (5.10,1.68) (4.28,1.91)
  };
\foreach \p in {
  (3.70,2.78),(4.30,2.95),(4.98,2.95),(5.66,2.87),(6.42,2.58),
  (7.05,2.18),(6.72,1.82),(5.92,1.63),(5.10,1.68),(4.28,1.91)
}
  \node[clothnode] at \p {};
\node[font=\scriptsize, text=black!82] at (3.72,3.42)
  {$\dim(x)=\mathcal O(N)$};
\draw[->, >=stealth, draw=black!68, line width=0.6pt]
  (3.10,3.38) arc[start angle=180,end angle=270,radius=0.60];
\draw[strand] (5.90,3.52) -- (6.66,3.23);
\node[clothnode] at (5.90,3.52) {};
\node[clothnode] at (6.66,3.23) {};
\draw[<->, >=stealth, draw=black!68, line width=0.6pt]
  (6.04,2.725) -- (6.28,3.375)
  node[midway, right, font=\scriptsize] {$A_i x=y_i,\ i\in\IAll$};
\draw[->, >=stealth, draw=black!68, line width=0.6pt]
  (5.66,2.87) -- (5.586,1.650)
  node[midway, left, font=\scriptsize] {$f_i(A_i x),\ i\in\IC$};
\draw[->, >=stealth, draw=black!68, line width=0.6pt]
  (4.28,1.91) -- (4.28,1.32)
  node[midway, left, font=\scriptsize] {$\varphi_i(A_i x),\ i\in\IO$};
\end{scope}

\begin{scope}[xshift=0.50cm]
\draw[strand]
  plot coordinates {
    (9.03,2.78) (9.35,2.86) (9.70,2.91) (10.07,2.94)
    (10.45,2.94) (10.83,2.90) (11.21,2.80) (11.58,2.64)
    (11.92,2.42) (12.20,2.18) (12.18,1.96) (11.91,1.80)
    (11.52,1.69) (11.09,1.64) (10.64,1.65) (10.19,1.71)
    (9.77,1.80) (9.40,1.92)
  };
\foreach \p in {
  (9.03,2.78),(9.35,2.86),(9.70,2.91),(10.07,2.94),(10.45,2.94),
  (10.83,2.90),(11.21,2.80),(11.58,2.64),(11.92,2.42),(12.20,2.18),
  (12.18,1.96),(11.91,1.80),(11.52,1.69),(11.09,1.64),(10.64,1.65),
  (10.19,1.71),(9.77,1.80),(9.40,1.92)
}
  \node[clothnode, inner sep=0.90pt] at \p {};
\end{scope}

\draw[ground] (-0.10,1.15) -- (13.50,1.15);
\draw[<->, >=stealth, draw=black!68, line width=0.6pt]
  (6.62,1.15) -- (6.62,1.63)
  node[midway, right, font=\scriptsize] {$m$};
\end{tikzpicture}
\caption{A family of problem instances simulating a 2D cloth falling on the ground.}
\label{fig:block-scaling-cloth}
\vspace{-5px}
\end{figure}
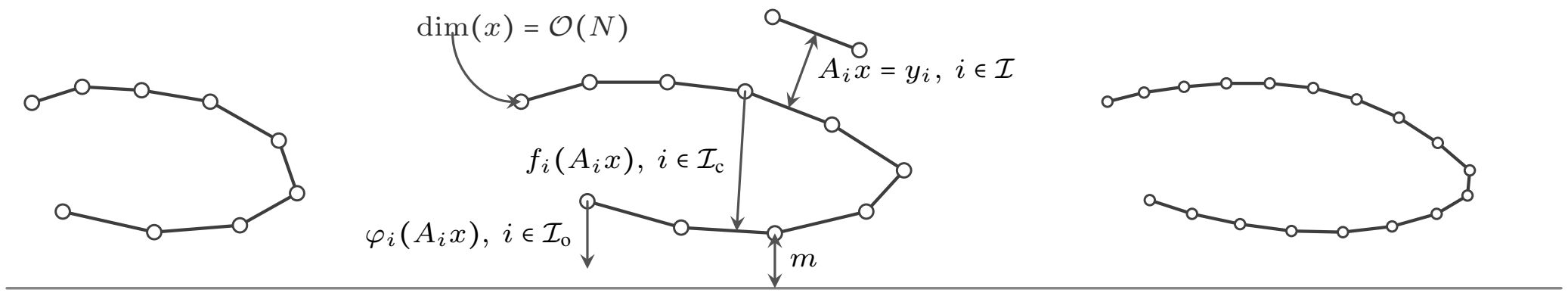

\prettyref{fig:block-scaling-cloth} illustrates \prettyref{ass:block-scaling} on a piece of cloth falling to the ground, discretized into $N$ elements of $2$ vertices each. Each element carries a fixed number of objective terms, such as gravitational and internal potential energy, so $\dim(x)=\mathcal O(N)$. Pairwise self-collision constraints give $M=\mathcal O(N^2)$, each $A_i$ gathers a bounded local stencil so $\|A_i\|_2=\Theta(1)$, and the other $\mathcal O(1)$ bounds encode fixed local energy, force, and curvature scales. The lower bounds on $\widehat C_i$ and $M_i$ normalize valid majorants, and that on $\varphi_i(A_ix^1)$ fixes an energy gauge, pinning each block's initial energy to a positive $\Theta(1)$ scale; the lower bound on $m$ is the substantive geometric condition, requiring a uniform initial self-collision and boundary margin. This gauge is a normalization choice: because the initial energy is defined only up to an additive constant, the family bounds use a fixed common per-block energy normalization and are not invariant to arbitrary additive energy shifts, such as the aggregate constant shift used to meet \prettyref{ass:obj} in \prettyref{appen:applications}. These conditions are natural for consistently normalized refinements that preserve clearance, but remain conditional under arbitrary mesh refinement~\cite{narain2016admm,overby2017admm,liu2026distributed}.

\subsection{Convergence of IP-ADMM}
Progress is measured by the reconstructed original-problem residuals $\CritWKKT$, $\CritbFJ$, and $\CritKKT$ defined in \prettyref{def:criteria}--\prettyref{eq:crit}.  These are nearby-point certificates: the objective and constraint (sub)gradients are evaluated at reconstructed block points $\breve y_i$ of~\prettyref{eq:reconstruction} that are close to, but not identical with, $A_ix$, with the gap controlled by the consensus term $\|Ax-\breve y\|$.  Evaluating stationarity at such nearby reconstructed points is standard for splitting and proximal methods on nonsmooth composite problems~\cite{rockafellar2009variational,yuan2025admm,jiagrimmer2025}.  Ordinary IP-ADMM returns the best certified candidate in the stated iteration window with the following guarantee.  Problem aggregates appearing below, such as $\widehat G_0$, are fixed in~\prettyref{eq:parameters}.

\begin{theorem}[Original-problem rates for IP-ADMM]
\label{thm:rate}
Run the ordinary IP-ADMM mode of~\prettyref{alg:IP-ADMM} with the parameter choice~\prettyref{eq:parameters} and polynomial schedule~\prettyref{eq:IP-schedules}.  Fix $p\in(0,1)$ and $\xi\in(0,1]$.  For the window $j=3,\ldots,K+1$ and $K\geq3$,
\begin{eqaligned}
\min_j\CritWKKT^j
&=\mathcal O(K^{(p-1)/2})+\mathcal O(K^{-p}),\\
\min_j\CritbFJ^j
&=\mathcal O(K^{(p-1)/2})+\mathcal O(K^{-p})
+\mathcal O(K^{-p/2})+\mathcal O(K^{-1}).
\end{eqaligned}
If, in addition,~\prettyref{ass:block-scaling} holds, the fixed-$\xi$ bounds specialize to $\mathcal O(\|A\|_2N\sqrt M\,K^{-1/3})$ for w-KKT at $p=1/3$ and $\mathcal O(\|A\|_2N\sqrt M\,K^{-1/4})$ for b-FJ at $p=1/2$.  For the optional tuned $p=1/2$ branch, define
\begin{eqaligned}
\label{eq:block-xi-star}
\xi_{\rm nom}\defeq m/(\|A\|_2N),
\qquad
\xi\defeq\min\{1,\xi_{\rm nom}\}.
\end{eqaligned}
On the interior branch $\xi_{\rm nom}<1$, write $\xi_\star\defeq\xi_{\rm nom}$.  Choosing $\xi=\xi_\star$ gives
\begin{eqaligned}
\min_j\CritbFJ^j
=\mathcal O\left(
\left[\|A\|_2N+(\|A\|_2^2\sqrt N/m)\right]K^{-1/2}
+\sqrt{\|A\|_2Nm}\,K^{-1/4}+K^{-1}
\right).
\end{eqaligned}
If $\xi_{\rm nom}\geq1$, the rule clips to $\xi=1$ and the fixed-$\xi$ bound applies.
\end{theorem}
\begin{assume}[Uniform Band-MFCQ]
\label{ass:MFCQ}
There are $\gamma_{\rm cq},\delta_{\rm cq}>0$ such that, whenever $Cx=d$, $f_i(y_i)\geq0$, and $\max_{i\in\IC}\|A_ix-y_i\|\leq\delta_{\rm cq}$, there is a unit $d_{\rm cq}\in\ker C$ satisfying
\begin{eqaligned}
\label{eq:cq}
\nabla f_i(y_i)^TA_id_{\rm cq}\geq\gamma_{\rm cq}\quad\text{whenever }0\leq f_i(y_i)\leq\delta_{\rm cq}.
\end{eqaligned}
The condition is monotone in the band radius, so replacing $\delta_{\rm cq}$ by $\min\{\delta_{\rm cq},1\}$ preserves it.  We henceforth use this normalized choice, so $0<\delta_{\rm cq}\leq1$.
\end{assume}
\begin{corollary}[KKT rates for IP-ADMM under Uniform Band-MFCQ]
\label{cor:KKT-MFCQ}
Run the ordinary IP-ADMM mode of~\prettyref{alg:IP-ADMM} with the parameter choice~\prettyref{eq:parameters}, $p=1/2$, $\xi\in(0,1]$, and polynomial schedule~\prettyref{eq:IP-schedules}, and suppose~\prettyref{ass:MFCQ} holds.  There is a finite, instance-dependent qualification threshold $K_{\rm cq}(1/2,\xi)\geq3$ such that, for every $K\geq K_{\rm cq}(1/2,\xi)$, $\min_{j=3,\ldots,K+1}\CritKKT^j=\mathcal O(K^{-1/4})$. If, in addition,~\prettyref{ass:block-scaling} holds, define $D_{\rm cq}\defeq2(1+\|A\|_2\sqrt{\widehat G_0}/\gamma_{\rm cq})=\Theta(1+\|A\|_2\sqrt N/\gamma_{\rm cq})$; then the fixed-$\xi$ upper bound is $\mathcal O(D_{\rm cq}\|A\|_2N\sqrt M\,K^{-1/4})$.
\end{corollary}

\paragraph{Proof sketch.}
Our schedule increases the penalty and decreases the smoothing, following the continuation idea of IPDS-ADMM~\cite{yuan2025admm}, whose analysis builds a Lyapunov function. Our nonlinear wall, however, raises a new issue: the continuation argument holds only while the split point stays strictly feasible. To handle this, \prettyref{lem:terminal-bound} lower-bounds the Lyapunov function without assuming that the split point is strictly feasible. \prettyref{thm:feasibility} then combines this bound with the descent of the Lyapunov candidate to prove strict feasibility by induction. The drift-corrected lower bound \prettyref{lem:lower-bound} gives the movement budget \prettyref{cor:movement-budget}. Exact objective-force reconstruction turns this budget into original-problem residuals in \prettyref{thm:rate-master}, and \prettyref{cor:KKT-MFCQ-exact} converts the same residual to KKT under Uniform Band-MFCQ. The family-level and tuned-$\xi$ bounds follow as specializations. Full constants and the dependency chain are in \prettyref{appen:convergence}.

\paragraph{Interpretation of the continuation balance.}
IPDS-ADMM allows problem-dependent continuation amplitudes $\xi$, including a condition-number-dependent choice in one regime~\cite{yuan2025admm}. Our family-level result instead identifies the explicit nominal balance $\xi_{\rm nom}=m/(\|A\|_2N)$ for the nonlinear-inequality split. It changes the leading $K^{-1/4}$ coefficient from the fixed-$\xi$ family scale $\|A\|_2N\sqrt M$ to $\sqrt{\|A\|_2Nm}$, improving the block-scale dependence of that leading high-accuracy term. It does not improve the finite-$K$ bound: the $K^{-1/2}$ terms and the clipping condition remain.

\subsection{Convergence of IP-ADMM-Ckpt}
Under the CQ condition, the checkpoint mechanism yields an improved asymptotic certificate rate after the first completed checkpoint, as stated below. This is a statement about the latest-checkpoint certificate as $K\to\infty$, not a claim of faster progress at a fixed budget: no certified output exists before the first checkpoint completes at iteration $K_0$, and the family bound on $K_0$ can be large.
\begin{theorem}[Checkpoint KKT rates for IP-ADMM-Ckpt under Uniform Band-MFCQ]
\label{thm:checkpoint-rate}
Run the $\mathsf{Ckpt}$ mode of~\prettyref{alg:IP-ADMM} with checkpoint test $E^{k+1}\leq\tau_{\rm cp}/\rho^k$, the explicit initialization \prettyref{eq:checkpoint-parameters}, and the geometric checkpoint promotion in \prettyref{alg:IP-ADMM}, and suppose \prettyref{ass:MFCQ} holds. Let $K_0$ be the iteration budget at which checkpoint $0$ completes, and for $K\geq K_0$ let ${\CritKKT}_{\rm out}(K)$ denote the residual of the latest-checkpoint output.  Then ${\CritKKT}_{\rm out}(K)=\mathcal O(K^{-1/3})$ whenever $K\geq K_0$, and no certified output exists for $K<K_0$. If, in addition, \prettyref{ass:block-scaling} holds across a family, the Uniform Band-MFCQ constants are uniformly bounded away from zero, and $\chi$, $\ell_{\rm cp}$, and $\tau_{\rm cp}$ are fixed, then ${\CritKKT}_{\rm out}(K)=\mathcal O\left(\|A\|_2N^{4/3}K^{-1/3}\right)$ whenever $K\geq K_0$. In this family the startup budget obeys $K_0=\mathcal O\left(N\left(N+\|A\|_2^2/\gamma_{\rm cq}^2\right)^3\right)$, so the improved exponent describes only post-startup behavior and need not reduce the total iteration budget.
\end{theorem}

\paragraph{Proof sketch.}
The checkpoint gate is loosely motivated by the progress-tested multiplier updates of adaptive proximal ADMM~\cite{maia2026adaptive}, but our proof uses a pointwise level induction. We condition each level on two caps: a cap on the stored multiplier and a cap on the entry merit. Under these caps, \prettyref{lem:ckpt-arrival} charges every rejection against the frozen-level merit and proves finite arrival at each level $j$ within $\mathcal O(\rho_j^3)$ iterations, where $\rho_j$ is the level-$j$ penalty. Independently, \prettyref{lem:ckpt-accepted} uses split feasibility, support control, and Uniform Band-MFCQ to bound the wall forces, the accepted virtual multiplier, and the consensus by $\mathcal O(\rho_j^{-1})$. The promotion charge is then summable, so \prettyref{lem:ckpt-transition} propagates both caps and closes the induction. Each completed level thus carries the pointwise $\mathcal O(\rho_j^{-1})$ KKT certificate of \prettyref{lem:ckpt-certificate}, and \prettyref{cor:ckpt-complexity} combines it with geometric level growth to give the latest-checkpoint rate $\mathcal O(K^{-1/3})$. Family and multiplier scaling are developed in \prettyref{appen:checkpoint}.

\paragraph{Multiplier bounds.}
The two analyses provide complementary multiplier controls.  For ordinary IP-ADMM, without any constraint qualification, the rescaled invariant $\sup_j\|\lambda^j\|^2/\rho^j<\infty$ holds, so $\lambda^j/\rho^j$ is bounded, as in IPDS-ADMM~\cite{yuan2025admm}, and tends to zero for $p>0$; the invariant bounds the raw multipliers only at the constant-penalty endpoint $p=0$.  For IP-ADMM-Ckpt, Uniform Band-MFCQ bounds the raw stored and accepted consensus multipliers and the accepted wall forces by constants fixed at initialization, so no rescaling is needed; rejected virtual multipliers are not covered.  See \prettyref{lem:scaled-multiplier} and \prettyref{cor:ckpt-multiplier-scaling}.

\section{Numerical Evaluation}
We evaluate capsule stacking and cloth resting for $n\in\{8,16,24\}$, with one deterministic run per family--size--method case. We compare against three baselines. SciPy's SLSQP and trust-constr, the standard active-set SQP and trust-region interior-point methods, are the general-purpose baselines used in every solver comparison (\prettyref{appen:experiments}), and BC-ADMM~\cite{pan2026bcadmm} is the closest published optimizer to ours, since it is also a parallel ADMM that starts from a feasible point and keeps every iterate original-space feasible. Every method receives the same strictly feasible initialization and the same analytic first derivatives, and trust-constr additionally receives exact Hessians. Across the $27$ runs we test every certificate that each method can produce: w-KKT, b-FJ, and KKT for both IP-ADMM variants; w-KKT and KKT for SLSQP and trust-constr, which have no b-FJ certificate; and, for BC-ADMM, only the w-KKT residual of its own relaxation. On the small instances ($n=8$), SLSQP, trust-constr, and both IP-ADMM variants reach the KKT threshold of $10^{-3}$, although on capsule stacking trust-constr does so only at a configuration with crossing medial axes, in which two capsules pass through each other. On the larger instances ($n\in\{16,24\}$), SLSQP, trust-constr, and both IP-ADMM variants still reach the w-KKT threshold, and the IP-ADMM variants also reach b-FJ, but neither of them reaches KKT, and SLSQP reaches KKT only on cloth resting, not on capsule stacking. \rev{trust-constr reaches KKT on cloth resting and on capsule stacking at $n=24$ but not $n=16$. Its capsule outputs at $n=8$ and $n=24$ satisfy every constraint yet contain crossing medial axes, which the endpoint--segment collision constraints cannot detect once a long step has jumped into the crossing (\prettyref{appen:applications}); every other output passes the exact-geometry checks.} This is consistent with CQ degeneracy or difficulty of KKT certification, but our experiments do not distinguish these explanations. Either way, it agrees with prior works~\cite{pan2026bcadmm,li2020ipc} that certify only w-KKT or b-FJ solutions. We further observe that the configuration settles early: the normalized state-space $L_2$ distance of each iterate to the final output (\prettyref{fig:visual-difference}) falls below $5\%$ of its initial value well before KKT would be certified, so an early stop, or a stopping rule based on the relative change in the physical configuration, is a useful practical heuristic in these experiments. \prettyref{appen:experiments} gives the protocol and the per-run numbers, and \prettyref{appen:ablation} connects the theorem-sufficient parameters to the tuned values used here. We further show that our method has applications beyond physics simulation, where \prettyref{appen:control} optimizes a UAV trajectory whose objective carries a sharp $\ell_1$ goal term, which our proximal oracle handles exactly. Finally, \prettyref{appen:scaling} reports real wall-clock time on the cloth testbed.

\paragraph{\rev{Learned contact constraints.}} \rev{\prettyref{appen:neural} tests a setting in which a trained network defines the constraints. A frozen $\tanh$ network, trained as the SDF of a nonconvex cross, supplies every body--body contact constraint while $n\in\{4,8,16\}$ crosses settle under gravity in a hexagonal container, which gives $69{,}120$ to $645{,}120$ constraints on $3n$ pose variables. Running on a GPU, IP-ADMM reaches a common target of b-FJ, feasibility, and force and torque balance in $174$, $1742$, and $1568$\,s, respectively. SciPy's SLSQP diverges in every case, and trust-constr needs $787$\,s at $n=4$ and does not reach it within the $3600$\,s budget at $n=8$ and $n=16$. Every certified configuration passes exact polygon checks.}

\section{Limitations and Conclusion}
This work lays a theoretical foundation for parallel physics optimization. To our knowledge, IP-ADMM is the first single-loop method that updates nonconvex inequality constraints locally and in parallel with finite-time guarantees, combining independent nonlinear-inequality updates, one equality-constrained quadratic consensus solve per iteration, a CQ-explicit hierarchy of finite-time certificates, and conditional family-aware bounds in $N$, $M$, and $\|A\|_2$ under a family-scaling convention. Together these give finite-time, original-problem guarantees that respect the local--global splitting used in practice.
Several limitations temper these guarantees. First, the KKT results rely on Uniform Band-MFCQ, whose constants $(\gamma_{\rm cq},\delta_{\rm cq})$ we cannot compute or even lower-bound in our applications. Second, our experiments use practical parameters outside the theorem-sufficient construction (\prettyref{eq:parameters}), as is common for first-order nonconvex methods~\cite{yuan2025admm,jiagrimmer2025,pan2026bcadmm}; \prettyref{appen:ablation} reports the sufficient construction itself. Third, the family-aware rates are conditional on the scaling convention and its energy-gauge normalization rather than intrinsic properties of the algorithm. Finally, our main convergence benchmarks remain small, static, single-threaded, and 2D; the additional multicore scaling study and the GPU study with learned constraints (\prettyref{appen:neural}) provide only limited scalability evidence.
These gaps point to concrete next steps. Distributed implementations and larger GPU benchmarks would test the parallel structure at scale; and richer benchmarks would stress the family-scaling assumptions. Integrating the solver into end-to-end learning pipelines and world models is a natural path toward large-scale Physics+AI.

\clearpage
\section*{AI Use Statement}

(This section is required by the ICLR 2027 AI Policy for Authors and does not count toward the page limit.)

We used generative AI tools in preparing this manuscript and disclose their use here in accordance with the ICLR 2027 AI Policy for Authors, which covers both recommended- and required-disclosure activities. Concretely, we used generative AI to (i)~search for and identify related literature; (ii)~edit and condense author-written text for clarity and length; (iii)~perform adversarial checking of the mathematical arguments and proofs, flagging gaps, unstated assumptions, and imprecise claims for the authors to resolve. The authors independently derived, checked, and verified all mathematical claims and proofs, and designed, implemented, ran, and verified all experiments, including data generation and analysis. Every AI-assisted contribution---references, edits, arguments, and interpretations---was reviewed and validated by the authors against the primary sources and the underlying mathematics before inclusion. The authors take full responsibility for the final content of this work, including all text, claims, proofs, and artifacts.

\bibliographystyle{iclr2027_conference}
\bibliography{template}

\clearpage
\appendix
\section{\texorpdfstring{\rev{Comparison with Prior Methods}}{Comparison with Prior Methods and Detailed Rates}}
\label{appen:prior}
\rev{\prettyref{table:prior} compares IP-ADMM with prior first-order methods for nonconvex constrained optimization, and \prettyref{table:rates} lists the detailed bounds of both IP-ADMM variants. A method is \emph{single-loop} if each iteration costs a fixed number of oracle calls and linear solves, with no inner loop whose length depends on the target accuracy $\epsilon$. Its constraint update is \emph{local} if each constraint is updated independently of the others from its own block, so all constraint updates can run in parallel, and \emph{global} if all constraints enter one joint subproblem, penalty, or linearization. One IP-ADMM sweep costs one proximal call per objective block, one closed-form update per constraint block, and one consensus solve. Dividing the first block row of the optimality system of the $x$-update by $\rho^k$ shows that its matrix, $(1+\theta_{2,\mathsf m})A^TA$ bordered by $C$, does not depend on $k$, so one factorization serves the whole run in both modes; a rejected IP-ADMM-Ckpt step also costs exactly one sweep, so both variants are single-loop. Since certificates, regularity conditions, and cost units differ across rows, \prettyref{table:prior} supports comparing exponents only; for instance, the certificate of~\cite{jiagrimmer2025} is stronger than our b-FJ, with exact feasibility and $\epsilon^2$ complementarity. Local--global physics solvers are listed for their architecture; we found no finite-time original-problem rate reported for them (N/A). Among the methods in \prettyref{table:prior}, IP-ADMM is the only one that combines a single loop, local updates of nonconvex constraints, and a finite-time rate. Local--global physics solvers report no rate. BC-ADMM solves its local subproblems iteratively and certifies a fixed-barrier relaxation. The other single-loop methods update their constraints globally. Here a nonconvex constraint is a smooth function accessed through its values and gradients, each updated independently in closed form. IPDS-ADMM also admits nonconvex sets, but only through an exact proximal (projection) oracle and with sequential block updates.}

\begin{table}[ht]
\centering
\caption{\rev{Comparison with prior first-order methods for nonconvex constrained optimization. Loop: S = single loop, D = double loop, SI = subproblem per iteration. Constr. (constraints): A = affine, NE = nonconvex equalities, NI = nonconvex inequalities, NEI = both; $^\dagger$ stochastic setting. Constr. Update (constraint update): L = local, G = global. Bold marks desirable properties: a single loop, nonconvex constraints, local constraint updates, and no CQ. Each certificate is listed on its own line with its CQ and the iteration or oracle complexity to reach accuracy $\epsilon$ ($\mu$ in~\cite{hinder2024worst}); relaxed stationary refers to the fixed-barrier biconvex relaxation of~\cite{pan2026bcadmm}, whose local subproblems are solved iteratively; cost units differ across methods; $\tilde{\mathcal O}$ hides logarithmic factors.}}
\label{table:prior}
\begin{revision}\scriptsize
\setlength{\tabcolsep}{3pt}
\renewcommand{\arraystretch}{1.1}
\begin{tabular}{@{}>{\raggedright\arraybackslash}p{0.30\linewidth}cccll@{}}
\toprule
Method & Loop & Constr. & Constr. Update & Certificate & CQ / regularity\\
\midrule
\multirow{2}{=}{Jia--Grimmer~\cite{jiagrimmer2025}} & \multirow{2}{*}{D} & \multirow{2}{*}{\textbf{NI}} & \multirow{2}{*}{G} & FJ: $\mathcal O(\epsilon^{-4})$ & \textbf{none}\\
&&&& KKT: $\mathcal O(\epsilon^{-4})$ & CQ\\
\midrule
\multirow{2}{=}{Lin--Ma--Xu~\cite{lin2019ippp}} & \multirow{2}{*}{D} & \multirow{2}{*}{\textbf{NEI}} & \multirow{2}{*}{G} & stationary: $\tilde{\mathcal O}(\epsilon^{-3})$ & non-singularity\\
&&&& weak stationary: $\tilde{\mathcal O}(\epsilon^{-4})$ & feasible start\\
\midrule
Li et al.~\cite{li2021ialm} & D & \textbf{NE} & G & KKT: $\tilde{\mathcal O}(\epsilon^{-3})$ & regularity along iterates\\
\midrule
Hinder--Ye~\cite{hinder2024worst} & SI & \textbf{NI} & G & FJ: $\mathcal O(\epsilon^{-7/4})$ & \textbf{none}; strictly feasible start\\
\midrule
BC-ADMM~\cite{pan2026bcadmm} & D & \textbf{NI} & \textbf{L} & relaxed stationary: $\mathcal O(\epsilon^{-2})$ & \textbf{none}; biconvex relaxation\\
\midrule
IPDS-ADMM~\cite{yuan2025admm} & \textbf{S} & A & G & critical: $\mathcal O(\epsilon^{-3})$ & \textbf{none}\\
\midrule
Yang et al.~\cite{yang2025singleloop} & \textbf{S} & \textbf{NI}$^\dagger$ & G & nearly KKT: $\mathcal O(\epsilon^{-6})$ & CQ-type regularity\\
\midrule
El Bourkhissi--Necoara~\cite{elbourkhissi2025} & \textbf{S} & \textbf{NE} & G & KKT: $\mathcal O(\epsilon^{-2})$ & $\sigma_{\min}(G)\geq\sigma$; bounded iterates\\
\midrule
Liu et al.~\cite{liu2026nrlalm} & \textbf{S} & \textbf{NE} & G & KKT: $\mathcal O(\epsilon^{-2})$ & uniform LICQ\\
\midrule
\multirow{3}{=}{IP-ADMM (ours)} & \multirow{3}{*}{\textbf{S}} & \multirow{3}{*}{\textbf{NI}, A} & \multirow{3}{*}{\textbf{L}} & w-KKT: $\mathcal O(\epsilon^{-3})$ & \textbf{none}\\
&&&& b-FJ: $\mathcal O(\epsilon^{-4})$ & \textbf{none}\\
&&&& KKT: $\mathcal O(\epsilon^{-4})$ & Uniform Band-MFCQ\\
\midrule
IP-ADMM-Ckpt (ours) & \textbf{S} & \textbf{NI}, A & \textbf{L} & KKT: $\mathcal O(\epsilon^{-3})$ after $K_0$ & Uniform Band-MFCQ\\
\midrule
Local--global physics ADMM~\cite{narain2016admm,overby2017admm,lee2023subsystem,zhu2024frictional,liu2026distributed} & \textbf{S} & \textbf{NEI} & \textbf{L} & N/A & --\\
\bottomrule
\end{tabular}
\end{revision}
\end{table}

\begin{table}[ht]
\centering
\caption{Convergence rate upper bounds of the IP-ADMM variants. MFCQ means Uniform Band-MFCQ. Family bounds assume \prettyref{ass:block-scaling}, under which $m=\mathcal O(1)$ and $M=\mathcal O(N^2)$; we omit $m$ for brevity but keep $M$ explicit, since practical instances often have $M\ll N^2$.  For IP-ADMM-Ckpt, $K_0$ is the startup budget at which the first checkpoint completes; no certified output exists for $K<K_0$, and the family bound on $K_0$ can be large (\prettyref{thm:checkpoint-rate}). The checkpoint family rate carries no separate $M$ outside $\|A\|_2$ only conditional on family-uniform qualification margins (\prettyref{ass:MFCQ}).  The w-KKT, b-FJ, and KKT certificates in the Cert.\ column are nearby-point residuals (\prettyref{def:criteria}).}
\label{table:rates}
\resizebox{\linewidth}{!}{%
\begin{tabular}{@{}llllll@{}}
\toprule
Mode & Cert. & CQ & Validity & General & Family \\
\midrule
\multirow{4}{*}{IP-ADMM} & w-KKT & none & $p=1/3$, fixed $\xi$, $K\geq3$ & $\mathcal O(K^{-1/3})$ & $\mathcal O(\|A\|_2N\sqrt M K^{-1/3})$ \\
& b-FJ & none & $p=1/2$, fixed $\xi$, $K\geq3$ & $\mathcal O(K^{-1/4})$ & $\mathcal O(\|A\|_2N\sqrt M K^{-1/4})$ \\
& b-FJ & none & $p=1/2$, $\xi=\xi_\star$, $\xi_{\rm nom}<1$, $K\geq3$ & $\mathcal O(K^{-1/4})$ & $\begin{gathered}\mathcal O\bigl((\|A\|_2N+\|A\|_2^2\sqrt N)K^{-1/2}\\{}+\sqrt{\|A\|_2N}\,K^{-1/4}+K^{-1}\bigr)\end{gathered}$ \\
& KKT & MFCQ & $p=1/2$, fixed $\xi$, $K\geq K_{\rm cq}(1/2,\xi)$ & $\mathcal O(K^{-1/4})$ & $\mathcal O(D_{\rm cq}\|A\|_2N\sqrt M K^{-1/4})$ \\
\midrule
IP-ADMM-Ckpt & KKT & MFCQ & $K\geq K_0$ & $\mathcal O(K^{-1/3})$ & $\mathcal O(\|A\|_2N^{4/3}K^{-1/3})$ \\
\bottomrule
\end{tabular}}
\end{table}

\begin{revision}
\paragraph{Finite-band contact solutions as b-FJ points.}
The following observation connects fixed-band contact models with the b-FJ certificate of \prettyref{def:criteria}.
\begin{proposition}[Finite-band solutions are b-FJ points]
\label{prop:finite-band}
Let $Cx=d$, set $\breve y_i=y_i=A_ix$ for all $i\in\IAll$, and suppose $f_i(y_i)\geq0$ for $i\in\IC$. Suppose $v_i\in\partial\varphi_i(y_i)$, $\iota_i\leq0$, and $\zeta$ satisfy
\begin{equation*}
\Bigl\|\sum_{i\in\IO}A_i^Tv_i+\sum_{i\in\IC}\iota_iA_i^T\nabla f_i(y_i)+C^T\zeta\Bigr\|\leq r,
\qquad
\iota_i=0\ \text{whenever}\ f_i(y_i)>b.
\end{equation*}
Then $x$ is an $\epsilon$-b-FJ point for every $\epsilon\geq\max\{r,b\}$.
\end{proposition}
\begin{proof}
Set $\alpha=1/(1+\|\iota\|_1)$, $\bar\iota=\alpha\iota$, and $\bar\zeta=\alpha\zeta$. Then $\alpha>0$, $\bar\iota\leq0$, and $\alpha+\|\bar\iota\|_1=1$. The consensus residual is zero, and the normalized stationarity residual is $\alpha$ times the original one, hence at most $\alpha r\leq\epsilon$. If $f_i(y_i)>b$, then $\bar\iota_i=0$; otherwise $|\bar\iota_if_i(y_i)|\leq|\bar\iota_i|\,b\leq b\leq\epsilon$, since $|\bar\iota_i|\leq1$. Finally, $f_i(y_i)>\epsilon\geq b$ implies $\bar\iota_i=0$, which is the support condition in~\prettyref{eq:eps-bFJ}.
\end{proof}
For frictionless IPC~\cite{li2020ipc}, the momentum balance of its Section~4.5, Eqs.~(8)--(9), has the form $\nabla E(x)-\sum_k\lambda_k\nabla d_k(x)=r$ for the incremental potential $E$, with contact forces $\lambda_k\geq0$ given by barrier derivatives that vanish for $d_k(x)\geq\hat d$; this is its $\epsilon$-complementarity, Section~2.1, Eq.~(3), with $\epsilon=\hat d$. Taking $f_k=d_k$ and $\iota_k=-\lambda_k$, \prettyref{prop:finite-band} shows that an IPC solution with force residual $r$ is a $\max\{\|r\|,\hat d\}$-b-FJ point. This implication is ours; IPC does not state a Fritz--John result, and its piecewise-smooth distances need not satisfy our standing assumptions, so the statement concerns the certificate only. A fixed band $\hat d$ does not converge to zero-band stationarity, whereas the b-FJ rate of \prettyref{thm:rate} drives the band to zero.
\end{revision}

\section{Application Formulations and Assumption Justifications}
\label{appen:applications}
This appendix records the exact block structure of the two benchmark families.  In both cases $C$ and $d$ are empty and $\psi_i\equiv0$.

\paragraph{Compact-domain realization.}
Suppose exact boundary blocks confine every feasible configuration to a compact set $\Omega$ (modulo exact periodic coordinates), and each remaining raw $f_i$ and $h_i$ is $C^{1,1}$ on an open neighborhood of $A_i\Omega$. Then \prettyref{ass:f}--\prettyref{ass:bg} hold automatically. Indeed, on this compact set each block has bounded value and gradient and a finite gradient-Lipschitz constant. A Whitney extension of this value--gradient jet, followed by a smooth cutoff outside a larger neighborhood, yields a globally defined block that agrees in value and gradient on $A_i\Omega$. The constraint extensions then have globally Lipschitz and bounded gradients, which supply \prettyref{ass:f}--\prettyref{ass:bg}. The objective extensions are globally Lipschitz with Lipschitz gradients, and an aggregate constant shift supplies the lower bound in \prettyref{ass:obj} without changing derivatives or minimizers. This shift is a gauge choice: it leaves derivatives and minimizers unchanged, but it does change absolute energy aggregates such as those entering the family-level constants of \prettyref{ass:block-scaling}, so those family bounds are stated relative to a fixed common energy normalization. Because the exact boundary blocks force feasible points into $\Omega$, the raw and extended constrained objectives agree on their common feasible envelope. This is the BC-ADMM extension argument of~\cite{pan2026bcadmm}. We use it to show that our two applications satisfy all the problem assumptions.

\paragraph{Distance regularization.}
To satisfy our problem assumptions, we locally mollify the objective and constraint functions using the helper functions defined below. Similar techniques have been used in~\cite{pan2026bcadmm,li2020ipc}. For the spring scale $\epsilon>0$, define the pseudo-Huber length modulation
\begin{eqaligned}
\mathscr M^{\rm len}_{\epsilon}(r)
=\sqrt{\|r\|^2+\epsilon^2}-\epsilon,
\end{eqaligned}
and for a squared distance $t\geq0$ define
\begin{eqaligned}
\mathscr M^{\rm rad}_{\delta}(t)
=\delta^2\left(\sqrt{1+t/\delta^2}-1\right).
\end{eqaligned}
The latter equals $\delta\mathscr M^{\rm len}_{\delta}(r)$ when $t=\|r\|^2$.  It is smooth and strictly increasing, its residual-space gradient is globally bounded by $\delta$, and that gradient is $1$-Lipschitz.  Consequently,
\begin{eqaligned}
\mathscr M^{\rm rad}_{\delta}(D)-\mathscr M^{\rm rad}_{\delta}(s^2)\geq0
\quad\Longleftrightarrow\quad D\geq s^2.
\end{eqaligned}
Thus radial modulation preserves the threshold set of the squared distance supplied to it.
For a point $p$ and a segment with endpoints $u,v$, and activation margin constants $\ell_{\rm act}>0$, let
\begin{eqaligned}
D_{\rm seg}(p;u,v)=\min_{t\in[0,1]}\|p-(1-t)u-tv\|^2,\;
D_0(p;u,v)=\left\|p-\frac{u+v}{2}\right\|^2,\;
\varsigma=\|v-u\|^2.
\end{eqaligned}
With $S(z)=6z^5-15z^4+10z^3$, define the blend
\begin{eqaligned}
\operatorname{blend}_{\ell_{\rm act}}(\varsigma)=
\begin{cases}
0,&\varsigma\leq\ell_{\rm act}^2/4,\\
S\!\left(\dfrac{\varsigma-\ell_{\rm act}^2/4}{3\ell_{\rm act}^2/4}\right),
&\ell_{\rm act}^2/4<\varsigma<\ell_{\rm act}^2,\\
1,&\varsigma\geq\ell_{\rm act}^2,
\end{cases}
\end{eqaligned}
and the short-segment modulation
\begin{eqaligned}
\mathscr M^{\rm seg}_{\ell_{\rm act}}(p;u,v)
=(1-\operatorname{blend}_{\ell_{\rm act}}(\varsigma))D_0
+\operatorname{blend}_{\ell_{\rm act}}(\varsigma)D_{\rm seg}.
\end{eqaligned}
The quintic blend has zero first and second derivatives at both junctions.  Since $D_{\rm seg}$ is jointly $C^{1,1}$ wherever $u\ne v$, including the endpoint/interior projection transitions, $\mathscr M^{\rm seg}_{\ell_{\rm act}}$ is endpoint-symmetric and $C^{1,1}$ through $u=v$.  It equals the standard squared point--segment distance exactly when $\|v-u\|\geq\ell_{\rm act}$.  Because $D_{\rm seg}\leq D_0$, the positive-margin modulation can overestimate the exact distance inside its short activation band; there it is a local model regularization rather than an exact reparameterization.  Thus $\mathscr M^{\rm seg}_{\ell_{\rm act}}$ provides a globally smooth surrogate collision constraint.  The exact distance $D_{\rm seg}$ is nonsmooth when $\varsigma=0$; the smooth blending keeps $\mathscr M^{\rm seg}_{\ell_{\rm act}}$ away from this nonsmooth configuration as $\varsigma\to0$.

\paragraph{Capsule stacking.}
Capsule $i$ has pose $z_i=(c_i,\vartheta_i)\in\R^2\times\R$, and $x=(z_1,\ldots,z_n)\in\R^{3n}$.  With planar rotation $R(\vartheta_i)$ and fixed body-coordinate medial-axis endpoints $a_e$, $e\in\{1,2\}$, its world endpoint is
\begin{eqaligned}
p_{ie}(z_i)=c_i+R(\vartheta_i)a_e.
\end{eqaligned}
There is one gravity block per capsule, $h_i(z_i)=e_2^Tc_i$, so $N=n$.  For each ordered source endpoint $(i,e)$ and target capsule $j\ne i$, let $D_{iej}(z_i,z_j)$ be its squared distance to capsule $j$'s medial segment.  With capsule radius $r_{\rm cap}$ and clearance $s=2r_{\rm cap}$, the capsule--capsule collision block is
\begin{eqaligned}
f_{iej}=\mathscr M^{\rm rad}_{\delta}(D_{iej})
-\mathscr M^{\rm rad}_{\delta}(s^2).
\end{eqaligned}
Each endpoint also gives the capsule--container boundary block
\begin{eqaligned}
f_{ie}^{\rm circ}=\mathscr M^{\rm rad}_{\delta}(R_{\rm eff}^2)
-\mathscr M^{\rm rad}_{\delta}(\|p_{ie}-c_{\rm circ}\|^2),
\end{eqaligned}
where $c_{\rm circ}$ is the container center, $R_{\rm circ}$ the container radius, and $R_{\rm eff}=R_{\rm circ}-r_{\rm cap}$.  Hence the dimension is $3n$, $N=n$, and $M=2n^2$.

\paragraph{Assumption justification.}
The exact boundary blocks place every medial endpoint in a closed disk and therefore bound every center.  Orientation is periodic, so one-period representatives give a compact center--orientation envelope; the raw trigonometric blocks are its periodic pullback.  Every target medial segment has fixed positive length, hence $D_{iej}$ is $C^{1,1}$ on a neighborhood of each compact local-pose image, and radial composition preserves this property.  The compact-domain realization above therefore supplies \prettyref{ass:f}--\prettyref{ass:bg}. Gravity is smooth on the compact center domain, its extension supplies the smooth clauses and lower bound of \prettyref{ass:obj}, and $\psi_i\equiv0$ is closed, convex, globally Lipschitz, and has the identity proximal map. The separated initial poses have strict pairwise and boundary clearance, verifying \prettyref{ass:feasible}. Thus the extended capsule formulation satisfies all standing assumptions.

\paragraph{Cloth resting.}
Let $x=(q_1,\ldots,q_n)\in\R^{2n}$ collect the planar node positions. The gravity energy is $h_i^{\rm grav}(q_i)=m_ge_2^Tq_i$ with gravity coefficient $m_g>0$ (node mass times gravity), and a stretch or bending energy on pair $(i,j)$ with stiffness $k_{ij}>0$ and rest length $\ell_{ij}$ is defined as
\begin{eqaligned}
h_{ij}^{\rm spring}(q_i-q_j)
=\frac{k_{ij}}{2}
\left(\mathscr M^{\rm len}_{\epsilon}(q_i-q_j)-\ell_{ij}\right)^2.
\end{eqaligned}
We add one stretch spring energy for each segment, and one bending spring energy for each pair of endpoints of two consecutive segments. Thus $N=n+(n-1)+(n-2)=3n-3$. For every node $i$ and nonincident segment $(q_j,q_{j+1})$, the theoretical regularized self-collision block takes $\ell_{\rm act}>0$ and uses
\begin{eqaligned}
f_{ij}^{\rm self}
&=\mathscr M^{\rm rad}_{\delta}
\!\left(\mathscr M^{\rm seg}_{\ell_{\rm act}}(q_i;q_j,q_{j+1})\right)
-\mathscr M^{\rm rad}_{\delta}(s^2),\\
f_i^{\rm circ}
&=\mathscr M^{\rm rad}_{\delta}(R_{\rm circ}^2)
-\mathscr M^{\rm rad}_{\delta}(\|q_i-c_{\rm circ}\|^2).
\end{eqaligned}
Here $s>0$ is the cloth node--segment clearance, and $R_{\rm circ}$ and $c_{\rm circ}$ are the container radius and center. The ordered nonincident construction has $(n-1)(n-2)$ blocks and the boundary conditions add $n$, so $M=(n-1)(n-2)+n$.

\paragraph{Assumption justification.}
The exact circle blocks confine every node to a closed disk, making the global configuration envelope compact. Positive $\ell_{\rm act}$ makes every short-segment squared distance $C^{1,1}$ even at collapsed adjacent edges; radial composition and compact-domain realization then supply \prettyref{ass:f}--\prettyref{ass:bg}. Gravity and the positive-$\epsilon$ spring blocks are $C^{1,1}$ on the compact objective images, so their extensions supply the smooth clauses and lower bound of \prettyref{ass:obj}; again $\psi_i\equiv0$ has the required convex proximal structure. The initialized strand has strict boundary and nonincident-clearance margins. Moreover, $\mathscr M^{\rm seg}_{\ell_{\rm act}}\geq D_{\rm seg}$, so the regularization preserves that strict clearance and verifies \prettyref{ass:feasible}. Thus the globally extended, positive-margin cloth formulation satisfies all standing assumptions.

\paragraph{What the collision blocks encode.} The capsule and cloth collision blocks are point-to-segment clearances---endpoint-to-medial-segment for capsules and node-to-nonincident-segment for cloth---and the construction includes both directions of every nonincident pair. These are surrogate contact conditions. In the plane, the distance between two disjoint segments is attained at an endpoint of one of them, so with both directions included, the unsmoothed point-to-segment clearances keep every non-crossing pair at least the clearance apart. They cannot, however, detect a pair whose segments already cross, since all four endpoint-to-opposite-segment distances can then remain positive. A continuous motion cannot reach such a crossing without first violating an endpoint clearance, but a long discrete step can jump into one. We therefore treat these blocks as surrogate benchmarks rather than a formal no-crossing guarantee, and certify the absence of crossings on every reported output with an independent exact geometric check that rejects capsule overlap or medial-axis crossing and cloth self-crossing or degenerate segments (\prettyref{appen:experiments}).

\section{Numerical Protocol and Additional Results}
\label{appen:experiments}
The evaluation contains capsule stacking and cloth resting at sizes $8$, $16$, and $24$, with one deterministic SLSQP, trust-constr, IP-ADMM, and IP-ADMM-Ckpt run per case, plus one BC-ADMM~\cite{pan2026bcadmm} run per cloth case, since BC-ADMM does not apply to capsule stacking (explained below). All methods receive the same strictly feasible initial configuration and analytic first derivatives; trust-constr additionally receives exact Hessians, as described next. Our empirical experimental setup and parameter choices are summarized in \prettyref{table:experiment-parameters}. Our software environment is Python with NumPy, SciPy for SLSQP and trust-constr, and PyTorch; CPU execution with a single numerical thread.

\paragraph{\rev{SciPy baselines.}} \rev{Every solver comparison in this paper, here and for the UAV problem of \prettyref{appen:control} and the learned-constraint problem of \prettyref{appen:neural}, compares against the two general-purpose constrained solvers of SciPy (the parameter study of \prettyref{appen:ablation} and the timing study of \prettyref{appen:scaling} run IP-ADMM alone), run in float64 on the CPU on the original nonlinear inequalities. SLSQP~\cite{kraft1988} is an active-set sequential quadratic programming method: each iteration solves a dense quadratic subproblem built from the constraint Jacobian and a BFGS approximation of the Lagrangian Hessian. trust-constr~\cite{byrd1999} is a trust-region interior-point method: it adds a logarithmic barrier on slack variables, takes trust-region SQP steps on each barrier subproblem, and decreases the barrier parameter until it meets its tolerance. Both receive the analytic objective gradient and constraint Jacobian; trust-constr also receives the exact Hessians of the objective and of the multiplier-weighted constraints, computed from the same PyTorch model (by automatic differentiation here and in \prettyref{appen:control}, and by differentiating the network twice in \prettyref{appen:neural}). We set \texttt{ftol} of SLSQP and \texttt{gtol}, \texttt{xtol} of trust-constr to $10^{-12}$, leave all other options at their SciPy defaults, and tune neither solver, so each stops on its own convergence test. Neither solver keeps the iterates original-space feasible or produces a native b-FJ certificate, so we score both post hoc: here at their final iterate, with the KKT diagnostic and the exact-geometry checks described next; in \prettyref{appen:control} at their final iterate, by goal gap, exact obstacle clearance, and path length, since the KKT reconstruction is unreliable at the $\ell_1$ kink; and in \prettyref{appen:neural} along their stored iterates, with the same certification target as IP-ADMM.}

All tables and images are generated from canonical stored run records. Completion means that a solver produced a final output; feasibility evaluates the original nonlinear inequalities. Physical validity applies independent exact-geometry checks for capsule overlap and medial-axis crossing or cloth self-crossing and segment degeneracy. Verification is method-specific. Both IP-ADMM variants report their native b-FJ certificate together with split and original-space feasibility. For SLSQP, trust-constr, and both IP-ADMM variants, the w-KKT and KKT values come from a post-hoc same-$x$ KKT diagnostic: constraints with $f_i\leq2\times10^{-4}$ are declared active, nonnegative inequality multipliers are fitted by bounded least squares, and stationarity and complementarity are then evaluated. These reconstructed multipliers are not native solver multipliers. BC-ADMM is scored by the residual of the problem it actually solves, as explained next.

BC-ADMM~\cite{pan2026bcadmm} is the closest published optimizer to ours, so we add it as a baseline. It is also a parallel ADMM that starts from a strictly feasible point and keeps every iterate original-space feasible, but it gets there a different way. It adds a helper variable to each nonconvex constraint and turns that constraint into a biconvex one. For a cloth node and a segment that have to stay apart, the helper variable is a separating plane: the node has to sit on one side of it and the segment on the other. With the plane held fixed, both sides of the constraint are convex, so every block update is a small convex problem, and the plane itself has a closed form once the blocks are fixed. Although BC-ADMM solves these local subproblems iteratively in general (\prettyref{appen:prior}), in this cloth setup each block update has a closed form, so one BC-ADMM iteration costs one effective sweep and its counts in \prettyref{table:all-experiments} are in the same unit as those of IP-ADMM. The springs get the same treatment, with the spring direction as the helper variable, and this reproduces the exact spring energy. The circular container is already convex and needs no helper variable, and gravity is linear, so it stays in the consensus solve. We keep the settings from that paper ($\kappa_y^{\rm BC}=2$, $\kappa^{\rm BC}=2.1$, $\gamma^{\rm BC}=0.95$, $\eta^{\rm BC}=1-10^{-5}$) and tune only the two penalties and the barrier margin, which is what it prescribes.
Unfortunately, BC-ADMM does not apply to capsule stacking. A capsule pose carries an angle, and the world position of a capsule tip is not a convex function of that angle, even with the separating plane held fixed. Therefore, we do not know a BC-ADMM-compatible biconvex reformulation for this capsule formulation, and the BC-ADMM paper lists this reliance on a biconvex relaxation as a limitation of the method. We therefore mark every capsule entry for BC-ADMM with --, meaning not achievable and not computed. BC-ADMM also does not solve our problem. Its barrier keeps a fixed margin that never shrinks, so its answer stops a fixed distance short of the constraint boundary, and its stationarity condition belongs to that relaxed problem rather than to ours. We therefore score it by the residual it actually drives to zero, its own w-KKT residual: the larger of the relaxed stationarity $\|\nabla f+A^{T}\nabla_{y}g\|_{\infty}$ and the consensus gap $\|Ax-y\|_{\infty}$. It has no b-FJ certificate at all, and its barrier takes the place of multipliers, so it has no complementarity residual for our constraints either. Both of those columns are therefore marked -- for BC-ADMM.

\prettyref{table:all-experiments} reports, for every run, the earliest effective sweep at which each certificate first reaches the shared $10^{-3}$ target; an effective sweep is one value-and-gradient evaluation of every block, the oracle work of one ADMM iteration, so the solvers are compared in one cost unit. The three ADMM methods have no built-in termination test, so each runs under a wall-clock budget, and we record the first sweep at which each certificate is established; N/A means it was not established within this budget. SLSQP and trust-constr instead stop on their own convergence tests (\prettyref{table:experiment-parameters}). All $27$ runs complete, and every reported output passes the original-space feasibility check within the numerical tolerances of \prettyref{table:experiment-parameters}: the $24$ runs of SLSQP, trust-constr, and the two IP-ADMM variants, plus the $3$ BC-ADMM cloth runs. \rev{Every output also passes the physical-validity checks except the trust-constr capsule outputs at $n=8$ and $n=24$, marked $\ddagger$. They satisfy every endpoint--segment constraint, yet one and three capsule pairs, respectively, end with crossing medial axes. A continuous motion cannot reach such a crossing without an endpoint first touching the other capsule, but a long interior-point step can jump past that contact. No other run ends in such a configuration.} For the two smallest instances ($n=8$) both IP-ADMM variants establish the w-KKT, b-FJ, and KKT certificates; for the larger instances they still establish w-KKT and b-FJ, but KKT is not reached within the budget. The two variants trade the certificates off: IP-ADMM-Ckpt establishes b-FJ earlier on every capsule instance and on cloth-8 ($10$k--$60$k versus $40$k--$100$k sweeps), equally early on cloth-16, and later on cloth-24, whereas it never establishes w-KKT earlier than ordinary IP-ADMM. \rev{Even the SciPy baselines do not robustly certify KKT on capsule stacking: both miss it at size $16$, SLSQP also misses it at size $24$, and trust-constr reaches it at sizes $8$ and $24$ only at configurations with crossing medial axes. On cloth resting both reach KKT at every size, and trust-constr ends at a lower objective than every other method, that is, at a different local minimum.} This is consistent with CQ degeneracy or difficulty of KKT certification, but our experiments do not distinguish these explanations. It suggests that requiring KKT is a fragile stopping criterion in practice. \rev{SLSQP and trust-constr have no b-FJ certificate, so we score them only at their final iterate}.  BC-ADMM is scored against its own relaxed target, which it reaches only on the smallest cloth instance; at the larger sizes its residual settles above the target and the run stops without a certificate. Its b-FJ and KKT columns are marked -- for the structural reasons given above. This is the trade we want to show. BC-ADMM assumes the nonconvex constraints are biconvex and exploits that structure, which we never assume, and it only ever has to certify a relaxed problem, one whose constraints are all strictly inactive by construction. The convergence of BC-ADMM is not cheap: on cloth-8 it needs $178$k iterations to bring its relaxed residual to $10^{-3}$, while ordinary IP-ADMM certifies b-FJ, a statement about the original constraints, at $60$k, and IP-ADMM-Ckpt at $10$k. BC-ADMM stalls because it refreshes the separating planes only when the merit strictly decreases. Once the iterate reaches the optimum for the current planes, the merit stops decreasing, so the planes are frozen for good and the residual that also asks them to be optimal stays far above the target. In contrast, our algorithm assumes less and certifies more: our b-FJ certificate is a statement about the original constraints, and \prettyref{thm:rate} bounds how many iterations it takes to reach a given level.

\prettyref{fig:experiment-trajectories} shows ordinary IP-ADMM snapshots from the start to the final output. The configuration settles early (\prettyref{fig:visual-difference}): the distance to the final output falls below $5\%$ of its initial value early on. The later snapshots therefore differ only slightly from the final output, most visibly for cloth at $n\in\{16,24\}$.  \prettyref{fig:certificates} plots the w-KKT, b-FJ, and KKT certificates against the total solver iteration for every instance and mode, with the shared $10^{-3}$ target marked; at the larger sizes the KKT curve stalls above the target while the other two drop below it.  \prettyref{fig:visual-difference} plots how far each iterate still is from the final output, measured by the $L_2$ distance between the current and final solution vectors and scaled so the start is $1$.  This distance drops below $5\%$ early and then flattens near zero, long before KKT is reached (among the plotted IP-ADMM methods, KKT is reached only at $n=8$).  The configuration already looks like the final output well before stationarity is certified, so an early stop is a reasonable heuristic in practice.  The catch is that this distance can only be measured after the run, because it needs the final solution.  On the same axes we therefore also plot a quantity available during the run: the relative step $r_{\rm rel}^{k}=\|x^{k}-x^{k-\Delta}\|/\max\{1,\|x^{k}\|\}$, where $\Delta$ is the gap between stored samples.  It uses only the current and previous iterate.  The plot shows that $r_{\rm rel}^{k}$ falls fastest while the distance to the final output drops, and then keeps decaying steadily on the log scale after the distance has flattened near zero.  It therefore does not mark the settling moment sharply, but a threshold on $r_{\rm rel}^{k}$ still gives a practical online stopping rule, as in earlier physics solvers~\cite{narain2016admm,overby2017admm,pan2026bcadmm}.  In short, a full KKT certificate can take much longer than the moment the motion has effectively stopped, and an online threshold on $r_{\rm rel}^{k}$ can approximate that moment.

\begin{table*}[t]
\centering\small
\caption{Empirical nondimensional benchmark and solver settings used in the canonical runs; these numerical choices are not claimed to satisfy the theorem-sufficient parameter construction.  Values written $a/b/c$ are the per-size settings for $n=8/16/24$.}
\label{table:experiment-parameters}
\resizebox{\textwidth}{!}{%
\begin{tabular}{@{}lll@{}}
\toprule
Group & Parameter & Value \\
\midrule
Capsule model & $n$ & 8, 16, 24 \\
 & length, diameter, container radius & 1.05, 0.2, 2 \\
 & constraint smoothing $\delta$ & 0.25 \\
Cloth model & $n$ & 8, 16, 24 \\
 & mass, gravity, container radius & 1, 1, 2 \\
 & stretch rest value, stiffness & 0.32, 100 \\
 & bending stiffness, spring $\epsilon$ & 4, 0.01 \\
 & clearance, constraint smoothing $\delta$ & 0.055, 0.25 \\
SLSQP & \texttt{ftol} & $10^{-12}$ \\
trust-constr & \texttt{gtol}, \texttt{xtol} & $10^{-12}$, $10^{-12}$ \\
Ordinary capsule & $\rho^0,L^0,G,\eta,\sigma,\xi,p$ & 80/80/160, $10^{6}$, 0.005, $3\times10^{5}$, 1, 0.1/0.04/0.04, 0.5 \\
 & $\theta_1,\theta_2$ & 1.01, 0.01 \\
Ordinary cloth & $\rho^0,L^0,G,\eta,\sigma,\xi,p$ & 40/40/160, $10^{6}$, 0.006, 15.5, 1, 0.1, 0.6/0.5/0.5 \\
 & $\theta_1,\theta_2$ & 1.01, 0.1 \\
Ckpt capsule & $\rho_{\rm cp}^0,G^{\rm cp},\eta^{\rm cp},\chi,\tau_{\rm cp},\ell_{\rm cp}$ & 320, 0.005, 1, 1.000003, 1000, 9.765625 \\
 & $\theta_1^{\rm cp},\theta_2^{\rm cp}$ & 1.01, 0.01 \\
Ckpt cloth & $\rho_{\rm cp}^0,G^{\rm cp},\eta^{\rm cp},\chi,\tau_{\rm cp},\ell_{\rm cp}$ & 40/40/160, 0.006, 1, 1.000003, 2, 625/625/156.25 \\
 & $\theta_1^{\rm cp},\theta_2^{\rm cp}$ & 1.01, 0.1 \\
BC-ADMM cloth & $\beta^{\rm BC},\beta_y^{\rm BC},\gamma^{\rm BC},\eta^{\rm BC},\kappa^{\rm BC},\kappa_y^{\rm BC}$ & 100, 100, 0.95, 0.99999, 2.1, 2 \\
 & margin $\varepsilon^{\rm BC}$, $\sigma^{\rm BC}$ & $10^{-3}$, $10^{-6}$ \\
Verification & active, violation, b-FJ, Stat., Comp. & $2\times10^{-4}$, $10^{-7}$, $10^{-3}$, $10^{-3}$, $10^{-3}$ \\
Physical gate & capsule, cloth tolerance & $10^{-10}$, $10^{-9}$ \\
Runtime & precision, device, CPU threads & float64, CPU, 1 \\
\bottomrule
\end{tabular}}
\end{table*}

\begin{table*}[t]
\centering\scriptsize
\caption{Counts are effective sweeps (k: thousands): one sweep evaluates the value and gradient of every block once, as one ADMM iteration does.  SciPy calls are converted by the number of blocks they evaluate, and an exact trust-constr Hessian counts as $\dim x$ calls; SciPy linear algebra and the per-sweep IP-ADMM certificate check are not counted.  w-KKT, b-FJ, and KKT give the first sweep at which each certificate reaches $10^{-3}$ (N/A: not reached; --: not achievable by the method, so not computed).  SLSQP and trust-constr have no b-FJ certificate and are scored at their final iterate.  Viol. is the maximum inequality violation.  $\dagger$: residual of the BC-ADMM relaxation, not comparable with the other rows; BC-ADMM~\cite{pan2026bcadmm} does not apply to capsules and has no b-FJ or KKT certificate for the original problem.  $^{\ddagger}$: feasible but with crossing medial axes.}
\label{table:all-experiments}
\setlength{\tabcolsep}{2pt}
\begin{minipage}[t]{0.495\textwidth}\centering\tiny
\textbf{Capsule stacking}\\[2pt]
\begin{tabular}{@{}clrrrrr@{}}
\toprule
$n$ & Method & Objective & Viol. & w-KKT & b-FJ & KKT \\
\midrule
\multirow{5}{*}{8} & SLSQP & -11.90 & $5.55\!\times\!10^{-16}$ & 46 & -- & 46 \\
 & trust-constr$^{\ddagger}$ & -11.71 & 0 & 2.2k & -- & 2.2k \\
 & IP-ADMM & -11.92 & 0 & 139k & 40k & 1329k \\
 & IP-ADMM-Ckpt & -11.92 & 0 & 220k & 10k & 995k \\
 & BC-ADMM & -- & -- & -- & -- & -- \\
\midrule
\multirow{5}{*}{16} & SLSQP & -19.20 & $8.88\!\times\!10^{-16}$ & 86 & -- & N/A \\
 & trust-constr & -19.04 & 0 & 3.9k & -- & N/A \\
 & IP-ADMM & -18.74 & 0 & 379k & 60k & N/A \\
 & IP-ADMM-Ckpt & -18.74 & 0 & 680k & 20k & N/A \\
 & BC-ADMM & -- & -- & -- & -- & -- \\
\midrule
\multirow{5}{*}{24} & SLSQP & -22.45 & $2.90\!\times\!10^{-15}$ & 84 & -- & N/A \\
 & trust-constr$^{\ddagger}$ & -22.63 & 0 & 6.6k & -- & 6.6k \\
 & IP-ADMM & -22.40 & 0 & 800k & 100k & N/A \\
 & IP-ADMM-Ckpt & -22.40 & 0 & 1080k & 60k & N/A \\
 & BC-ADMM & -- & -- & -- & -- & -- \\
\bottomrule
\end{tabular}
\end{minipage}\hfill
\begin{minipage}[t]{0.495\textwidth}\centering\tiny
\textbf{Cloth resting}\\[2pt]
\begin{tabular}{@{}clrrrrr@{}}
\toprule
$n$ & Method & Objective & Viol. & w-KKT & b-FJ & KKT \\
\midrule
\multirow{5}{*}{8} & SLSQP & -14.25 & $1.22\!\times\!10^{-12}$ & 183 & -- & 183 \\
 & trust-constr & -14.91 & 0 & 1.2k & -- & 1.2k \\
 & IP-ADMM & -14.38 & 0 & 160k & 60k & 261k \\
 & IP-ADMM-Ckpt & -14.38 & 0 & 160k & 10k & 841k \\
 & BC-ADMM & -14.35 & 0 & 178k$^{\dagger}$ & -- & -- \\
\midrule
\multirow{5}{*}{16} & SLSQP & -27.08 & $6.22\!\times\!10^{-15}$ & 372 & -- & 372 \\
 & trust-constr & -28.65 & 0 & 2.2k & -- & 2.2k \\
 & IP-ADMM & -27.87 & 0 & 560k & 320k & N/A \\
 & IP-ADMM-Ckpt & -27.87 & 0 & 595k & 320k & N/A \\
 & BC-ADMM & -28.09 & 0 & N/A & -- & -- \\
\midrule
\multirow{5}{*}{24} & SLSQP & -39.77 & $7.44\!\times\!10^{-15}$ & 403 & -- & 403 \\
 & trust-constr & -41.36 & 0 & 11k & -- & 11k \\
 & IP-ADMM & -39.93 & 0 & 960k & 280k & N/A \\
 & IP-ADMM-Ckpt & -39.93 & 0 & 1720k & 440k & N/A \\
 & BC-ADMM & -38.80 & 0 & N/A & -- & -- \\
\bottomrule
\end{tabular}
\end{minipage}
\end{table*}

\clearpage
\begin{figure*}[p]
\centering
\setlength{\tabcolsep}{1pt}
\begin{tabular}{@{}ccccc@{}}
\toprule
Initial & 25\% & 50\% & 75\% & Final \\
\midrule
\addlinespace[2pt]
\multicolumn{5}{c}{\textbf{Capsule 8}} \\
\includegraphics[width=0.194\textwidth,trim=34 53 34 53,clip]{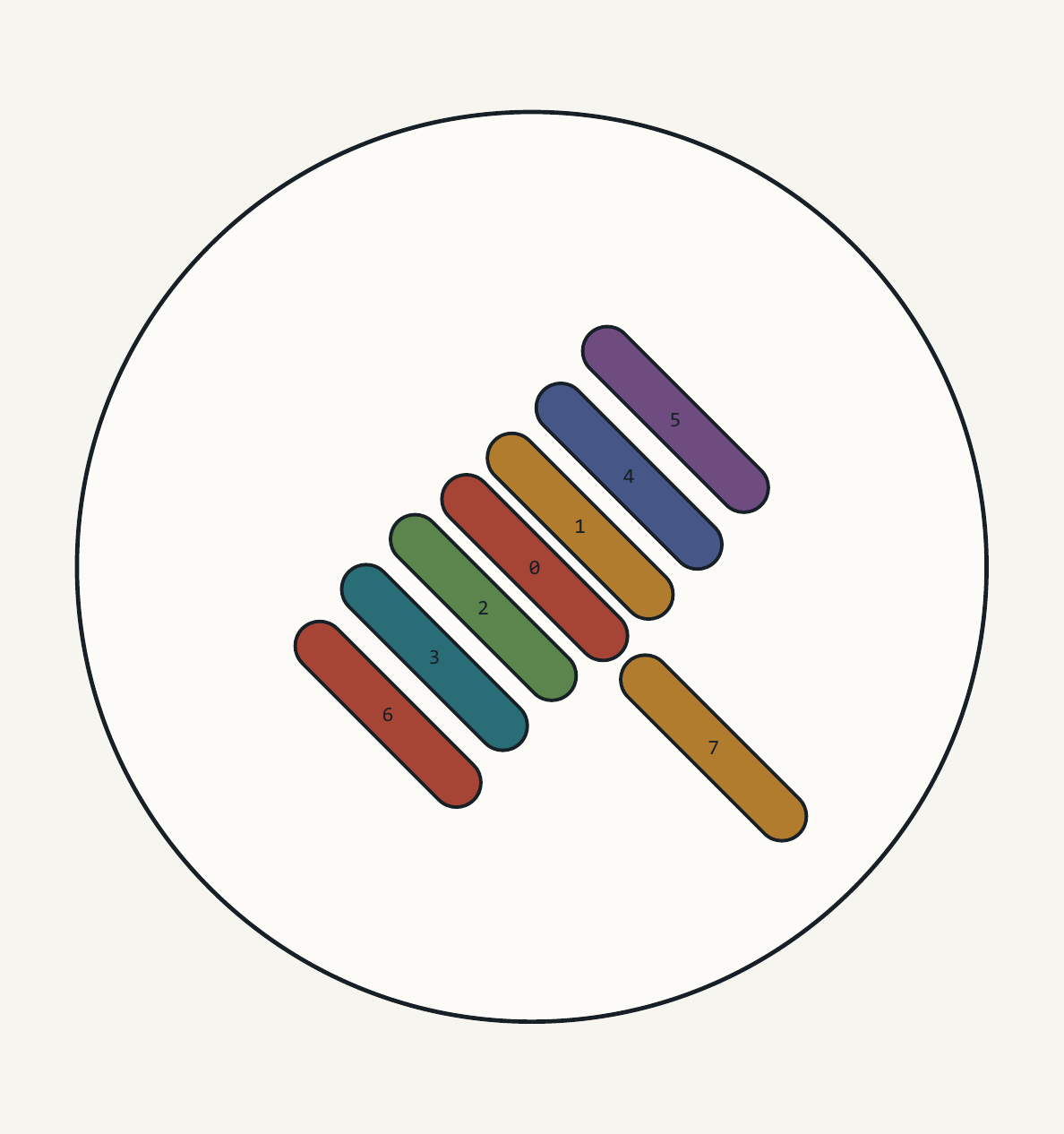} & \includegraphics[width=0.194\textwidth,trim=34 53 34 53,clip]{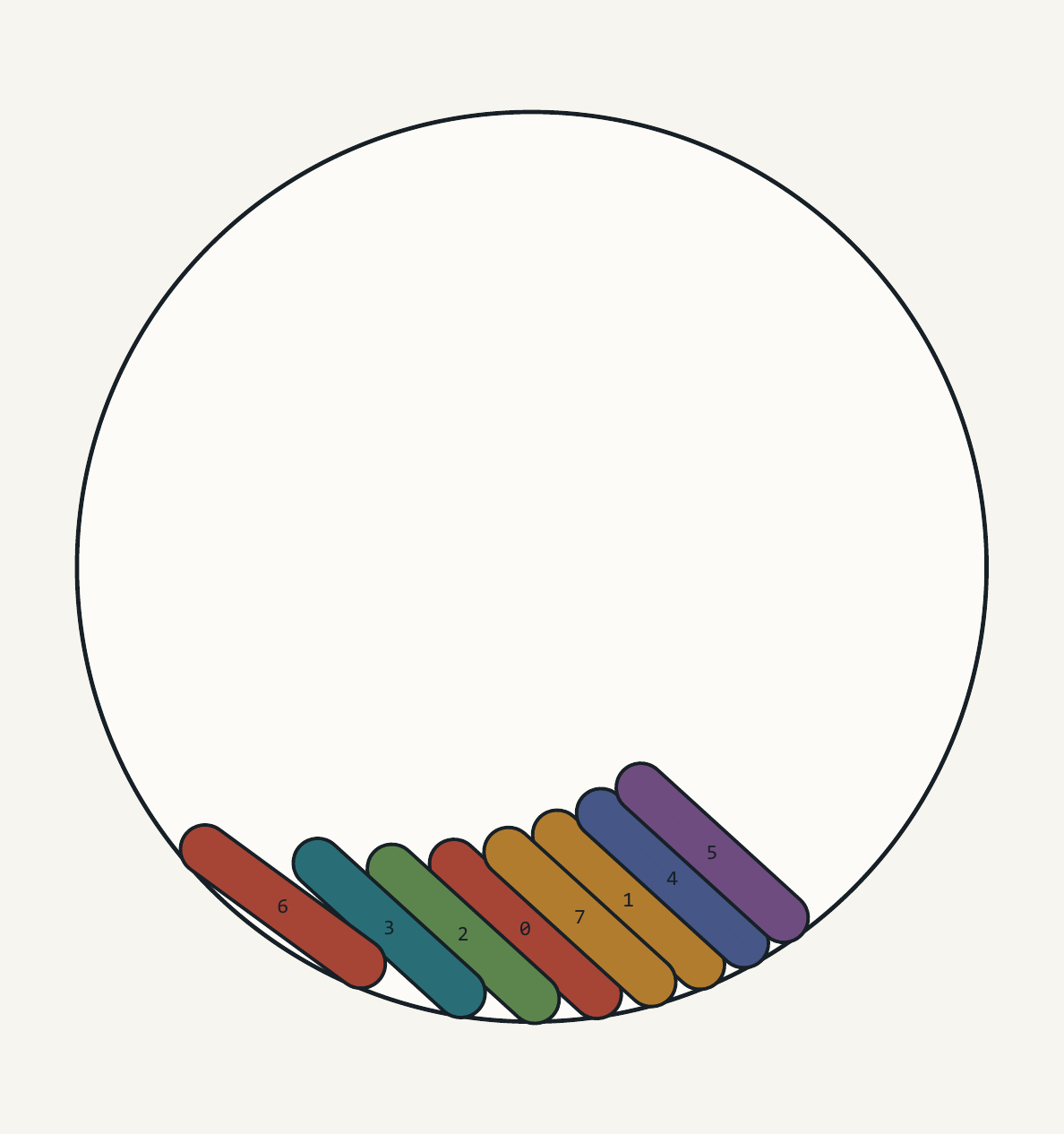} & \includegraphics[width=0.194\textwidth,trim=34 53 34 53,clip]{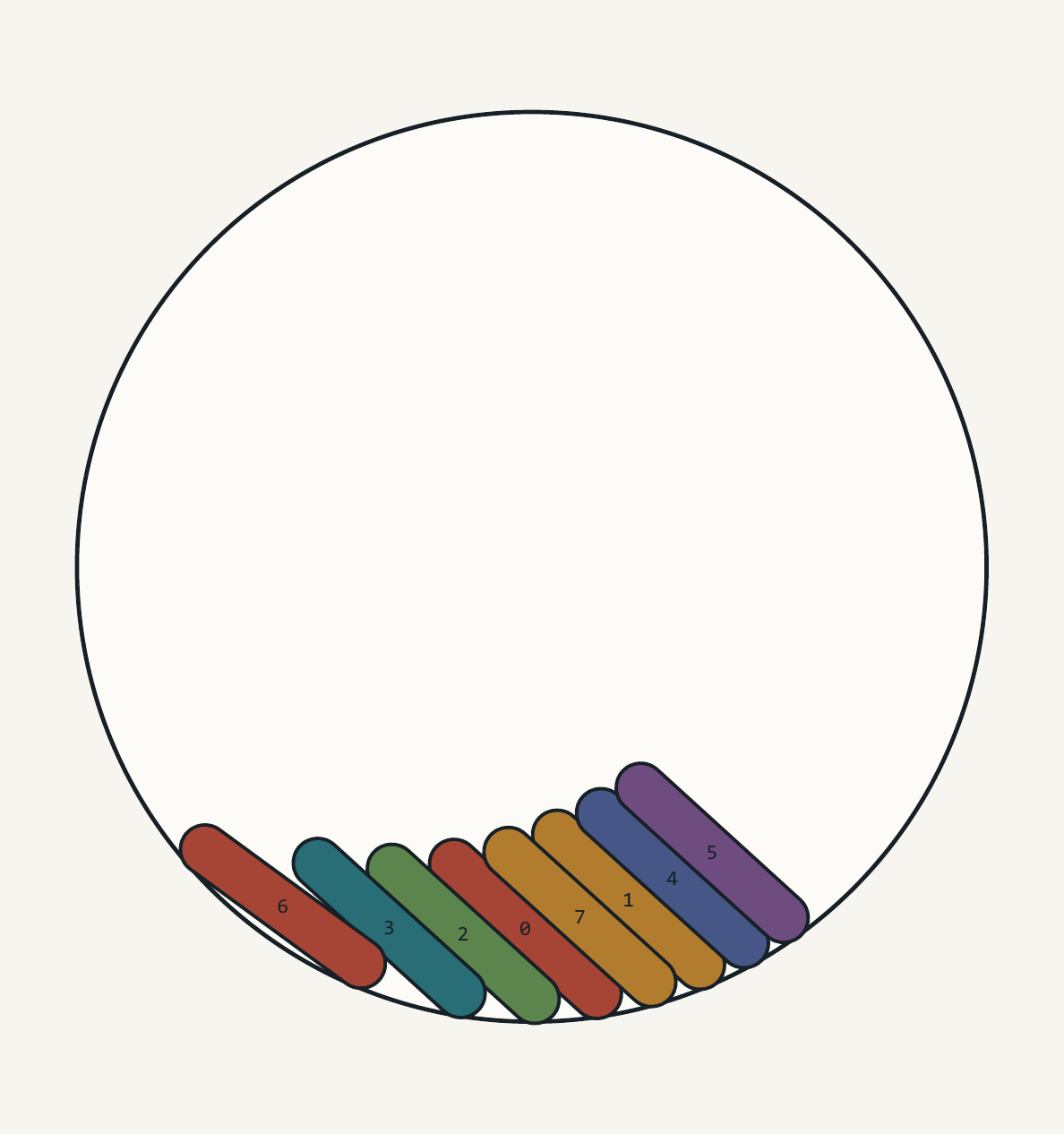} & \includegraphics[width=0.194\textwidth,trim=34 53 34 53,clip]{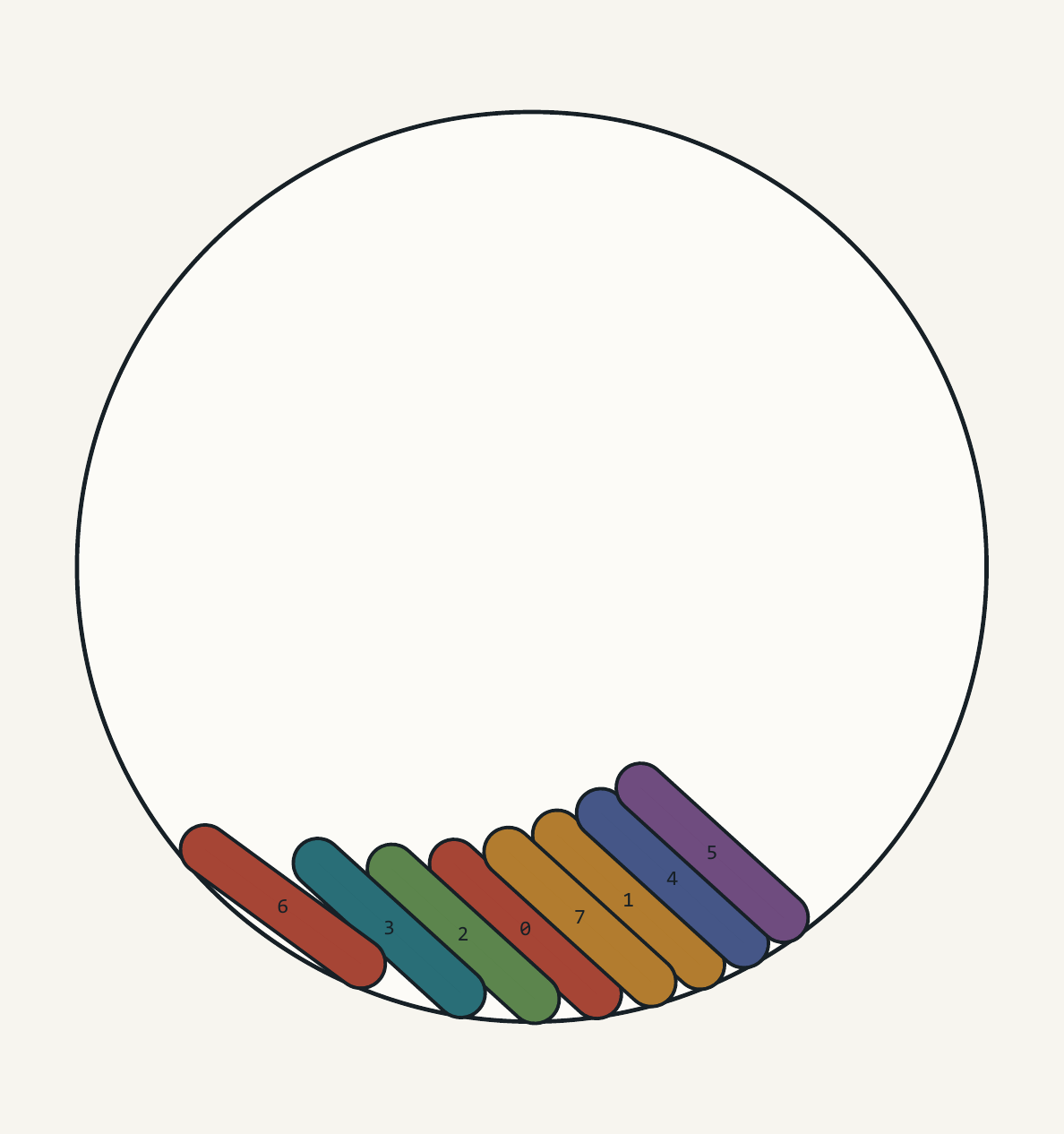} & \includegraphics[width=0.194\textwidth,trim=34 53 34 53,clip]{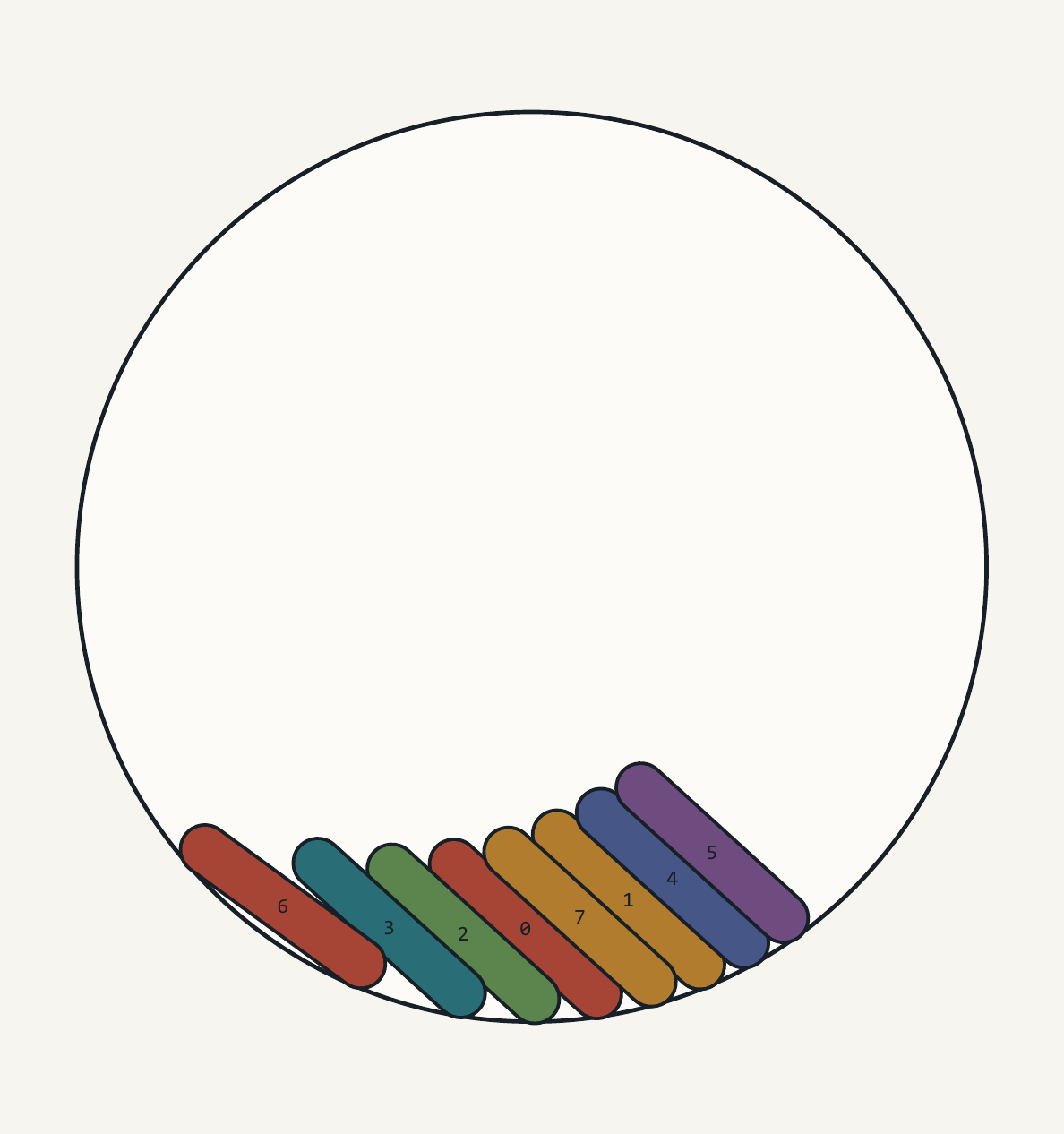} \\
\midrule
\addlinespace[2pt]
\multicolumn{5}{c}{\textbf{Capsule 16}} \\
\includegraphics[width=0.194\textwidth,trim=34 53 34 53,clip]{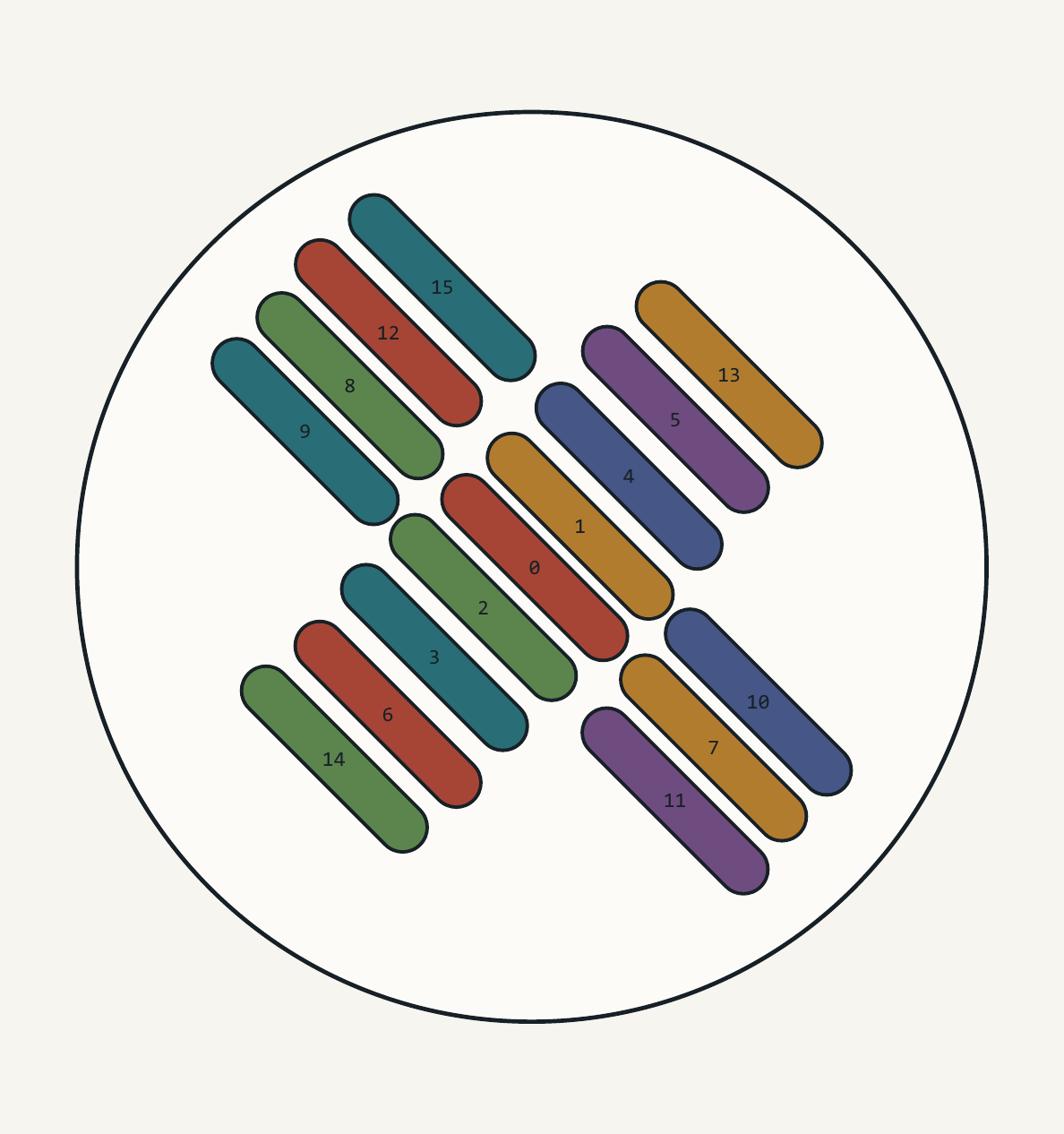} & \includegraphics[width=0.194\textwidth,trim=34 53 34 53,clip]{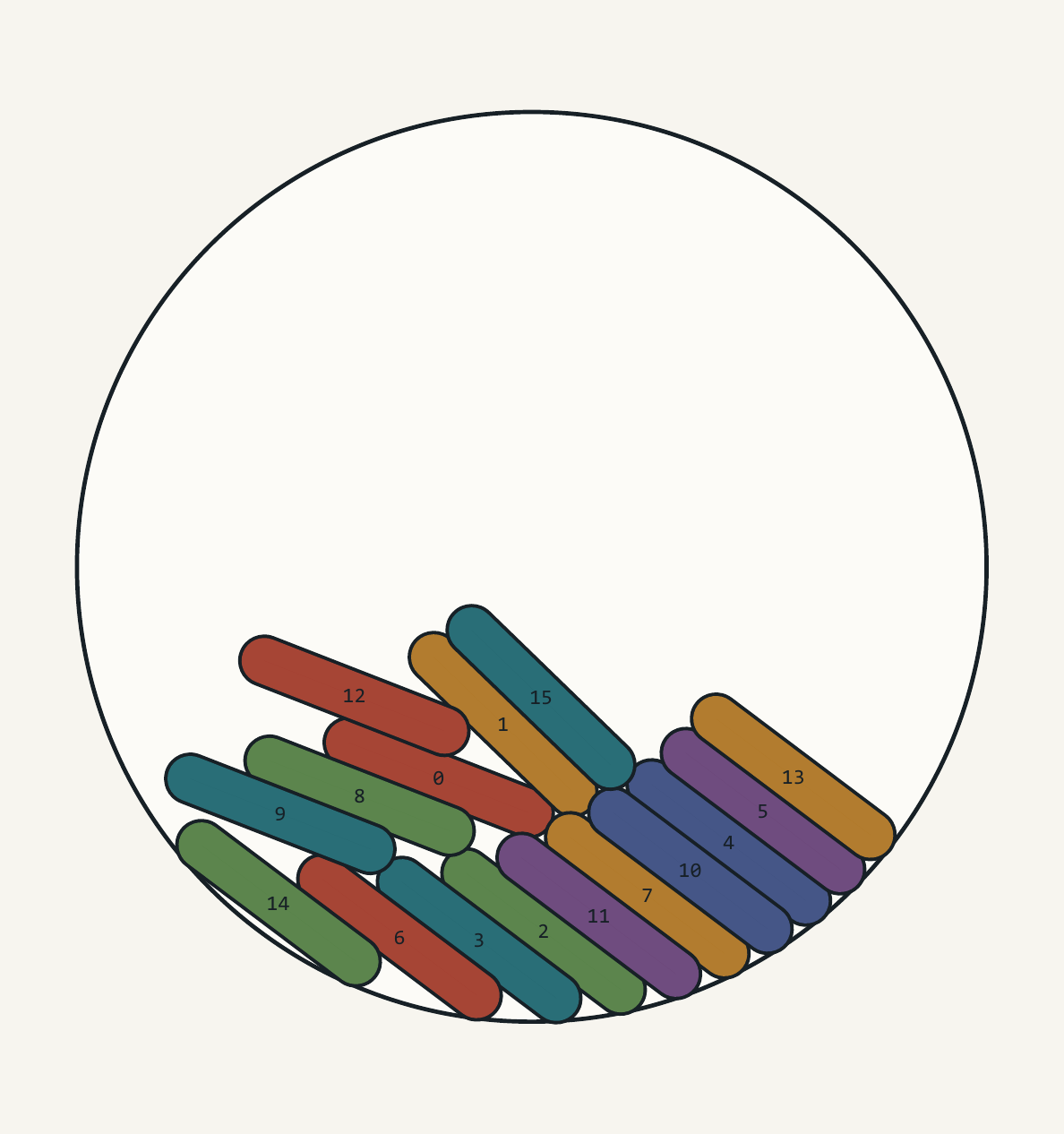} & \includegraphics[width=0.194\textwidth,trim=34 53 34 53,clip]{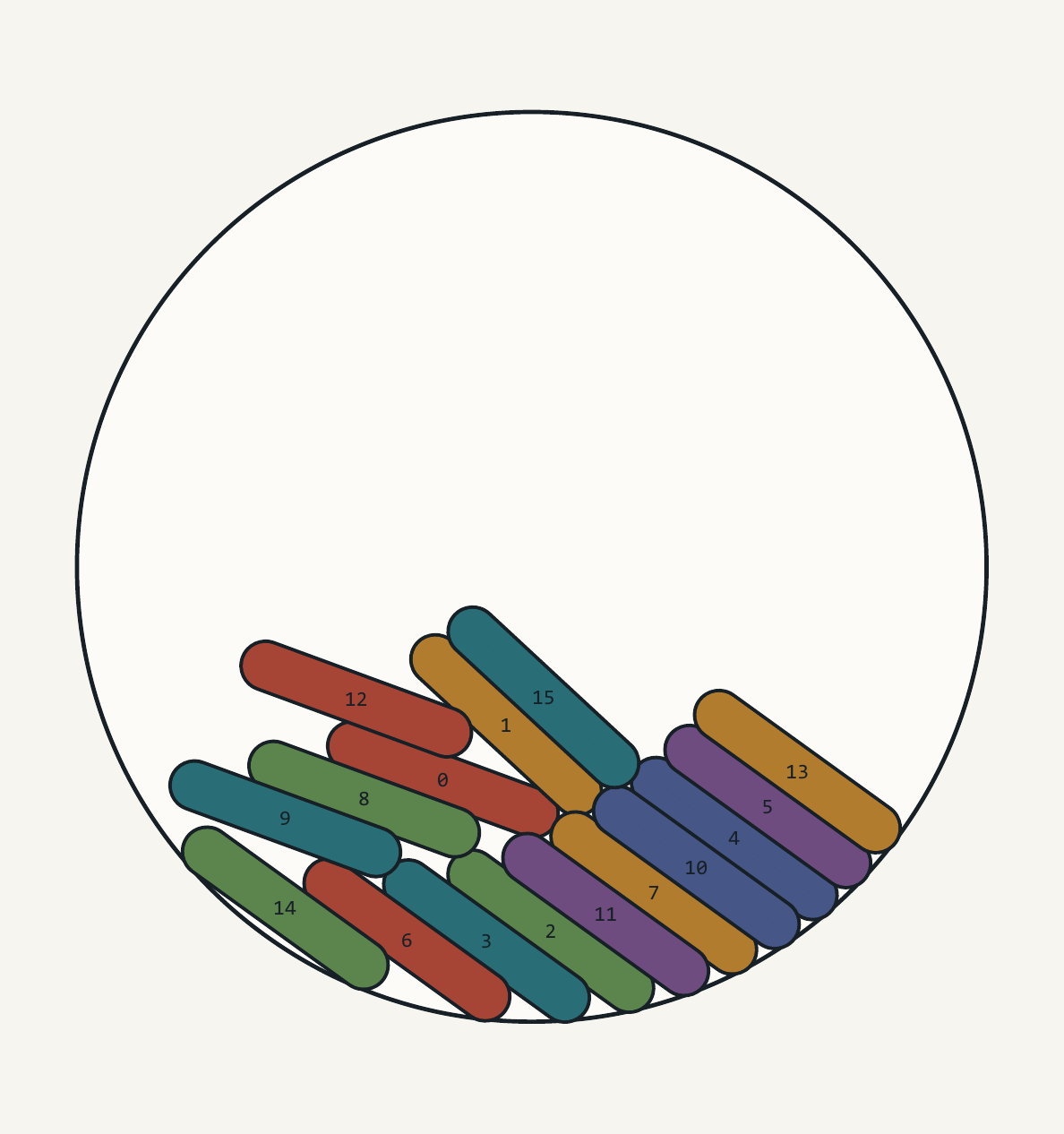} & \includegraphics[width=0.194\textwidth,trim=34 53 34 53,clip]{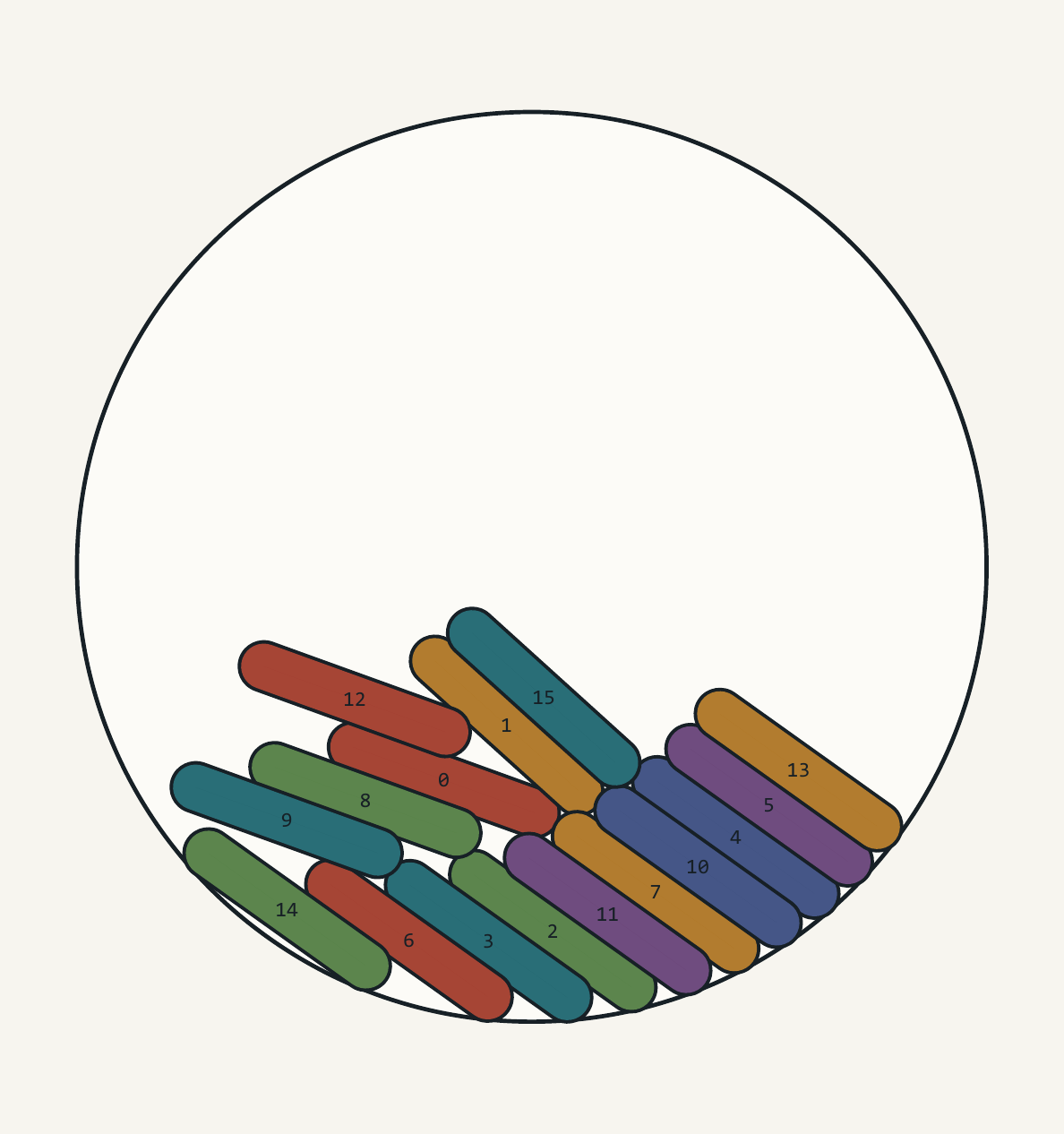} & \includegraphics[width=0.194\textwidth,trim=34 53 34 53,clip]{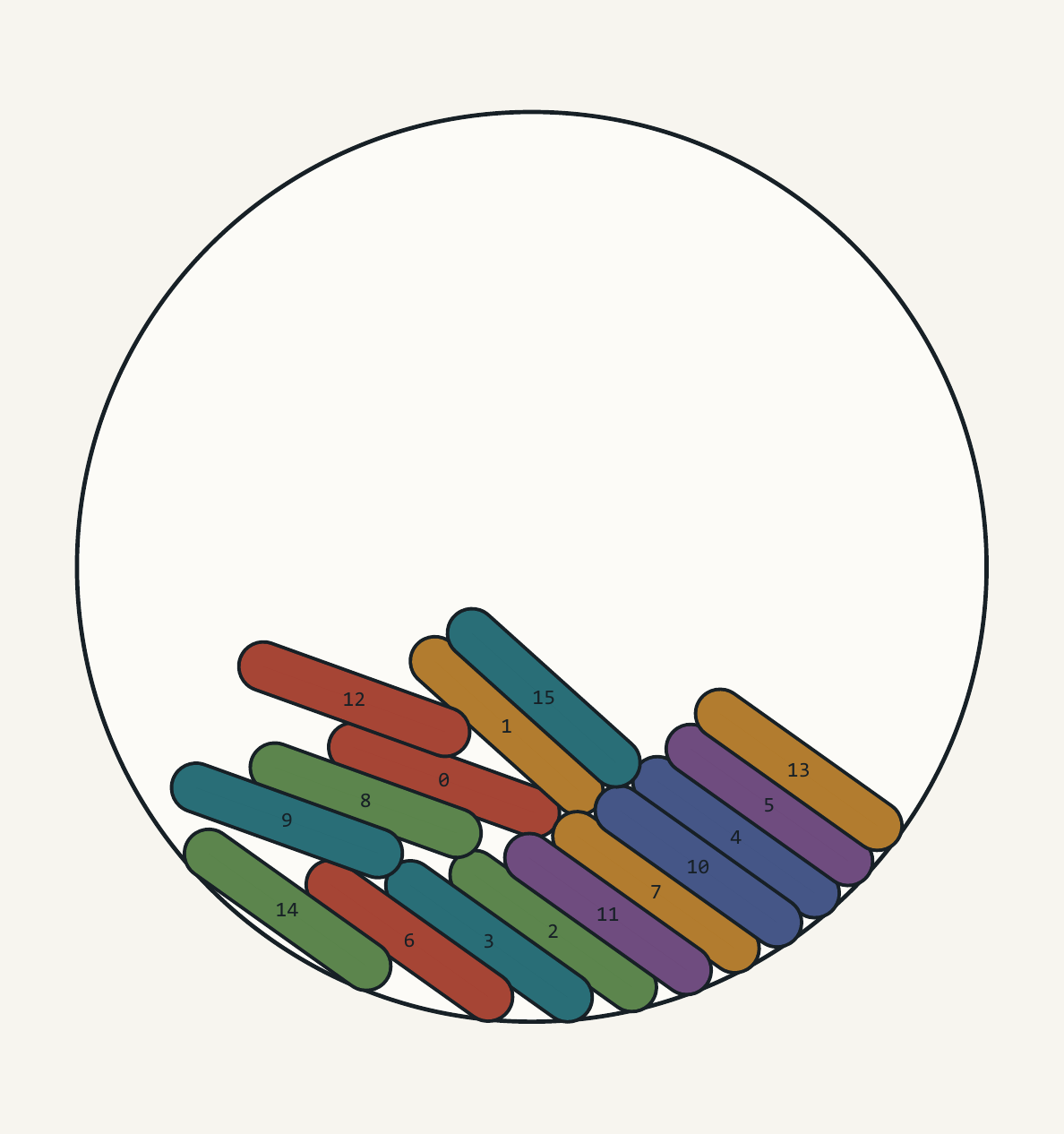} \\
\midrule
\addlinespace[2pt]
\multicolumn{5}{c}{\textbf{Capsule 24}} \\
\includegraphics[width=0.194\textwidth,trim=34 53 34 53,clip]{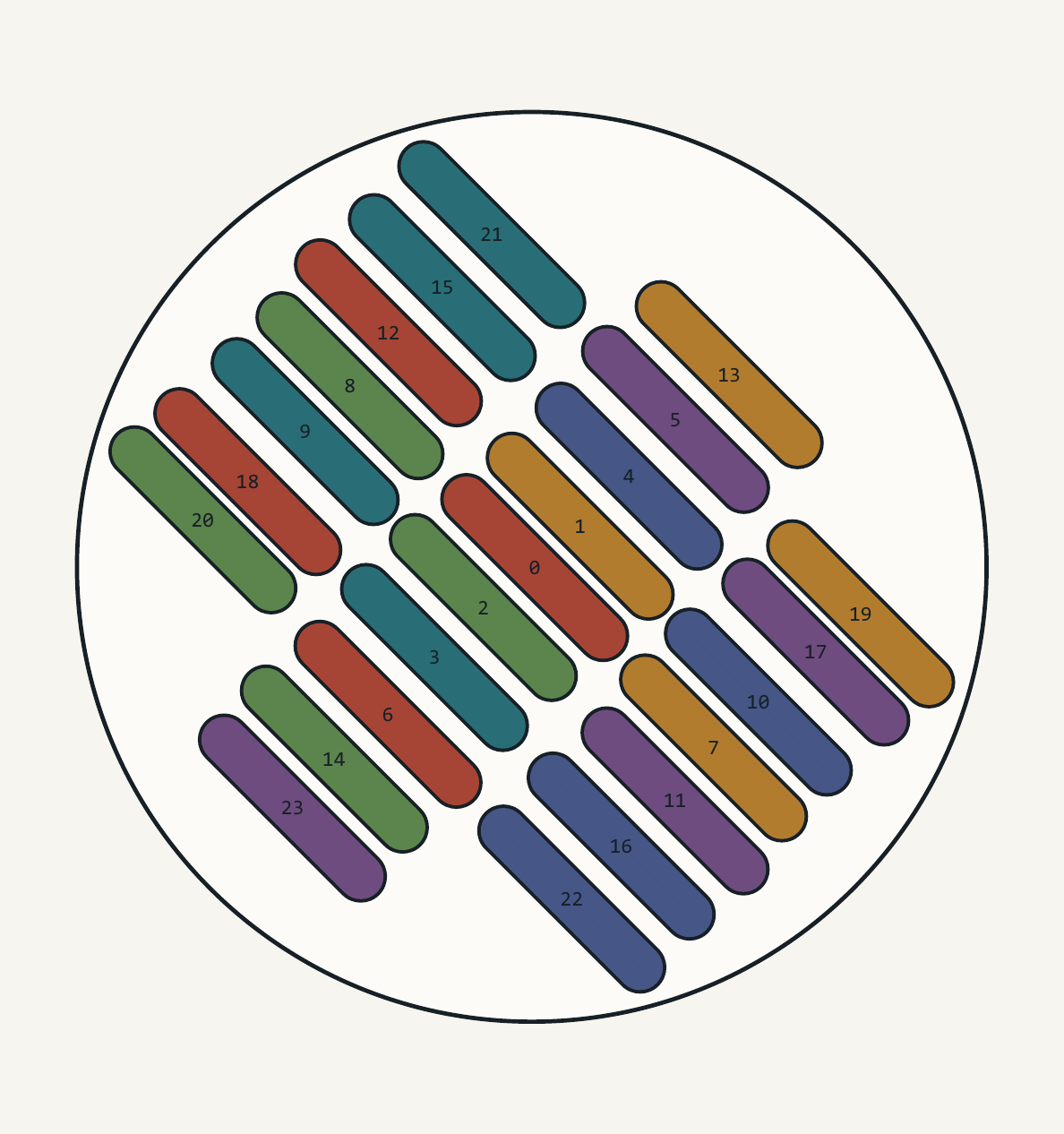} & \includegraphics[width=0.194\textwidth,trim=34 53 34 53,clip]{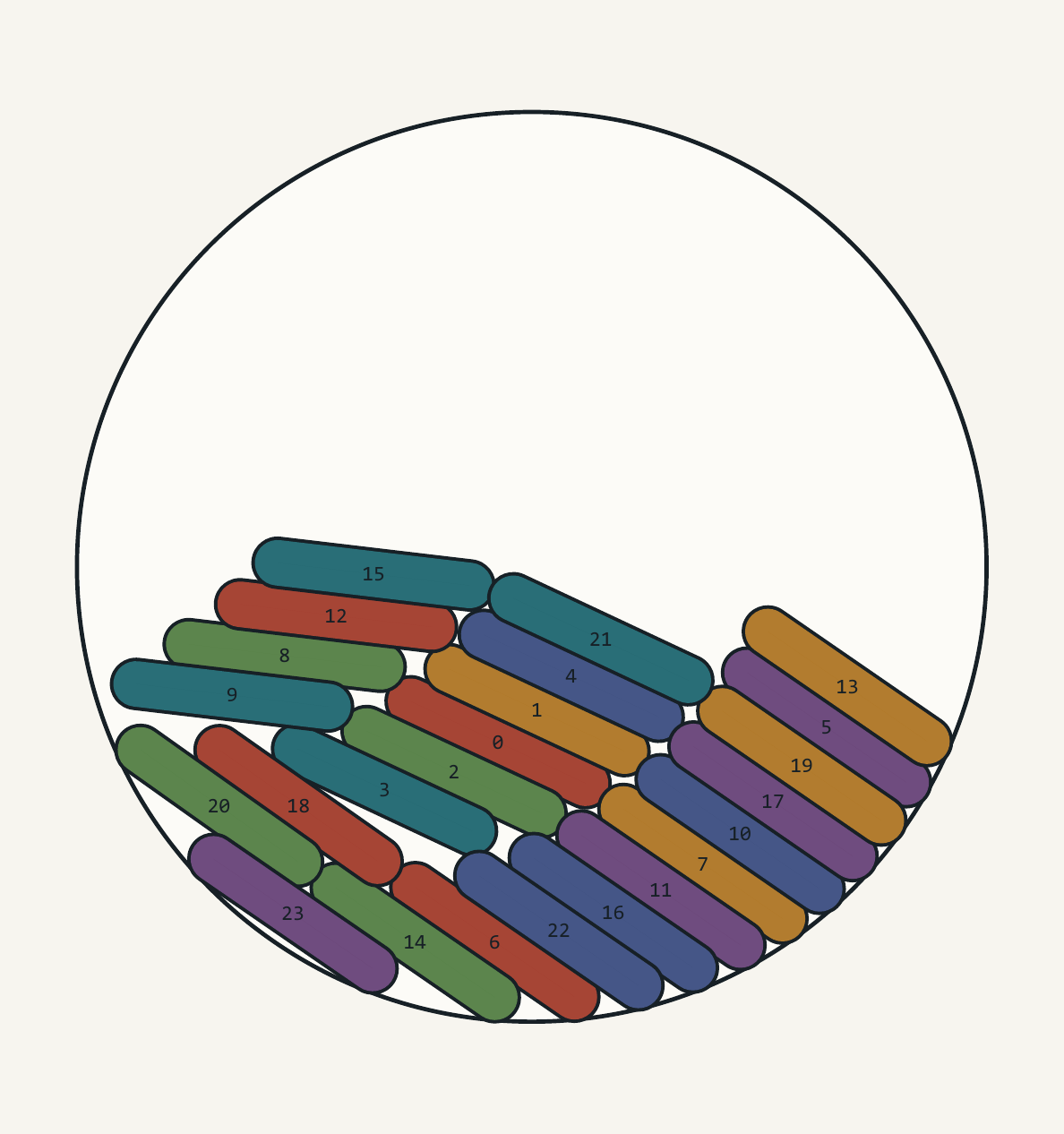} & \includegraphics[width=0.194\textwidth,trim=34 53 34 53,clip]{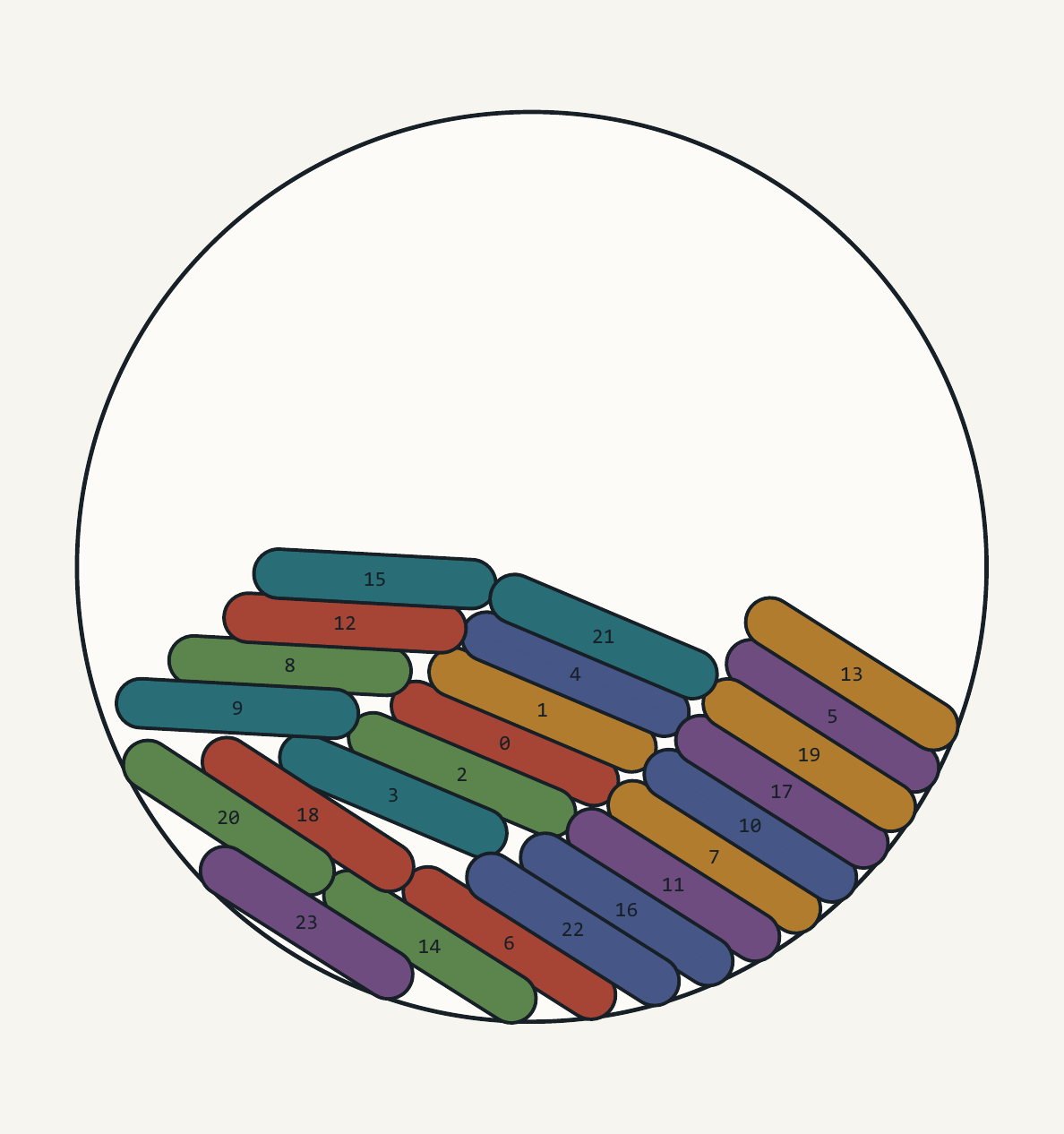} & \includegraphics[width=0.194\textwidth,trim=34 53 34 53,clip]{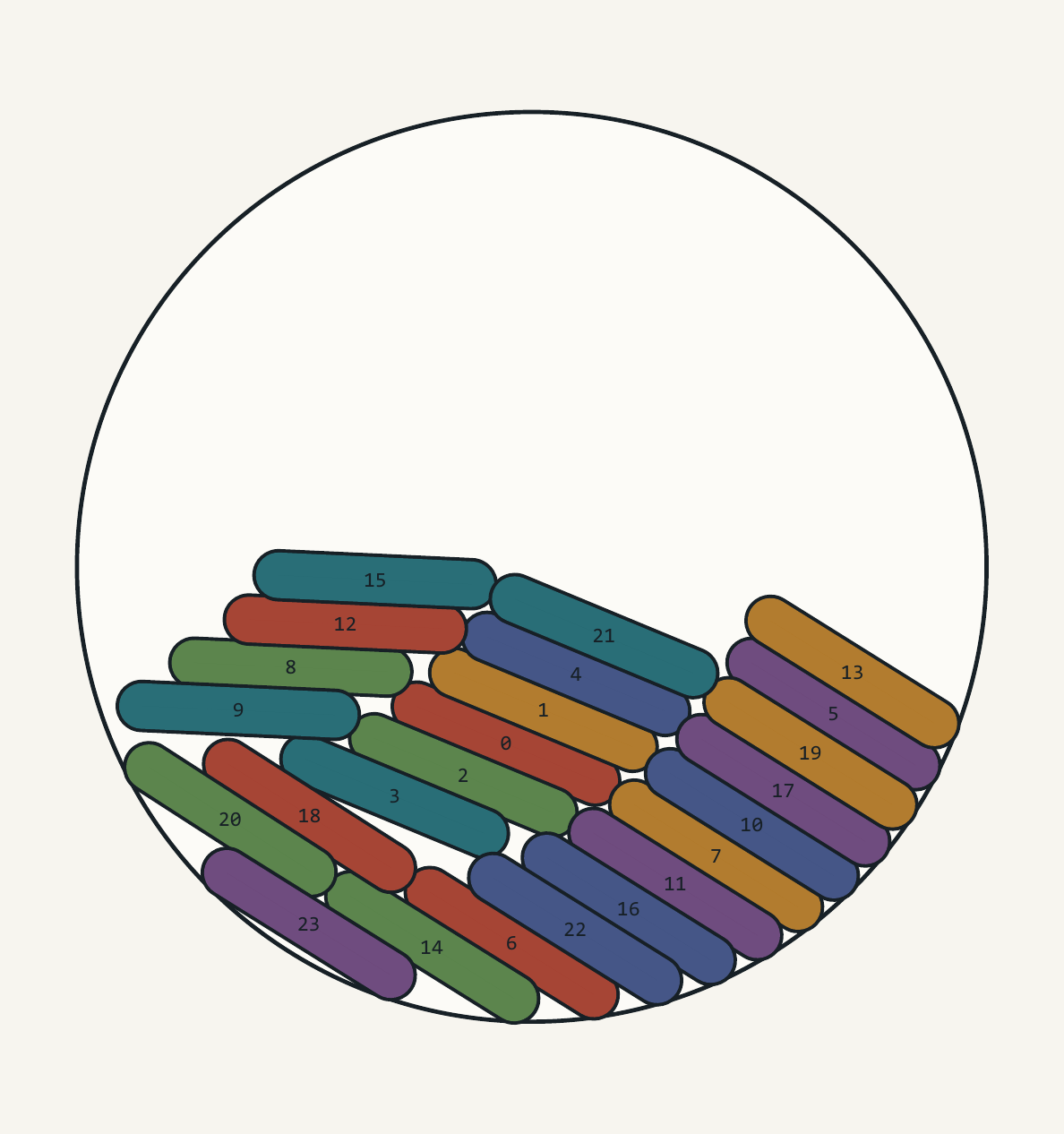} & \includegraphics[width=0.194\textwidth,trim=34 53 34 53,clip]{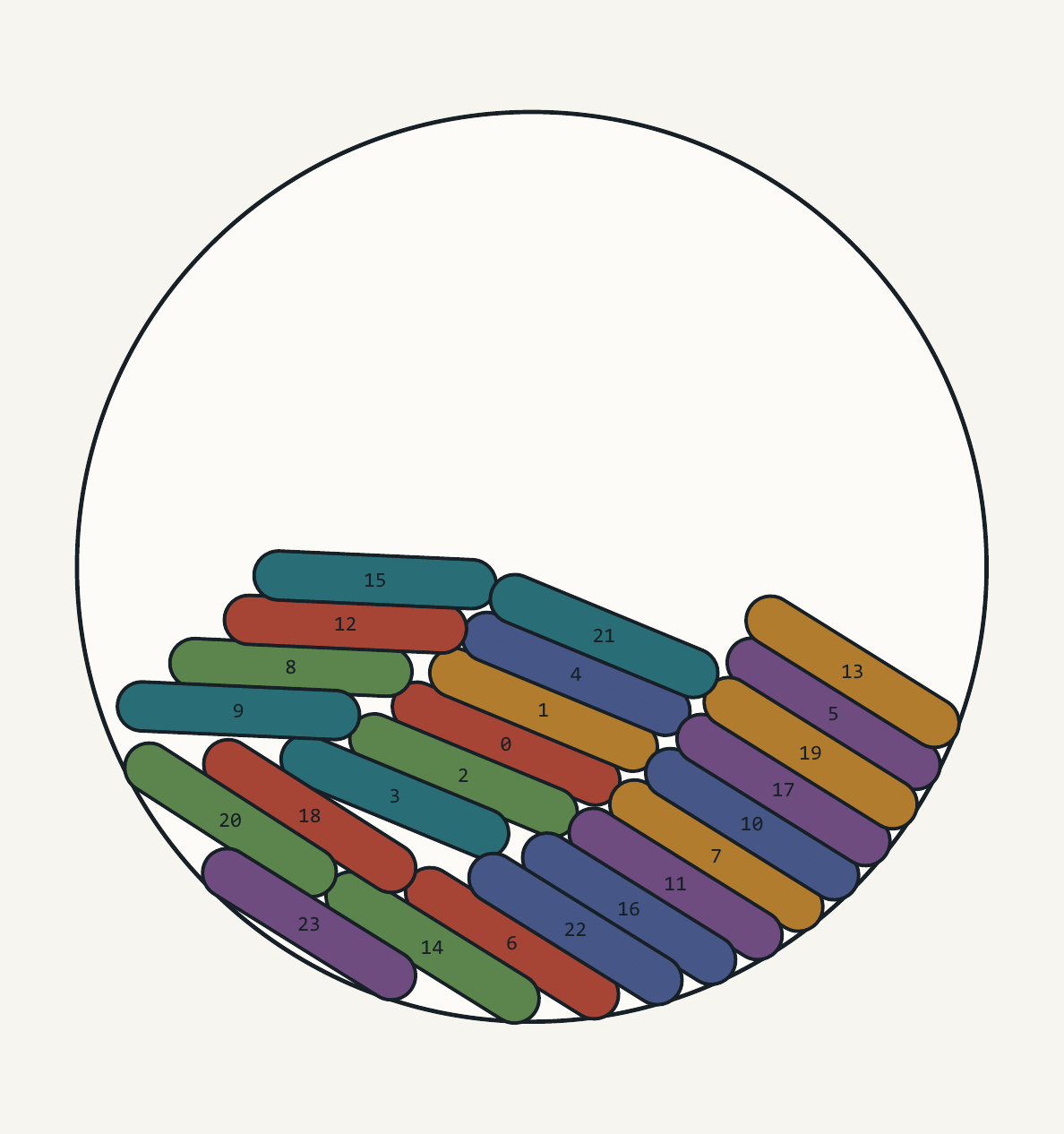} \\
\midrule
\addlinespace[2pt]
\multicolumn{5}{c}{\textbf{Cloth 8}} \\
\includegraphics[width=0.194\textwidth,trim=34 53 34 53,clip]{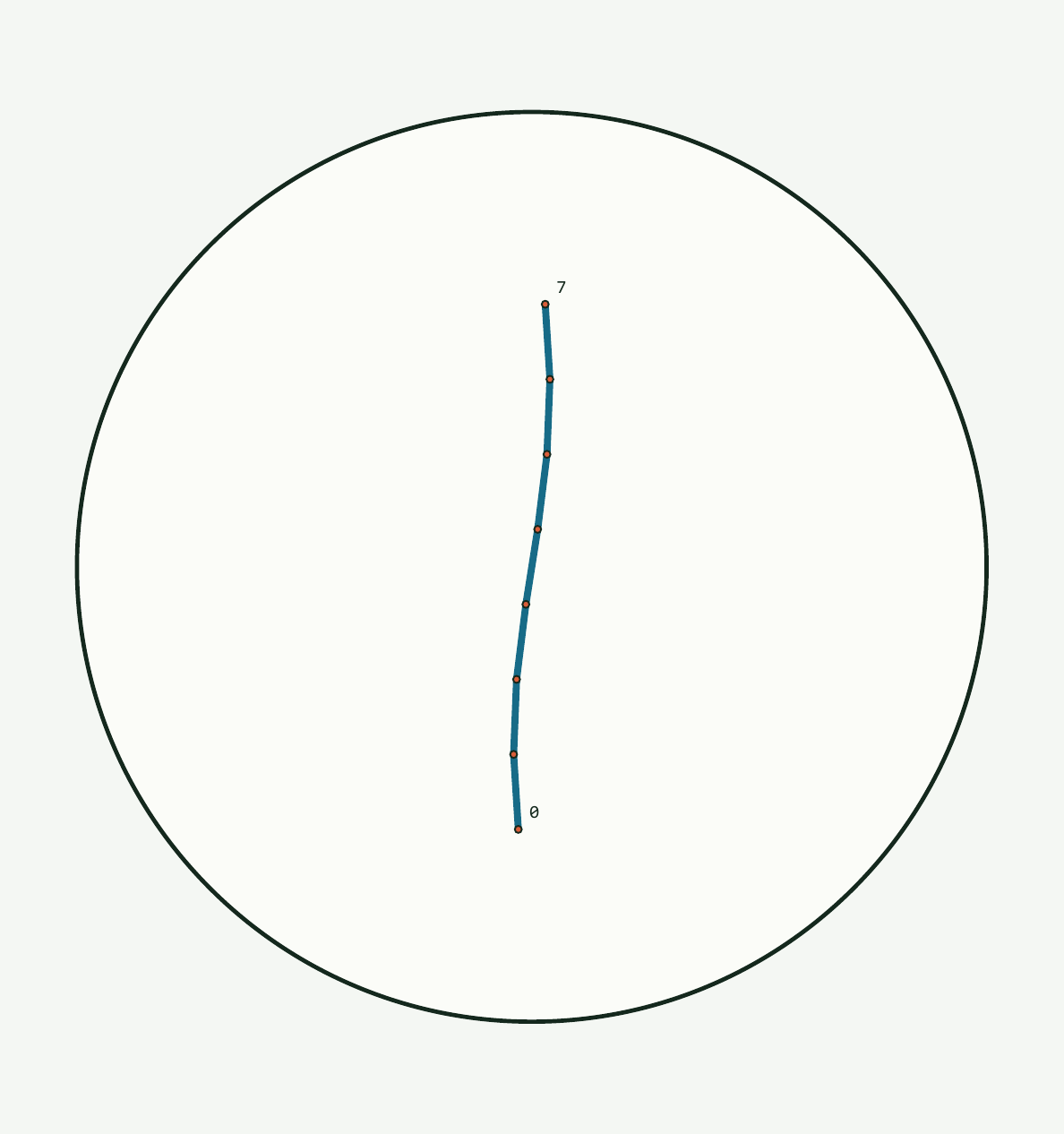} & \includegraphics[width=0.194\textwidth,trim=34 53 34 53,clip]{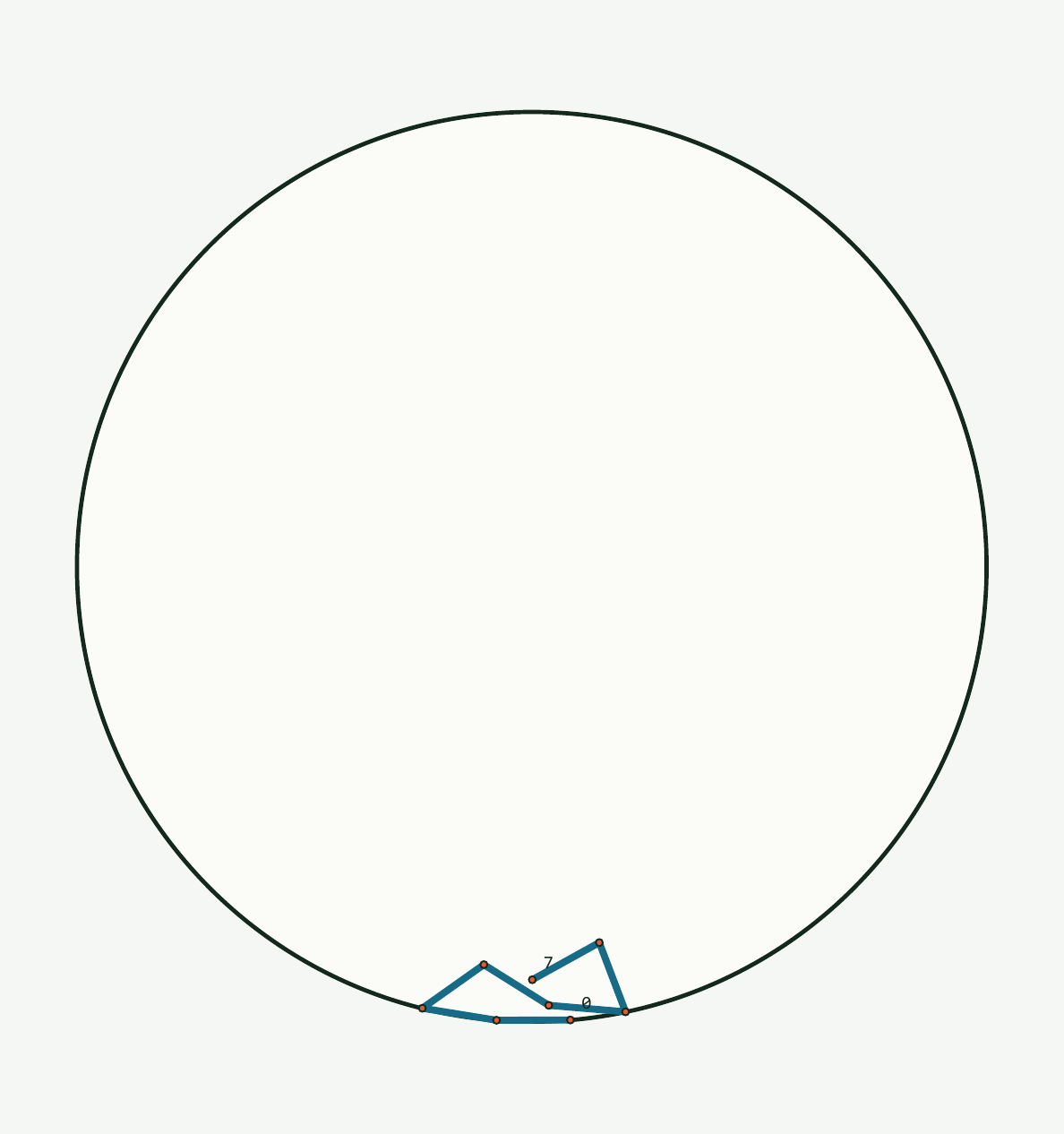} & \includegraphics[width=0.194\textwidth,trim=34 53 34 53,clip]{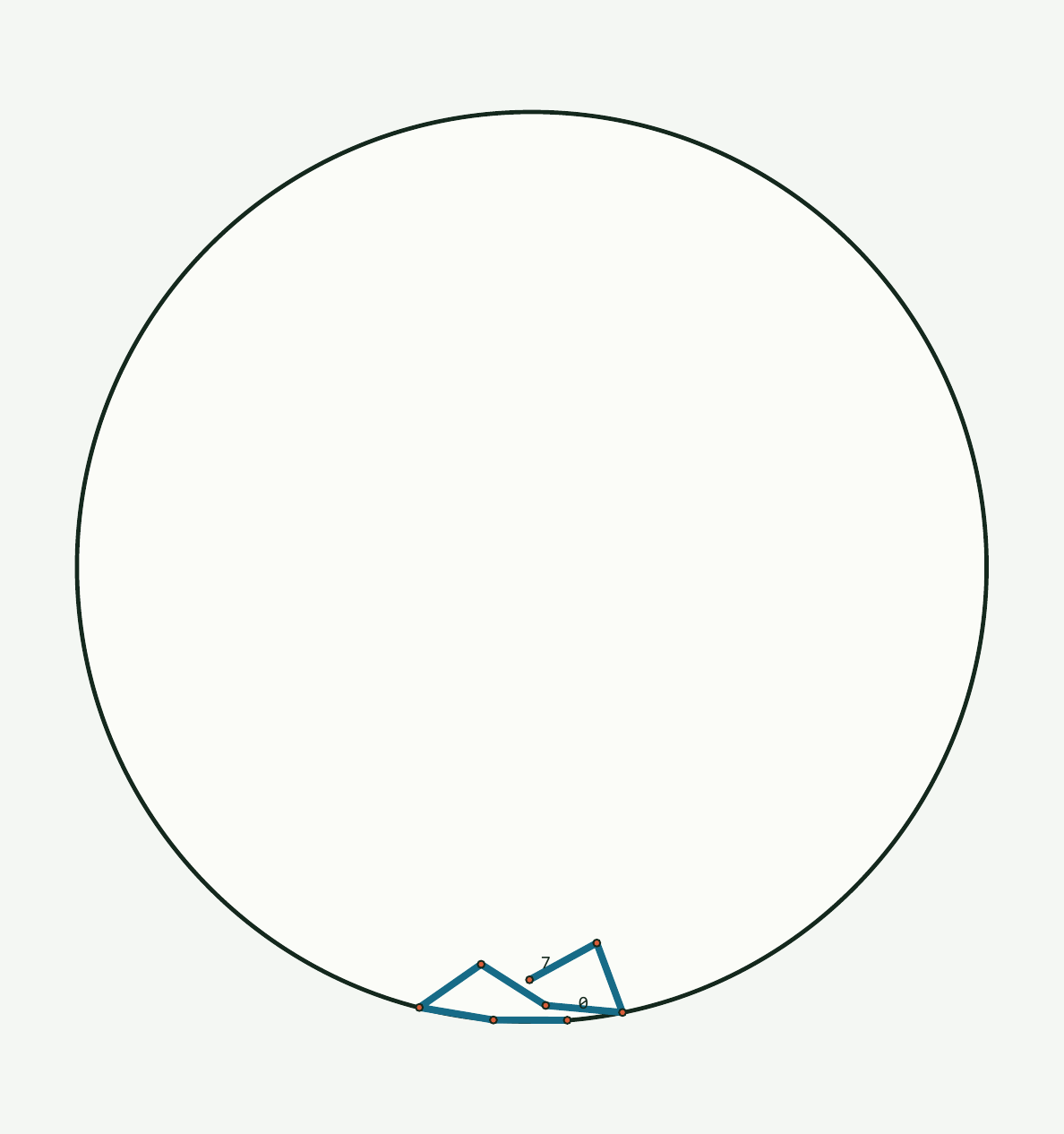} & \includegraphics[width=0.194\textwidth,trim=34 53 34 53,clip]{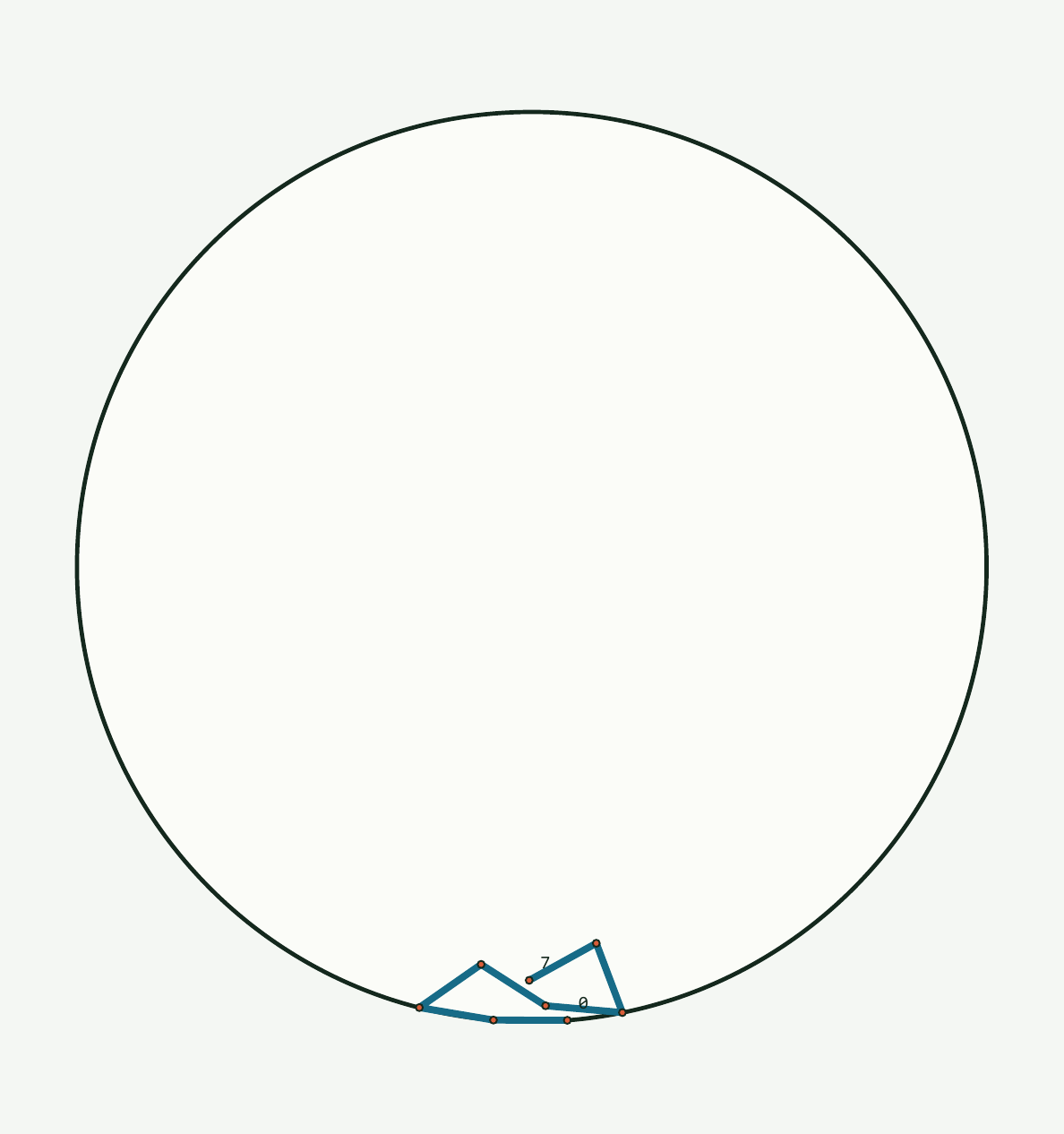} & \includegraphics[width=0.194\textwidth,trim=34 53 34 53,clip]{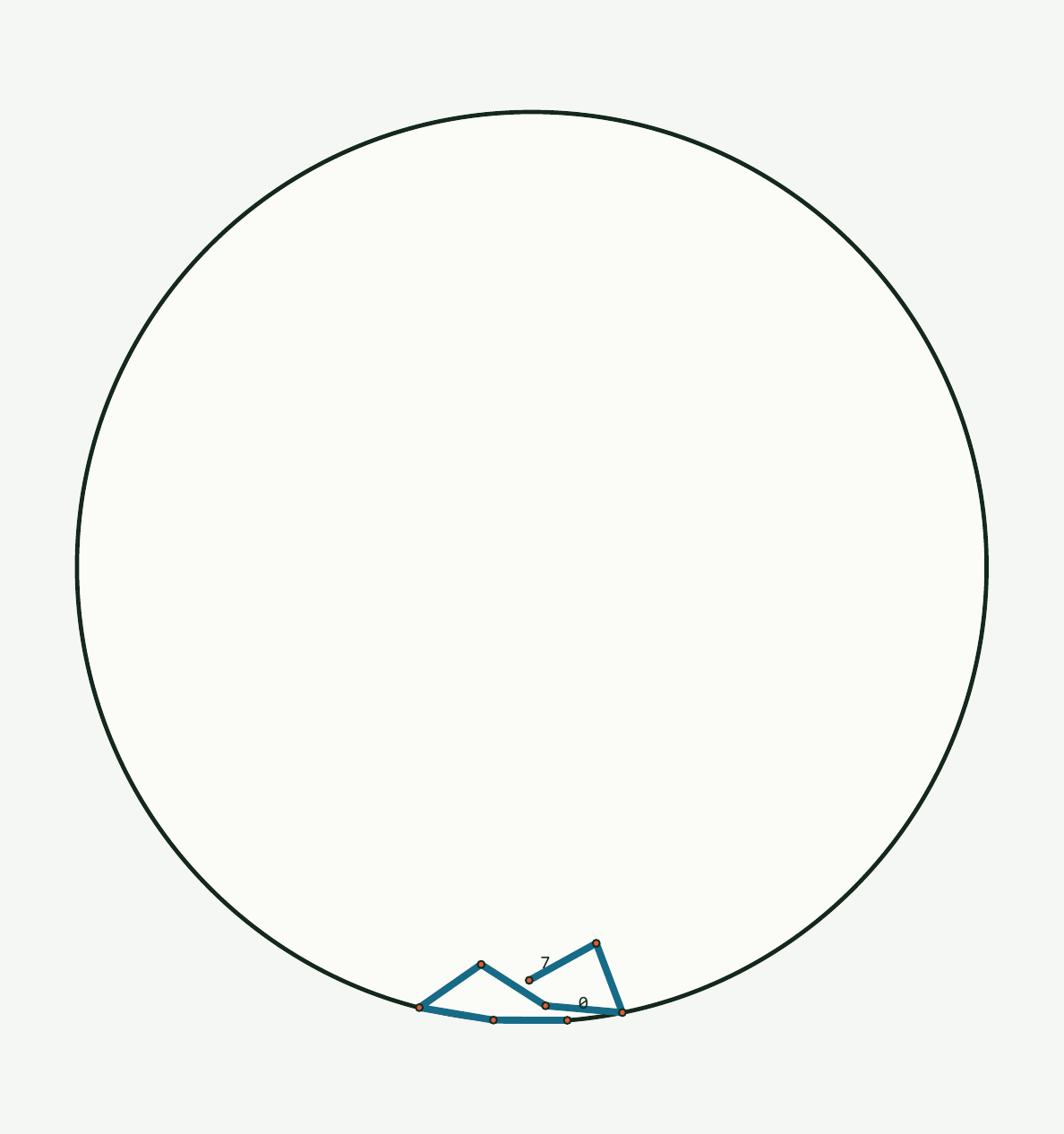} \\
\midrule
\addlinespace[2pt]
\multicolumn{5}{c}{\textbf{Cloth 16}} \\
\includegraphics[width=0.194\textwidth,trim=34 53 34 53,clip]{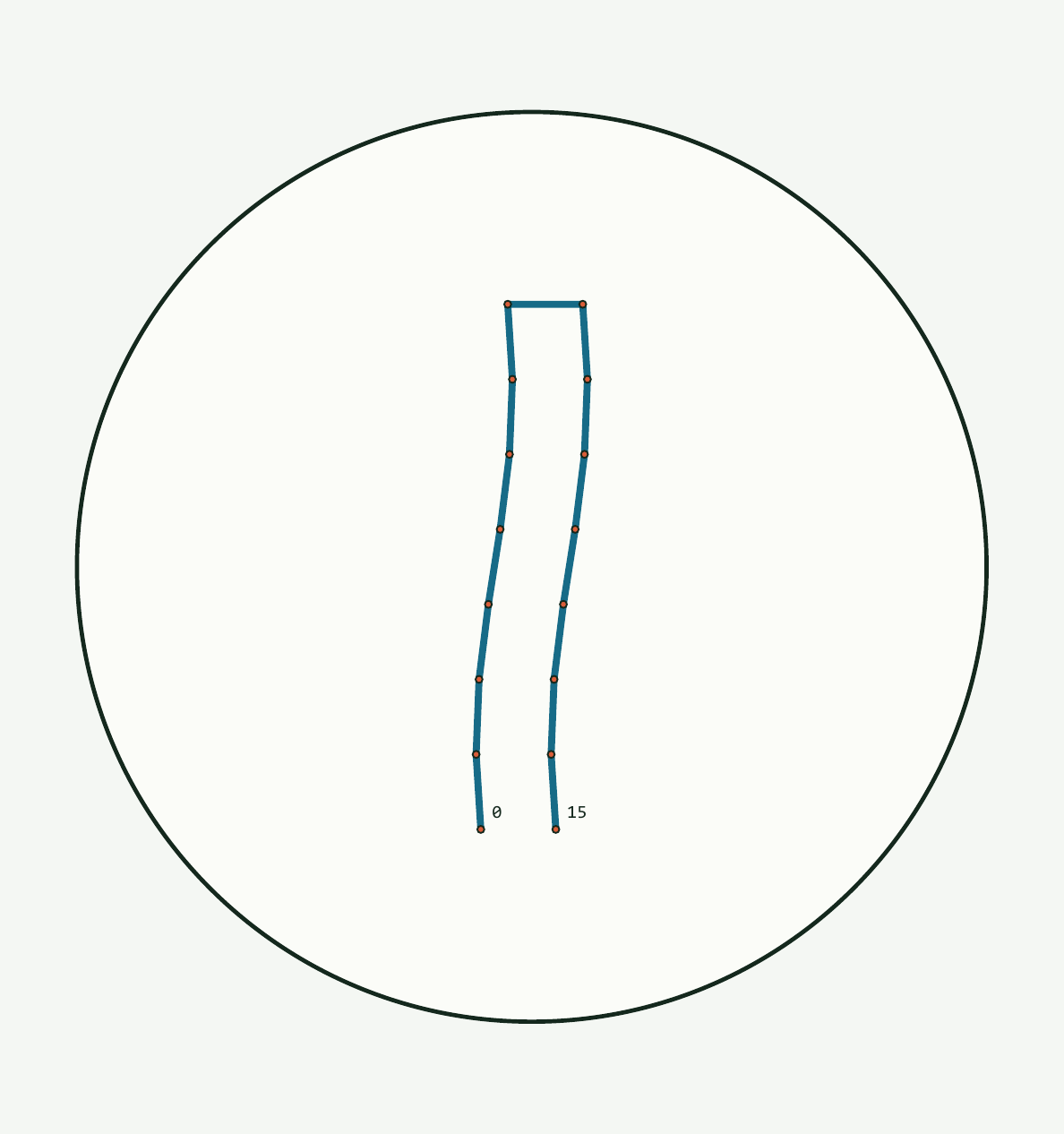} & \includegraphics[width=0.194\textwidth,trim=34 53 34 53,clip]{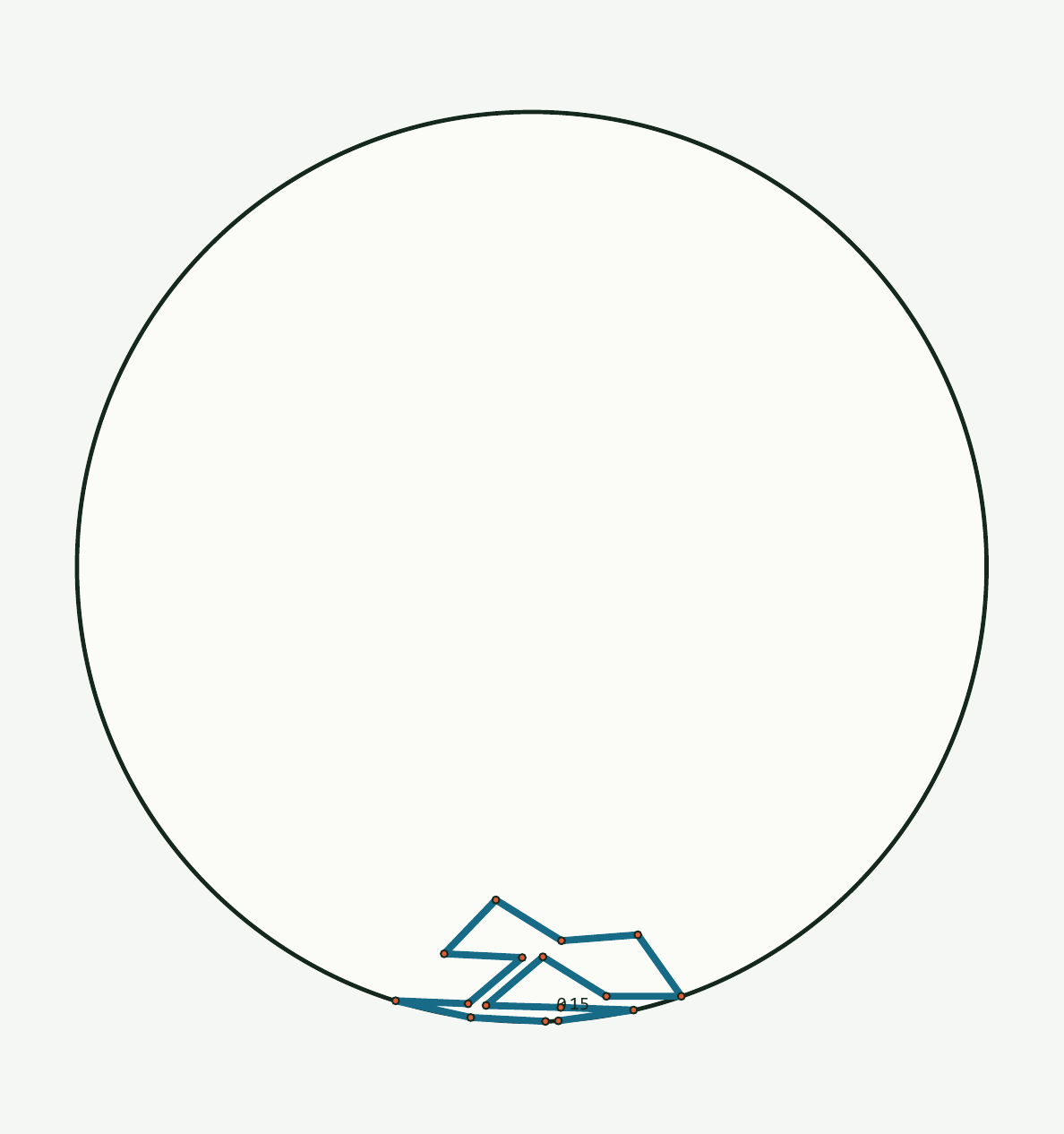} & \includegraphics[width=0.194\textwidth,trim=34 53 34 53,clip]{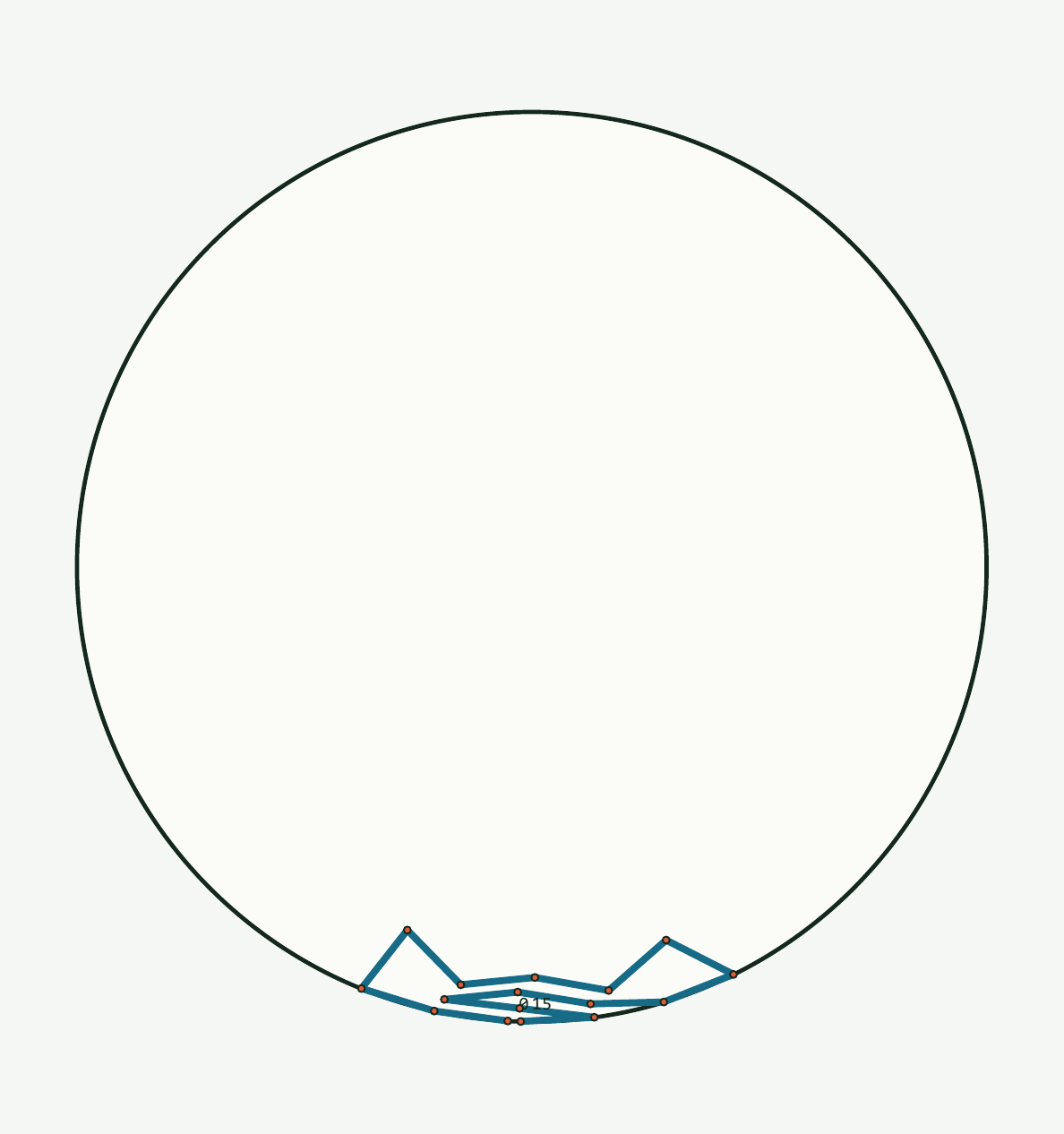} & \includegraphics[width=0.194\textwidth,trim=34 53 34 53,clip]{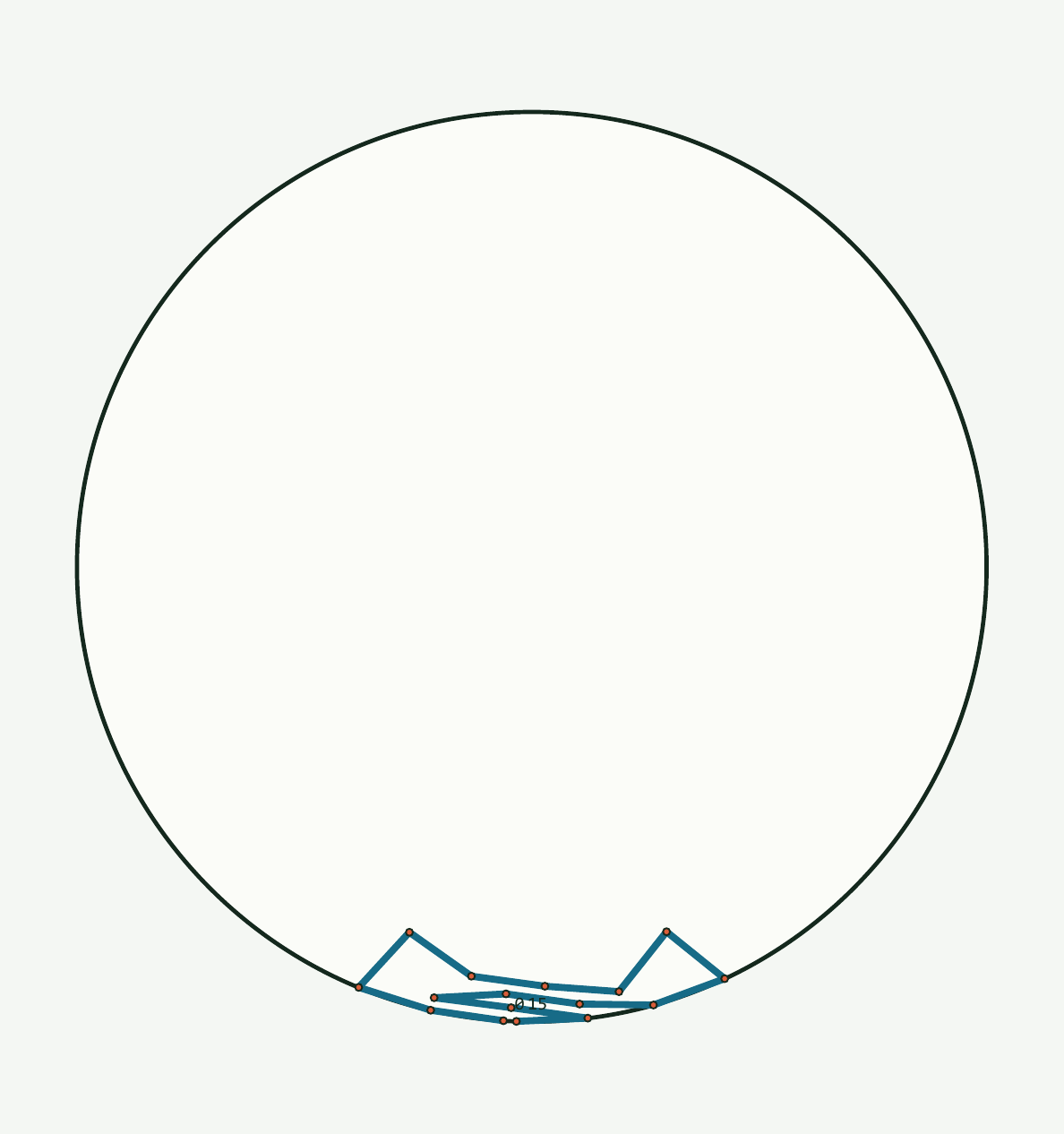} & \includegraphics[width=0.194\textwidth,trim=34 53 34 53,clip]{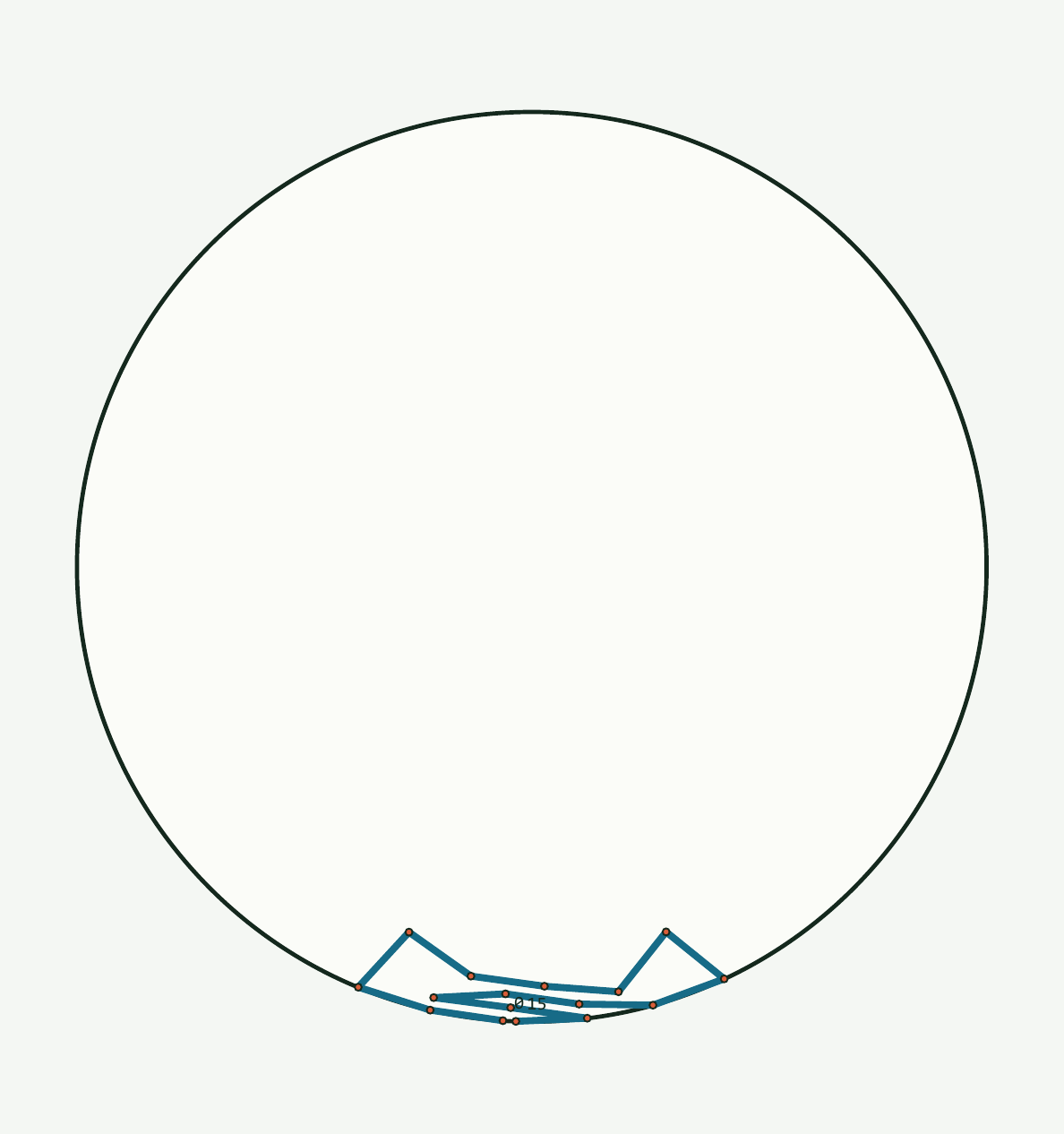} \\
\midrule
\addlinespace[2pt]
\multicolumn{5}{c}{\textbf{Cloth 24}} \\
\includegraphics[width=0.194\textwidth,trim=34 53 34 53,clip]{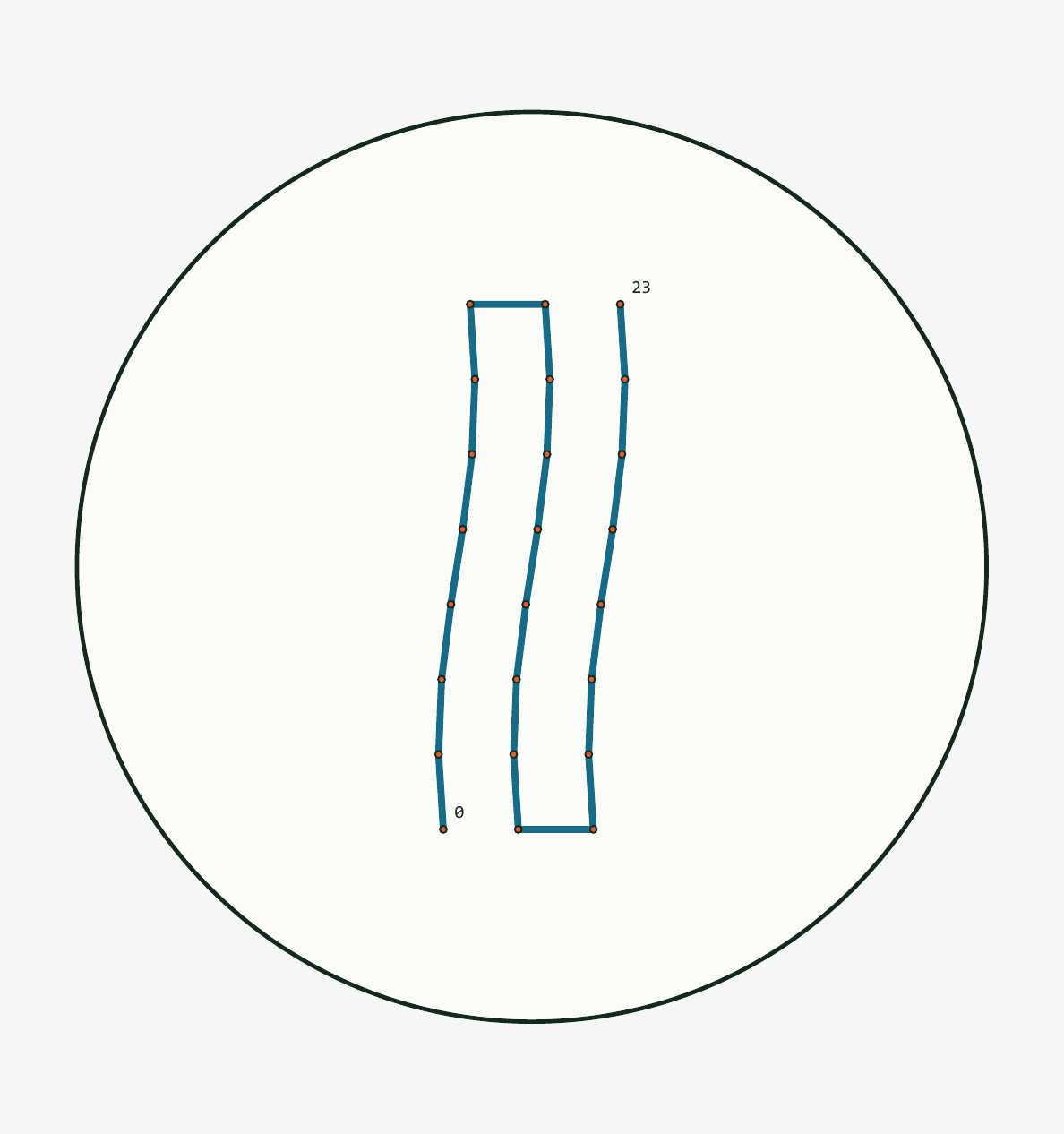} & \includegraphics[width=0.194\textwidth,trim=34 53 34 53,clip]{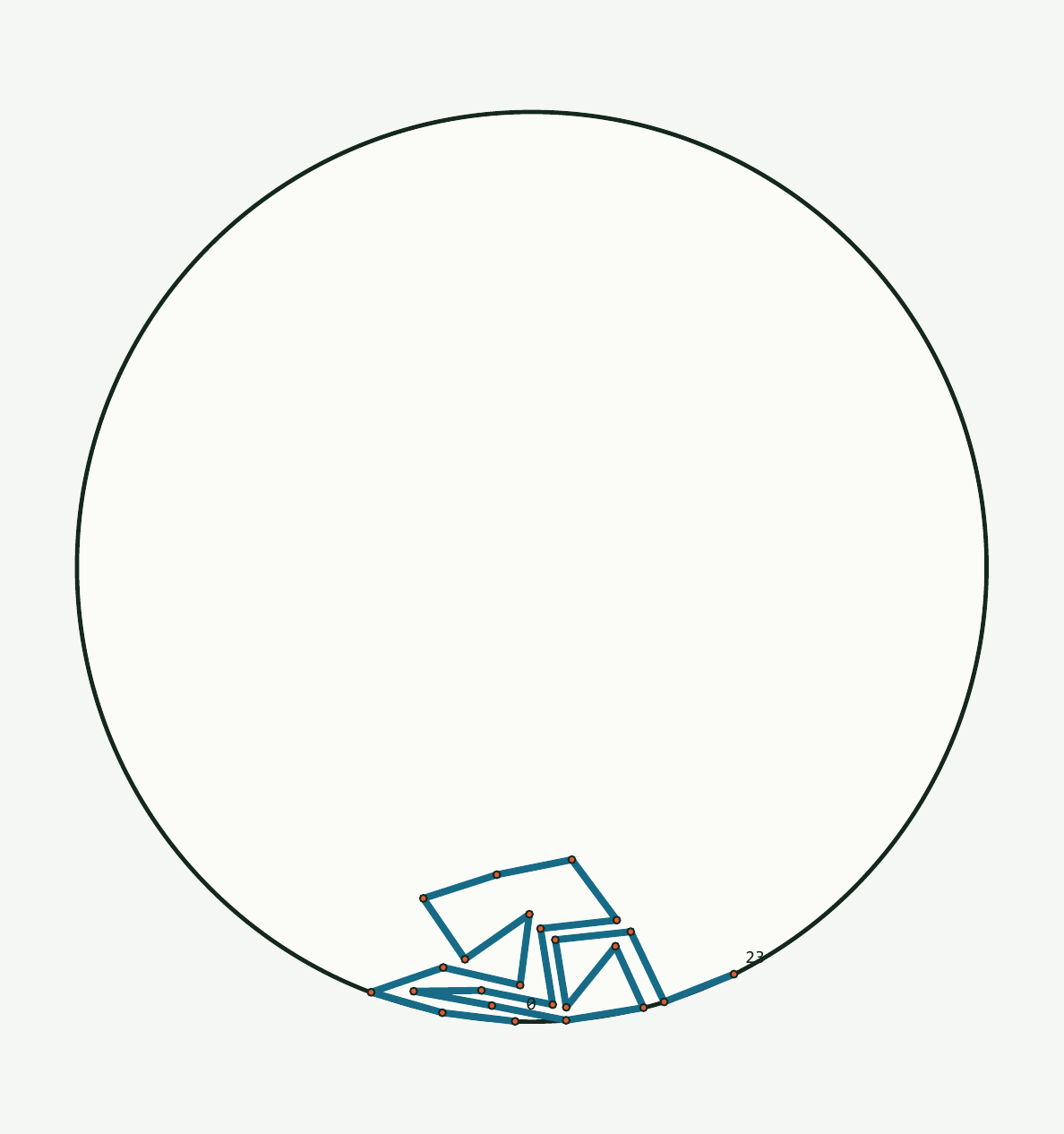} & \includegraphics[width=0.194\textwidth,trim=34 53 34 53,clip]{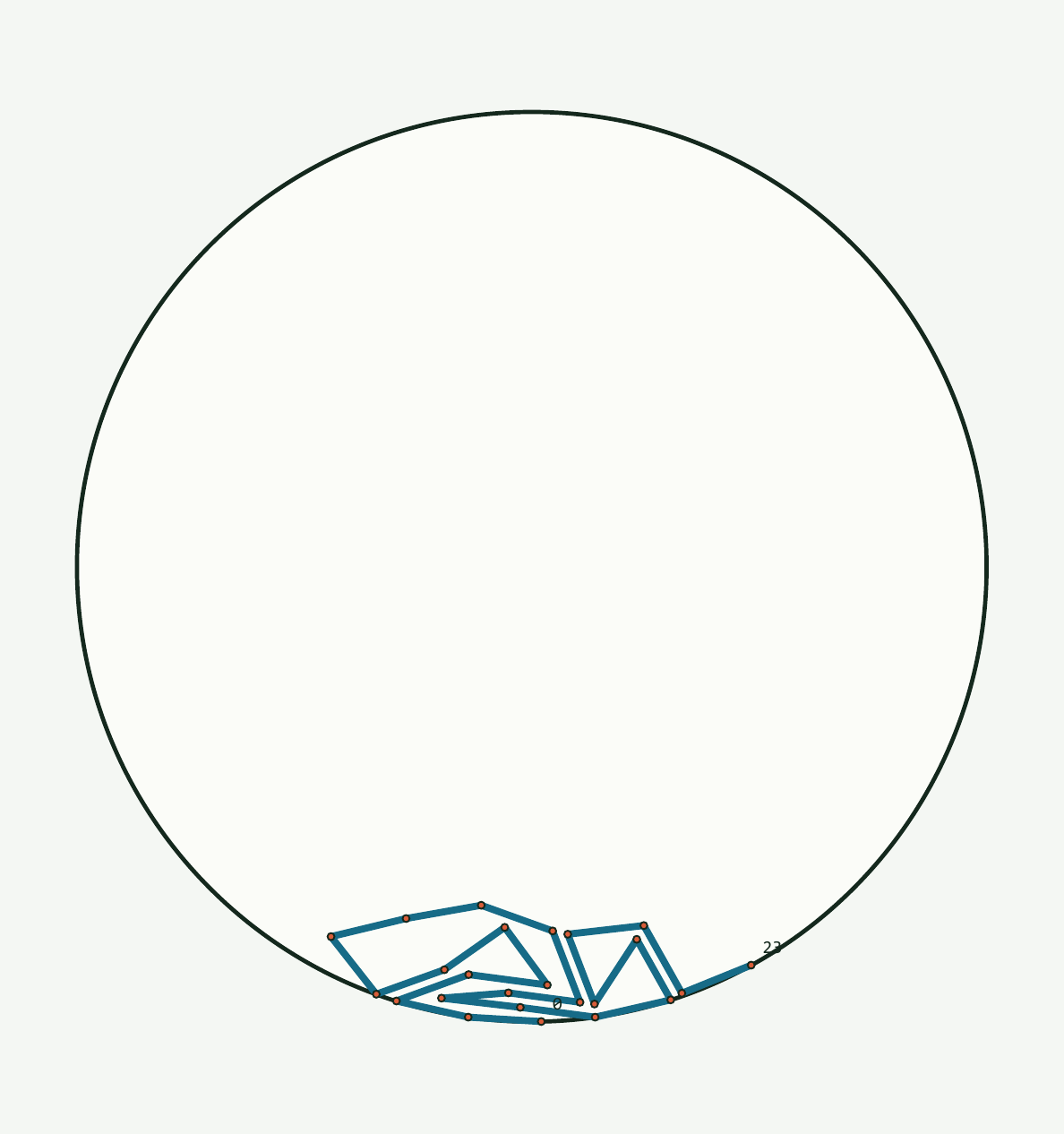} & \includegraphics[width=0.194\textwidth,trim=34 53 34 53,clip]{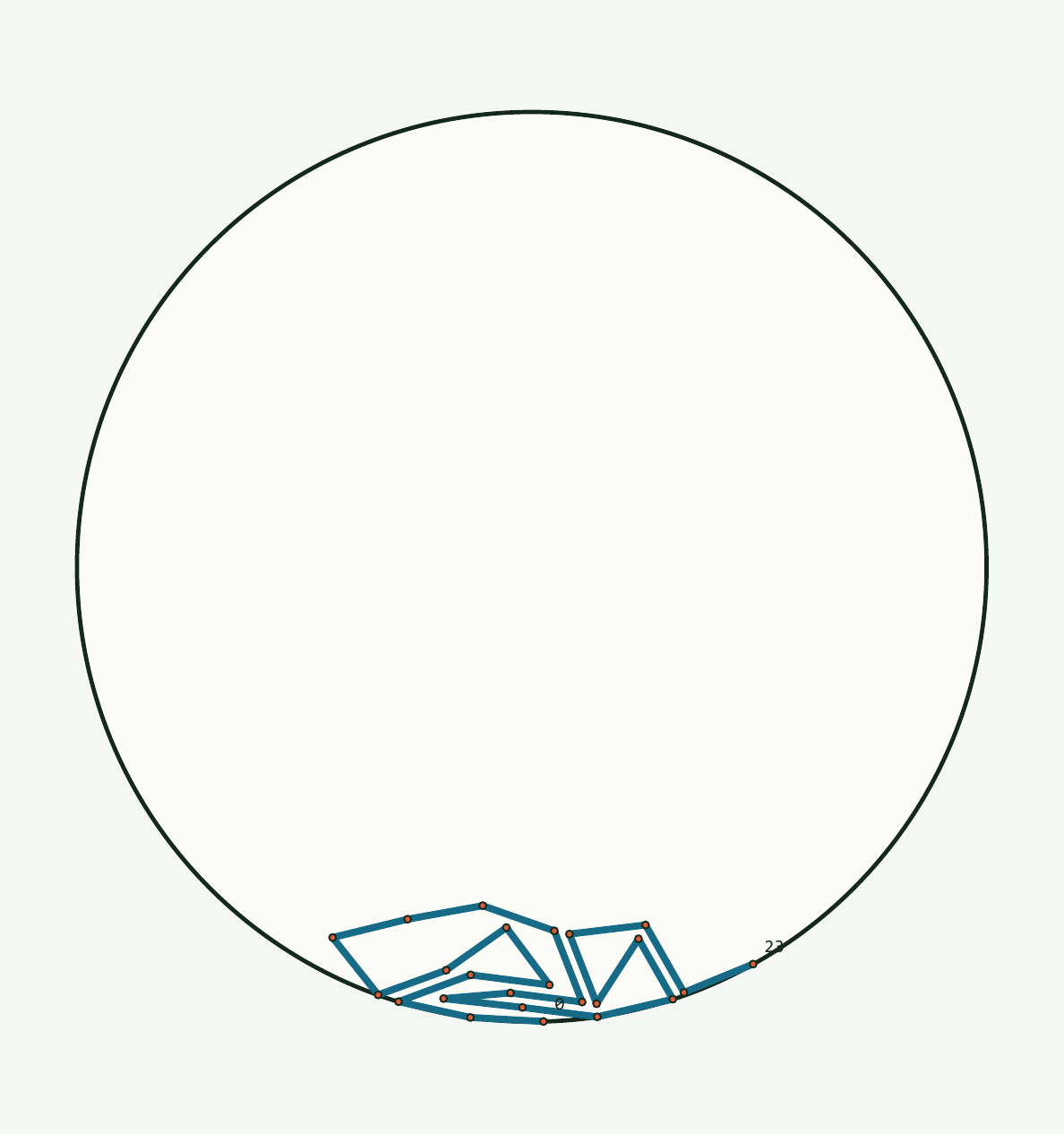} & \includegraphics[width=0.194\textwidth,trim=34 53 34 53,clip]{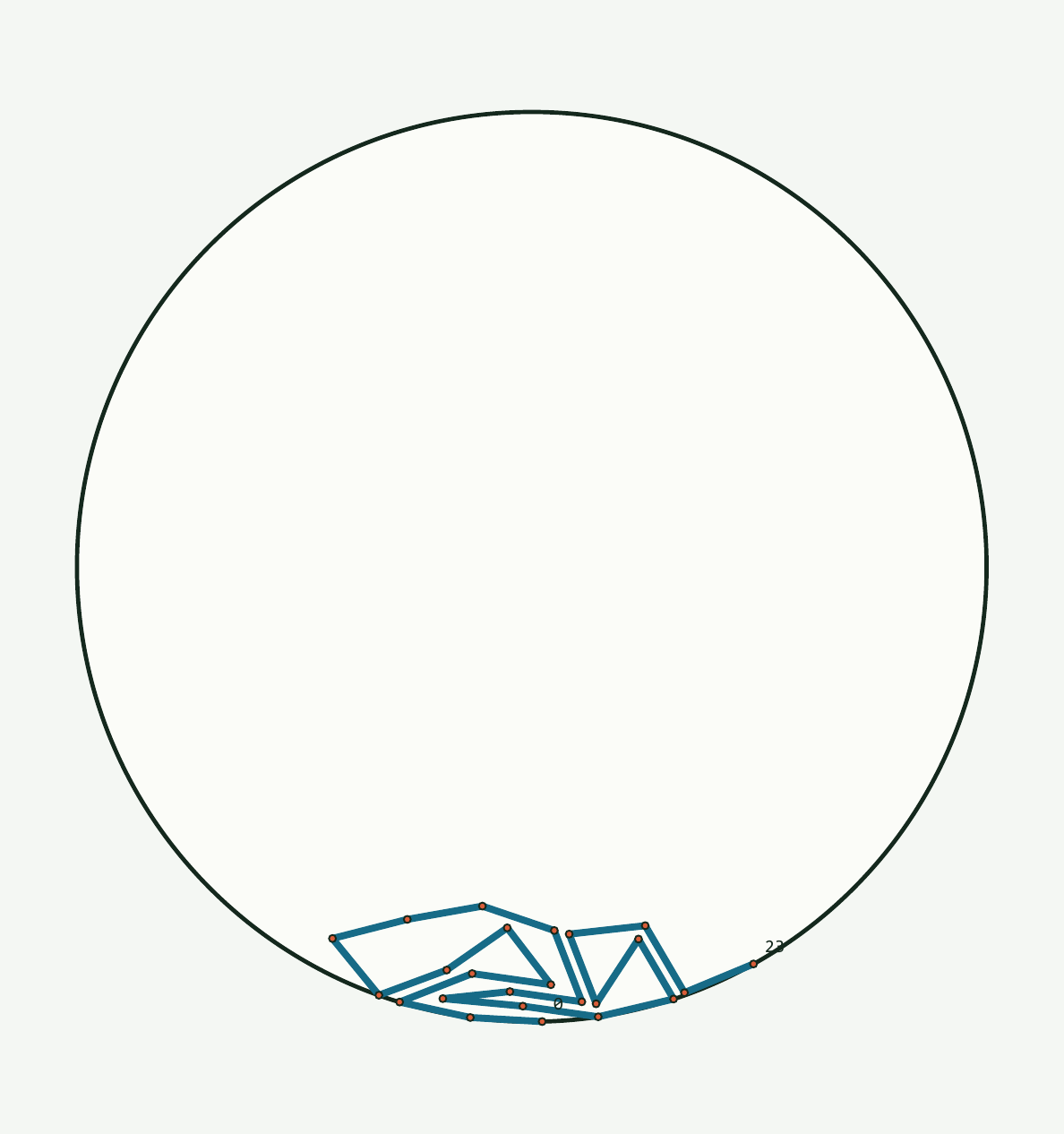} \\
\bottomrule
\end{tabular}
\caption{Ordinary IP-ADMM trajectories for all six instances.  Intermediate columns are the stored samples nearest 25\%, 50\%, and 75\% of total execution; the last column is the final output.}
\label{fig:experiment-trajectories}
\end{figure*}

\begin{figure*}[p]
\centering
\setlength{\tabcolsep}{1pt}
\begin{tabular}{@{}ccc@{}}
\multicolumn{3}{c}{\textbf{Capsule stacking, IP-ADMM}} \\
\includegraphics[width=0.32\textwidth]{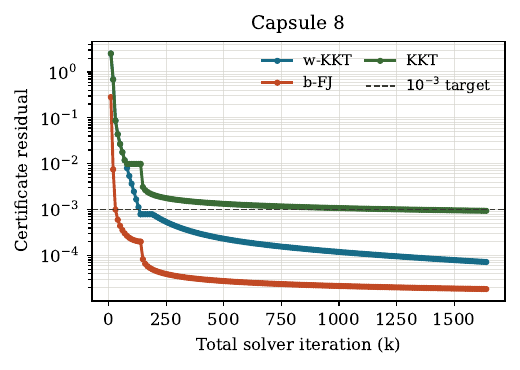} &
\includegraphics[width=0.32\textwidth]{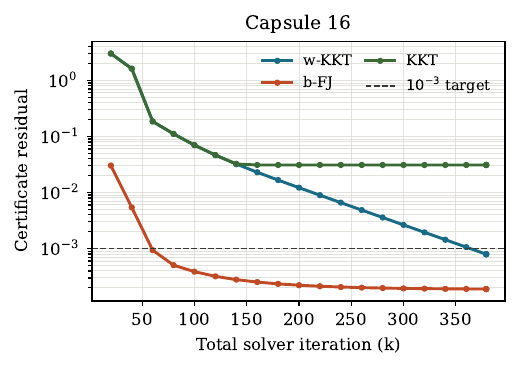} &
\includegraphics[width=0.32\textwidth]{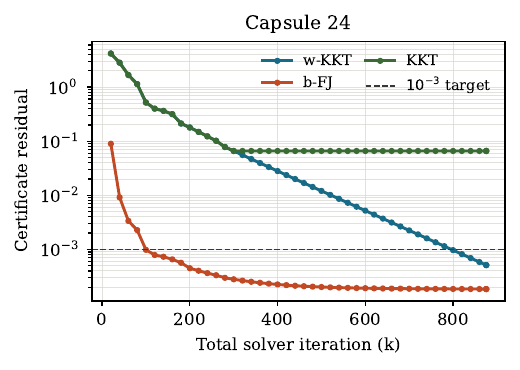} \\
\multicolumn{3}{c}{\textbf{Capsule stacking, IP-ADMM-Ckpt}} \\
\includegraphics[width=0.32\textwidth]{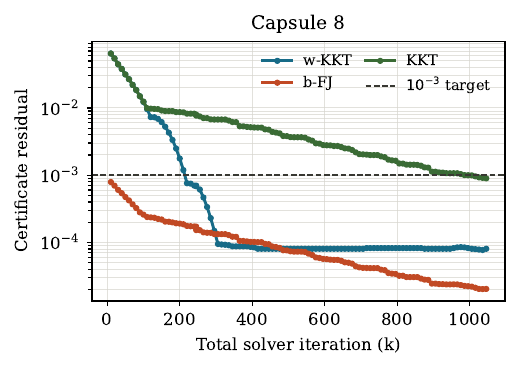} &
\includegraphics[width=0.32\textwidth]{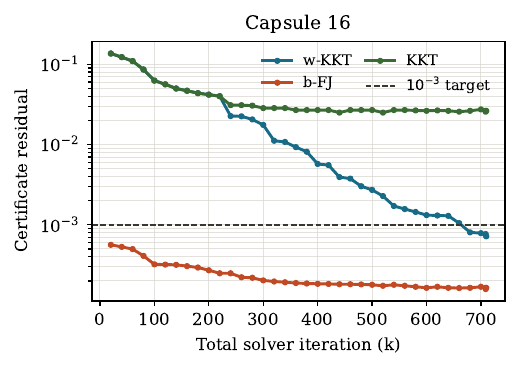} &
\includegraphics[width=0.32\textwidth]{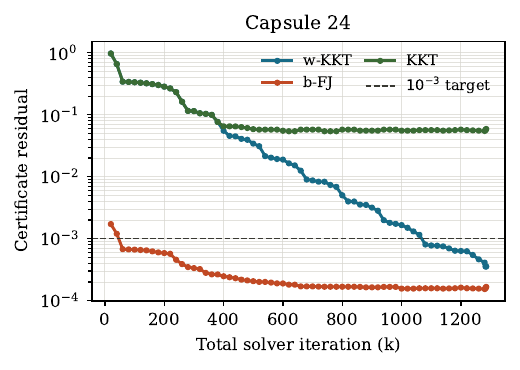} \\
\multicolumn{3}{c}{\textbf{Cloth resting, IP-ADMM}} \\
\includegraphics[width=0.32\textwidth]{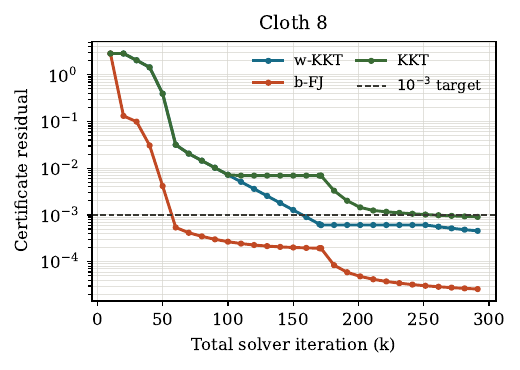} &
\includegraphics[width=0.32\textwidth]{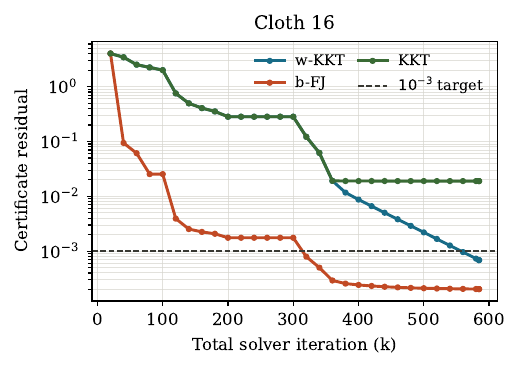} &
\includegraphics[width=0.32\textwidth]{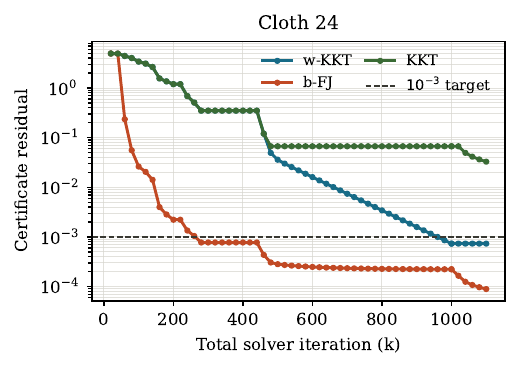} \\
\multicolumn{3}{c}{\textbf{Cloth resting, IP-ADMM-Ckpt}} \\
\includegraphics[width=0.32\textwidth]{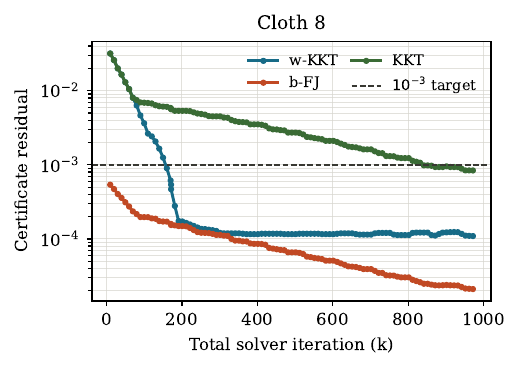} &
\includegraphics[width=0.32\textwidth]{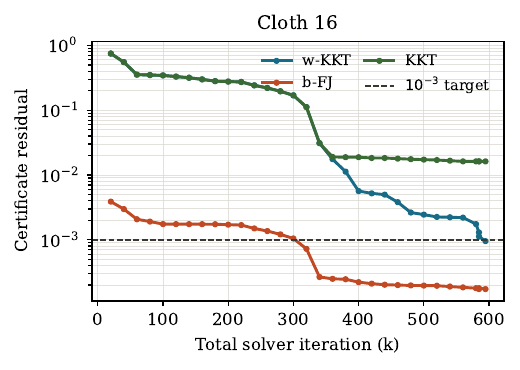} &
\includegraphics[width=0.32\textwidth]{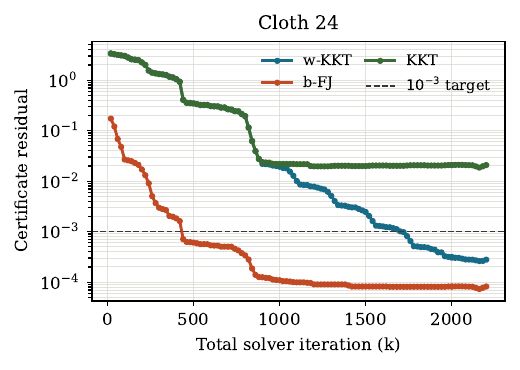} \\
\end{tabular}
\caption{w-KKT, b-FJ, and KKT certificates versus total solver iteration (in thousands) for every instance and ADMM mode. The dashed line is the shared $10^{-3}$ target. KKT reaches the target only at $n=8$; at larger sizes it stalls above the target while w-KKT and b-FJ still descend below it.}
\label{fig:certificates}
\end{figure*}

\begin{figure*}[p]
\centering
\setlength{\tabcolsep}{1pt}
\begin{tabular}{@{}ccc@{}}
\multicolumn{3}{c}{\textbf{Capsule stacking, IP-ADMM}} \\
\includegraphics[width=0.32\textwidth]{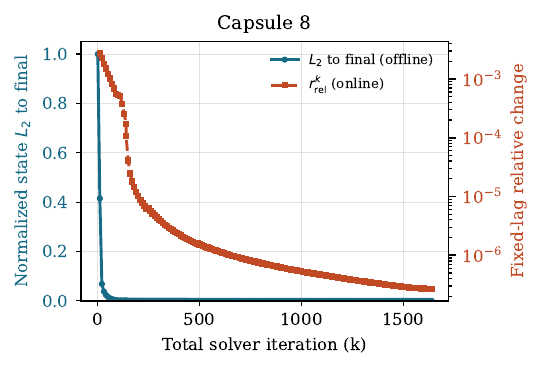} &
\includegraphics[width=0.32\textwidth]{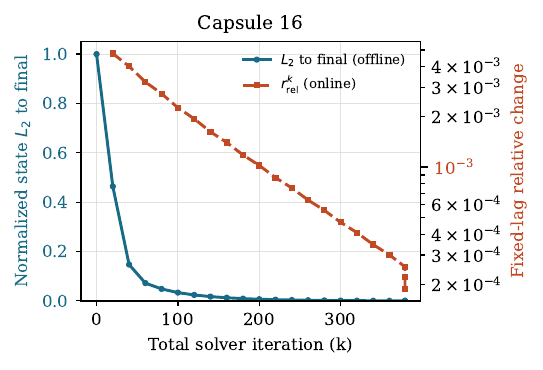} &
\includegraphics[width=0.32\textwidth]{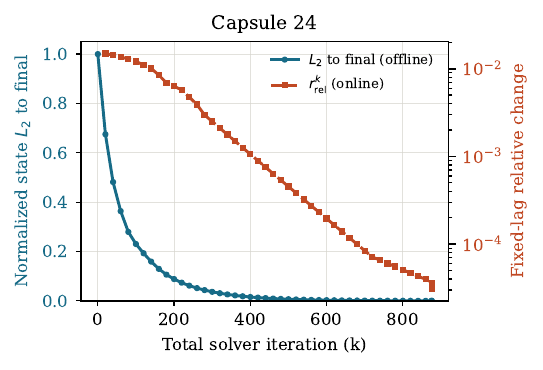} \\
\multicolumn{3}{c}{\textbf{Capsule stacking, IP-ADMM-Ckpt}} \\
\includegraphics[width=0.32\textwidth]{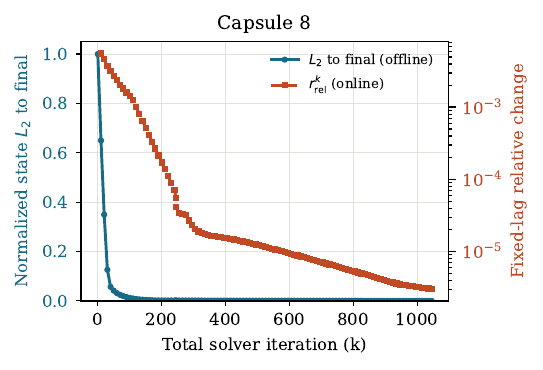} &
\includegraphics[width=0.32\textwidth]{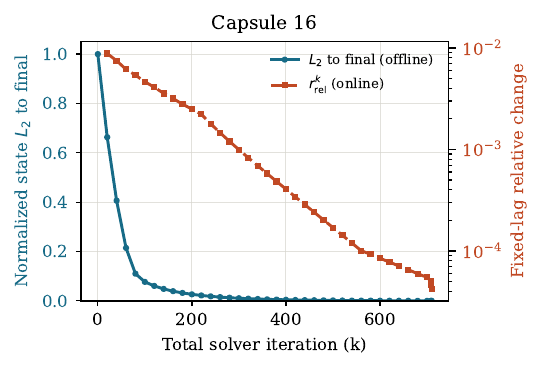} &
\includegraphics[width=0.32\textwidth]{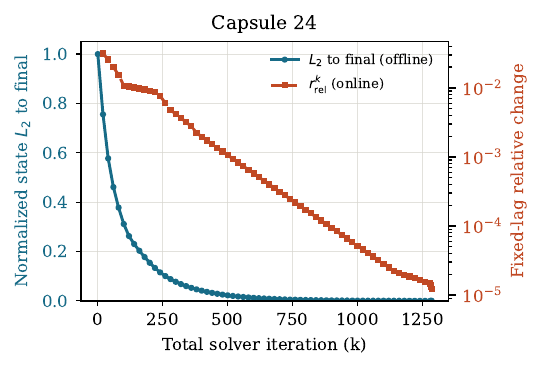} \\
\multicolumn{3}{c}{\textbf{Cloth resting, IP-ADMM}} \\
\includegraphics[width=0.32\textwidth]{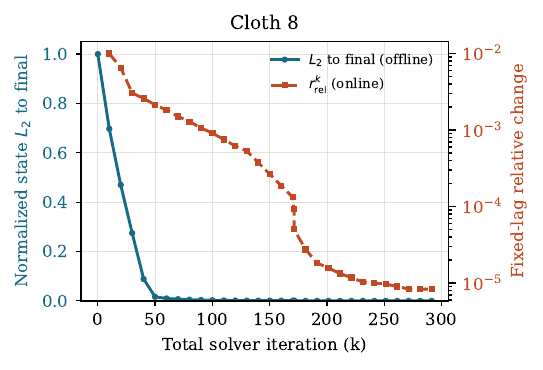} &
\includegraphics[width=0.32\textwidth]{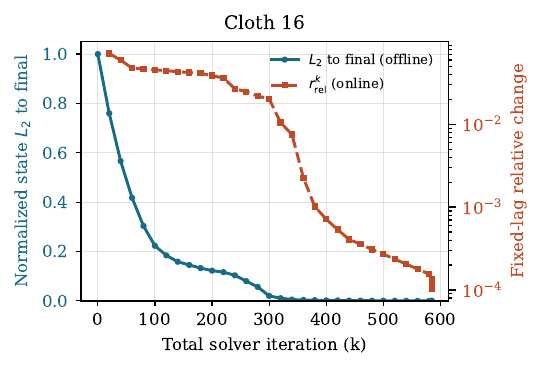} &
\includegraphics[width=0.32\textwidth]{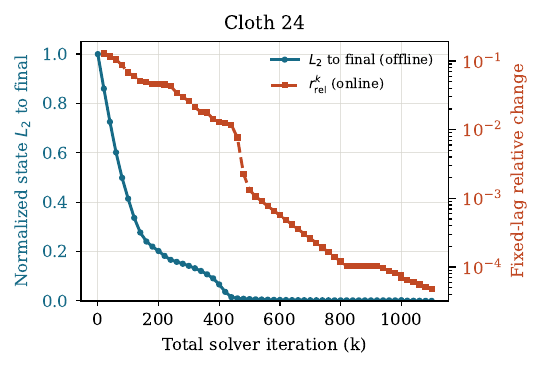} \\
\multicolumn{3}{c}{\textbf{Cloth resting, IP-ADMM-Ckpt}} \\
\includegraphics[width=0.32\textwidth]{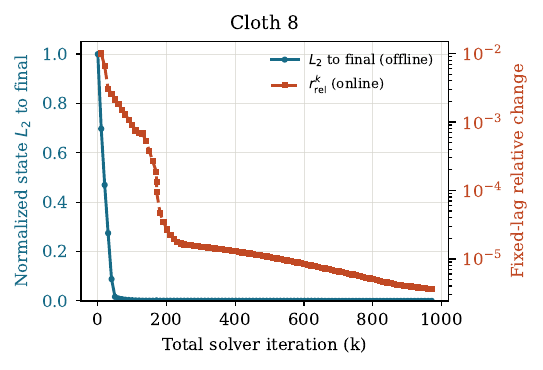} &
\includegraphics[width=0.32\textwidth]{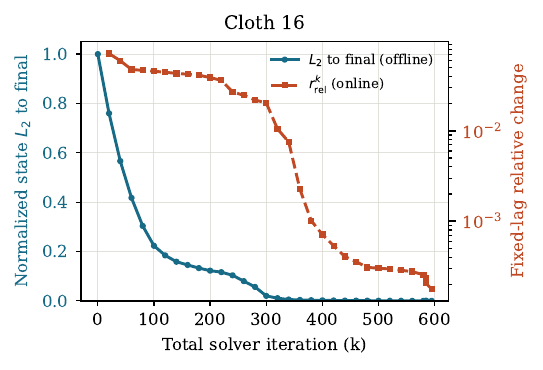} &
\includegraphics[width=0.32\textwidth]{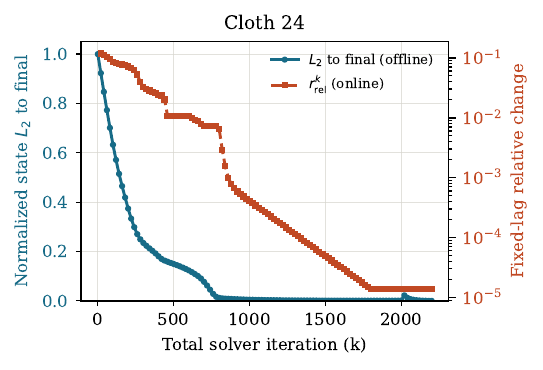} \\
\end{tabular}
\caption{Left axis (blue): normalized state-space $L_2$ distance (Euclidean distance between the current and final solution vectors) of each iterate to the final output, versus total solver iteration (in thousands), for every instance and ADMM mode.  Right axis (orange, log): the online relative-change diagnostic $\|x^{k}-x^{k-\Delta}\|/\max\{1,\|x^{k}\|\}$ over consecutive samples ($\Delta$ the history sampling interval).  The distance flattens near zero well before KKT is established (cf.\ \prettyref{fig:certificates}); the online diagnostic falls fastest during the same early phase and then keeps decaying steadily, so a threshold on it approximates the settling point without marking it sharply.}
\label{fig:visual-difference}
\end{figure*}
\clearpage

\section{Applications in Control and Motion Planning}
\label{appen:control}
Although our experiments are mostly focused on physical simulation, our method can be extended to many other applications, e.g. differentiable trajectory optimization, model-predictive control with obstacles. This section illustrates one of such example, where we optimize a trajectory subject to contact or collision constraints. Specifically, we assume a planar UAV follows a piecewise-linear trajectory $x_0,x_1,\ldots,x_n\in\R^2$ with a fixed start $x_0$.  We optimize the free nodes $x=(x_1,\ldots,x_n)$ to reach a target $x_{\rm goal}$ while avoiding circular obstacles. The objective is a small smoothness (squared-acceleration) term plus a sharp goal term,
\begin{eqaligned}
\sum_{j=2}^{n}\frac{w}{2}\bigl\|x_j-2x_{j-1}+x_{j-2}\bigr\|^2
\;+\;\gamma\,\bigl\|x_n-x_{\rm goal}\bigr\|_1,
\end{eqaligned}
with weights $w>0$ and $\gamma>0$, and each obstacle with center $c$ and radius $r$ contributes, for every segment $(x_{s-1},x_s)$, the clearance constraint $\mathscr M^{\rm rad}_{\delta}(D_{\rm seg}(c;x_{s-1},x_s))-\mathscr M^{\rm rad}_{\delta}(r^2)\geq0$, i.e.\ the same radially mollified point-to-segment distance used for capsule stacking (\prettyref{appen:applications}). Goal reaching is more important than smoothness, so we make it sharp: the $\ell_1$ term is exactly proximable, and its proximal map is soft-thresholding. This is exactly an instance of \prettyref{pb:prob}: the smoothness terms are the smooth blocks $h_i$, the goal term is a convex nonsmooth block $\psi=\gamma\|\cdot-x_{\rm goal}\|_1$, and the clearance blocks are the nonlinear inequalities.  We use $n=20$ nodes and two obstacles placed on the straight start--goal line, so the smoothness-optimal straight path is infeasible and the obstacle constraints are active at the optimum.

\begin{figure*}[t]
\centering
\begin{minipage}[c]{0.48\textwidth}%
\centering
\includegraphics[width=\textwidth]{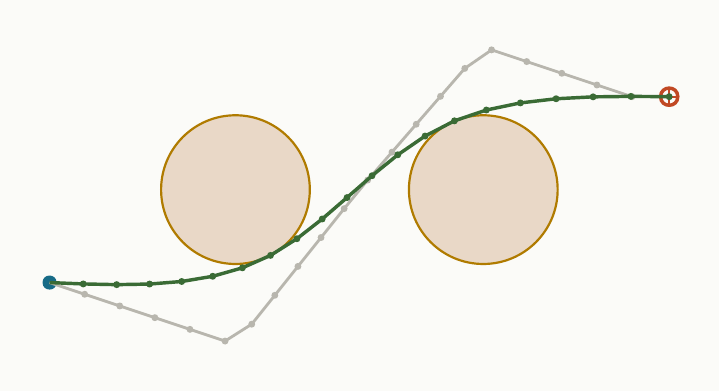}%
\end{minipage}\hfill
\begin{minipage}[c]{0.48\textwidth}%
\centering
\includegraphics[width=\textwidth]{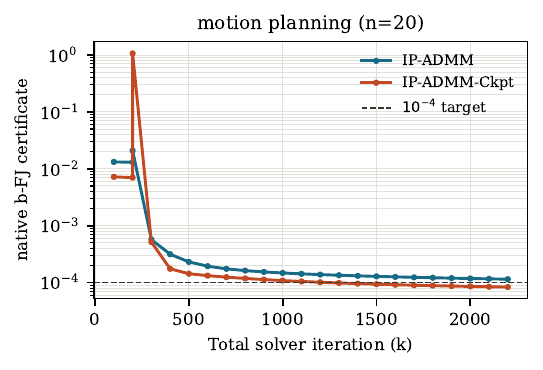}%
\end{minipage}
\caption{UAV trajectory planning.  Left: the initial detour (gray) and the optimized IP-ADMM trajectory (green) with the two circular obstacles; the start is the filled dot and the target is the crossed circle.  Right: the native banded Fritz--John certificate versus total solver iteration for the two IP-ADMM variants, with the shared $10^{-4}$ target marked.}
\label{fig:control}
\end{figure*}

\prettyref{table:control} and \prettyref{fig:control} summarize one deterministic run of SLSQP, trust-constr (both as described in \prettyref{appen:experiments}), BC-ADMM~\cite{pan2026bcadmm}, IP-ADMM, and IP-ADMM-Ckpt from the same strictly feasible initial detour. Except for trust-constr, all converge to essentially the same taut, collision-free path (\prettyref{fig:control} draws the IP-ADMM path, and the SLSQP, BC-ADMM, and both IP-ADMM path lengths agree to two decimals in \prettyref{table:control}), and the $\ell_1$ term drives the final node onto the target. \rev{trust-constr also places the final node on the target, but its trust region shrinks below \texttt{xtol} after $537$ iterations while the path is still $14\%$ longer and $0.39$ away from the obstacles, so it stops at an objective of $8.28$ against $0.90$ for the others. Its model assumes a twice-differentiable objective, which the sharp goal term is not.} The checkpoint variant certifies the solution below the $10^{-4}$ banded Fritz--John level, while the ordinary variant is still approaching it at the same iteration budget. \prettyref{table:rates} gives no b-FJ rate for IP-ADMM-Ckpt, so this is an empirical observation; it parallels the improved post-checkpoint KKT exponent listed there. BC-ADMM can also be run on this instance, since each obstacle-clearance constraint becomes a separating plane between the obstacle disk and one trajectory segment, and the $\ell_1$ goal term is convex with a soft-thresholding proximal map.  BC-ADMM lands on the same taut path and reaches the target exactly, but it certifies nothing for our problem: its barrier holds a fixed margin, so it stops short of the active obstacle constraints (\prettyref{table:control}) and, as in \prettyref{appen:experiments}, has no b-FJ certificate for the original constraints. Its convergence analysis also asks for a twice-differentiable $g$, which this objective is not, so the sharp term we handle exactly lies outside the class that analysis covers.

\begin{table}[t]
\centering\small
\caption{UAV trajectory planning ($n=20$ nodes, two obstacles).  The nonsmooth goal term $\gamma\|x_n-x_{\mathrm{goal}}\|_1$ makes goal reaching exact.  native b-FJ is the solver's banded Fritz--John certificate (IP-ADMM variants only); SLSQP, trust-constr, and BC-ADMM~\cite{pan2026bcadmm} have no native certificate, and the post-hoc KKT reconstruction is unreliable at the $\ell_1$ kink.  trust-constr stops on its \texttt{xtol} test with a collapsed trust region on a longer, non-taut path.  Clearance is the exact minimum obstacle margin (nonnegative means collision-free).  BC-ADMM's clearance is at least two orders of magnitude larger than that of SLSQP and the IP-ADMM variants because its barrier holds a fixed margin, so it stops short of the active constraint rather than resolving it; that standoff is why it has no certificate for the original constraints.  trust-constr's larger clearance comes from its early stop on a non-taut path.  \rev{Objective is the objective value at the final iterate.  Iters counts SciPy iterations for SLSQP and trust-constr and ADMM iterations otherwise; unlike \prettyref{table:all-experiments}, the SciPy counts are not converted to effective sweeps.}}
\label{table:control}
\begin{tabular}{@{}lrrrrrr@{}}
\toprule
Method & Objective & Goal gap & Clearance & Path len. & native b-FJ & Iters \\
\midrule
SLSQP & 0.90 & $2.17\!\times\!10^{-14}$ & $2.44\!\times\!10^{-15}$ & 10.84 & -- & 79 \\
trust-constr & 8.28 & $2.91\!\times\!10^{-13}$ & 0.39 & 12.34 & -- & 537 \\
IP-ADMM & 0.90 & $2.13\!\times\!10^{-5}$ & $1.82\!\times\!10^{-5}$ & 10.84 & $1.15\!\times\!10^{-4}$ & 2200k \\
IP-ADMM-Ckpt & 0.90 & $8.08\!\times\!10^{-6}$ & $2.63\!\times\!10^{-7}$ & 10.84 & $8.35\!\times\!10^{-5}$ & 2200k \\
BC-ADMM & 0.90 & $7.32\!\times\!10^{-15}$ & $1.99\!\times\!10^{-3}$ & 10.84 & -- & 118k \\
\bottomrule
\end{tabular}
\end{table}

\section{\texorpdfstring{\rev{Applications in Learned Contact Constraints}}{Learned Contact Constraints}}
\label{appen:neural}
\begin{revision}
In this appendix, we introduce a problem setup in which the constraints are formulated using a trained neural network: a frozen neural signed distance field (SDF) of a nonconvex cross supplies every body--body contact constraint while $n\in\{4,8,16\}$ crosses settle under gravity in a hexagonal container. As illustrated in \prettyref{fig:neural-setup}, every boundary sample of cross~B is expressed in the body frame of cross~A, where the learned SDF of A is evaluated, and the constraint requires this value to be at least a clearance~$\varepsilon$. This problem setup maps directly to prior works. Testing the boundary samples of one body against the distance field of the other is the classical treatment of nonconvex rigid-body contact~\cite{guendelman2003}; representing this field by a neural network that is evaluated together with its input gradient follows differentiable rigid-body simulation with learned shapes~\cite{strecke2021diffsdfsim} and neural-field collision avoidance in robotics~\cite{koptev2023njsdf,bukhari2025composite}; the smooth $\tanh$ activation is standard for such fields~\cite{gropp2020igr,sitzmann2020siren}; and the clearance offset follows~\cite{strecke2021diffsdfsim}. Since point sampling can miss contact between samples~\cite{macklin2020}, every reported configuration is checked against the exact polygons.
\begin{figure}[t]
\centering
\includegraphics[width=\linewidth]{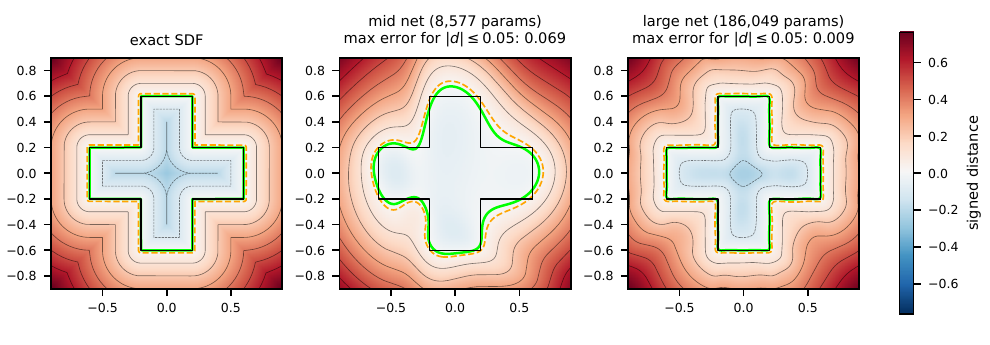}
\caption{Exact SDF of the cross and two trained $\tanh$ fields (green: zero level; orange dashed: level $\varepsilon=0.02$; black: exact polygon). The mid-size field rounds and distorts the zero level set, and its near-surface error ($0.069$) exceeds $\varepsilon$; the large field, used in all runs, stays below $\varepsilon$ ($0.009$).}
\label{fig:neural-sdf}
\end{figure}

\begin{figure}[t]
\centering
\includegraphics[width=\linewidth]{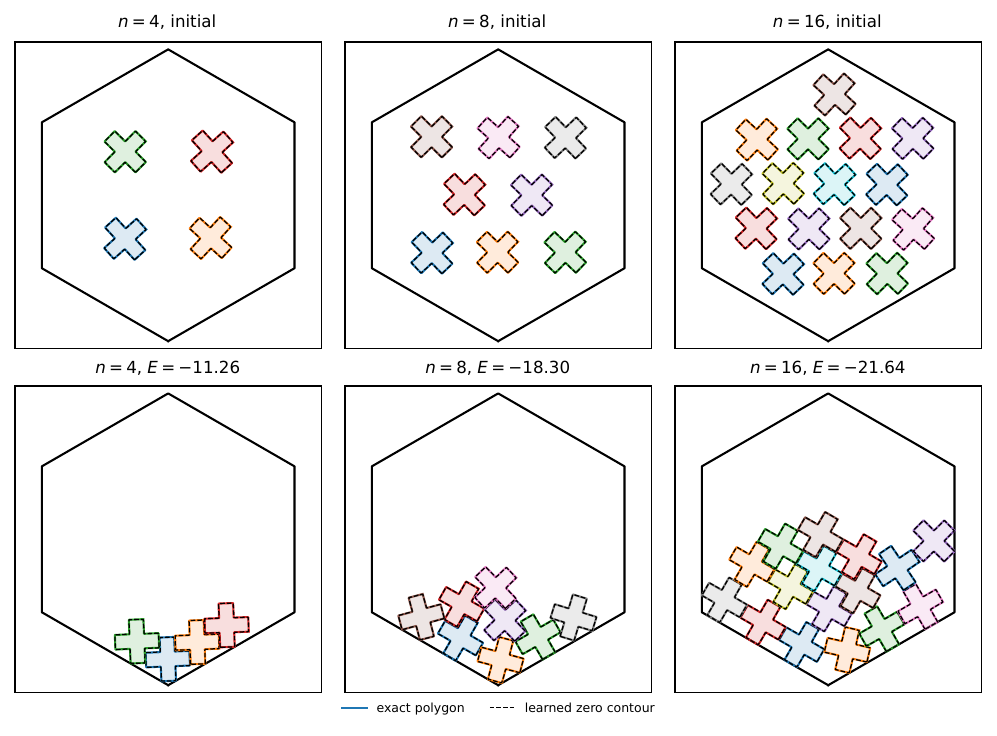}
\caption{Initial (top) and certified (bottom) configurations of the three certified IP-ADMM runs in \prettyref{table:neural-results}, with exact polygons (solid) and the learned zero level set (dashed). Each certified configuration passes the exact polygon checks.}
\label{fig:neural-scenes}
\end{figure}

\begin{figure}[t]
\centering
\includegraphics[width=\linewidth]{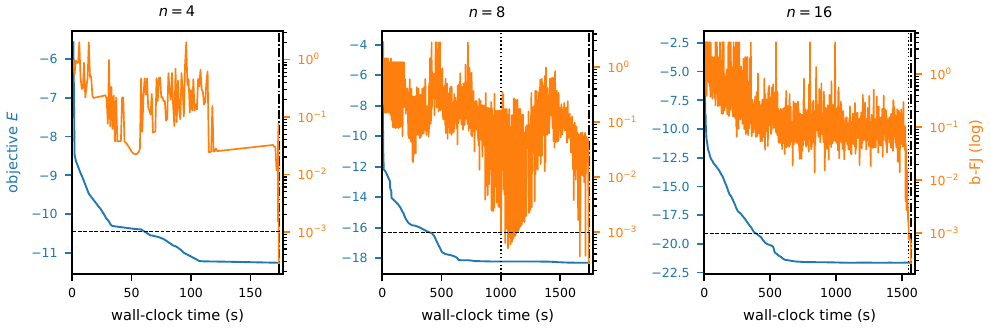}
\caption{Objective (left axis) and b-FJ residual (right axis, log scale; dashed: $10^{-3}$) versus wall-clock time for the three certified IP-ADMM runs, plotted up to the certified iterate. Dotted and dash-dotted lines mark the first iterates meeting the b-FJ target and certified.}
\label{fig:neural-curves}
\end{figure}
\paragraph{Formulation.}
The cross is $([-0.6,0.6]\times[-0.2,0.2])\cup([-0.2,0.2]\times[-0.6,0.6])$. Its SDF $s_\theta$ is a $\tanh$ network with six hidden layers of width $192$ ($186{,}049$ parameters), trained on exact signed distances over the body-frame disk of radius $9$ and then frozen; \prettyref{fig:neural-sdf} shows that a smaller network does not represent the shape accurately enough. Cross $i$ has pose $q_i=(t_i,\vartheta_i)$, the objective is the gravitational energy $\sum_i(t_i)_2$, and $S$ boundary samples $p_k$ are placed at density $\rho_{\rm s}$. Every ordered pair $i\neq j$ and sample $k$ gives the constraint
\begin{eqaligned}
f_{ij,k}(q_i,q_j)=s_\theta\bigl(R(\vartheta_j)^T(R(\vartheta_i)p_k+t_i-t_j)\bigr)-\varepsilon\geq0,\qquad\varepsilon=0.02,
\label{eq:neural-pair}
\end{eqaligned}
and every sample also gives six exact half-plane constraints of the hexagon (circumradius $4$). Each $A_\ell$ selects $(q_i,q_j)$ or $q_i$ from $x$, so the coupling stays affine. \prettyref{table:neural-settings} lists the sizes: $69{,}120$ to $645{,}120$ constraints on $12$ to $48$ pose variables. To keep the number of active constraints small, we use an active-set strategy that progressively inserts a constraint once its value falls within twice the current wall band, $f_\ell\leq2\beta^k$. This follows the broadphase collision detection of prior contact solvers, which generate contact constraints only for primitive pairs within a distance threshold~\cite{bridson2002cloth,li2020ipc}.

\paragraph{Assumption justification.}
The compact-domain realization of \prettyref{appen:applications} applies as for capsule stacking. The exact wall constraints keep every vertex of every cross in the hexagon, and orientations are periodic, so the pose envelope is compact. A $\tanh$ network with fixed finite weights is $C^\infty$ with globally bounded input Jacobian and Hessian (by propagating weight norms through $|\tanh'|\leq1$ and $|\tanh''|\leq4/(3\sqrt3)$), and the relative-position map in \prettyref{eq:neural-pair} is smooth, so every constraint is $C^{1,1}$ near the compact envelope; this supplies \prettyref{ass:f}--\prettyref{ass:bg}. Linear gravity with a constant shift supplies \prettyref{ass:obj}, and the strictly feasible initial poses supply \prettyref{ass:feasible}. Uniform Band-MFCQ cannot be verified for such densely sampled learned constraints, so the CQ-free b-FJ certificate is the guarantee that can be expected in general.

\input{artifacts_neural_collision_illustration}
\paragraph{Baseline comparison.}
We compare IP-ADMM (GPU, float32) with the SLSQP and trust-constr baselines described in \prettyref{appen:experiments}. All solvers use the same frozen field on the GPU, the same initial poses, and a budget of $3600$\,s. A configuration is certified when its b-FJ residual is at most $10^{-3}$, its violation at most $10^{-6}$, its forces and torques balance to $10^{-3}$, and its exact polygons are disjoint and inside the container. \prettyref{table:neural-results} shows that IP-ADMM is certified first in the three reported cases: in $174$\,s at $n=4$, against $787$\,s for trust-constr, and as the only certified solver at $n=8$ and $n=16$. IP-ADMM certifies $n=16$ slightly faster than $n=8$ because the sixteen crosses already fill most of the container, so they fall a shorter distance before coming to rest. SLSQP diverges in every case, and trust-constr runs out of time for $n\geq8$. Both SciPy solvers treat the $M$ rows as one block, so each iteration solves a quadratic or barrier subproblem over all $M$ rows on the CPU, and SLSQP's linearized steps of the $\tanh$ field leave the feasible set and diverge. IP-ADMM instead updates each row independently in parallel on the GPU, couples the rows only through one $3n\times3n$ linear system, and keeps every split state strictly feasible, so each constraint row evaluates the network only at points where that row holds; the consensus iterate $x$ itself need not be feasible, which is why certification also checks the violation. \prettyref{fig:neural-scenes} shows the certified configurations and \prettyref{fig:neural-curves} the convergence histories.

\begin{table}[t]
\centering
\caption{Learned-field experiment: sizes and settings. Values $a/b/c$ are for $n=4/8/16$. The IP-ADMM settings are practical, not the theorem-sufficient construction; the pair-row constant $L_\ell$ is a quantile of sampled local estimates.}
\label{table:neural-settings}
\begin{revision}\scriptsize
\setlength{\tabcolsep}{3pt}
\begin{tabular}{@{}lll@{}}
\toprule
Group & Parameter & Value \\
\midrule
Problem & crosses $n$; pose dimension $3n$ & $4/8/16$; $12/24/48$ \\
 & sample density $\rho_{\rm s}$; samples per cross $S$ & $400$; $1920$ \\
 & constraints $M$ (pair $+$ wall) & $69{,}120/199{,}680/645{,}120$ \\
Field & architecture; parameters; training domain & $\tanh$, $6\times192$; $186{,}049$; disk of radius $9$ \\
IP-ADMM & $\rho^0,L^0,G,\eta,\sigma,\xi,p,\theta_1,\theta_2$ & $20$, $10^6$, $0.006$, $3\times10^5$, $1$, $0.1$, $0.5$, $1.01$, $0.01$ \\
SLSQP & \texttt{ftol} & $10^{-12}$ \\
trust-constr & \texttt{gtol}, \texttt{xtol} & $10^{-12}$, $10^{-12}$ \\
Targets & active rows; violation; b-FJ & $f_\ell\leq2\times10^{-4}$; $10^{-6}$; $10^{-3}$ \\
Runtime & IP-ADMM; SciPy iterations; field & GPU float32; CPU float64; GPU float32 \\
\bottomrule
\end{tabular}
\end{revision}
\end{table}

\begin{table}[t]
\centering
\caption{Baseline comparison. b-FJ is the smallest b-FJ residual over feasible iterates, viol.\ the constraint violation of the certified iterate for certified runs and of the last iterate otherwise.}
\label{table:neural-results}
\begin{revision}\scriptsize
\setlength{\tabcolsep}{4pt}
\begin{tabular}{@{}rrllll@{}}
\toprule
$n$ & $M$ & Solver & b-FJ & viol. & Status \\
\midrule
\multirow{3}{*}{$4$} & \multirow{3}{*}{$69{,}120$} & IP-ADMM & $2.5\times10^{-4}$ & $0$ & \textbf{certified at 174\,s} \\
& & SLSQP & $0.15$ & $1.6\times10^{11}$ & diverged after 65 iterations \\
& & trust-constr & $2.0\times10^{-4}$ & $0$ & certified at 787\,s \\
\midrule
\multirow{3}{*}{$8$} & \multirow{3}{*}{$199{,}680$} & IP-ADMM & $2.6\times10^{-4}$ & $0$ & \textbf{certified at 1742\,s} \\
& & SLSQP & $3.7\times10^{-3}$ & $1.3\times10^{-4}$ & diverged after 41 iterations \\
& & trust-constr & $5.8\times10^{-2}$ & $0$ & time limit \\
\midrule
\multirow{3}{*}{$16$} & \multirow{3}{*}{$645{,}120$} & IP-ADMM & $2.7\times10^{-4}$ & $0$ & \textbf{certified at 1568\,s} \\
& & SLSQP & $0.67$ & $121.6$ & diverged after 13 iterations \\
& & trust-constr & $0.39$ & $2.7\times10^{-5}$ & time limit \\
\bottomrule
\end{tabular}
\end{revision}
\end{table}
\end{revision}

\section{\texorpdfstring{Parameter Sensitivity:\\Theorem-Sufficient versus Practical Tuning}{Parameter Sensitivity: Theorem-Sufficient versus Practical Tuning}}
\label{appen:ablation}
\begin{figure*}[t]
\centering
\includegraphics[width=\textwidth]{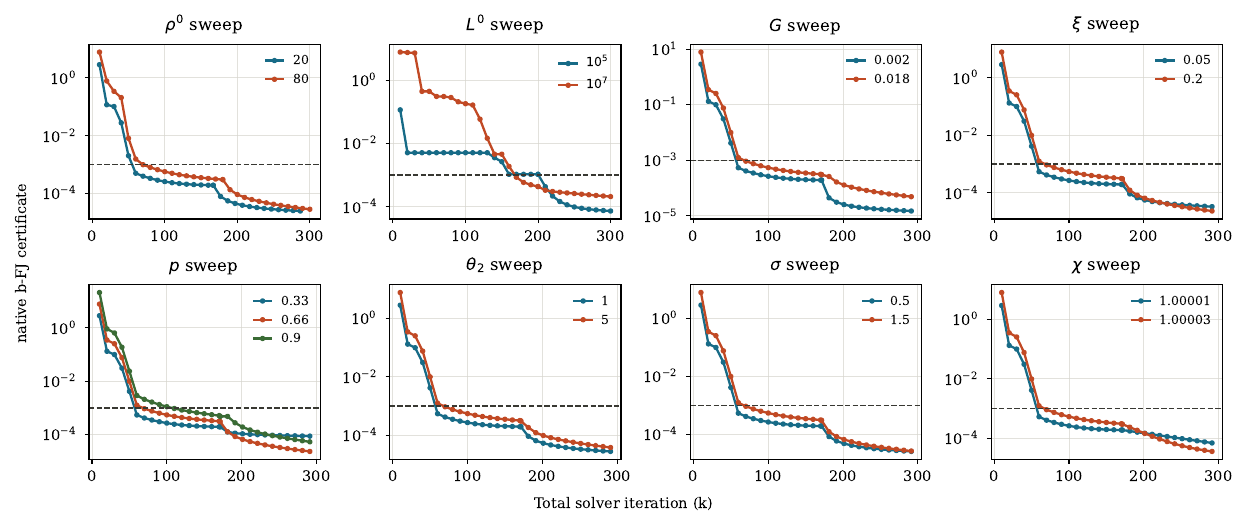}\\[2pt]
\includegraphics[width=\textwidth]{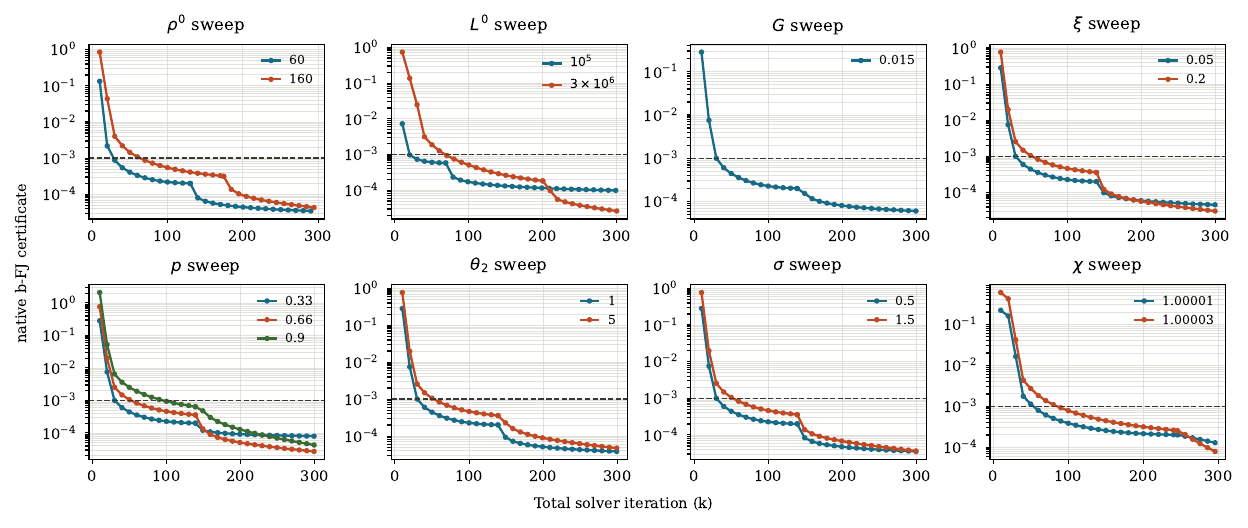}
\caption{Native b-FJ certificate versus total solver iteration (in thousands) for the cloth-8 (top) and capsule-8 (bottom) parameter sweeps. All curves stop at the common count of $3\times10^{5}$ total sweeps used in \prettyref{table:ablation}.}
\label{fig:ablation}
\end{figure*}

Our convergence analysis uses one specific, data-dependent choice of parameters (\prettyref{eq:parameters}). This choice is enough to prove the finite-time guarantees, but it is deliberately conservative. The main runs in \prettyref{appen:experiments} instead use values that we tuned aggressively for speed. This appendix connects the two settings on the $n=8$ instances. We report one run of the exact theorem-sufficient choice per family, and a sensitivity study around the practical baseline. Each ablation changes only one certified parameter, so we can see how that parameter alone affects the certified schedule. \prettyref{table:ablation} lists the certificates of the sensitivity study, and \prettyref{fig:ablation} plots the native b-FJ certificate against the total iteration count in each sweep. To compare runs like for like, every table row and every curve stops at a common count of $3\times10^{5}$ total sweeps.
\rev{Using practical rather than theorem-sufficient parameters follows common practice for first-order nonconvex methods: IPDS-ADMM fixes its continuation and proximal constants instead of evaluating its condition-number-dependent rule~\cite{yuan2025admm}, Jia--Grimmer and ACVI run fixed budgets or heuristic barrier parameters~\cite{jiagrimmer2025,yang2023acvi}, ABIP adds over-relaxation and adaptive penalties beyond its analysis~\cite{lin2021abip}, and BC-ADMM tunes penalties whose theoretical thresholds are difficult to estimate~\cite{pan2026bcadmm}; \citet{maia2026adaptive} further report that theorem-satisfying choices for an earlier proximal ADMM did not reach the target within their iteration limit.}

\paragraph{Theorem-sufficient construction.}
Evaluating \prettyref{eq:parameters} on cloth-8 gives $\rho^0\approx4.7\times10^{7}$, $L^0\approx1.2\times10^{5}$, and $G\approx96$, compared with the practical $\rho^0=40$, $L^0=10^{6}$, $G=6\times10^{-3}$; capsule-8 is even larger ($\rho^0\approx1.6\times10^{9}$). This run still executes and remains numerically feasible (zero original-space violation, positive split margin), but its native b-FJ certificate levels off near $7\times10^{-2}$ (cloth) and $4\times10^{-2}$ (capsule), with w-KKT and KKT of $5\times10^{-2}$ (cloth) and $7.9\times10^{-4}$ and $1\times10^{-2}$ (capsule). Its b-FJ is about three orders of magnitude looser than the practical baselines (\prettyref{table:ablation}) and misses the $10^{-3}$ target in both families.
This agrees with the analysis rather than contradicting it. The certificate the solver reports is exactly $\CritbFJ$ from \prettyref{eq:crit}, and \prettyref{thm:rate} only bounds $\min_{j\leq K}\CritbFJ^{j}\leq C\,K^{-1/4}$ for a fixed constant $C$ that does not depend on $K$. This run uses $K\approx2.9\times10^{5}$ sweeps in total, so $K^{-1/4}\approx0.043$, and the level $\approx6.8\times10^{-2}$ at which the certificate settles under the theorem-sufficient schedule gives $C\approx1.6$. So reaching a tolerance $\epsilon$ needs about $K\approx(C/\epsilon)^{4}$ iterations. The constant is large here for a structural reason. At the conservative $G$, the wall force on every near-active block saturates at $-\sqrt{2L^kG}$, which is far above the true multipliers, so the normalized stationarity term $\alpha^{k}\lVert\Xi^{k}\rVert$ of \prettyref{eq:crit} stays near a floor set by the geometry. The $G$ sweep in \prettyref{table:ablation} shows this directly: shrinking $G$ tightens the native b-FJ (cloth $G{=}2\times10^{-3}\!\to1.5\times10^{-5}$ versus $G{=}1.8\times10^{-2}\!\to4.8\times10^{-5}$). In short, the conservative $G$ inflates the rate constant, but it does not contradict the $\mathcal O(K^{-1/4})$ rate.

\paragraph{Practical sensitivity.}
Around the practical baseline the certified phase is robust. The exponent trades the certificates off: in both families $p=0.66$ gives the smallest b-FJ, below the baseline $p=0.6$ (cloth) and $p=0.5$ (capsule); $p=0.33$ gives the smallest w-KKT but two to four times the baseline b-FJ and KKT; and $p=0.9$ is worse on all three. The relaxation $\sigma$ leaves every certificate at $3\times10^{5}$ sweeps largely unchanged. Changes in $\rho^0$ and in $\theta_2$ (up to $500\times$ the baseline) move b-FJ and KKT by less than a factor of three (w-KKT by at most four). A larger $\chi$ lowers b-FJ (cloth $1.4\times10^{-4}\!\to3.6\times10^{-5}$) but raises w-KKT. The most sensitive knob is $L^0$: tenfold changes make b-FJ worse (cloth $L^0=10^{7}\!\to2.0\times10^{-4}$, capsule $L^0=10^{5}\!\to9.8\times10^{-5}$), while the threefold capsule increase to $3\times10^{6}$ does not ($2.5\times10^{-5}$ versus the baseline $3.5\times10^{-5}$ at the same $3\times10^{5}$ sweeps). Every completed ablation meets the b-FJ target with valid feasibility and physical checks (\prettyref{table:ablation}).
The capsule family adds one more check: it keeps every split state inside the compact domain that its certified contact-curvature bounds rely on. Two perturbations---reducing $G$ to $G/3$ and halving $\rho^0$ to $40$---push the iterates just outside that domain (margin $\approx-10^{-6}$) and abort. We therefore leave these two out of the capsule sweep; all other capsule perturbations in \prettyref{table:ablation}, including the threefold increase of $L^0$, complete normally. Cloth has no such guard and tolerates both directions.

\begin{table*}[t]
\centering\scriptsize
\caption{Size-8 parameter sensitivity.  Each ablation changes one certified parameter.  Every row, including the canonical Ordinary and Ckpt baselines, reports the last recorded iterate at or before a common budget of $300{,}000$ total sweeps, the same cutoff used in \prettyref{fig:ablation}.  native b-FJ is the solver's banded Fritz--John certificate; w-KKT and KKT are the reconstructed common residuals.  Every run meets the $10^{-3}$ b-FJ target with valid feasibility and physical checks.}
\label{table:ablation}
\setlength{\tabcolsep}{3pt}
\begin{minipage}[t]{0.49\textwidth}\centering
\textbf{Cloth resting ($n=8$)}\\[2pt]
\begin{tabular}{@{}lllll@{}}
\toprule
Parameter & Value & native b-FJ & w-KKT & KKT \\
\midrule
Ordinary & -- & $2.59\!\times\!10^{-5}$ & $4.56\!\times\!10^{-4}$ & $9.07\!\times\!10^{-4}$ \\
Ckpt & -- & $1.39\!\times\!10^{-4}$ & $1.40\!\times\!10^{-4}$ & $5.46\!\times\!10^{-3}$ \\
\addlinespace
$\rho^0$ & 20 & $2.50\!\times\!10^{-5}$ & $3.91\!\times\!10^{-4}$ & $9.07\!\times\!10^{-4}$ \\
 & 80 & $2.78\!\times\!10^{-5}$ & $6.12\!\times\!10^{-4}$ & $9.07\!\times\!10^{-4}$ \\
$L^0$ & $10^{5}$ & $7.30\!\times\!10^{-5}$ & $5.34\!\times\!10^{-4}$ & $2.82\!\times\!10^{-3}$ \\
 & $10^{7}$ & $2.02\!\times\!10^{-4}$ & $1.78\!\times\!10^{-2}$ & $1.78\!\times\!10^{-2}$ \\
$G$ & 0.002 & $1.50\!\times\!10^{-5}$ & $2.98\!\times\!10^{-4}$ & $5.13\!\times\!10^{-4}$ \\
 & 0.018 & $4.75\!\times\!10^{-5}$ & $9.40\!\times\!10^{-4}$ & $1.59\!\times\!10^{-3}$ \\
$\xi$ & 0.05 & $3.28\!\times\!10^{-5}$ & $3.56\!\times\!10^{-4}$ & $1.26\!\times\!10^{-3}$ \\
 & 0.2 & $2.24\!\times\!10^{-5}$ & $6.75\!\times\!10^{-4}$ & $6.75\!\times\!10^{-4}$ \\
$p$ & 0.33 & $8.68\!\times\!10^{-5}$ & $1.05\!\times\!10^{-4}$ & $3.60\!\times\!10^{-3}$ \\
 & 0.66 & $2.24\!\times\!10^{-5}$ & $6.70\!\times\!10^{-4}$ & $6.70\!\times\!10^{-4}$ \\
 & 0.9 & $4.97\!\times\!10^{-5}$ & $4.11\!\times\!10^{-3}$ & $4.11\!\times\!10^{-3}$ \\
$\theta_2$ & 1 & $2.74\!\times\!10^{-5}$ & $5.78\!\times\!10^{-4}$ & $9.07\!\times\!10^{-4}$ \\
 & 5 & $3.61\!\times\!10^{-5}$ & $1.31\!\times\!10^{-3}$ & $1.31\!\times\!10^{-3}$ \\
$\sigma$ & 0.5 & $2.59\!\times\!10^{-5}$ & $4.56\!\times\!10^{-4}$ & $9.07\!\times\!10^{-4}$ \\
 & 1.5 & $2.59\!\times\!10^{-5}$ & $4.56\!\times\!10^{-4}$ & $9.07\!\times\!10^{-4}$ \\
$\chi$ & 1.00001 & $7.15\!\times\!10^{-5}$ & $5.79\!\times\!10^{-4}$ & $2.71\!\times\!10^{-3}$ \\
 & 1.00003 & $3.61\!\times\!10^{-5}$ & $1.97\!\times\!10^{-3}$ & $1.97\!\times\!10^{-3}$ \\
\bottomrule
\end{tabular}
\end{minipage}\hfill
\begin{minipage}[t]{0.49\textwidth}\centering
\textbf{Capsule stacking ($n=8$)}\\[2pt]
\begin{tabular}{@{}lllll@{}}
\toprule
Parameter & Value & native b-FJ & w-KKT & KKT \\
\midrule
Ordinary & -- & $3.55\!\times\!10^{-5}$ & $4.25\!\times\!10^{-4}$ & $1.63\!\times\!10^{-3}$ \\
Ckpt & -- & $1.70\!\times\!10^{-4}$ & $1.97\!\times\!10^{-4}$ & $8.60\!\times\!10^{-3}$ \\
\addlinespace
$\rho^0$ & 60 & $3.47\!\times\!10^{-5}$ & $3.47\!\times\!10^{-4}$ & $1.63\!\times\!10^{-3}$ \\
 & 160 & $4.19\!\times\!10^{-5}$ & $7.92\!\times\!10^{-4}$ & $1.75\!\times\!10^{-3}$ \\
$L^0$ & $10^{5}$ & $9.85\!\times\!10^{-5}$ & $9.48\!\times\!10^{-4}$ & $4.63\!\times\!10^{-3}$ \\
 & $3\times10^{6}$ & $2.53\!\times\!10^{-5}$ & $5.04\!\times\!10^{-4}$ & $1.06\!\times\!10^{-3}$ \\
$G$ & 0.015 & $6.05\!\times\!10^{-5}$ & $5.29\!\times\!10^{-4}$ & $2.84\!\times\!10^{-3}$ \\
$\xi$ & 0.05 & $4.58\!\times\!10^{-5}$ & $2.21\!\times\!10^{-4}$ & $2.27\!\times\!10^{-3}$ \\
 & 0.2 & $2.98\!\times\!10^{-5}$ & $7.68\!\times\!10^{-4}$ & $1.16\!\times\!10^{-3}$ \\
$p$ & 0.33 & $7.97\!\times\!10^{-5}$ & $5.77\!\times\!10^{-5}$ & $4.13\!\times\!10^{-3}$ \\
 & 0.66 & $2.62\!\times\!10^{-5}$ & $1.36\!\times\!10^{-3}$ & $1.36\!\times\!10^{-3}$ \\
 & 0.9 & $4.00\!\times\!10^{-5}$ & $2.07\!\times\!10^{-3}$ & $2.07\!\times\!10^{-3}$ \\
$\theta_2$ & 1 & $3.84\!\times\!10^{-5}$ & $6.91\!\times\!10^{-4}$ & $1.63\!\times\!10^{-3}$ \\
 & 5 & $4.57\!\times\!10^{-5}$ & $1.21\!\times\!10^{-3}$ & $1.63\!\times\!10^{-3}$ \\
$\sigma$ & 0.5 & $3.55\!\times\!10^{-5}$ & $4.26\!\times\!10^{-4}$ & $1.63\!\times\!10^{-3}$ \\
 & 1.5 & $3.55\!\times\!10^{-5}$ & $4.25\!\times\!10^{-4}$ & $1.63\!\times\!10^{-3}$ \\
$\chi$ & 1.00001 & $1.29\!\times\!10^{-4}$ & $2.96\!\times\!10^{-4}$ & $6.16\!\times\!10^{-3}$ \\
 & 1.00003 & $7.65\!\times\!10^{-5}$ & $8.89\!\times\!10^{-4}$ & $2.32\!\times\!10^{-3}$ \\
\bottomrule
\end{tabular}
\end{minipage}
\end{table*}

\section{Parallel Wall-Clock Scaling}
\label{appen:scaling}
The main experiments are small and run on a single thread, so we use this appendix to measure real wall-clock time on the cloth-resting testbed of \prettyref{appen:experiments}: where each iteration spends its time, and how much faster it runs when we add CPU threads.  We do not solve to convergence here---we fix the penalty schedule, run a fixed number of sweeps, and split each sweep into the four steps of \prettyref{alg:IP-ADMM}: the local objective update, the constraint update, the global $x$-solve, and the rest (residuals, movement, and dual updates). In the cloth resting problem, we grow the cloth chain over $n\in\{16,24,32,48,64\}$ nodes ($2n$ variables, $3n-3$ objective blocks, and $(n-1)(n-2)+n$ collision constraints), going well past the $n\le24$ sizes used for the convergence study. For each size and each thread count in $\{1,2,4,8,16\}$ we set \texttt{torch.set\_num\_threads}, and run $300$ sweeps with \texttt{perf\_counter}. Each point runs in its own process, and we report the fastest of nine repeats.

\begin{figure*}[t]
\centering
\includegraphics[width=0.9\textwidth]{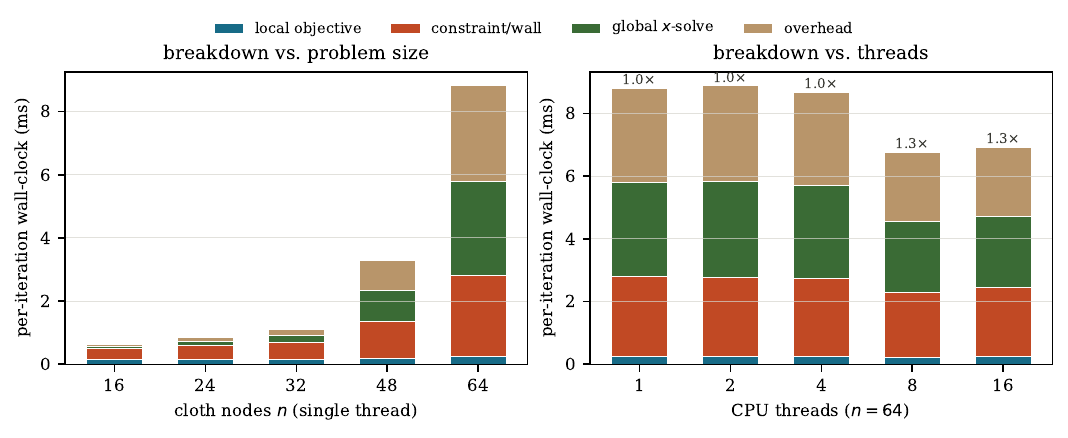}
\caption{Per-iteration wall-clock time split into the four sweep steps. Left: versus cloth size $n$ on a single thread; the batched local objective update stays small while the constraint update, global $x$-solve, and overhead grow. Right: versus CPU thread count at $n=64$, where adding threads cuts the per-iteration time by up to $1.3\times$ (annotated).}
\label{fig:scaling}
\end{figure*}

\prettyref{fig:scaling} (left) shows how one iteration divides among the four steps on a single thread. The local objective update is applied to all blocks at once, so it stays cheap: it falls from about a quarter of the time at $n=16$ to a thin sliver at $n=64$. As the problem grows, the constraint update, the global $x$-solve, and the overhead take up most of the time.
\prettyref{fig:scaling} (right) shows the same split at $n=64$ as we add CPU threads.  At $n=64$ the batched linear algebra is large enough to spread across cores, and the constraint, global, and overhead steps all get faster, cutting the per-iteration time by up to $1.3\times$ (at $8$ threads). This is a real speedup on a single machine, and it supports the parallelizability claim with wall-clock numbers rather than structure alone: at $n=64$, simply adding threads already makes each iteration faster.

\section{Explicit Parameter Constructions}
\label{appen:parameters}
This appendix records the sufficient data-dependent choices referenced by the main convergence statements.  The data aggregates in the ordinary display are shared by both modes, while the algorithmic tunings and wall-height constructions are mode-specific.

\subsection{Ordinary IP-ADMM}
For $p\in(0,1)$ define
\begin{eqaligned}
\label{eq:Hp}
H_p(\xi)=\xi^2\sum_{j=1}^\infty\frac{j^{2p-2}}{(1+\xi j^p)^2},
\end{eqaligned}
and set $H_0(\xi)=0$ by convention.
\begin{equation}
\label{eq:parameters}
\left\{
\begin{aligned}
H_{\max}&=\max_{i\in\IO}H_i,\quad
L_{\max}=\max_{i\in\IC}L_i,\quad M_{\max}=\max_{i\in\IC}M_i,
\quad\bar M=\sum_{i\in\IC}M_i^2,\\
S_0&=\sum_{i\in\IO}\varphi_i(A_ix^1),\quad
\widehat G_0=\sum_{i\in\IO}\widehat C_i^2,\quad
G_\psi=\sum_{i\in\IO}(C_i^\psi)^2,\\
M_F&=\|A\|_2\sqrt{\widehat G_0},\quad
J_H=\left(\sum_{i\in\IO}(H_iC_i^\psi)^2\right)^{1/2},\\
\sigma&=1,\quad
\theta_1=\frac{121909}{121000},\quad
\theta_2=25,\quad
\eta=1209,\quad d_0=\frac12,\\
L^0&=\max\left\{1,\frac1m,\frac{12S_0}{m^2}\right\},\qquad
G=\frac{m^2L^0}2,\\
\rho^0&=\max\left\{
\begin{gathered}
1,1209H_{\max},1209mL^0L_{\max},1209L^0M_{\max}^2,\\
\displaystyle 6\left(\frac\eta2+\frac{11}2\right)\frac{\widehat G_0}G,\\
\displaystyle 6p^2H_p(\xi)
\frac{\frac{25}2G_\psi+\frac{25}4GL^0\bar M}G
\end{gathered}
\right\}.
\end{aligned}
\right.
\end{equation}
Here $S_0\geq0$ by~\prettyref{ass:obj}, and $G>0$ because $m>0$.

\subsection{IP-ADMM-Ckpt}
\begingroup\small
\begin{equation}
\label{eq:checkpoint-parameters}
\left\{
\begin{aligned}
\theta_1^{\rm cp}&=\frac{101}{100},\quad
\theta_2^{\rm cp}=\frac1{100},\quad
\eta^{\rm cp}=\frac1{100},\\
G^{\rm cp}&=\frac12\max\{m^2,m,3S_0\},\\
\Delta_G&=\frac{G^{\rm cp}-S_0}{4},\quad
\overline V=S_0+2\Delta_G,\quad a_0=\frac{51}{100},\\
b_{\rm cp}&=\sqrt{\frac{2G^{\rm cp}}{\ell_{\rm cp}}},\quad
r_{\rm cp}=\max\left\{H_{\max},\sqrt{2G^{\rm cp}\ell_{\rm cp}}L_{\max}\right\},\\
c_{\upsilon,{\rm cp}}&=\theta_1^{\rm cp}(1+r_{\rm cp})-1,\quad
c_{\rm mov}=\sqrt{(r_{\rm cp}+c_{\upsilon,{\rm cp}})^2+(1+\theta_2^{\rm cp})^2},\\
S_\Xi&=\|A\|_2(c_{\rm mov}\tau_{\rm cp}+\eta^{\rm cp}J_H),\quad
I_{\rm cp}=\frac{M_F+S_\Xi}{\gamma_{\rm cq}},\\
Q_{\rm cp}&=\sqrt{\widehat G_0}+M_{\max}I_{\rm cp}+(c_{\upsilon,{\rm cp}}+1)\tau_{\rm cp},\\
D_{\rm cp}&=\frac{(Q_{\rm cp}+\sqrt{\widehat G_0})^2}{2},\\
\rho_{\rm cp}^0&=\max\left\{
\begin{gathered}
1,\frac{\eta^{\rm cp}G_\psi}{2\Delta_G},\frac{D_{\rm cp}}{\Delta_G},
\frac{b_{\rm cp}}m,\frac{2b_{\rm cp}}{\delta_{\rm cq}},
\frac{8\overline V}{\delta_{\rm cq}^2},\\
\frac{2(\chi+1)Q_{\rm cp}^2}{(1-\chi^{-1})\Delta_G},
\frac{4(Q_{\rm cp}+\sqrt{\widehat G_0})}{\delta_{\rm cq}},
\left(\frac{L_{\max}\tau_{\rm cp}^2}{\delta_{\rm cq}}\right)^{1/4}
\end{gathered}\right\}.
\end{aligned}\right.
\end{equation}
\endgroup
Since $S_0\geq0$ and $m>0$, $G^{\rm cp}\geq\max\{m^2,3S_0\}/2>S_0$, so $\Delta_G>0$.

\paragraph{Geometric checkpoint schedule.} If $j\geq0$ checkpoints have been committed, the level parameters are
\begin{eqaligned}
\label{eq:checkpoint-schedule}
\rho_j&\defeq\rho_{\rm cp}^0\chi^j,\qquad
L_j\defeq\ell_{\rm cp}\rho_j^2,\qquad
\mu_j\defeq\frac{\eta^{\rm cp}}{\rho_j},\qquad
\beta_j\defeq\frac{b_{\rm cp}}{\rho_j}.
\end{eqaligned}
The algorithm initializes $(\rho^1,L^1,\mu^1)\defeq(\rho_0,L_0,\mu_0)$.  The subscript $j$ indexes checkpoint levels and the superscript $k$ indexes outer iterations: iteration $k$ uses $(\rho_j,L_j,\mu_j)$ when it begins after exactly $j$ committed checkpoints, and an accepted checkpoint records its tuple with these pre-promotion values before promotion.

\section{Algebraic Lemmas}
\label{appen:algebraic}
We collect the componentwise envelope and wall calculus used below.  The objective assumptions immediately give the global estimates
\begin{eqaligned}
\label{eq:objective-derived-bounds}
\|\nabla h_i(y)\|\leq D_i,\qquad \|w\|\leq C_i^\psi\quad(w\in\partial\psi_i(y)),\qquad |\varphi_i(u)-\varphi_i(v)|\leq\widehat C_i\|u-v\|.
\end{eqaligned}
Indeed, differentiability and global Lipschitz continuity bound $\nabla h_i$, and every subgradient of the finite convex Lipschitz function $\psi_i$ has norm at most $C_i^\psi$.  Adding the two Lipschitz bounds gives the last estimate, with the main-text definition $\widehat C_i=D_i+C_i^\psi$.  The Moreau gradient bound below also yields $\|\nabla\bar f_i(y,L,\mu)\|\leq D_i+C_i^\psi=\widehat C_i$ for $i\in\IO$, so the objective surrogate is globally $\widehat C_i$-Lipschitz, uniformly in $(L,\mu)$.

\begin{lemma}[Convex Moreau calculus]
\label{lem:moreau-basic}
For $i\in\IO$, $\mu>0$, and every $y$, $P_{i,\mu}^\psi(y)$ is single valued, $(\psi_i)_\mu$ is convex and continuously differentiable, and
\begin{eqaligned}
\label{eq:moreau-grad}
\nabla(\psi_i)_\mu(y)
=\frac{y-P_{i,\mu}^\psi(y)}{\mu}
\in\partial\psi_i(P_{i,\mu}^\psi(y)).
\end{eqaligned}
Furthermore,
\begin{eqaligned}
\label{eq:moreau-bound}
\|\nabla(\psi_i)_\mu(y)\|\leq C_i^\psi,\qquad\|y-P_{i,\mu}^\psi(y)\|\leq\mu C_i^\psi,
\end{eqaligned}
and
\begin{eqaligned}
\label{eq:moreau-gap}
0\leq\psi_i(y)-(\psi_i)_\mu(y)\leq\frac\mu2(C_i^\psi)^2.
\end{eqaligned}
\end{lemma}
\begin{proof}
Introduce the minimization variable explicitly and write
\[
p=P_{i,\mu}^\psi(y)
=\argmin{r}\left\{\psi_i(r)+\frac1{2\mu}\|r-y\|^2\right\}.
\]
The objective is $1/\mu$-strongly convex and coercive, so $p$ is unique, and its optimality condition gives $(y-p)/\mu\in\partial\psi_i(p)$.  The envelope $(\psi_i)_\mu$ is convex as the partial minimization over $r$ of the jointly convex function $(r,y)\mapsto\psi_i(r)+\|r-y\|^2/(2\mu)$.  If $p_d=P_{i,\mu}^\psi(y+d)$, comparison at the two minimizers gives
\[
\left\langle\frac{y+d-p_d}{\mu},d\right\rangle-\frac{\|d\|^2}{2\mu}
\leq(\psi_i)_\mu(y+d)-(\psi_i)_\mu(y)
\leq\left\langle\frac{y-p}{\mu},d\right\rangle+\frac{\|d\|^2}{2\mu}.
\]
Monotonicity of $\partial\psi_i$ applied to the two optimality conditions gives $\|p_d-p\|^2\leq\langle p_d-p,d\rangle$, hence $p_d\to p$ as $d\to0$.  Dividing the comparison inequalities by $\|d\|$ proves differentiability and the gradient formula; the same proximal Lipschitz bound makes the gradient continuous.  Every subgradient of a finite convex $C_i^\psi$-Lipschitz function has norm at most $C_i^\psi$, proving the two bounds.  Finally, evaluating the proximal objective at $r=y$ gives $(\psi_i)_\mu(y)\leq\psi_i(y)$, while $\psi_i(p)\geq\psi_i(y)-C_i^\psi\|p-y\|$.  With the scalar variable $r=\|p-y\|\geq0$, minimizing $-C_i^\psi r+r^2/(2\mu)$ gives the right gap inequality.
\end{proof}

\begin{lemma}[Moreau gradient and parameter drift]
For $i\in\IO$, every $\mu>0$, and all $u,v,y$,
\begin{eqaligned}
\label{eq:env-smooth}
\|\nabla(\psi_i)_\mu(u)-\nabla(\psi_i)_\mu(v)\|\leq\frac1\mu\|u-v\|.
\end{eqaligned}
For $0<\mu_2\leq\mu_1$,
\begin{eqaligned}
\label{eq:env-value-drift}
(\psi_i)_{\mu_2}(y)-(\psi_i)_{\mu_1}(y)\leq\frac{(C_i^\psi)^2}{2}(\mu_1-\mu_2),
\end{eqaligned}
and
\begin{eqaligned}
\label{eq:env-grad-drift}
\|\nabla(\psi_i)_{\mu_2}(y)-\nabla(\psi_i)_{\mu_1}(y)\|\leq C_i^\psi(\mu_1/\mu_2-1).
\end{eqaligned}
\end{lemma}
\begin{proof}
For $P=P_{i,\mu}^\psi$, monotonicity of $\partial\psi_i$ in the optimality conditions at $u$ and $v$ gives
\[
\|P(u)-P(v)\|^2\leq\langle P(u)-P(v),u-v\rangle.
\]
Thus $P$ is firmly nonexpansive.  Expanding this inequality shows that $I-P$ satisfies the same inequality and is therefore nonexpansive; together with~\prettyref{eq:moreau-grad}, this proves~\prettyref{eq:env-smooth}.  Write $p_\alpha=P_{i,\alpha}^\psi(y)$.  For $\nu>\mu$, comparison at $p_\mu$ and $p_\nu$ gives
\[
\left(\frac1{2\nu}-\frac1{2\mu}\right)\|y-p_\nu\|^2
\leq(\psi_i)_\nu(y)-(\psi_i)_\mu(y)
\leq\left(\frac1{2\nu}-\frac1{2\mu}\right)\|y-p_\mu\|^2.
\]
Moreover, monotonicity of $\partial\psi_i$ in the two optimality conditions and~\prettyref{eq:moreau-bound} give $\|p_\nu-p_\mu\|\leq C_i^\psi|\nu-\mu|$.  Dividing the comparison inequalities by $\nu-\mu$ and taking either one-sided limit therefore yields
\[
\begin{aligned}
\frac{d}{d\mu}(\psi_i)_\mu(y)
=-\frac{\|y-P_{i,\mu}^\psi(y)\|^2}{2\mu^2}
=-\frac12\|\nabla(\psi_i)_\mu(y)\|^2.
\end{aligned}
\]
The right side is continuous on $(0,\infty)$.  Thus, integrating over $[\mu_2,\mu_1]$ and applying~\prettyref{eq:moreau-bound} yields
\[
\begin{aligned}
(\psi_i)_{\mu_2}(y)-(\psi_i)_{\mu_1}(y)
=\int_{\mu_2}^{\mu_1}\frac12\|\nabla(\psi_i)_\mu(y)\|^2\,d\mu
\leq\frac{(C_i^\psi)^2}{2}(\mu_1-\mu_2).
\end{aligned}
\]
For the gradient estimate put $g_j=\nabla(\psi_i)_{\mu_j}(y)$ and $p_j=y-\mu_jg_j$.  Monotonicity of $\partial\psi_i$ gives $0\leq\langle g_2-g_1,p_2-p_1\rangle$, hence $\mu_2\|g_2-g_1\|^2\leq (\mu_1-\mu_2)\langle g_2-g_1,g_1\rangle$.  If $g_2=g_1$ the claim is immediate; otherwise divide by $\mu_2\|g_2-g_1\|$ and use $\|g_1\|\leq C_i^\psi$.
\end{proof}

\begin{lemma}[One-call envelope prox]
\label{lem:one-call-envelope-prox}
For $i\in\IO$, $t,\mu>0$, and every $z$, let $y^+\defeq\prox_{t(\psi_i)_\mu}(z)$ and $p\defeq\prox_{(\mu+t)\psi_i}(z)$.  Then
\[
y^+=\frac\mu{\mu+t}z+\frac t{\mu+t}p,
\]
and
\begin{eqaligned}
\label{eq:one-call-outputs}
p&=P_{i,\mu}^\psi(y^+),\\
\nabla(\psi_i)_\mu(y^+)&=\frac{z-p}{\mu+t}
=\frac{y^+-p}{\mu}
\in\partial\psi_i(p).
\end{eqaligned}
In particular, with $t_i^k=1/m_i^k$ and $z_{o,i}^k$ defined in the main text after the block models~\prettyref{eq:block-models}, the objective-block quantities satisfy the exact one-call identities
\begin{eqaligned}
\label{eq:prox-oracle}
p_i^{k+1}&=\prox_{(\mu^k+t_i^k)\psi_i}(z_{o,i}^k),\\
y_i^{k+1}&=\frac{\mu^k}{\mu^k+t_i^k}z_{o,i}^k
+\frac{t_i^k}{\mu^k+t_i^k}p_i^{k+1},\\
w_i^{k+1}&=\frac{z_{o,i}^k-p_i^{k+1}}{\mu^k+t_i^k}.
\end{eqaligned}
\end{lemma}
\begin{proof}
Let $\hat y\defeq\mu z/(\mu+t)+tp/(\mu+t)$.  The optimality condition for $p$ is $(z-p)/(\mu+t)\in\partial\psi_i(p)$.  The affine combination gives both $\hat y-p=\mu(z-p)/(\mu+t)$ and $z-\hat y=t(z-p)/(\mu+t)$.  The first identity and $(z-p)/(\mu+t)\in\partial\psi_i(p)$ give $(\hat y-p)/\mu\in\partial\psi_i(p)$, hence $p=\prox_{\mu\psi_i}(\hat y)$.  By~\prettyref{eq:moreau-grad}, the second identity then gives $0=\nabla(\psi_i)_\mu(\hat y)+(\hat y-z)/t$.  Since $(\psi_i)_\mu$ is convex by~\prettyref{lem:moreau-basic}, this condition is necessary and sufficient for $\hat y$ to be the unique minimizer $y^+$ of $(\psi_i)_\mu(u)+\|u-z\|^2/(2t)$.  Hence $\hat y=y^+$, which gives the displayed identity and~\prettyref{eq:one-call-outputs}.

For the main-text objective-block subproblem, completing the square and dropping terms independent of $u$ gives
\[
(\psi_i)_{\mu^k}(u)+\frac1{2t_i^k}\|u-z_{o,i}^k\|^2.
\]
Thus its minimizer is $\prox_{t_i^k(\psi_i)_{\mu^k}}(z_{o,i}^k)$, and substituting $(t,\mu,z)=(t_i^k,\mu^k,z_{o,i}^k)$ gives~\prettyref{eq:prox-oracle}.
\end{proof}

\paragraph{Original-problem reconstruction.}
For the original-problem candidate certificate tuple generated by step $k$, let $p_i^{k+1}$ and $w_i^{k+1}$ be the one-call outputs in~\prettyref{eq:prox-oracle}, and define
\begin{eqaligned}
\label{eq:reconstruction}
\breve y_i^{k+1}&\defeq
\begin{cases}p_i^{k+1},&i\in\IO,\\y_i^{k+1},&i\in\IC,\end{cases}&
v_i^{k+1}&\defeq\nabla h_i(p_i^{k+1})+w_i^{k+1}\quad(i\in\IO).
\end{eqaligned}
Then
\begin{eqaligned}
\label{eq:recon-error}
v_i^{k+1}&\in\partial\varphi_i(\breve y_i^{k+1}),
\qquad \|y^{k+1}-\breve y^{k+1}\|\leq\mu^k C_\psi,\\
C_\psi&\defeq\left(\sum_{i\in\IO}(C_i^\psi)^2\right)^{1/2},
\qquad \|v_i^{k+1}-q_i^k\|\leq H_i(\|d_i^k\|+\mu^kC_i^\psi).
\end{eqaligned}
Indeed,~\prettyref{eq:one-call-outputs} gives $w_i^{k+1}\in\partial\psi_i(p_i^{k+1})$ and $y_i^{k+1}-p_i^{k+1}=\mu^kw_i^{k+1}$.  The exact smooth-plus-convex sum rule gives the subgradient membership, the bound $\|w_i^{k+1}\|\leq C_i^\psi$ gives the aggregate reconstruction error, and the $H_i$-Lipschitz continuity of $\nabla h_i$ gives the final force-transfer bound.

\begin{lemma}[Moreau continuation corrections]
\label{lem:Rpsi-discharge}
Recall $G_\psi$ from~\prettyref{eq:parameters} and define $R_\psi(\mu)\defeq\mu G_\psi/2$.  For $0<\mu'\leq\mu$, arbitrary $L,L'>0$, and every $y=(y_i)_{i\in\IO}$,
\begin{eqaligned}
\label{eq:Rpsi-monotone}
\sum_{i\in\IO}\bar f_i(y_i,L',\mu')+R_\psi(\mu')
\leq\sum_{i\in\IO}\bar f_i(y_i,L,\mu)+R_\psi(\mu).
\end{eqaligned}
Moreover, for every $L,\mu>0$ and $x$,
\begin{eqaligned}
\label{eq:Rpsi-lower}
\sum_{i\in\IO}\bar f_i(A_ix,L,\mu)+R_\psi(\mu)\geq0.
\end{eqaligned}
For the ordinary IP schedule define $\mu_{\rm IP}(L)\defeq\eta L^0/(\rho^0L)$ and $R(L)\defeq R_\psi(\mu_{\rm IP}(L))$.  Then, for $L'\geq L>0$ and every $y=(y_i)_{i\in\IO}$,
\begin{eqaligned}
\label{eq:R2}
\sum_{i\in\IO}\bar f_i(y_i,L',\mu_{\rm IP}(L'))+R(L')\leq\sum_{i\in\IO}\bar f_i(y_i,L,\mu_{\rm IP}(L))+R(L),
\end{eqaligned}
and
\begin{eqaligned}
\label{eq:R3}
\sum_{i\in\IO}\bar f_i(A_ix,L,\mu_{\rm IP}(L))+R(L)\geq0\qquad\forall x.
\end{eqaligned}
\end{lemma}
\begin{proof}
For $i\in\IO$, the branch $\bar f_i(\cdot,L,\mu)=h_i+(\psi_i)_\mu$ is independent of $L$.  Hence the $L$-argument is inert on $\IO$, the $h_i$ terms cancel, and only the envelopes change.  Applying~\prettyref{eq:env-value-drift} blockwise shows that their total increase when $\mu$ decreases to $\mu'$ is at most $(\mu-\mu')G_\psi/2=R_\psi(\mu)-R_\psi(\mu')$, proving~\prettyref{eq:Rpsi-monotone}.  By~\prettyref{eq:moreau-gap}, for every $i\in\IO$ and $x$,
\[
h_i(A_ix)+(\psi_i)_\mu(A_ix)+\frac\mu2(C_i^\psi)^2
\geq\varphi_i(A_ix).
\]
Sum over $i\in\IO$ and use the aggregate lower bound in~\prettyref{ass:obj} to obtain~\prettyref{eq:Rpsi-lower}.  Since $\mu_{\rm IP}$ is decreasing, $L'\geq L$ gives $\mu_{\rm IP}(L')\leq\mu_{\rm IP}(L)$; hence~\prettyref{eq:R2} and~\prettyref{eq:R3} are~\prettyref{eq:Rpsi-monotone} and~\prettyref{eq:Rpsi-lower} with $(\mu,\mu')=(\mu_{\rm IP}(L),\mu_{\rm IP}(L'))$.
\end{proof}

\begin{lemma}[Wall calculus, model error, and exact affine prox]
\label{lem:added_algebraic_3}
Let $\mathsf m$ be the active mode, with wall $g_{\mathsf m}$ and height $G_{\mathsf m}$ from the main text; $G_{\mathsf m}>0$ in both modes by~\prettyref{appen:parameters}.  Recall the activation width
\begin{eqaligned}
\label{eq:band}
\beta(L)=\sqrt{2G_{\mathsf m}/L}.
\end{eqaligned}
For every $L>0$, the wall function $g_{\mathsf m}(\cdot,L)$ is nonnegative, continuously differentiable, convex, and $L$-smooth, with explicit force
\begin{eqaligned}
\label{eq:band-force}
g_{\mathsf m}'(s,L)=\max\{\min\{Ls-\sqrt{2LG_{\mathsf m}},0\},-\sqrt{2LG_{\mathsf m}}\}.
\end{eqaligned}
In particular, $|g_{\mathsf m}'(s,L)|\leq\sqrt{2LG_{\mathsf m}}$ and $g_{\mathsf m}(s,L)>G_{\mathsf m}$ whenever $s<0$. For all $L,L'>0$ and $s\in\R$,
\begin{eqaligned}
\label{eq:wall-parameter-drift}
|g_{\mathsf m}'(s,L)-g_{\mathsf m}'(s,L')|\leq\sqrt{2G_{\mathsf m}}\,|\sqrt L-\sqrt{L'}|.
\end{eqaligned}
For $i\in\IC$, every iteration index $k$, and every $u$,
\begin{eqaligned}
|g_{\mathsf m}(f_i(u),L^k)-g_{\mathsf m}(\ell_i^k(u),L^k)|\leq\frac{\sqrt{2L^kG_{\mathsf m}}\,L_i}{2}\|u-y_i^k\|^2.
\end{eqaligned}
Moreover, consider
\[
\min_u\left\{g_{\mathsf m}(c+a^Tu,L)+\frac Q2\|u-z\|^2\right\},\qquad Q>0.
\]
Put $s_z=c+a^Tz$, $a_L=\sqrt{2LG_{\mathsf m}}$, and $\kappa=\|a\|^2/Q$.  Its unique minimizer is
\begin{eqaligned}
\label{eq:affine-prox}
u^+=z-\frac\iota Q a,
\qquad
\iota=
\begin{cases}
-a_L,&s_z<-\kappa a_L,\\
\displaystyle L\frac{s_z-\beta(L)}{1+\kappa L},
&-\kappa a_L\leq s_z<\beta(L),\\
0,&s_z\geq\beta(L).
\end{cases}
\end{eqaligned}
The formulas agree at both boundaries, and $\iota=g_{\mathsf m}'(c+a^Tu^+,L)$.
\end{lemma}
\begin{proof}
The definition of $g_{\mathsf m}$, together with $G_{\mathsf m},L>0$, directly gives nonnegativity and $g_{\mathsf m}(s,L)>G_{\mathsf m}$ for $s<0$.  At $s=0$, the negative and middle branches have the same value and one-sided derivative; within the $s\geq0$ formula, the middle and zero branches likewise match in value and derivative at $s=\beta(L)$.  These branchwise derivatives and matching boundary values prove continuous differentiability and give~\prettyref{eq:band-force}.  The derivative is nondecreasing and $L$-Lipschitz in $s$, proving convexity and $L$-smoothness; its explicit form also gives $|g_{\mathsf m}'(s,L)|\leq\sqrt{2LG_{\mathsf m}}$.

It remains to prove the uniform parameter bound.  Put $t=\sqrt L$ and $c_G=\sqrt{2G_{\mathsf m}}$.  For fixed $s$,
\[
\begin{aligned}
g_{\mathsf m}'(s,t^2)=
\begin{cases}
-c_Gt,&s\leq0,\\
t^2s-c_Gt,&s>0,\ 0<t<c_G/s,\\
0,&s>0,\ t\geq c_G/s.
\end{cases}
\end{aligned}
\]
This map is continuous and piecewise $C^1$.  On its smooth pieces its derivative is respectively $-c_G$, $2ts-c_G$, and $0$; on the middle piece, $0<ts<c_G$ gives $|2ts-c_G|<c_G$, so every derivative has absolute value at most $c_G$.  Integrating piecewise between $\sqrt L$ and $\sqrt{L'}$ proves~\prettyref{eq:wall-parameter-drift}.  The Lipschitz-gradient remainder bound gives $|f_i(u)-\ell_i^k(u)|\leq L_i\|u-y_i^k\|^2/2$; composing with $|g_{\mathsf m}'|\leq\sqrt{2L^kG_{\mathsf m}}$ proves the model-error bound.

For the prox, the objective is $Q$-strongly convex, so its minimizer is unique; when $a=0$, the formulas reduce directly to $u^+=z$. Optimality has the form $Q(u-z)+\iota a=0$ with $\iota=g_{\mathsf m}'(c+a^Tu,L)$, hence $u=z-\iota a/Q$ and $s=c+a^Tu=s_z-\kappa\iota$.  On the first branch, $\iota=-a_L$ and $s=s_z+\kappa a_L<0$ exactly when $s_z<-\kappa a_L$.  On the middle branch, using $L\beta(L)=a_L$ gives
\[
\iota=L(s-\beta(L)),\qquad
s=\frac{s_z+\kappa a_L}{1+\kappa L},
\qquad
\iota=L\frac{s_z-\beta(L)}{1+\kappa L},
\]
and $0\leq s<\beta(L)$ is equivalent to $-\kappa a_L\leq s_z<\beta(L)$.  On the last branch, $\iota=0$ and $s=s_z\geq\beta(L)$.  At $s_z=-\kappa a_L$ the first two formulas both give $\iota=-a_L$; at $s_z=\beta(L)$ the latter two both give zero.  Thus \prettyref{eq:affine-prox} satisfies the unique optimality condition and $\iota=g_{\mathsf m}'(c+a^Tu^+,L)$.
\end{proof}

The two true/model discrepancies share the notation $\tau_i^k$ from the main text.  If
\begin{eqaligned}
\widetilde f_i^k(u)=
\begin{cases}
h_i(y_i^k)+\langle\nabla h_i(y_i^k),u-y_i^k\rangle+(\psi_i)_{\mu^k}(u),&i\in\IO,\\
g_{\mathsf m}(\ell_i^k(u),L^k),&i\in\IC,
\end{cases}
\end{eqaligned}
then, globally,
\begin{eqaligned}
\label{eq:model-two-sided}
|\bar f_i(u,L^k,\mu^k)-\widetilde f_i^k(u)|
\leq\frac{\tau_i^k}{2}\|u-y_i^k\|^2,
\end{eqaligned}
and $\widetilde f_i^k$ is $\zeta_i^k$-smooth.  Consequently, adding $\langle e_i^k,u-y_i^k\rangle+\mathsf B_i^k\|u-y_i^k\|^2/2$, where $\mathsf B_i^k=\rho^k+\tau_i^k$, to $\widetilde f_i^k$ yields, up to an additive constant, a majorizer of the true block Lagrangian at $x^k$.

\paragraph{Schedule consequences.}
The ordinary IP schedules in~\prettyref{eq:IP-schedules} are $\rho^k=\rho^0(1+\xi k^p)$, $L^k=L^0(1+\xi k^p)$, and $\mu^k=\eta/\rho^k=\mu_{\rm IP}(L^k)$.  Thus, for every $k\geq1$,
\begin{eqaligned}
\label{eq:rho-L-ratio}
\frac{\rho^k}{L^k}=\frac{\rho^0}{L^0}.
\end{eqaligned}
They are nondecreasing:
\begin{eqaligned}
\label{eq:schedule-monotone}
\rho^{k+1}\geq\rho^k,\qquad L^{k+1}\geq L^k.
\end{eqaligned}

\begin{lemma}[IP schedule ratios]
\label{lem:rho-increase}
In ordinary IP mode, for $k\geq1$, $1\leq\rho^{k+1}/\rho^k=L^{k+1}/L^k\leq1+\xi\leq2$.
\end{lemma}
\begin{proof}
The schedules share the factor $1+\xi k^p$, so
\[
\begin{aligned}
\frac{\rho^{k+1}}{\rho^k}=\frac{L^{k+1}}{L^k}
=\frac{1+\xi(k+1)^p}{1+\xi k^p}.
\end{aligned}
\]
For $p\in(0,1)$, the map $t\mapsto t^p$ is nondecreasing and concave on $(0,\infty)$, so the ratio is at least one and, for $k\geq1$,
\[
(k+1)^p-k^p\leq p k^{p-1}\leq1,
\]
where the last inequality uses $p\leq1$ and $k^{p-1}\leq1$.  Therefore, explicitly using $1+\xi k^p\geq1$,
\[
\begin{aligned}
\frac{1+\xi(k+1)^p}{1+\xi k^p}&=1+\frac{\xi((k+1)^p-k^p)}{1+\xi k^p}\\
&\leq1+\xi((k+1)^p-k^p)\\
&\leq1+\xi.
\end{aligned}
\]
When $p=0$, the ratio is exactly one.  Finally, $\xi\leq1$ gives $1+\xi\leq2$.
\end{proof}

\begin{lemma}
Under the ordinary polynomial schedule~\prettyref{eq:IP-schedules}, with $p\in[0,1)$ and $\xi\in(0,1]$, the schedule sums satisfy
\begin{eqaligned}
\label{eq:sequence-sum-a}
\sum_{k=2}^\infty\frac{(\sqrt{L^k}-\sqrt{L^{k-1}})^2}{\rho^k}
\leq\frac{L^0p^2}{4\rho^0}H_p(\xi),
\end{eqaligned}
and
\begin{eqaligned}
\label{eq:sequence-sum-c}
\sum_{k=2}^\infty\frac{(\rho^k/\rho^{k-1}-1)^2}{\rho^k}
\leq\frac{p^2}{\rho^0}H_p(\xi).
\end{eqaligned}
\end{lemma}
\begin{proof}
Let $p>0$ and $k\geq2$, and put $j=k-1\geq1$.  By the mean value theorem, $k^p-j^p=pc^{p-1}(k-j)=pc^{p-1}$ for some $c\in(j,k)$, since $k-j=1$.  Since $t\mapsto t^p$ is concave, its derivative is nonincreasing, and therefore
\begin{eqaligned}
k^p-j^p\leq pj^{p-1}.
\end{eqaligned}
For the ratio series, the common schedule gives the identity
\[
\begin{aligned}
\frac{\rho^k}{\rho^{k-1}}-1&=\frac{\xi(k^p-j^p)}{1+\xi j^p}
\leq\frac{\xi p j^{p-1}}{1+\xi j^p}.
\end{aligned}
\]
Both sides are nonnegative by schedule monotonicity, so squaring preserves the inequality.  Using $\rho^k\geq\rho^0$ and shifting from $k\geq2$ to $j=k-1\geq1$ yield
\[
\sum_{k=2}^\infty
\frac{(\rho^k/\rho^{k-1}-1)^2}{\rho^k}
\leq\frac{p^2}{\rho^0}\xi^2\sum_{j=1}^\infty
\frac{j^{2p-2}}{(1+\xi j^p)^2}
=\frac{p^2}{\rho^0}H_p(\xi).
\]

For the square-root series, rationalization gives
\[
\sqrt{L^k}-\sqrt{L^{k-1}}=\frac{L^0\xi(k^p-j^p)}{\sqrt{L^k}+\sqrt{L^{k-1}}}.
\]
Moreover,
\[
\begin{aligned}
\sqrt{L^k}+\sqrt{L^{k-1}}
\geq2\sqrt{L^{k-1}}
=2\sqrt{L^0}\sqrt{1+\xi j^p}.
\end{aligned}
\]
Consequently,
\[
\sqrt{L^k}-\sqrt{L^{k-1}}\leq\frac{\sqrt{L^0}\,\xi p j^{p-1}}{2\sqrt{1+\xi j^p}}.
\]
Both sides are nonnegative by schedule monotonicity, so this inequality may also be squared.  Use the sharper denominator bound $\rho^k=\rho^0(1+\xi k^p)\geq\rho^0(1+\xi j^p)$ and again set $j=k-1$ to obtain
\[
\sum_{k=2}^\infty
\frac{(\sqrt{L^k}-\sqrt{L^{k-1}})^2}{\rho^k}
\leq\frac{L^0p^2}{4\rho^0}\xi^2\sum_{j=1}^\infty
\frac{j^{2p-2}}{(1+\xi j^p)^2}
=\frac{L^0p^2}{4\rho^0}H_p(\xi).
\]
The defining series converges because, for fixed $\xi>0$, its full summand is at most $j^{-2}$. If $p=0$, then $k^p=1$ for every positive index, so both $L^k$ and $\rho^k$ are constant and every summand on the two left sides vanishes.
\end{proof}

\begin{lemma}[Small-$\xi$ regime of $H_{1/2}$]
\label{lem:Hp-small-xi}
As $\xi\downarrow0$, $H_{1/2}(\xi)=\Theta(\xi^2\log(e/\xi))$.
\end{lemma}
\begin{proof}
Because the claim is asymptotic as $\xi\downarrow0$, it suffices to take $0<\xi\leq1/2$.  Let $J=\lceil\xi^{-2}\rceil$ and split the series in~\prettyref{eq:Hp} with $p=1/2$ at $J$.  Since $\xi^{-2}\geq1$ and hence $J\leq2\xi^{-2}$, for $j\leq J$ one has $1\leq1+\xi j^{1/2}\leq1+\sqrt2$.  Thus the head denominator is bounded above and below by absolute positive constants, and, by comparison with $\int_1^J t^{-1}\,dt$, the head is comparable to
\[
\xi^2\sum_{j=1}^Jj^{-1}=\Theta(\xi^2\log J).
\]
On $j>J$, one has $\xi j^{1/2}>1$, so the full summand in~\prettyref{eq:Hp}, including the $\xi^2$ prefactor, lies between $\tfrac14j^{-2}$ and $j^{-2}$, and the tail is $\Theta(J^{-1})=\Theta(\xi^2)$.  Hence the head dominates, and $\log J=\Theta(\log(e/\xi))$ proves the claim.
\end{proof}

\begin{lemma}[Geometric checkpoint schedule identities]
\label{lem:checkpoint-schedule}
Under the checkpoint schedule~\prettyref{eq:checkpoint-schedule}, with $\mathsf m=\mathsf{Ckpt}$ (so $G_{\mathsf m}=G^{\rm cp}$), $\rho_{\rm cp}^0>0$, $\chi>1$, $\ell_{\rm cp}>0$, and $b_{\rm cp}$ from~\prettyref{eq:checkpoint-parameters}, the following wall scalings hold for every integer $j\geq0$ and every $i\in\IC$:
\begin{eqaligned}
\label{eq:checkpoint-wall-scaling}
\beta(L_j)=\frac{b_{\rm cp}}{\rho_j}=\beta_j,
\qquad
\frac{\sqrt{2L_jG^{\rm cp}}\,L_i}{\rho_j}=\sqrt{2G^{\rm cp}\ell_{\rm cp}}\,L_i.
\end{eqaligned}
The geometric sums satisfy
\begin{eqaligned}
\label{eq:checkpoint-geometric-sums}
\sum_{j=0}^\infty\frac1{\rho_j}
&=\frac1{\rho_{\rm cp}^0(1-\chi^{-1})},
\qquad
\sum_{j=0}^J\rho_j^3
&\leq\frac{\rho_J^3}{1-\chi^{-3}}\quad\text{for every integer }J\geq0.
\end{eqaligned}
For every integer $j\geq0$ and $s\geq0$, one has $g_{\mathsf m}(s,L_{j+1})\leq g_{\mathsf m}(s,L_j)$.
\end{lemma}
\begin{proof}
Since $\rho_j>0$, substituting $L_j=\ell_{\rm cp}\rho_j^2$ and $b_{\rm cp}=\sqrt{2G^{\rm cp}/\ell_{\rm cp}}$, and $\beta_j=b_{\rm cp}/\rho_j$ into $\beta(L_j)=\sqrt{2G^{\rm cp}/L_j}$ gives the first identity in~\prettyref{eq:checkpoint-wall-scaling}; the second follows from $\sqrt{2L_jG^{\rm cp}}=\rho_j\sqrt{2G^{\rm cp}\ell_{\rm cp}}$, using $\rho_j>0$.  Both sums in~\prettyref{eq:checkpoint-geometric-sums} are geometric; for the finite sum, write $\rho_j^3=\rho_J^3\chi^{-3(J-j)}$ and enlarge it to the infinite backward series.  Finally, for $s\geq0$ and $L>0$, moving $\sqrt L>0$ inside the max in the explicit wall formula gives
\[
g_{\mathsf m}(s,L)=\frac12\max\{\sqrt{2G^{\rm cp}}-\sqrt L\,s,0\}^2.
\]
Since $\sqrt{2G^{\rm cp}}-s\sqrt L$ is nonincreasing in $L$ for $s\geq0$ and $t\mapsto\max\{t,0\}^2$ is nondecreasing, the display is nonincreasing in $L$; use $L_{j+1}=\chi^2L_j\geq L_j$.
\end{proof}

The following force-drift quantities are used only in the ordinary mode ($\mathsf m=\mathsf{IP}$, so $G_{\mathsf m}=G$) and for $k\geq2$.  Set
\begin{eqaligned}
b_i^k\defeq
\begin{cases}H_i,&i\in\IO,\\2L^kM_i^2+\sqrt{2L^kG}\,L_i,&i\in\IC,
\end{cases}
\end{eqaligned}
and
\begin{eqaligned}
\delta_i^k\defeq
\begin{cases}
C_i^\psi(\rho^k/\rho^{k-1}-1),&i\in\IO,\\
\sqrt{2G}M_i(\sqrt{L^k}-\sqrt{L^{k-1}}),&i\in\IC.
\end{cases}
\end{eqaligned}
Then
\begin{eqaligned}
\label{eq:delta}
\|q_i^k-q_i^{k-1}\|\leq\zeta_i^k\|d_i^k\|+b_i^k\|d_i^{k-1}\|+\delta_i^k,
\end{eqaligned}
and define
\begin{eqaligned}
\label{eq:U}
U^k\defeq\frac{5}{\rho^k}\sum_{i\in\IAll}(\delta_i^k)^2.
\end{eqaligned}
\begin{lemma}[Uniform force drift]
\label{lem:fi-Smooth-L}
Under the polynomial schedule~\prettyref{eq:IP-schedules}, for every $i\in\IAll$ and $k\geq2$, the quantities $\zeta_i^k$, $b_i^k$, and $\delta_i^k$ defined above satisfy~\prettyref{eq:delta} at every state, without a feasibility hypothesis.  Moreover, $\delta_i^k\geq0$.
\end{lemma}
\begin{proof}
Schedule monotonicity first gives $\delta_i^k\geq0$ on both branches.  For objective blocks, split the force difference into a state change of the convex envelope at the common parameter $\mu^k$, a parameter change at the common point $y_i^k$, and the change of $\nabla h_i$.  Apply~\prettyref{eq:env-smooth}, $H_i$-smoothness, and~\prettyref{eq:env-grad-drift}; the identity $\mu^{k-1}/\mu^k=\rho^k/\rho^{k-1}$, immediate from $\mu^k=\eta/\rho^k$ in~\prettyref{eq:IP-schedules}, gives the objective branch of $\delta_i^k$.  For constraint blocks put $s_i^k=\ell_i^k(y_i^{k+1})$.  Bounded gradients give $|s_i^k-s_i^{k-1}|\leq M_i\|d_i^k\|+2M_i\|d_i^{k-1}\|$.  Since $L^{k-1}\leq L^k$ by~\prettyref{eq:schedule-monotone}, both scalar state-change terms may use the $L^k$ Lipschitz and force bounds.  The parameter estimate in~\prettyref{eq:wall-parameter-drift} contributes $\sqrt{2G}(\sqrt{L^k}-\sqrt{L^{k-1}})$, where the absolute value drops by the same monotonicity.  Multiplying by $\nabla f_i(y_i^k)$ and then using smoothness to transfer that gradient from $y_i^{k-1}$ yields the constraint branch.
\end{proof}

\begin{lemma}[Summable force drift]
\label{lem:drift-sum}
Under~\prettyref{eq:parameters} and the polynomial schedule~\prettyref{eq:IP-schedules} with $p\in[0,1)$ and $\xi\in(0,1]$,
\begin{eqaligned}
\label{eq:drift-sum-a}
\sum_{k=2}^\infty U^k\leq\Sigma_U\defeq
\frac{p^2H_p(\xi)}{\rho^0}
\left(5G_\psi+\frac52GL^0\bar M\right),
\end{eqaligned}
where $G_\psi$ and $\bar M$ are the quantities already defined in~\prettyref{eq:parameters}.
\end{lemma}
\begin{proof}
Because $\IAll$ is finite and every summand is nonnegative, the block sum in~\prettyref{eq:U} may be interchanged with the infinite $k$-sum.  The definition~\prettyref{eq:U} contributes the common factor $5$.  For objective blocks, inserting $\delta_i^k=C_i^\psi(\rho^k/\rho^{k-1}-1)$ into~\prettyref{eq:sequence-sum-c} and summing over $i\in\IO$ gives $p^2H_p(\xi)G_\psi/\rho^0$ before that factor, hence the objective constant $5$.  For constraint blocks, $(\delta_i^k)^2=2GM_i^2(\sqrt{L^k}-\sqrt{L^{k-1}})^2$; applying~\prettyref{eq:sequence-sum-a} and summing over $i\in\IC$ gives $p^2H_p(\xi)GL^0\bar M/(2\rho^0)$ before the common factor, hence the constraint constant $5/2$.  At $p=0$, the definition $H_0(\xi)=0$ makes $\Sigma_U=0$.
\end{proof}

\section{Convergence Analysis of IP-ADMM}
\label{appen:convergence}
This appendix proves the ordinary-mode results of \prettyref{thm:rate} and \prettyref{cor:KKT-MFCQ}.  In the $\mathsf{IP}$ mode of \prettyref{alg:IP-ADMM}, the multiplier is updated at every iteration with the relaxation $\sigma=1$ fixed in~\prettyref{eq:parameters}, and $(\rho^k,L^k,\mu^k)$ follow~\prettyref{eq:IP-schedules}.  The coupled relation $\mu^k=\mu_{\rm IP}(L^k)$ from~\prettyref{eq:IP-schedules} and~\prettyref{eq:rho-L-ratio} supplies the objective continuation correction~\prettyref{lem:Rpsi-discharge}.  The nonlinear wall adds the main difficulty: comparing successive wall parameters requires $f_i(y_i^{k+1})\geq0$ for every $i\in\IC$, which is not known when the step is first analyzed.  We resolve this dependence with an unconditional terminal estimate and explicit parameter margins before using the wall continuation comparison globally.  The checkpoint mode uses a different proof and is treated in \prettyref{appen:checkpoint}.

\paragraph{Organization.}
The local model and multiplier identities in \prettyref{sub:lyapunov-descent} give one-step descent and a Lyapunov bridge whose terminal continuation comparison is still conditional.  \prettyref{lem:terminal-bound} controls that terminal point unconditionally, and \prettyref{lem:initial-bound} and the margins in \prettyref{lem:parameter} supply the finite budget with which \prettyref{thm:feasibility} closes the wall-feasibility induction.  In \prettyref{sub:stationarity-rates}, \prettyref{lem:lower-bound} yields the movement budget \prettyref{cor:movement-budget}, which feeds \prettyref{thm:rate-master}; \prettyref{cor:KKT-MFCQ-exact} converts the same selected candidate under Uniform Band-MFCQ, and the family and tuned-$\xi$ results specialize the exact rate constants.

\subsection{Lagrangian and Lyapunov Descent}
\label{sub:lyapunov-descent}
Call a split iterate $y^k$ \emph{wall-feasible} when $f_i(y_i^k)\geq0$ for every $i\in\IC$.  Put $\Delta x^k=x^{k+1}-x^k$, $\Delta\lambda^{k+1}=\lambda^{k+1}-\lambda^k$, and
\begin{eqaligned}
\upsilon_i^k&\defeq m_i^k-\rho^k,&
\Gamma_i^k&\defeq|\upsilon_i^k|+\zeta_i^k,&
\widehat\Gamma_i^k&\defeq|\upsilon_i^{k-1}|+b_i^k.
\end{eqaligned}
Write $\Upsilon^k\defeq\operatorname{diag}(\upsilon_i^kI_{n_i})_{i\in\IAll}$ for the block-diagonal curvature matrix, so that $(\Upsilon^kd^k)_i=\upsilon_i^kd_i^k$. For $k\geq2$, define
\begin{eqaligned}
\Theta^k\defeq{}&\frac{5}{\rho^k}
\sum_{i\in\IAll}\left[(\widehat\Gamma_i^k)^2\|d_i^{k-1}\|^2
+(\rho^{k-1})^2\|A_i\Delta x^{k-1}\|^2\right],\\
\kappa_\Psi&\defeq\frac{\xi+4}{2},\\
\label{eq:Lyapunov-def}
\bar\Psi^k&\defeq\mathcal L(x^k,y^k,\lambda^k,L^{k-1},\mu^{k-1},\rho^k)
+R(L^{k-1})+\kappa_\Psi\Theta^k,\\
\Psi^k&\defeq\mathcal L(x^k,y^k,\lambda^k,L^k,\mu^k,\rho^k)
+R(L^k)+\kappa_\Psi\Theta^k.
\end{eqaligned}

For either mode $\mathsf m$, the exact optimality condition of the active block model is
\begin{eqaligned}
\label{eq:opt-yi}
0=q_i^k+\lambda_i^k+\rho^k(y_i^k-A_ix^k)+m_i^kd_i^k.
\end{eqaligned}
For ordinary IP-ADMM, using $\Delta\lambda_i^{k+1}=\rho^k(y_i^{k+1}-A_ix^{k+1})$ and $y_i^k-A_ix^k=y_i^{k+1}-A_ix^{k+1}-d_i^k+A_i\Delta x^k$ gives
\begin{eqaligned}
\label{eq:block-dual-identity}
\lambda_i^{k+1}=-q_i^k-\upsilon_i^kd_i^k-\rho^kA_i\Delta x^k.
\end{eqaligned}

For either mode $\mathsf m$, the $x$-subproblem has a multiplier $\varpi^{k+1}$ satisfying
\begin{eqaligned}
\label{eq:opt-x}
(1+\theta_{2,\mathsf m})A^TA\Delta x^k+C^T\varpi^{k+1}&=A^T(y^{k+1}-Ax^k+\lambda^k/\rho^k),\\
Cx^{k+1}&=d.
\end{eqaligned}

From this point through the end of the ordinary appendix, set $\mathsf m=\mathsf{IP}$, so $\theta_{1,\mathsf m}=\theta_1$, $\theta_{2,\mathsf m}=\theta_2$, $G_{\mathsf m}=G$, and $g\defeq g_{\mathsf m}$.

\begin{lemma}[One-step descent]
\label{lem:descendent}
For $k\geq1$,
\begin{eqaligned}
&\mathcal L(x^{k+1},y^{k+1},\lambda^{k+1},L^k,\mu^k,\rho^{k+1})
-\mathcal L(x^k,y^k,\lambda^k,L^k,\mu^k,\rho^k)\\
&\leq-\frac{2\theta_1-1}{2}\sum_{i\in\IAll}\mathsf B_i^k\|d_i^k\|^2
-\frac{1+2\theta_2}{2}\rho^k\|A\Delta x^k\|^2
+\frac{\xi+2}{2\rho^k}\|\Delta\lambda^{k+1}\|^2.
\end{eqaligned}
If $y^{k+1}$ is wall-feasible, then
\begin{eqaligned}
\label{eq:Lk-update}
&\mathcal L(x^{k+1},y^{k+1},\lambda^{k+1},L^{k+1},\mu^{k+1},\rho^{k+1})+R(L^{k+1})\\
\leq&\mathcal L(x^{k+1},y^{k+1},\lambda^{k+1},L^k,\mu^k,\rho^{k+1})+R(L^k).
\end{eqaligned}
\end{lemma}
\begin{proof}
For a block $i$, apply $\bar f_i(u,L^k,\mu^k)\leq\widetilde f_i^k(u)+\tau_i^k\|u-y_i^k\|^2/2$ from~\prettyref{eq:model-two-sided} at $u=y_i^{k+1}$; at $u=y_i^k$ the two functions agree exactly.  Together with exact expansion of the consensus quadratic, this shows that the true Lagrangian change of the block is at most its model change after removing the proximal term, plus $\mathsf B_i^k\|d_i^k\|^2/2$.  Since $\mathcal M_i^k$ is $m_i^k$-strongly convex and is minimized at $y_i^{k+1}$,
\[
\mathcal M_i^k(y_i^{k+1})-\mathcal M_i^k(y_i^k)
\leq-\frac{m_i^k}{2}\|d_i^k\|^2.
\]
Removing the proximal term contributes another $-m_i^k\|d_i^k\|^2/2$; therefore the true block change is at most $-(m_i^k-\mathsf B_i^k/2)\|d_i^k\|^2 =-(2\theta_1-1)\mathsf B_i^k\|d_i^k\|^2/2$. Summing gives the first term of the lemma.

Contracting~\prettyref{eq:opt-x} with $\Delta x^k$ and using $C\Delta x^k=0$ gives the $x$-change $-(1+2\theta_2)\rho^k\|A\Delta x^k\|^2/2$. The multiplier and penalty-parameter changes equal $(\rho^{k+1}+\rho^k)\|y^{k+1}-Ax^{k+1}\|^2/2$.  The multiplier update $\Delta\lambda^{k+1}=\rho^k(y^{k+1}-Ax^{k+1})$ substitutes the residual identity
\[
\|y^{k+1}-Ax^{k+1}\|^2=\frac{\|\Delta\lambda^{k+1}\|^2}{(\rho^k)^2},
\]
while~\prettyref{lem:rho-increase} gives $\rho^{k+1}+\rho^k\leq(\xi+2)\rho^k$ through $\rho^{k+1}\leq(1+\xi)\rho^k$.  Multiplying the two yields at most the last term $(\xi+2)\|\Delta\lambda^{k+1}\|^2/(2\rho^k)$ of the first assertion.

For the continuation comparison, the objective aggregate plus $R$ decreases by~\prettyref{eq:R2}.  Wall-feasibility of $y^{k+1}$ places every constraint block on the branch $s\geq0$ of the piecewise wall, where it admits the equivalent form
\[
g(s,L)=\frac12\max\{\sqrt{2G}-\sqrt L\,s,0\}^2,
\]
which is manifestly nonincreasing in $L$; with $L^{k+1}\geq L^k$ from~\prettyref{eq:schedule-monotone} this gives $g(s,L^{k+1})\leq g(s,L^k)$.  The multiplier and penalty terms do not depend on $L$, proving~\prettyref{eq:Lk-update}.
\end{proof}

Subtracting~\prettyref{eq:block-dual-identity} at consecutive indices yields
\begin{eqaligned}
\label{eq:lambda-a}
\Delta\lambda_i^{k+1}=-z_i^k,
\end{eqaligned}
where
\begin{eqaligned}
z_i^k=q_i^k-q_i^{k-1}+\upsilon_i^kd_i^k-\upsilon_i^{k-1}d_i^{k-1}
+\rho^kA_i\Delta x^k-\rho^{k-1}A_i\Delta x^{k-1}.
\end{eqaligned}
\begin{lemma}[Multiplier-difference estimate]
\label{lem:lambda}
For $k\geq2$,
\begin{eqaligned}
\frac{\|\Delta\lambda^{k+1}\|^2}{\rho^k}
\leq{}&U^k+5\sum_{i\in\IAll}
\left[\frac{(\Gamma_i^k)^2}{\rho^k}
+\frac{(\widehat\Gamma_i^{k+1})^2}{\rho^{k+1}}\right]\|d_i^k\|^2
+10\rho^k\|A\Delta x^k\|^2+\Theta^k-\Theta^{k+1}.
\end{eqaligned}
\end{lemma}
\begin{proof}
By~\prettyref{eq:lambda-a}, $\|\Delta\lambda_i^{k+1}\|=\|z_i^k\|$.  By~\prettyref{eq:delta}, which~\prettyref{lem:fi-Smooth-L} makes valid at every state, wall-feasible or not,
\[
\|z_i^k\|\leq
\Gamma_i^k\|d_i^k\|+\widehat\Gamma_i^k\|d_i^{k-1}\|
+\rho^k\|A_i\Delta x^k\|+\rho^{k-1}\|A_i\Delta x^{k-1}\|
+\delta_i^k.
\]
Square with $(a_1+\cdots+a_5)^2\leq5\sum_ja_j^2$, sum over blocks, and divide by $\rho^k$.  Add and subtract the next-index terms $5(\widehat\Gamma_i^{k+1})^2\|d_i^k\|^2/\rho^{k+1}$ and $5(\rho^k)^2\|A_i\Delta x^k\|^2/\rho^{k+1}$. The two terms carrying the current displacement $A\Delta x^k$ are the squared term $5\rho^k\|A_i\Delta x^k\|^2$ produced by the third summand of the bound on $\|z_i^k\|$ and the added next-index term $5(\rho^k)^2\|A_i\Delta x^k\|^2/\rho^{k+1}$; schedule monotonicity~\prettyref{eq:schedule-monotone} bounds the latter by the former, so the two together are at most $10\rho^k\|A\Delta x^k\|^2$.  The lagged terms and the subtracted next-index terms give exactly $\Theta^k-\Theta^{k+1}$.  The drift squares are exactly $U^k$.
\end{proof}

Define the positive candidate coefficients
\begin{eqaligned}
D_{y,i}^k&\defeq\frac{2\theta_1-1}{2}\mathsf B_i^k
-5\kappa_\Psi\left[
\frac{(\Gamma_i^k)^2}{\rho^k}
+\frac{(\widehat\Gamma_i^{k+1})^2}{\rho^{k+1}}\right],\\
D_x^k&\defeq\left(\frac{1+2\theta_2}{2}
-10\kappa_\Psi\right)\rho^k.
\end{eqaligned}
Combining~\prettyref{lem:descendent} and~\prettyref{lem:lambda}, with $(\xi+2)/2-\kappa_\Psi=-1$ on the dual increment, gives
\begin{eqaligned}
\label{eq:Lyapunov}
\bar\Psi^{k+1}-\Psi^k
\leq-\frac{1}{\rho^k}\|\Delta\lambda^{k+1}\|^2
-\sum_{i\in\IAll}D_{y,i}^k\|d_i^k\|^2
-D_x^k\|A\Delta x^k\|^2+\kappa_\Psi U^k.
\end{eqaligned}
If $y^{k+1}$ is wall-feasible,~\prettyref{eq:Lk-update} gives $\Psi^{k+1}\leq\bar\Psi^{k+1}$.

\subsection{Split-Feasibility Bounds}
For later use define
\begin{eqaligned}
C_{d,i}^k&\defeq\frac{m_i^k}{2}+\frac{3\tau_i^k}{2},\\
\widetilde C_{r,i}^k&\defeq
\begin{cases}
\dfrac{m_i^k}{2}+\tau_i^k+\dfrac{\zeta_i^k}{2},&i\in\IO,\\[1mm]
\dfrac{m_i^k}{2}+\dfrac{5\tau_i^k}{2}+2\zeta_i^k,&i\in\IC.
\end{cases}
\end{eqaligned}

\begin{lemma}[Unconditional terminal lower bound]
\label{lem:terminal-bound}
For every $k\geq1$, without assuming that $y^{k+1}$ is wall-feasible,
\begin{eqaligned}
\label{eq:terminal-bound}
&\mathcal L(x^{k+1},y^{k+1},\lambda^{k+1},L^k,\mu^k,\rho^{k+1})\\
&\geq\sum_{i\in\IO}\bar f_i(A_ix^{k+1},L^k,\mu^k)
+\frac12\sum_{i\in\IC}\bar f_i(y_i^{k+1},L^k,\mu^k)\\
&\quad-\sum_{i\in\IAll}\left[
\widetilde C_{r,i}^k\|y_i^{k+1}-A_ix^{k+1}\|^2
+C_{d,i}^k\|d_i^k\|^2
+\frac{\rho^k}{4}\|y_i^k-A_ix^k\|^2\right].
\end{eqaligned}
\end{lemma}
\begin{proof}
Write $r_i^{k+1}=y_i^{k+1}-A_ix^{k+1}$ and let $\widetilde f_i^k$ be the model-only function from \prettyref{eq:model-two-sided}.  It is $\zeta_i^k$-smooth and $q_i^k=\nabla\widetilde f_i^k(y_i^{k+1})$.  Apply \prettyref{eq:model-two-sided} at $y_i^{k+1}$ and $A_ix^{k+1}$ and apply the descent lemma to $\widetilde f_i^k$ between those points.  Since $A_ix^{k+1}-y_i^k=d_i^k-r_i^{k+1}$, expanding and using the Young inequality $\|d_i^k-r_i^{k+1}\|^2\leq2\|d_i^k\|^2+2\|r_i^{k+1}\|^2$ gives
\begin{eqaligned}
\bar f_i(y_i^{k+1},L^k,\mu^k)-\langle q_i^k,r_i^{k+1}\rangle
\geq\bar f_i(A_ix^{k+1},L^k,\mu^k)
-\frac{3\tau_i^k}{2}\|d_i^k\|^2
-\left(\tau_i^k+\frac{\zeta_i^k}{2}\right)\|r_i^{k+1}\|^2.
\end{eqaligned}
Model optimality and the multiplier update give
\[
\lambda_i^{k+1}=-q_i^k-m_i^kd_i^k-\rho^k(y_i^k-A_ix^k) +\rho^kr_i^{k+1}.
\]
Insert this into the terminal Lagrangian.  The inequalities
\[
-m_i^k\langle d_i^k,r_i^{k+1}\rangle \geq-\frac{m_i^k}{2}(\|d_i^k\|^2+\|r_i^{k+1}\|^2)
\]
and
\[
-\rho^k\langle y_i^k-A_ix^k,r_i^{k+1}\rangle \geq-\rho^k\|r_i^{k+1}\|^2 -\frac{\rho^k}{4}\|y_i^k-A_ix^k\|^2
\]
bound the cross terms, and the term $-\rho^k\|r_i^{k+1}\|^2$ produced by the second cancels the term $+\rho^k\|r_i^{k+1}\|^2$ contributed by $\lambda_i^{k+1}$.  Discarding the remaining positive terminal penalty gives the objective-block coefficients above and the pre-wall constraint coefficients.

For $i\in\IC$ and $F_i^k(y)=g(f_i(y),L^k)$, the chain rule gives $\nabla F_i^k(y)=g'(f_i(y),L^k)\nabla f_i(y)$, so $F_i^k$ is $(\tau_i^k+\zeta_i^k)$-smooth: the product-difference bound contributes $\sqrt{2L^kG}L_i=\tau_i^k$ through $|g'|\leq\sqrt{2L^kG}$ and the $L_i$-smoothness of $f_i$, and $L^kM_i^2=\zeta_i^k$ through the $L^k$-Lipschitz continuity of $g'(\cdot,L^k)$ and $\|\nabla f_i\|\leq M_i$.  Moreover $|g'(s,L)|^2\leq2Lg(s,L)$ for every $s$: for $s\geq0$, \prettyref{eq:band-force} gives $|g'(s,L)|^2=L^2\max(\beta(L)-s,0)^2=2Lg(s,L)$, while for $s<0$ it gives $|g'(s,L)|^2=2LG\leq2Lg(s,L)$ because $g(s,L)>G$ there by~\prettyref{lem:added_algebraic_3}.  Hence $\|\nabla F_i^k(y)\|^2\leq|g'(f_i(y),L^k)|^2M_i^2\leq2L^kM_i^2F_i^k(y)=2\zeta_i^kF_i^k(y)$.  Combining the descent lemma for the $(\tau_i^k+\zeta_i^k)$-smooth $F_i^k$ with this gradient bound and the Young inequality $\sqrt{2\zeta_i^kF_i^k(y)}\|r\|\leq F_i^k(y)/2+\zeta_i^k\|r\|^2$, we obtain
\[
F_i^k(y-r)\geq\frac12F_i^k(y) -\left(\frac{\tau_i^k}{2}+\frac{3\zeta_i^k}{2}\right)\|r\|^2.
\]
Apply this with $y=y_i^{k+1}$ and $r=r_i^{k+1}$.  Adding its loss to the pre-wall residual coefficient gives the sharper coefficient $m_i^k/2+3\tau_i^k/2+2\zeta_i^k$.  Since $\tau_i^k\geq0$, enlarging it to the displayed $m_i^k/2+5\tau_i^k/2+2\zeta_i^k$ weakens the lower bound in the valid direction.  This completes the proof.
\end{proof}

Combining~\prettyref{eq:terminal-bound} with~\prettyref{eq:Lyapunov} and $r_i^{k+1}=\Delta\lambda_i^{k+1}/\rho^k$ gives
\begin{eqaligned}
\label{eq:Lyapunov-lower-bound}
&W_{k+1}-\Psi^k\\
&\leq\sum_i\left[
\frac{\widetilde C_{r,i}^k}{(\rho^k)^2}
-\frac{1}{\rho^k}\right]\|\Delta\lambda_i^{k+1}\|^2
+\sum_i(C_{d,i}^k-D_{y,i}^k)\|d_i^k\|^2\\
&\quad-D_x^k\|A\Delta x^k\|^2+\kappa_\Psi U^k
+\frac{\rho^k}{4(\rho^{k-1})^2}\|\Delta\lambda^k\|^2,
\end{eqaligned}
where
\begin{eqaligned}
W_{k+1}\defeq\sum_{i\in\IO}\bar f_i(A_ix^{k+1},L^k,\mu^k)+R(L^k)
+\frac12\sum_{i\in\IC}\bar f_i(y_i^{k+1},L^k,\mu^k).
\end{eqaligned}
This estimate is valid before $y^{k+1}$ is known to be wall-feasible.

\begin{lemma}[Initialization]
\label{lem:initial-bound}
Under~\prettyref{eq:parameters}, $d_i^1=0$ for $i\in\IC$ and $\|d_i^1\|\leq\widehat C_i/m_i^1$ for $i\in\IO$.  Moreover, $\|y^2-Ax^2\|\leq\|d^1\|$, and
\begin{eqaligned}
\label{eq:Psi-2-bound}
\Psi^2\leq{}&S_0+R(L^1)+\sum_{i\in\IO}\frac{\widehat C_i^2}{m_i^1}\\
&+\left(\rho^1+\frac{\rho^2}{2}\right)
\sum_{i\in\IO}\frac{\widehat C_i^2}{(m_i^1)^2}\\
&+\frac{25}{2\rho^2}
\sum_{i\in\IO}\left[(\widehat\Gamma_i^2)^2
+\frac{(\rho^1)^2}{26^2}\right]
\frac{\widehat C_i^2}{(m_i^1)^2}.
\end{eqaligned}
\end{lemma}
\begin{proof}
By~\prettyref{eq:band} and~\prettyref{eq:p7}, $\beta^1=\sqrt{2G/L^1}\leq m$.  At initialization $e_i^1=0$.  If $i\in\IC$, then $f_i(y_i^1)\geq m\geq\beta^1$ and hence $g(\ell_i^1(y_i^1),L^1)=0$.  The constraint model is a sum of two nonnegative terms and has value zero at $y_i^1$; strong convexity makes this its unique minimizer, so $d_i^1=0$.  For $i\in\IO$, the block optimality condition~\prettyref{eq:opt-yi} at $k=1$ collapses to $0=q_i^1+m_i^1d_i^1$ with $q_i^1=\nabla h_i(y_i^1)+w_i^2$, because $\lambda^1=0$ and $y_i^1=A_ix^1$ annihilate its two consensus terms.  By~\prettyref{eq:objective-derived-bounds}, $\|\nabla h_i(y_i^1)\|\leq D_i$, while $w_i^2\in\partial\psi_i(p_i^2)$ by~\prettyref{eq:one-call-outputs} gives $\|w_i^2\|\leq C_i^\psi$; hence $\|d_i^1\|=\|q_i^1\|/m_i^1\leq(D_i+C_i^\psi)/m_i^1=\widehat C_i/m_i^1$.

Let $Q$ be the orthogonal projector onto $A(\ker C)$.  At $k=1$ the right side of~\prettyref{eq:opt-x} is $A^Td^1$: its multiplier term drops because $\lambda^1=0$, and $y^1=Ax^1$ turns $y^2-Ax^1$ into $d^1$.  Contract \prettyref{eq:opt-x} with arbitrary directions in $\ker C$.  Since $C\Delta x^1=0$, this shows that $d^1-(1+\theta_2)A\Delta x^1$ is orthogonal to $A(\ker C)$, while $A\Delta x^1\in A(\ker C)$.  Applying $Q$ therefore gives $A\Delta x^1=Qd^1/(1+\theta_2)$.  Hence $\|A\Delta x^1\|\leq\|d^1\|/(1+\theta_2)$.  Moreover, $y^2-Ax^2=d^1-A\Delta x^1=(I-Q)d^1+[\theta_2/(1+\theta_2)]Qd^1$; orthogonality of $Qd^1$ and $(I-Q)d^1$ gives $\|y^2-Ax^2\|\leq\|d^1\|$. The corrected objective aggregate at $y^2$ is at most
\[
S_0+R(L^1)+\sum_{i\in\IO}\widehat C_i\|d_i^1\| \leq S_0+R(L^1)+\sum_{i\in\IO}\widehat C_i^2/m_i^1,
\]
using~\prettyref{eq:R2}, the $\widehat C_i$-Lipschitz property, and $(\psi_i)_{\mu^1}\leq\psi_i$.  For $i\in\IC$, $d_i^1=0$, hence $y_i^2=y_i^1$.  Since $\beta^2\leq\beta^1\leq m\leq f_i(y_i^2)$, one has $g(f_i(y_i^2),L^2)=0$.  Thus every constraint-wall value $g(f_i(y_i^2),L^2)$ appearing in $\Psi^2$ is exactly zero.  Expanding $\Psi^2$ and inserting the two projector bounds yields \prettyref{eq:Psi-2-bound}.  Here $\lambda^2=\rho^1(y^2-Ax^2)$ is the multiplier update at $k=1$ read with $\lambda^1=0$, while $\kappa_\Psi\leq5/2$ and $1+\theta_2=26$ bound $\kappa_\Psi\Theta^2$ by the last term displayed in~\prettyref{eq:Psi-2-bound}.
\end{proof}

\subsection{Parameter Verification and Split-Feasibility Induction}
\begin{lemma}[Explicit parameter margins]
\label{lem:parameter}
Under~\prettyref{eq:parameters}, for every block and every admissible index,
\begin{align}
\label{eq:p2}
\frac{\rho^k}{4(\rho^{k-1})^2}-\frac{1}{\rho^{k-1}}&\leq0,\\
\label{eq:p3}
\frac{\widetilde C_{r,i}^k}{(\rho^k)^2}-\frac{1}{\rho^k}&\leq0,\\
\label{eq:p4}
D_x^k&\geq d_0\rho^k,\\
\label{eq:p5}
D_{y,i}^k&\geq d_0\rho^k,\\
\label{eq:p6}
C_{d,i}^k-D_{y,i}^k&\leq0,\\
\label{eq:p7}
\sqrt{2G/L^1}&\leq m.
\end{align}
In addition, $0\leq\tau_i^k/\rho^k\leq1/1209$, and the finite-wall budget holds:
\begin{eqaligned}
\label{eq:p8}
&S_0+R(L^1)+\sum_{i\in\IO}\frac{\widehat C_i^2}{m_i^1}
+\left(\rho^1+\frac{\rho^2}{2}+\frac{\rho^2}{4}\right)
\sum_{i\in\IO}\frac{\widehat C_i^2}{(m_i^1)^2}\\
&\quad+\frac{25}{2\rho^2}
\sum_{i\in\IO}\left[(\widehat\Gamma_i^2)^2
+\frac{(\rho^1)^2}{26^2}\right]
\frac{\widehat C_i^2}{(m_i^1)^2}
+\kappa_\Psi\Sigma_U\leq\frac G2.
\end{eqaligned}
\end{lemma}
\begin{proof}
Write $\varepsilon\defeq1/1209$, and recall $\theta_1=121909/121000$, $\theta_2=25$, $\eta=1209$, and $d_0=1/2$ from~\prettyref{eq:parameters}.  The first four clauses in the definition of $\rho^0$, namely the constant $1$ together with the three data clauses $1209H_{\max}$, $1209mL^0L_{\max}$, and $1209L^0M_{\max}^2$, together with $\eta=1209$, imply, for every $i,k$,
\begin{eqaligned}
0\leq\frac{\tau_i^k}{\rho^k}\leq\varepsilon,
\qquad0\leq\frac{\zeta_i^k}{\rho^k}\leq\varepsilon.
\end{eqaligned}
Indeed, this is immediate for objective blocks.  For constraint blocks use $\sqrt{2L^kG}=mL^0\sqrt{1+\xi k^p}$ and the constant ratio $L^k/\rho^k=L^0/\rho^0$.  Since
\[
\theta_1-1=\frac{909}{121000},
\qquad
\upsilon_i^k=(\theta_1-1)\rho^k+\theta_1\tau_i^k,
\]
the displayed bounds give the exact normalized estimates
\begin{eqaligned}
\frac{|\upsilon_i^k|}{\rho^k}
&\leq\theta_1(1+\varepsilon)-1=\frac{1009}{120900},\\
\frac{\Gamma_i^k}{\rho^k}
    &=\frac{|\upsilon_i^k|+\zeta_i^k}{\rho^k}
    \leq\frac{1009}{120900}+\frac1{1209}
    =\frac{1109}{120900}.
\end{eqaligned}
The same clauses bound the lagged drift constant by $b_i^k/\rho^k\leq3\varepsilon=1/403$: for $i\in\IO$ it equals $H_i/\rho^k\leq\varepsilon$, and for $i\in\IC$ it equals $(2\zeta_i^k+\tau_i^k)/\rho^k\leq3\varepsilon$.  Put
\[
\vartheta=\frac1{403},\qquad
\varrho_k=\frac{\rho^{k+1}}{\rho^k}\in[1,2].
\]
Then $\widehat\Gamma_i^{k+1}=|\upsilon_i^k|+b_i^{k+1}$ gives
\[
\frac{\widehat\Gamma_i^{k+1}}{\rho^{k+1}}\leq\frac{1009}{120900\varrho_k}+\vartheta.
\]
The exact inequality
\[
\left(\frac{1009}{120900}\right)^2-2\vartheta^2=\frac{838081}{14616810000}>0
\]
implies $(1009/120900)^2/\varrho_k\geq(1009/120900)^2/2>\vartheta^2$.  Hence the elementary identity
\[
\left(\frac{1009}{120900}+\vartheta\right)^2
-\varrho_k\left(\frac{1009}{120900\varrho_k}+\vartheta\right)^2
=(\varrho_k-1)\left[\frac{(1009/120900)^2}{\varrho_k}-\vartheta^2\right]\geq0
\]
gives, with $1009/120900+\vartheta=1309/120900$,
\[
\frac{(\widehat\Gamma_i^{k+1})^2}{\rho^k\rho^{k+1}}
\leq\left(\frac{1309}{120900}\right)^2.
\]
Since $\kappa_\Psi\leq5/2$, the definitions of $D_{y,i}^k$ and $D_x^k$ now give
\begin{eqaligned}
\frac{D_{y,i}^k}{\rho^k}
&\geq\frac{178630370053}{353726802000}>\frac12,\\
\frac{D_x^k}{\rho^k}
&\geq\frac{1+2\theta_2}{2}-25=\frac12.
\end{eqaligned}
Thus~\prettyref{eq:p4}--\prettyref{eq:p5} hold with $d_0=1/2$.  For the terminal residual coefficient,
\[
\frac{\widetilde C_{r,i}^k}{\rho^k}
\leq\frac{\theta_1(1+\varepsilon)}2+\frac{9\varepsilon}{2}
=\frac{122809}{241800}<1,
\]
which proves~\prettyref{eq:p3}.  To prove~\prettyref{eq:p6} without separately rounding $C_{d,i}^k$ and $D_{y,i}^k$, put $s=\tau_i^k/\rho^k\in[0,\varepsilon]$.  The exact bounds above yield
\[
\frac{D_{y,i}^k-C_{d,i}^k}{\rho^k}
\geq\frac{(\theta_1-1)(1+s)}2-\frac{3s}{2}
-\frac{25}{2}\left[
\left(\frac{1109}{120900}\right)^2+
\left(\frac{1309}{120900}\right)^2
\right].
\]
The right side is affine in $s$, and at its two endpoints $s=0$ and $s=\varepsilon$ it equals, respectively,
\[
\frac{27393814}{22107925125}>0,
\qquad
\frac{881}{584672400}>0.
\]
It is therefore nonnegative throughout the interval, proving~\prettyref{eq:p6}.  Since $\rho^k/\rho^{k-1}\leq2$, the left side of~\prettyref{eq:p2} is at most $-1/(2\rho^{k-1})$.  Condition~\prettyref{eq:p7} follows from $G=m^2L^0/2$ and $L^1=L^0(1+\xi)\geq L^0$: indeed $\sqrt{2G/L^1}=m\sqrt{L^0/L^1}\leq m$.

It remains to prove~\prettyref{eq:p8}.  The branch $L^0\geq12S_0/m^2$ gives
\[
S_0\leq\frac{m^2L^0}{12}=\frac G6.
\]
Next, $R(L^1)=\eta G_\psi/(2\rho^1)\leq\eta\widehat G_0/(2\rho^0)$, using $\rho^1\geq\rho^0$ and $G_\psi=\sum_{i\in\IO}(C_i^\psi)^2\leq\sum_{i\in\IO}\widehat C_i^2=\widehat G_0$, which holds because $\widehat C_i=D_i+C_i^\psi\geq C_i^\psi$.  Since $m_i^1\geq\theta_1\rho^1\geq\theta_1\rho^0$, $\rho^1\leq2\rho^0$, and $\rho^2\leq3\rho^0$ (using $\xi\leq1$ and $2^p\leq2$), the linear and three quadratic initialization terms are bounded by
\[
\left(\frac1{\theta_1}+\frac{17}{4\theta_1^2}\right)
\frac{\widehat G_0}{\rho^0}.
\]
Furthermore, $5\kappa_\Psi\leq25/2$, $(\widehat\Gamma_i^2)^2\leq(1309/120900)^2\rho^1\rho^2$, and $1+\theta_2=26$, so the initial $\Theta^2$ term is at most
\[
\frac{25}{2\theta_1^2}
\left[\left(\frac{1309}{120900}\right)^2+\frac1{26^2}\right]
\frac{\widehat G_0}{\rho^0}.
\]
The exact initialization coefficient satisfies
\[
\frac1{\theta_1}+\frac{17}{4\theta_1^2}
+\frac{25}{2\theta_1^2}
\left[\left(\frac{1309}{120900}\right)^2+\frac1{26^2}\right]
=\frac{112940321939035250}{21723216943256361}<\frac{11}{2}.
\]
Consequently the initialization and smoothing contribution is at most
\[
\left(\frac\eta2+\frac{11}{2}\right)
\frac{\widehat G_0}{\rho^0}\leq\frac G6
\]
by the penultimate clause in the definition of $\rho^0$.  Finally,~\prettyref{eq:drift-sum-a} and $\kappa_\Psi\leq5/2$ give the exact bound
\[
\kappa_\Psi\Sigma_U\leq
\frac{p^2H_p(\xi)}{\rho^0}
\left(\frac{25}{2}G_\psi+\frac{25}{4}GL^0\bar M\right)
\leq\frac G6
\]
by the last clause in $\rho^0$.  Thus the three allocated contributions on the left side of~\prettyref{eq:p8} are each at most $G/6$, and their sum is at most $G/2$.
\end{proof}

\begin{theorem}[Feasibility and one-step Lyapunov decrease]
\label{thm:feasibility}
Under~\prettyref{eq:parameters}, $f_i(y_i^k)\geq0$ for every $i\in\IC$ and $k\geq1$.  Moreover, for $k\geq2$,
\begin{eqaligned}
\label{eq:Lyapunov-decrease}
\Psi^{k+1}-\Psi^k
\leq-\frac{d_0}{\rho^k}\|\Delta\lambda^{k+1}\|^2
-d_0\rho^k\|d^k\|^2-d_0\rho^k\|A\Delta x^k\|^2
+\kappa_\Psi U^k.
\end{eqaligned}
\end{theorem}
\begin{proof}
The initial split iterate is wall-feasible: $y_i^1=A_ix^1$ and $f_i(A_ix^1)\geq m>0$ by~\prettyref{ass:feasible}.  By~\prettyref{lem:initial-bound}, every constraint block satisfies $y_i^2=y_i^1$, so $y^2$ is wall-feasible as well.

For the first unverified point, set $k=2$ in \prettyref{eq:Lyapunov-lower-bound}.  Conditions~\prettyref{eq:p3}, \prettyref{eq:p4}, and~\prettyref{eq:p6} discard its current dual, $x$, and $y$ terms.  The unmatched lagged term is
\[
\frac{\rho^2}{4(\rho^1)^2}\|\Delta\lambda^2\|^2 =\frac{\rho^2}{4}\|y^2-Ax^2\|^2 \leq\frac{\rho^2}{4}\sum_{i\in\IO}\frac{\widehat C_i^2}{(m_i^1)^2},
\]
where the first equality uses $\Delta\lambda^2=\lambda^2=\rho^1(y^2-Ax^2)$, i.e.\ the multiplier update at $k=1$ together with $\lambda^1=0$.  Combine this with~\prettyref{eq:Psi-2-bound}, $U^2\leq\Sigma_U$, and~\prettyref{eq:p8}; this yields $W_3\leq G/2$. By~\prettyref{eq:R3}, the corrected objective part of $W_3$ is nonnegative, so $\tfrac12\sum_{i\in\IC}g(f_i(y_i^3),L^2)\leq G/2$; since every wall term is nonnegative, each one individually satisfies $g(f_i(y_i^3),L^2)\leq G$.  Since $g(s,L)>G$ for $s<0$ by~\prettyref{lem:added_algebraic_3}, all $f_i(y_i^3)\geq0$.

Suppose now that $y^2,\ldots,y^K$ are wall-feasible, $K\geq3$.  For $k=2,\ldots,K-1$, continuation comparison is legal and \prettyref{eq:Lyapunov}, \prettyref{eq:p4}, and \prettyref{eq:p5} give
\[
\Psi^{k+1}-\Psi^k\leq -\frac{1}{\rho^k}\|\Delta\lambda^{k+1}\|^2+\kappa_\Psi U^k.
\]
Apply~\prettyref{eq:Lyapunov-lower-bound} at the unverified terminal index $K$.  Conditions~\prettyref{eq:p3}, \prettyref{eq:p4}, and \prettyref{eq:p6} absorb its current terms, while~\prettyref{eq:p2} absorbs its lagged dual term $\rho^K\|\Delta\lambda^K\|^2/(4(\rho^{K-1})^2)$ against the matching term $-\|\Delta\lambda^K\|^2/\rho^{K-1}$, which the telescoped decrease above contributes at $k=K-1$. Consequently
\[
W_{K+1}\leq\Psi^2+\kappa_\Psi\sum_{k=2}^KU^k\leq G/2,
\]
because, by~\prettyref{eq:Psi-2-bound} and $\sum_{k=2}^KU^k\leq\Sigma_U$, the middle expression is dominated by the left side of~\prettyref{eq:p8}.  Applying the argument of the previous paragraph to $W_{K+1}\leq G/2$ then proves $f_i(y_i^{K+1})\geq0$ for every $i\in\IC$, that is, $y^{K+1}$ is wall-feasible.  Since $y^1$, $y^2$, and $y^3$ are wall-feasible and this step is available for every $K\geq3$, strong induction shows that $y^k$ is wall-feasible for every $k\geq1$, which is the first assertion of the theorem.

Once wall-feasibility is known, $\Psi^{k+1}\leq\bar\Psi^{k+1}$ is always legal. Combining~\prettyref{eq:Lyapunov} with $d_0=1/2<1$, \prettyref{eq:p4}, and~\prettyref{eq:p5} gives \prettyref{eq:Lyapunov-decrease}.
\end{proof}

\subsection{Stationarity and Rates}
\label{sub:stationarity-rates}
\paragraph{Original-problem certificates.}
\begin{define}[Original-problem residuals]
\label{def:criteria}
In all three notions, $\breve y_i=y_i$ for every constraint block $i\in\IC$. These are nearby-point certificates: stationarity is evaluated at the reconstructed block points $\breve y_i$ rather than at the common original point $A_ix$, and the two agree only up to the consensus term $\|Ax-\breve y\|$. The weak $\epsilon$-KKT and $\epsilon$-KKT conditions are, respectively,
\begin{equation}
\begin{aligned}
\text{weak $\epsilon$-KKT:}\quad
&\left\{
\begin{aligned}
&\|Ax-\breve y\|\leq\epsilon,\\
&\left\|\sum_{i\in\IO}A_i^Tv_i+
\sum_{i\in\IC}\iota_iA_i^T\nabla f_i(y_i)+C^T\zeta\right\|\leq\epsilon,\\
&v_i\in\partial\varphi_i(\breve y_i),\ Cx=d,\ f_i(y_i)\geq0,\ \iota_i\leq0;
\end{aligned}
\right.\\[1ex]
\text{$\epsilon$-KKT:}\quad
&\left\{
\begin{aligned}
&\|Ax-\breve y\|\leq\epsilon,\\
&\left\|\sum_{i\in\IO}A_i^Tv_i+
\sum_{i\in\IC}\iota_iA_i^T\nabla f_i(y_i)+C^T\zeta\right\|\leq\epsilon,\\
&v_i\in\partial\varphi_i(\breve y_i),\ Cx=d,\ f_i(y_i)\geq0,\ \iota_i\leq0,\\
&|\iota_if_i(y_i)|\leq\epsilon\quad(i\in\IC).
\end{aligned}
\right.
\end{aligned}
\end{equation}
The $\epsilon$-b-FJ conditions are
\begin{equation}
\label{eq:eps-bFJ}
\left\{
\begin{aligned}
&\|Ax-\breve y\|\leq\epsilon,\\
&\left\|\alpha\sum_{i\in\IO}A_i^Tv_i+
\sum_{i\in\IC}\bar\iota_iA_i^T\nabla f_i(y_i)+C^T\bar\zeta\right\|\leq\epsilon,\\
&v_i\in\partial\varphi_i(\breve y_i),\ Cx=d,\ f_i(y_i)\geq0,\
\alpha>0,\ \bar\iota_i\leq0,\quad\alpha+\|\bar\iota\|_1=1,\\
&|\bar\iota_if_i(y_i)|\leq\epsilon,\qquad
\bar\iota_i\max\{f_i(y_i)-\epsilon,0\}=0\quad(i\in\IC).
\end{aligned}
\right.
\end{equation}
\end{define}

For the candidate certificate tuple generated by step $k$, let $\varpi^{k+1}$ be the multiplier of the $x^{k+1}$ subproblem and set
\begin{eqaligned}
\zeta^{k+1}&=\rho^k\varpi^{k+1},&
\alpha^{k+1}&=\frac1{1+\|\iota^{k+1}\|_1},\\
\bar\iota^{k+1}&=\alpha^{k+1}\iota^{k+1},&
\bar\zeta^{k+1}&=\alpha^{k+1}\zeta^{k+1},\\
\omega^{k+1}&\defeq\beta^k+\frac12\max_{i\in\IC}L_i\|d_i^k\|^2.
\end{eqaligned}
By its definition, $\alpha^{k+1}\in(0,1]$. Define
\begin{eqaligned}
\Xi^{k+1}&\defeq\sum_{i\in\IO}A_i^Tv_i^{k+1}
+\sum_{i\in\IC}\iota_i^{k+1}A_i^T\nabla f_i(y_i^{k+1})+C^T\zeta^{k+1},\\
\label{eq:crit}
\CritWKKT^{k+1}&\defeq\|Ax^{k+1}-\breve y^{k+1}\|+\|\Xi^{k+1}\|,\\
\CritKKT^{k+1}&\defeq\CritWKKT^{k+1}+
\sum_{i\in\IC}|\iota_i^{k+1}f_i(y_i^{k+1})|,\\
\CritbFJ^{k+1}&\defeq\|Ax^{k+1}-\breve y^{k+1}\|
+\alpha^{k+1}\|\Xi^{k+1}\|
+\sum_{i\in\IC}|\bar\iota_i^{k+1}f_i(y_i^{k+1})|+\omega^{k+1}.
\end{eqaligned}

\begin{lemma}[Lyapunov lower bound]
\label{lem:lower-bound}
Under~\prettyref{eq:parameters} and the ordinary polynomial schedule~\prettyref{eq:IP-schedules}, for every $k\geq2$,
\begin{eqaligned}
\Psi^k\geq\underline\Psi\defeq
-\frac{\widehat G_0}{2\rho^0}-\kappa_\Psi\Sigma_U> -\infty.
\end{eqaligned}
\end{lemma}
\begin{proof}
By~\prettyref{eq:R3}, with $\mu^k=\mu_{\rm IP}(L^k)$ by~\prettyref{eq:IP-schedules} and~\prettyref{eq:rho-L-ratio}, and the $\widehat C_i$-Lipschitz property of the objective surrogate,
\[
\sum_{i\in\IO}\bar f_i(y_i^k,L^k,\mu^k)+R(L^k) \geq-\sum_{i\in\IO}\widehat C_i\|y_i^k-A_ix^k\|.
\]
For each $i\in\IO$, Young's inequality completes the square between $-\widehat C_i\|y_i^k-A_ix^k\|$ and the corresponding objective-block quadratic penalty, contributing at worst $-\widehat C_i^2/(2\rho^k)\geq-\widehat C_i^2/(2\rho^0)$.  The constraint-block quadratic penalties are nonnegative and are discarded, as are the nonnegative wall and $\Theta^k$ terms.  Using the multiplier update then gives
\[
\Psi^k\geq(\lambda^k)^T(y^k-Ax^k)-\frac{\widehat G_0}{2\rho^0} =\frac{(\lambda^k)^T\Delta\lambda^k}{\rho^{k-1}} -\frac{\widehat G_0}{2\rho^0}.
\]
By~\prettyref{eq:Lyapunov-def}, $\kappa_\Psi=(\xi+4)/2\in(2,5/2]$.  Moreover, $U^j\geq0$ by~\prettyref{eq:U}, and~\prettyref{lem:drift-sum} gives the finite tail $0\leq\sum_{j=k}^\infty U^j\leq\Sigma_U$. Define the exact nonnegative drift-tail correction $\widehat\Psi^k=\Psi^k+\kappa_\Psi\sum_{j=k}^\infty U^j$. By~\prettyref{eq:Lyapunov-decrease}, $\widehat\Psi^k$ is nonincreasing. If it ever fell below $-\widehat G_0/(2\rho^0)$, its tail sum relative to that level would diverge to $-\infty$.  On the other hand, the polarization identity $(\lambda^k)^T\Delta\lambda^k=(\|\lambda^k\|^2-\|\lambda^{k-1}\|^2+\|\Delta\lambda^k\|^2)/2$ and $\rho^{k-1}\leq\rho^k$ give
\[
\widehat\Psi^k+\frac{\widehat G_0}{2\rho^0} \geq\frac{\|\lambda^k\|^2-\|\lambda^{k-1}\|^2} {2\rho^{k-1}} \geq\frac{\|\lambda^k\|^2}{2\rho^k} -\frac{\|\lambda^{k-1}\|^2}{2\rho^{k-1}}.
\]
The partial sums telescope; discarding the nonnegative terminal term leaves the finite initial endpoint as a uniform lower bound, a contradiction.  Hence $\widehat\Psi^k\geq-\widehat G_0/(2\rho^0)$, and subtracting its nonnegative drift tail, bounded above by $\kappa_\Psi\Sigma_U$, proves the claim.
\end{proof}

\begin{corollary}[Global movement budget]
\label{cor:movement-budget}
Under~\prettyref{eq:parameters} and the ordinary polynomial schedule~\prettyref{eq:IP-schedules}, define
\begin{eqaligned}
\label{eq:global-budget}
\mathcal B\defeq\frac{\Psi^2-\underline\Psi+\kappa_\Psi\Sigma_U}{d_0}.
\end{eqaligned}
Then $\mathcal B<\infty$ and
\begin{eqaligned}
\label{eq:movement-budget}
\sum_{k=2}^\infty\left[
\frac{\|\Delta\lambda^{k+1}\|^2}{\rho^k}
+\rho^k\|d^k\|^2+\rho^k\|A\Delta x^k\|^2
\right]\leq\mathcal B.
\end{eqaligned}
\end{corollary}
\begin{proof}
The constant $\mathcal B$ is finite because~\prettyref{eq:Psi-2-bound} bounds $\Psi^2$ above, $\underline\Psi$ and $\Sigma_U$ are finite, and $d_0>0$.  Sum~\prettyref{eq:Lyapunov-decrease} from $k=2$ to $K$, use $\Psi^{K+1}\geq\underline\Psi$ from~\prettyref{lem:lower-bound}, and bound $\sum_{k=2}^KU^k\leq\Sigma_U$ by~\prettyref{lem:drift-sum}.  The result is
\[
d_0\sum_{k=2}^K\left[ \frac{\|\Delta\lambda^{k+1}\|^2}{\rho^k} +\rho^k\|d^k\|^2+\rho^k\|A\Delta x^k\|^2 \right] \leq\Psi^2-\underline\Psi+\kappa_\Psi\Sigma_U.
\]
After division by $d_0>0$, the partial sums on the left of~\prettyref{eq:movement-budget} are nondecreasing because every summand is nonnegative, and they are uniformly bounded above by the finite constant $\mathcal B$.  They therefore converge to a finite limit; letting $K\to\infty$ proves~\prettyref{eq:movement-budget}.
\end{proof}

\begin{lemma}[Bounded rescaled multipliers]
\label{lem:scaled-multiplier}
Under~\prettyref{eq:parameters} and the ordinary polynomial schedule~\prettyref{eq:IP-schedules}, define
\begin{eqaligned}
c_\upsilon&\defeq\theta_1-1+\frac{\theta_1}{1209}
=\frac{1009}{120900},\\
C_\lambda&\defeq
\frac{\sqrt{\widehat G_0}}{\sqrt{\rho^0}}
+\sqrt{\frac{2GL^0\bar M}{\rho^0}}
+(1+c_\upsilon)\sqrt{\mathcal B}.
\end{eqaligned}
Then, for every $j\geq1$,
\begin{eqaligned}
\label{eq:scaled-multiplier-same-index}
\frac{\|\lambda^j\|}{\sqrt{\rho^j}}&\leq C_\lambda,&
\frac{\|\lambda^j\|^2}{\rho^j}&\leq C_\lambda^2,&
\frac{\|\lambda^j\|}{\rho^j}&\leq
\frac{C_\lambda}{\sqrt{\rho^j}}.
\end{eqaligned}
If $p>0$, then
\begin{eqaligned}
\label{eq:scaled-multiplier-decay}
\frac{\|\lambda^j\|}{\rho^j}
=\mathcal O((\rho^j)^{-1/2})
=\mathcal O(j^{-p/2})\longrightarrow0.
\end{eqaligned}
When $p=0$, the sequence remains bounded, but no decay is asserted.
\end{lemma}
\begin{proof}
We first bound the stacked force $q^k$.  For $i\in\IO$, the current force formula, the bound $\|\nabla h_i(y_i^k)\|\leq D_i$, and \prettyref{eq:moreau-bound} give $\|q_i^k\|\leq D_i+C_i^\psi=\widehat C_i$.  For $i\in\IC$, \prettyref{eq:band-force} and $\|\nabla f_i(y_i^k)\|\leq M_i$ give $\|q_i^k\|\leq\sqrt{2L^kG}\,M_i$.  Since $q^k$ is the Euclidean stacking of these blocks,
\[
\|q^k\|^2=\sum_{i\in\IAll}\|q_i^k\|^2
\leq\widehat G_0+2L^kG\bar M.
\]
Thus $\sqrt{a+b}\leq\sqrt a+\sqrt b$ for $a,b\geq0$ gives
\[
\|q^k\|\leq\sqrt{\widehat G_0}+\sqrt{2L^kG\bar M}.
\]

By definition, $\upsilon_i^k=m_i^k-\rho^k=(\theta_1-1)\rho^k+\theta_1\tau_i^k$.  The bound $0\leq\tau_i^k/\rho^k\leq1/1209$ from~\prettyref{lem:parameter} therefore gives $|\upsilon_i^k|/\rho^k\leq\theta_1-1+\theta_1/1209=c_\upsilon$ for every block.  For the block-diagonal matrix $\Upsilon^k$ of~\prettyref{sub:lyapunov-descent}, this yields
\[
\|\Upsilon^kd^k\|^2 =\sum_{i\in\IAll}|\upsilon_i^k|^2\|d_i^k\|^2 \leq c_\upsilon^2(\rho^k)^2\|d^k\|^2.
\]
Stacking \prettyref{eq:block-dual-identity} over the blocks gives
\[
\lambda^{k+1}=-q^k-\Upsilon^kd^k-\rho^kA\Delta x^k.
\]
For $k\geq2$, each nonnegative summand in \prettyref{eq:movement-budget} is at most $\mathcal B$.  Hence $\sqrt{\rho^k}\|d^k\|\leq\sqrt{\mathcal B}$ and $\sqrt{\rho^k}\|A\Delta x^k\|\leq\sqrt{\mathcal B}$.  Dividing the preceding identity by $\sqrt{\rho^k}$, applying the triangle inequality and the block-diagonal estimate, and then using~\prettyref{eq:rho-L-ratio} gives
\begin{align*}
\frac{\|\lambda^{k+1}\|}{\sqrt{\rho^k}}
&\leq\frac{\|q^k\|}{\sqrt{\rho^k}}
+c_\upsilon\sqrt{\rho^k}\|d^k\|
+\sqrt{\rho^k}\|A\Delta x^k\|\\
&\leq\frac{\sqrt{\widehat G_0}}{\sqrt{\rho^k}}
+\sqrt{\frac{2GL^k\bar M}{\rho^k}}
+(1+c_\upsilon)\sqrt{\mathcal B}\\
&\leq\frac{\sqrt{\widehat G_0}}{\sqrt{\rho^0}}
+\sqrt{\frac{2GL^0\bar M}{\rho^0}}
+(1+c_\upsilon)\sqrt{\mathcal B}=C_\lambda.
\end{align*}
Here the last line also uses $\rho^k\geq\rho^0$.

It remains to treat $k=1$, because~\prettyref{eq:movement-budget} starts at $k=2$.  \prettyref{lem:initial-bound} gives $d_i^1=0$ for $i\in\IC$ and $\|d_i^1\|\leq\widehat C_i/m_i^1$ for $i\in\IO$.  Since $m_i^1\geq\theta_1\rho^1$,
\[
\|d^1\|\leq\frac{\sqrt{\widehat G_0}}{\theta_1\rho^1}.
\]
The second endpoint estimate $\|y^2-Ax^2\|\leq\|d^1\|$ is also exactly the conclusion of~\prettyref{lem:initial-bound}; the multiplier update together with the initialization $\lambda^1=0$ then gives $\lambda^2=\rho^1(y^2-Ax^2)$.  Therefore
\[
\frac{\|\lambda^2\|}{\sqrt{\rho^1}} \leq\frac{\sqrt{\widehat G_0}}{\theta_1\sqrt{\rho^1}} \leq\frac{\sqrt{\widehat G_0}}{\sqrt{\rho^0}} \leq C_\lambda,
\]
where the middle inequality uses $\theta_1>1$ and $\rho^1\geq\rho^0$.  Hence, for every $k\geq1$,
\begin{eqaligned}
\label{eq:scaled-multiplier-step}
\frac{\|\lambda^{k+1}\|}{\sqrt{\rho^k}}\leq C_\lambda.
\end{eqaligned}

For $j\geq2$, apply~\prettyref{eq:scaled-multiplier-step} with $k=j-1$ and use schedule monotonicity~\prettyref{eq:schedule-monotone}:
\[
\frac{\|\lambda^j\|}{\sqrt{\rho^j}} \leq\frac{\|\lambda^j\|}{\sqrt{\rho^{j-1}}}\leq C_\lambda.
\]
The same inequality is immediate for $j=1$ because $\lambda^1=0$.  Squaring the first of the three same-index bounds in~\prettyref{eq:scaled-multiplier-same-index} gives $\|\lambda^j\|^2/\rho^j\leq C_\lambda^2$; dividing that first bound by $\sqrt{\rho^j}$ gives $\|\lambda^j\|/\rho^j\leq C_\lambda/\sqrt{\rho^j}$.  Finally, for $p>0$, the schedule $\rho^j=\rho^0(1+\xi j^p)$ proves \prettyref{eq:scaled-multiplier-decay}.  At $p=0$ it is constant on positive indices, so the proved bound remains uniform without implying decay.
\end{proof}

For the step-$k$ candidate certificate tuple, the objective-block reconstruction and its transfer bounds are given in~\prettyref{eq:reconstruction}--\prettyref{eq:recon-error}.  For constraint blocks,
\begin{eqaligned}
\|\iota_i^{k+1}\nabla f_i(y_i^{k+1})-q_i^k\|
\leq\sqrt{2L^kG}\,L_i\|d_i^k\|.
\end{eqaligned}
Also, $\iota_i^{k+1}\neq0$ implies $\ell_i^k(y_i^{k+1})<\beta^k$, and smoothness of $f_i$ gives
\begin{eqaligned}
\label{eq:band-support}
\iota_i^{k+1}\neq0\quad\Longrightarrow\quad
0\leq f_i(y_i^{k+1})<\beta^k+\frac{L_i}{2}\|d_i^k\|^2
\leq\omega^{k+1}.
\end{eqaligned}

Define the reconstructed-force gaps $u_i^k\defeq v_i^{k+1}-q_i^k$ for $i\in\IO$ and $u_i^k\defeq\iota_i^{k+1}\nabla f_i(y_i^{k+1})-q_i^k$ for $i\in\IC$, stacked as $u^k\defeq(u_i^k)_{i\in\IAll}$.  Combining~\prettyref{eq:opt-yi} and \prettyref{eq:opt-x} gives the exact identity
\begin{eqaligned}
\label{eq:residual-stationary}
\Xi^{k+1}=A^T[u^k-\Upsilon^kd^k]
-(1+\theta_2)\rho^kA^TA\Delta x^k.
\end{eqaligned}
Set
\[
r_*\defeq\max\left\{H_{\max}/\rho^0,
\sqrt{2GL^0}L_{\max}/\rho^0\right\},
\qquad c_y\defeq\frac{1009}{120900}
+\frac{242909}{121000}r_*.
\]
The transfer bounds and $c_\upsilon=1009/120900$ imply
\begin{eqaligned}
\|\Xi^{k+1}\|
\leq\|A\|_2\left[c_y\rho^k\|d^k\|
+26\rho^k\|A\Delta x^k\|+\mu^kJ_H\right].
\end{eqaligned}
The multiplier update and reconstruction similarly give
\begin{eqaligned}
\|Ax^{k+1}-\breve y^{k+1}\|
\leq\frac{\|\Delta\lambda^{k+1}\|}{\rho^k}+\mu^kC_\psi.
\end{eqaligned}
With $Q_\psi\defeq C_\psi+\|A\|_2J_H$, where $J_H$ is fixed in~\prettyref{eq:parameters}, it follows that
\begin{eqaligned}
\label{eq:wKKT-pointwise}
\CritWKKT^{k+1}
\leq\frac{\|\Delta\lambda^{k+1}\|}{\rho^k}
+\|A\|_2c_y\rho^k\|d^k\|
+26\|A\|_2\rho^k\|A\Delta x^k\|+\mu^kQ_\psi.
\end{eqaligned}

Normalization gives $\iota_i^{k+1},\bar\iota_i^{k+1}\leq0$ and $\alpha^{k+1}+\|\bar\iota^{k+1}\|_1=1$.  From \prettyref{eq:band-support} and feasibility,
\begin{eqaligned}
\sum_{i\in\IC}|\bar\iota_i^{k+1}f_i(y_i^{k+1})|
\leq\omega^{k+1}\|\bar\iota^{k+1}\|_1\leq\omega^{k+1}.
\end{eqaligned}
Thus
\begin{eqaligned}
\CritbFJ^{k+1}\leq\CritWKKT^{k+1}+2\omega^{k+1}.
\end{eqaligned}
If $\CritbFJ^{k+1}\leq\epsilon$, then $\omega^{k+1}\leq\epsilon$ and \prettyref{eq:band-support} proves the exact final support line of \prettyref{eq:eps-bFJ}.  Therefore the residual certifies the unchanged original b-FJ definition.

Recall $M_F=\|A\|_2\sqrt{\widehat G_0}$ from~\prettyref{eq:parameters}.
\begin{lemma}[Uniform Band-MFCQ bounds the objective multiplier]
\label{lem:alpha-lower}
Under~\prettyref{ass:MFCQ}, let $(x,y,\breve y,v,\alpha,\bar\iota,\bar\zeta)$ be an $\epsilon$-b-FJ point in the sense of~\prettyref{def:criteria}, with $\epsilon\leq\min\{\delta_{\rm cq},\gamma_{\rm cq}/2\}$.  Then
\begin{eqaligned}
\label{eq:alpha-lower}
\alpha\geq\frac{\gamma_{\rm cq}}{2(\gamma_{\rm cq}+M_F)}.
\end{eqaligned}
\end{lemma}
\begin{proof}
Let $a_i=-\bar\iota_i\geq0$.  Since $\bar\iota_i\leq0$, the b-FJ normalization gives $\alpha+\sum_{i\in\IC}a_i=1$.  Thus $\alpha\in(0,1]$ and $\sum_{i\in\IC}a_i=1-\alpha$.  By~\prettyref{def:criteria}, $\breve y_i=y_i$ for every $i\in\IC$, and therefore
\[
\max_{i\in\IC}\|A_ix-y_i\|\leq\|Ax-\breve y\| \leq\epsilon\leq\delta_{\rm cq}.
\]
For every index with $a_i>0$, the sign and support conditions in~\prettyref{def:criteria} also give $0\leq f_i(y_i)\leq\epsilon\leq\delta_{\rm cq}$.  Hence~\prettyref{ass:MFCQ} supplies a single unit $d_{\rm cq}\in\ker C$ for which~\prettyref{eq:cq} holds at every such index.  Write the b-FJ stationarity residual as
\[
s=\alpha\sum_{i\in\IO}A_i^Tv_i-
\sum_{i\in\IC}a_iA_i^T\nabla f_i(y_i)+C^T\bar\zeta,
\qquad \|s\|\leq\epsilon.
\]
For this direction, $Cd_{\rm cq}=0$, so the equality-multiplier term vanishes.  Moreover, the exact sum rule $\partial\varphi_i=\nabla h_i+\partial\psi_i$ and~\prettyref{eq:objective-derived-bounds} give $\|v_i\|\leq D_i+C_i^\psi=\widehat C_i$.  Zero-padding the vector $(v_i)_{i\in\IO}$ on $\IC$ therefore gives
\[
\left\|\sum_{i\in\IO}A_i^Tv_i\right\|
\leq\|A\|_2\left(\sum_{i\in\IO}\|v_i\|^2\right)^{1/2}
\leq\|A\|_2\sqrt{\widehat G_0}=M_F.
\]
On the other hand, $\nabla f_i(y_i)^TA_id_{\rm cq}\geq\gamma_{\rm cq}$ whenever $a_i>0$; the terms with $a_i=0$ vanish.  Contracting with $d_{\rm cq}$ and using $-d_{\rm cq}^Ts\leq\epsilon$ therefore yields
\[
\gamma_{\rm cq}(1-\alpha)
\leq\sum_{i\in\IC}a_i\nabla f_i(y_i)^TA_id_{\rm cq}
=\alpha d_{\rm cq}^T\sum_{i\in\IO}A_i^Tv_i-d_{\rm cq}^Ts
\leq\alpha M_F+\epsilon.
\]
Rearrangement gives $\alpha\geq(\gamma_{\rm cq}-\epsilon)/(\gamma_{\rm cq}+M_F)$, and $\epsilon\leq\gamma_{\rm cq}/2$ proves~\prettyref{eq:alpha-lower}.
\end{proof}

\begin{theorem}[Exact finite-time original-problem rates]
\label{thm:rate-master}
Under~\prettyref{eq:parameters} and the ordinary polynomial schedule~\prettyref{eq:IP-schedules}, let $\mathcal B<\infty$ be the movement constant supplied by~\prettyref{cor:movement-budget}.  For every $K\geq3$, let
\begin{eqaligned}
R_K&=\sum_{k=2}^K\rho^k,&I_K&=\sum_{k=2}^K(\rho^k)^{-1},&
B_K&=\sum_{k=2}^K\beta^k,\\
\mathcal Q&=\|A\|_2^2(c_y^2+26^2).
\end{eqaligned}
Then
\begin{eqaligned}
\label{eq:rate-explicit}
F_{\rm W}(K)&\defeq\frac{\sqrt{\mathcal B}}{K-1}
\sqrt{I_K+\mathcal Q R_K}+\frac{\eta Q_\psi}{K-1}I_K,\\
F_{\rm bFJ}(K)&\defeq F_{\rm W}(K)+\frac{2B_K}{K-1}
+\frac{L_{\max}\mathcal B}{\rho^2(K-1)},\\
\min_{j=3,\ldots,K+1}\CritWKKT^j&\leq F_{\rm W}(K),\\
\min_{j=3,\ldots,K+1}\CritbFJ^j&\leq F_{\rm bFJ}(K).
\end{eqaligned}
The second minimum is attained at one selected certificate-bearing iterate whose tuple is an $F_{\rm bFJ}(K)$-b-FJ point in the sense of \prettyref{eq:eps-bFJ}.  For $0<p<1$,
\begin{eqaligned}
\label{eq:rate}
F_{\rm W}(K)&=\mathcal O(K^{(p-1)/2})+\mathcal O(K^{-p}),\\
F_{\rm bFJ}(K)&=\mathcal O(K^{(p-1)/2})+\mathcal O(K^{-p})
+\mathcal O(K^{-p/2})+\mathcal O(K^{-1}).
\end{eqaligned}
Here the hidden constants may depend on the fixed $p\in(0,1)$, $\xi$, and problem and algorithm data, but not on $K$; these estimates are not uniform as $p$ or $\xi$ approaches an endpoint. Thus $p=1/3$ gives w-KKT rate $\mathcal O(K^{-1/3})$, while $p=1/2$ gives b-FJ rate $\mathcal O(K^{-1/4})$.
\end{theorem}
\begin{proof}
Sum~\prettyref{eq:wKKT-pointwise} for $k=2,\ldots,K$. For each $k$, write the three movement terms as the inner product of
\[
\left(\frac{\|\Delta\lambda^{k+1}\|}{\sqrt{\rho^k}},
\sqrt{\rho^k}\|d^k\|,
\sqrt{\rho^k}\|A\Delta x^k\|\right)
\quad\hbox{and}\quad
\left(\frac1{\sqrt{\rho^k}},
\|A\|_2c_y\sqrt{\rho^k},
26\|A\|_2\sqrt{\rho^k}\right).
\]
A single Cauchy--Schwarz inequality over all indices and all three channels, followed by~\prettyref{eq:movement-budget}, gives
\[
\sum_{k=2}^K\left[ \frac{\|\Delta\lambda^{k+1}\|}{\rho^k} +\|A\|_2c_y\rho^k\|d^k\| +26\|A\|_2\rho^k\|A\Delta x^k\|\right]
\leq\sqrt{\mathcal B}\sqrt{I_K+\mathcal Q R_K}.
\]
Since $\sum_{k=2}^K\mu^k=\eta I_K$, averaging yields the bound $F_{\rm W}(K)$.

Next, $\beta^k\geq0$ by~\prettyref{eq:band}, and $2\omega^{k+1}\leq2\beta^k+L_{\max}\|d^k\|^2$.  Since $\rho^k\geq\rho^2$ for $k\geq2$, the movement budget gives $\sum_{k=2}^K\|d^k\|^2\leq\mathcal B/\rho^2$.  Consequently,
\[
\sum_{k=2}^K\CritbFJ^{k+1}
\leq\sum_{k=2}^K\CritWKKT^{k+1}
+2\sum_{k=2}^K\omega^{k+1}
\leq(K-1)F_{\rm bFJ}(K).
\]
Choose $j_*\in\{3,\ldots,K+1\}$ minimizing the entire b-FJ residual. Then $\CritbFJ^{j_*}\leq F_{\rm bFJ}(K)$ at one common selected certificate-bearing iterate; its consensus, normalized stationarity, complementarity, and nonnegative $\omega^{j_*}$ terms are all bounded there, while~\prettyref{eq:band-support} supplies the exact support clause.

For $0<p<1$,
\[
R_K=\mathcal O(\rho^0[K+\xi K^{p+1}]),\qquad I_K=\mathcal O((\rho^0\xi)^{-1}K^{1-p}),\qquad B_K=\mathcal O(m\xi^{-1/2}K^{1-p/2}).
\]
Substitution gives~\prettyref{eq:rate}; balancing its exponents yields $p=1/3$ for w-KKT and $p=1/2$ for b-FJ.  The widened-support term is $\mathcal O(K^{-1})$.
\end{proof}

\begin{corollary}[Exact KKT conversion under Uniform Band-MFCQ]
\label{cor:KKT-MFCQ-exact}
Under~\prettyref{eq:parameters}, the ordinary polynomial schedule~\prettyref{eq:IP-schedules}, and~\prettyref{ass:MFCQ}, for each fixed $p\in(0,1)$ and $\xi\in(0,1]$, there is a finite, instance-dependent threshold $K_{\rm cq}(p,\xi)\geq3$, independent of $K$, such that for $K\geq K_{\rm cq}(p,\xi)$,
\begin{eqaligned}
\label{eq:kkt-rate}
\min_{j=3,\ldots,K+1}\CritKKT^j
\leq\frac{2(\gamma_{\rm cq}+M_F)}{\gamma_{\rm cq}}F_{\rm bFJ}(K).
\end{eqaligned}
In particular, for fixed $\xi$ and $p=1/2$, the same range $j\in\{3,\ldots,K+1\}$ gives $\mathcal O(K^{-1/4})$ for every $K\geq K_{\rm cq}(1/2,\xi)$.
\end{corollary}
\begin{proof}
Pointwise, the definitions and $\bar\iota^j=\alpha^j\iota^j$, $\alpha^j\in(0,1]$, and $\omega^j\geq0$ give the explicit chain
\begin{eqaligned}
\label{eq:kkt-vs-bFJ}
\CritKKT^j
&=\|Ax^j-\breve y^j\|+\|\Xi^j\|
+\sum_{i\in\IC}|\iota_i^jf_i(y_i^j)|\\
&\leq\frac1{\alpha^j}\left[
\|Ax^j-\breve y^j\|+\alpha^j\|\Xi^j\|
+\sum_{i\in\IC}|\bar\iota_i^jf_i(y_i^j)|+\omega^j\right]
=\frac{\CritbFJ^j}{\alpha^j}.
\end{eqaligned}
For fixed $p\in(0,1)$ and $\xi\in(0,1]$, $F_{\rm bFJ}(K)\to0$.  Choose the stated $K_{\rm cq}\geq3$ so that $F_{\rm bFJ}(K)\leq\min\{\delta_{\rm cq},\gamma_{\rm cq}/2\}$ for every $K\geq K_{\rm cq}(p,\xi)$.  Let $j_*\in\{3,\ldots,K+1\}$ minimize $\CritbFJ^j$; by~\prettyref{thm:rate-master}, the tuple at $j_*$ is an $F_{\rm bFJ}(K)$-b-FJ point.  The two threshold clauses permit~\prettyref{lem:alpha-lower} to be applied at $j_*$ and give $1/\alpha^{j_*}\leq2(\gamma_{\rm cq}+M_F)/\gamma_{\rm cq}$.  Therefore
\[
\min_{j=3,\ldots,K+1}\CritKKT^j
\leq\CritKKT^{j_*}
\leq\frac{\CritbFJ^{j_*}}{\alpha^{j_*}}
\leq\frac{2(\gamma_{\rm cq}+M_F)}{\gamma_{\rm cq}}F_{\rm bFJ}(K),
\]
which proves~\prettyref{eq:kkt-rate}; the $p=1/2$ rate follows from~\prettyref{eq:rate}, with the $K$-independent factor $2(\gamma_{\rm cq}+M_F)/\gamma_{\rm cq}$ absorbed into the $\mathcal O(\cdot)$ constant.
\end{proof}

\paragraph{Same-$x$ KKT residual in the smooth case.}
The residuals of~\prettyref{def:criteria} are nearby-point certificates: the objective forces $v_i$ use the proximal points $\breve y_i=p_i$ and the constraint gradients use the split points $y_i$, both of which approach $A_ix$ only as the consensus term $\|Ax-\breve y\|$ vanishes. When $\psi_i\equiv0$, the Moreau force vanishes ($w_i=0$), so $p_i=y_i$ and $\breve y_i=y_i$ for every block, and each $\varphi_i=h_i$ is smooth. Replacing every $y_i$ by $A_ix$ then changes the stationarity vector by at most $\bigl(\|A\|_2\sum_{i\in\IO}H_i+\|\iota\|_1\|A\|_2\max_{i\in\IC}L_i\bigr)\,\|Ax-y\|$, the feasibility slack by at most $\max_{i\in\IC}M_i\,\|Ax-y\|$, and each complementarity term $|\iota_if_i(A_ix)|$ by at most $|\iota_i|M_i\|Ax-y\|$, using the Lipschitz gradients $H_i,L_i$ and the bounded constraint gradients $M_i$. On the certified iterate, \prettyref{ass:MFCQ} and~\prettyref{lem:alpha-lower} bound $\|\iota\|_1\leq2(\gamma_{\rm cq}+M_F)/\gamma_{\rm cq}$, while $\|Ax-y\|\leq\CritWKKT\leq\CritKKT$. Hence the conventional same-$x$ KKT residual, with all objective and constraint terms evaluated at $A_ix$, is bounded by a fixed multiple of $\CritKKT$ and obeys the same $\mathcal O(K^{-1/4})$ rate as~\prettyref{cor:KKT-MFCQ-exact}. This is the standard transfer from a near-point to a same-point stationarity certificate~\cite{rockafellar2009variational}.

For the optional family analysis set
\begin{eqaligned}
\label{eq:block-W}
W_p(\xi)\defeq\max\{1,MH_p(\xi)\}.
\end{eqaligned}

\begin{lemma}[Aggregate block scaling]
\label{lem:block-aggregate-scaling}
Under~\prettyref{ass:block-scaling}, uniformly for $\xi\in(0,1]$,
\begin{eqaligned}
&\widehat G_0=\Theta(N),\quad
H_{\max},L_{\max}=\mathcal O(1),\quad
m=\Theta(1),\quad L^0,G=\Theta(N),\\
&\rho^0=\Theta(NW_p(\xi)),\qquad \Sigma_U=\mathcal O(N).
\end{eqaligned}
Moreover $H_p(\xi)=\mathcal O(1)$ uniformly on $(0,1]$.  For fixed $\xi\in(0,1]$, $W_p(\xi)=\Theta(M)$.
\end{lemma}
\begin{proof}
The per-block and clearance upper bounds in~\prettyref{ass:block-scaling} give $S_0,\widehat G_0,G_\psi=\mathcal O(N)$, $\bar M=\mathcal O(M)$, and $H_{\max},L_{\max},M_{\max},m=\mathcal O(1)$.  The uniform positive lower bound on $\varphi_i(A_ix^1)$ gives $S_0=\Omega(N)$; the positive lower bounds on $\widehat C_i$ and $M_i$ give $\widehat G_0=\Omega(N)$ and $\bar M=\Omega(M)$; and those on $M_i$ and $m$ give $M_{\max},m=\Omega(1)$. Hence $S_0,\widehat G_0=\Theta(N)$, $\bar M=\Theta(M)$, and $M_{\max},m=\Theta(1)$.  Therefore~\prettyref{eq:parameters} gives $L^0,G=\Theta(N)$.

First,~\prettyref{eq:Hp} implies the uniform bound
\[
0\leq H_p(\xi) \leq\sum_{j=1}^\infty j^{-2},
\]
because $(1+\xi j^p)^2\geq\xi^2j^{2p}$.  In the maximum defining $\rho^0$, the term $1209L^0M_{\max}^2$ is $\Theta(N)$, while the final term is
\[
\mathcal O\bigl(H_p(\xi)[1+NM]\bigr).
\]
The remaining entries in that maximum, namely $1$, $1209H_{\max}$, $1209mL^0L_{\max}$, and $6(\eta/2+11/2)\widehat G_0/G$, are each $\mathcal O(N)$.  Conversely, the $L^0\bar M$ part of the final entry is $\Theta(NMH_p(\xi))$.  Since $N,M\geq1$, taking the maximum proves $\rho^0=\Theta(N\max\{1,MH_p(\xi)\})$, with constants allowed to depend on the fixed $p$.

Finally, the definition~\prettyref{eq:drift-sum-a} gives
\[
\Sigma_U =\mathcal O\left( \frac{H_p(\xi)}{NW_p(\xi)}[N+N^2M]\right) =\mathcal O(N),
\]
using $MH_p(\xi)/W_p(\xi)\leq1$ and the uniform bound on $H_p$. For fixed $\xi>0$, the positive series in~\prettyref{eq:Hp} is a fixed constant, so $W_p(\xi)=\Theta(M)$.
\end{proof}

\begin{lemma}[Matrix scaling and automatic tuned saturation]
\label{lem:auto-saturation}
Under~\prettyref{ass:block-scaling}, $\|A\|_2=\Omega(1)$.  For $p=1/2$ on the interior branch $\xi_{\rm nom}<1$, recall $\xi_\star$ from~\prettyref{thm:rate}; then $W_{1/2}(\xi_\star)=\Theta(1)$ uniformly over the family.
\end{lemma}
\begin{proof}
Let $s_A\in\mathbb N$ and $0<c_A\leq C_A<\infty$ be the family-uniform constants from~\prettyref{ass:block-scaling}.  Norm equivalence on an $s_A$-sparse vector gives
\[
\frac{c_A}{\sqrt{s_A}}
\leq\|A_i\|_F
\leq C_A
\qquad(i\in\IAll).
\]
Since $A$ vertically stacks the $N+M$ blocks,
\[
\|A\|_F^2
=\sum_{i\in\IAll}\|A_i\|_F^2
=\Theta(N+M).
\]
Using $\operatorname{rank}(A)\leq\dim(x)=\mathcal O(N)$ yields
\[
\|A\|_2^2
\geq\frac{\|A\|_F^2}{\operatorname{rank}(A)}
=\Omega\left(\frac{N+M}{N}\right),
\]
while $\|A\|_2^2\leq\|A\|_F^2=\mathcal O(N+M)$.  In particular, $\|A\|_2=\Omega(1)$.

By~\prettyref{lem:Hp-small-xi}, there are $\xi_0\in(0,1)$ and $C_0>0$ such that
\[
H_{1/2}(\xi)\leq C_0\xi^2\log\frac e\xi
\qquad(0<\xi\leq\xi_0).
\]
For $\xi\in[\xi_0,1]$, the uniform bound in~\prettyref{lem:block-aggregate-scaling} and the positive minimum of $\xi^2\log(e/\xi)$ on this compact interval give the same inequality after increasing $C_0$.  Hence it holds for every $\xi\in(0,1]$.  On the interior branch, substituting $\xi_\star=m/(\|A\|_2N)$ gives
\[
M H_{1/2}(\xi_\star)
=\mathcal O\left(
\frac{Mm^2}{\|A\|_2^2N^2}
\log\frac{e\|A\|_2N}{m}
\right).
\]
The upper bound $\|A\|_2^2=\mathcal O(N+M)=\mathcal O(N^2)$, together with $m=\Theta(1)$, gives
\[
\log\frac{e\|A\|_2N}{m}=\mathcal O(\log(eN)).
\]
Combining this with $\|A\|_2^2=\Omega((N+M)/N)$ yields
\[
M H_{1/2}(\xi_\star)
=\mathcal O\left(
\frac{M}{N(N+M)}\log(eN)
\right)
=\mathcal O\left(\frac{\log(eN)}{N}\right).
\]
The right side is uniformly bounded for $N\geq1$.  Since $W_{1/2}(\xi)=\max\{1,MH_{1/2}(\xi)\}$ by~\prettyref{eq:block-W}, this proves $W_{1/2}(\xi_\star)=\Theta(1)$ uniformly.
\end{proof}

\begin{lemma}[Block scaling of exact rate constants]
\label{lem:block-rate-scaling}
Under~\prettyref{ass:block-scaling}, uniformly for $\xi\in(0,1]$,
\[
\mathcal B=\mathcal O(N),\qquad
Q_\psi=\mathcal O(\|A\|_2\sqrt N),\qquad
\mathcal Q=\Theta(\|A\|_2^2).
\]
\end{lemma}
\begin{proof}
The left side of~\prettyref{eq:p8} equals the initialization upper bound in~\prettyref{eq:Psi-2-bound} plus two nonnegative extra terms.  Thus~\prettyref{eq:p8} and~\prettyref{lem:initial-bound} imply $\Psi^2\leq G/2=\mathcal O(N)$.  The exact formula in \prettyref{lem:lower-bound} and~\prettyref{lem:block-aggregate-scaling} give $|\underline\Psi|=\mathcal O(N)$.  Since $\Psi^2\geq\underline\Psi$, also $|\Psi^2|=\mathcal O(N)$. Substitution in~\prettyref{eq:global-budget} from \prettyref{cor:movement-budget} proves $\mathcal B=\mathcal O(N)$.

The per-block upper bounds give $J_H,C_\psi=\mathcal O(\sqrt N)$.  The lower bound $\|A\|_2=\Omega(1)$ from~\prettyref{lem:auto-saturation} is used here to write $C_\psi=\mathcal O(\|A\|_2\sqrt N)$, while $\|A\|_2J_H=\mathcal O(\|A\|_2\sqrt N)$ follows directly.  Thus the definition of $Q_\psi$ before~\prettyref{eq:wKKT-pointwise} gives $Q_\psi=\mathcal O(\|A\|_2\sqrt N)$.

Moreover,
\[
\frac{H_{\max}}{\rho^0}=\mathcal O(W_p(\xi)^{-1}), \qquad \frac{\sqrt{2GL^0}L_{\max}}{\rho^0} =\mathcal O(W_p(\xi)^{-1}),
\]
so $r_*=\mathcal O(W_p(\xi)^{-1})$.  Since $W_p\geq1$, the positive floor in $c_y=1009/120900+(242909/121000)r_*$ gives $c_y=\Theta(1)$.  The exact definition of $\mathcal Q$ in~\prettyref{thm:rate-master} then gives $\mathcal Q=\Theta(\|A\|_2^2)$: its lower bound is supplied by the summand $26^2\|A\|_2^2$, while the matching upper bound uses $c_y=\Theta(1)$.
\end{proof}

Define the two block-explicit upper-bound shapes
\begin{eqaligned}
\label{eq:block-rate-shapes}
\mathfrak F_{\rm W}(K)&\defeq
\|A\|_2N\sqrt{W_p(\xi)}
\left[K^{-1/2}+\sqrt\xi K^{(p-1)/2}\right]
+\frac{\|A\|_2}{\sqrt N W_p(\xi)\xi}K^{-p},\\
\mathfrak F_{\rm bFJ}(K)&\defeq\mathfrak F_{\rm W}(K)
+m\xi^{-1/2}K^{-p/2}+\frac1{W_p(\xi)K}.
\end{eqaligned}

\begin{corollary}[Block-explicit finite-time rates]
\label{cor:block-explicit-rate}
Under~\prettyref{ass:block-scaling}, there is a constant $C_{\rm bs}>0$, possibly depending on the fixed $p$ and the family-uniform constants in~\prettyref{ass:block-scaling} but independent of $N,M,\xi,K$, and $\|A\|_2$, such that every $K\geq3$ satisfies
\begin{eqaligned}
\label{eq:block-explicit-rate}
\min_{j=3,\ldots,K+1}\CritWKKT^j&\leq C_{\rm bs}\mathfrak F_{\rm W}(K),\\
\min_{j=3,\ldots,K+1}\CritbFJ^j&\leq C_{\rm bs}\mathfrak F_{\rm bFJ}(K).
\end{eqaligned}
The second minimum retains the common-iterate certificate of \prettyref{thm:rate-master}.  For fixed $\xi\in(0,1]$,
\begin{eqaligned}
\label{eq:block-fixed-rate}
p=\frac13:&\quad
\min_{j=3,\ldots,K+1}\CritWKKT^j
=\mathcal O(\|A\|_2N\sqrt M K^{-1/3}),\\
p=\frac12:&\quad
\min_{j=3,\ldots,K+1}\CritbFJ^j=\mathcal O(\|A\|_2N\sqrt M K^{-1/4}).
\end{eqaligned}
For fixed $\xi$, \prettyref{lem:block-aggregate-scaling} gives $W_p(\xi)=\Theta(M)$.  At $p=1/3$, the movement and reconstruction contributions both scale as $K^{-1/3}$, while the reconstruction coefficient is dominated by $\|A\|_2N\sqrt M$.  At $p=1/2$, the movement and activation-band contributions determine the displayed $K^{-1/4}$ b-FJ rate, and all remaining terms are no larger.
\end{corollary}
\begin{proof}
For $K\geq3$, the schedules give
\[
R_K=\mathcal O(\rho^0[K+\xi K^{p+1}]),\quad I_K=\mathcal O((\rho^0\xi)^{-1}K^{1-p}),
\]
and $B_K=\mathcal O(m\xi^{-1/2}K^{1-p/2})$.  By the definition in~\prettyref{thm:rate-master}, $\mathcal Q\geq26^2\|A\|_2^2$; the lower bound $\|A\|_2=\Omega(1)$ from~\prettyref{lem:auto-saturation} therefore makes $\mathcal Q$ family-uniformly positive.  Since $\rho^k\geq\rho^0\geq1$, each $(\rho^k)^{-1}\leq\rho^k$, so $I_K\leq R_K$ and hence $I_K+\mathcal Q R_K=\mathcal O(\mathcal Q R_K)$. Together with $\mathcal B=\mathcal O(N)$ and $\mathcal Q=\Theta(\|A\|_2^2)$ from~\prettyref{lem:block-rate-scaling}, the movement part of \prettyref{eq:rate-explicit} is at most
\[
\mathcal O\left( \|A\|_2N\sqrt{W_p(\xi)} [K^{-1/2}+\sqrt\xi K^{(p-1)/2}] \right).
\]
The sharper inverse-schedule estimate, together with $Q_\psi=\mathcal O(\|A\|_2\sqrt N)$ and $\rho^0=\Theta(NW_p(\xi))$, bounds the reconstruction part by
\[
\mathcal O\left( \frac{\|A\|_2}{\sqrt N W_p(\xi)\xi}K^{-p} \right).
\]
The band sum contributes $\mathcal O(m\xi^{-1/2}K^{-p/2})$.

For the final term, $\rho^2=\rho^0(1+\xi2^p)=\Theta(\rho^0)=\Theta(NW_p(\xi))$, so $L_{\max}\mathcal B/[\rho^2(K-1)] =\mathcal O([W_p(\xi)K]^{-1})$. Substitution in~\prettyref{eq:rate-explicit} proves \prettyref{eq:block-explicit-rate}, including the common-iterate certificate.

For fixed $\xi$,~\prettyref{lem:block-aggregate-scaling} gives $W_p(\xi)=\Theta(M)$.  At $p=1/3$, the movement terms have exponents $K^{-1/2}$ and $K^{-1/3}$, and the reconstruction term has exponent $K^{-1/3}$; moreover, since $N,M\geq1$, $\|A\|_2/(\sqrt N M)\leq\|A\|_2N\sqrt M$.  Thus all three are bounded by the displayed w-KKT monomial $\|A\|_2N\sqrt M K^{-1/3}$.  At $p=1/2$, the movement exponents are $K^{-1/2}$ and $K^{-1/4}$, reconstruction is $K^{-1/2}$, the activation band is $K^{-1/4}$, and widened support is $K^{-1}$; the family upper bound $m=\mathcal O(1)$ and the lower bound $\|A\|_2=\Omega(1)$ from~\prettyref{lem:auto-saturation}, together with $N,M\geq1$, let the displayed coefficient $\|A\|_2N\sqrt M$ absorb the remaining coefficients.  This proves~\prettyref{eq:block-fixed-rate}.
\end{proof}

\begin{corollary}[Tuned block rate on the interior branch]
\label{cor:tuned-block-rate}
Suppose~\prettyref{ass:block-scaling} holds with $p=1/2$ and the nominal value in~\prettyref{eq:block-xi-star} satisfies $\xi_{\rm nom}<1$.  Recall $\xi_\star$ from~\prettyref{thm:rate} and choose $\xi=\xi_\star$.  Then, for every $K\geq3$,
\begin{eqaligned}
\label{eq:tuned-block-rate}
\min_{j=3,\ldots,K+1}\CritbFJ^j
=\mathcal O\left(
\left[\|A\|_2N+\frac{\|A\|_2^2\sqrt N}{m}\right]K^{-1/2}
+\sqrt{\|A\|_2Nm}\,K^{-1/4}+K^{-1}
\right).
\end{eqaligned}
The minimum is attained at the same selected certificate-bearing iterate as in~\prettyref{thm:rate-master}.  The two $K^{-1/2}$ contributions are the baseline movement term and the inverse-schedule reconstruction term after substituting $\xi_\star$.  Balancing the movement and activation-band coefficients produces the $\sqrt{\|A\|_2Nm}\,K^{-1/4}$ term, while the widened-support correction contributes $K^{-1}$.  This shape applies on the interior tuned branch; when the nominal choice clips to $\xi=1$, the fixed-$\xi$ rate in \prettyref{eq:block-fixed-rate} applies instead.
\end{corollary}
\begin{proof}
\prettyref{lem:auto-saturation} gives $W_{1/2}(\xi_\star)=\Theta(1)$.  The two $K^{-1/4}$ coefficients in \prettyref{eq:block-rate-shapes} are then $\|A\|_2N\sqrt{\xi_\star}$, up to the factor $\sqrt{W_{1/2}(\xi_\star)}=\Theta(1)$, and $m\xi_\star^{-1/2}$, which are balanced by the exact identity $\xi_\star=\xi_{\rm nom}=m/(\|A\|_2N)$ from \prettyref{eq:block-xi-star}.  Substitution in \prettyref{cor:block-explicit-rate} gives \prettyref{eq:tuned-block-rate}; in particular the widened-support term is $\mathcal O(K^{-1})$.  The clipped case is the fixed-$\xi$ case $\xi=1$ of~\prettyref{eq:block-fixed-rate}.
\end{proof}

\begin{corollary}[Block scaling under Uniform Band-MFCQ]
\label{cor:block-KKT-scaling}
Suppose~\prettyref{eq:parameters}, the ordinary polynomial schedule~\prettyref{eq:IP-schedules},~\prettyref{ass:block-scaling}, and~\prettyref{ass:MFCQ} hold, with $D_{\rm cq}$ as in~\prettyref{cor:KKT-MFCQ}.  Then the conversion factor in~\prettyref{eq:kkt-rate} equals $D_{\rm cq}$, and for fixed $\xi\in(0,1]$, $p=1/2$, and every $K\geq K_{\rm cq}(1/2,\xi)$, with the instance-dependent threshold from~\prettyref{cor:KKT-MFCQ-exact},
\begin{eqaligned}
\label{eq:block-KKT-fixed}
\min_{j=3,\ldots,K+1}\CritKKT^j
=\mathcal O(D_{\rm cq}\|A\|_2N\sqrt M K^{-1/4}).
\end{eqaligned}
\end{corollary}
\begin{proof}
By the definition of $M_F$ in~\prettyref{eq:parameters} and~\prettyref{lem:block-aggregate-scaling}, the exact factor in~\prettyref{eq:kkt-rate} is
\[
\frac{2(\gamma_{\rm cq}+M_F)}{\gamma_{\rm cq}} =2\left(1+\frac{\|A\|_2\sqrt{\widehat G_0}} {\gamma_{\rm cq}}\right)=D_{\rm cq} =\Theta\left(1+\frac{\|A\|_2\sqrt N}{\gamma_{\rm cq}}\right).
\]
Fix $\xi\in(0,1]$, set $p=1/2$, and let $K\geq K_{\rm cq}(1/2,\xi)$.  Let $j_*\in\{3,\ldots,K+1\}$ minimize $\CritbFJ^j$; by~\prettyref{thm:rate-master}, the tuple at $j_*$ is an $F_{\rm bFJ}(K)$-b-FJ point.  As in the proof of~\prettyref{cor:KKT-MFCQ-exact}, \prettyref{lem:alpha-lower} applies at $j_*$ and gives $1/\alpha^{j_*}\leq D_{\rm cq}$.  The pointwise conversion~\prettyref{eq:kkt-vs-bFJ} and the $p=1/2$ case of~\prettyref{eq:block-fixed-rate} then give
\[
\min_{j=3,\ldots,K+1}\CritKKT^j\leq\CritKKT^{j_*}\leq D_{\rm cq}\CritbFJ^{j_*}
=D_{\rm cq}\min_{j=3,\ldots,K+1}\CritbFJ^j=\mathcal O(D_{\rm cq}\|A\|_2N\sqrt M K^{-1/4}),
\]
which is~\prettyref{eq:block-KKT-fixed}.
\end{proof}

\begin{proof}[Proof of \prettyref{thm:rate}]
The instance-wise asymptotic bounds are the corresponding specializations of \prettyref{thm:rate-master}.  The fixed-$\xi$ family statements follow from \prettyref{cor:block-explicit-rate}, and the interior tuned statement follows from \prettyref{cor:tuned-block-rate} and \prettyref{lem:auto-saturation} with the exact rule \prettyref{eq:block-xi-star}.  These results give every assertion of \prettyref{thm:rate}.
\end{proof}

\begin{proof}[Proof of \prettyref{cor:KKT-MFCQ}]
Apply \prettyref{cor:KKT-MFCQ-exact} with $p=1/2$.  Its fixed-family specialization is \prettyref{cor:block-KKT-scaling}, which gives the stated $D_{\rm cq}$ dependence.
\end{proof}

\section{Convergence Analysis of IP-ADMM-Ckpt}
\label{appen:checkpoint}
This appendix proves \prettyref{thm:checkpoint-rate} for the $\mathsf{Ckpt}$ mode of \prettyref{alg:IP-ADMM}.  Throughout this appendix, set $\mathsf m=\mathsf{Ckpt}$, so that $\theta_{1,\mathsf m}=\theta_1^{\rm cp}$, $\theta_{2,\mathsf m}=\theta_2^{\rm cp}$, and $G_{\mathsf m}=G^{\rm cp}$ with the values in~\prettyref{eq:checkpoint-parameters}, and write $g\defeq g_{\mathsf m}$.  The residual objects $\varpi^{k+1}$, $\zeta^{k+1}$, $\Xi^{k+1}$, and $\CritKKT^{k+1}$ of~\prettyref{def:criteria}--\prettyref{eq:crit} and the reconstructed-force gaps $u^k$ defined just before~\prettyref{eq:residual-stationary} are introduced in \prettyref{appen:convergence}; they are mode-neutral and are read here with $\mathsf m=\mathsf{Ckpt}$.  The mode uses the unit virtual-multiplier update shown in the algorithm, the checkpoint test $E^{k+1}\leq\tau_{\rm cp}/\rho^k$, the explicit parameters in~\prettyref{eq:checkpoint-parameters}, and the geometric schedule~\prettyref{eq:checkpoint-schedule}.  Its proof conditions each level on two caps, derives both finite arrival and accepted-step control under those caps, and then propagates the caps across checkpoint transitions.  This produces pointwise certificates for checkpoints.  The ordinary $\mathsf{IP}$ mode instead uses a global continuation argument and is analyzed in \prettyref{appen:convergence}.

\paragraph{Levels, frozen state, and checkpoint outputs.}
Level $j\geq0$ begins after exactly $j$ checkpoints have been committed.  While that level is active, the level-defining state $(\lambda^{[j]},\rho_j,L_j,\mu_j)$ is fixed, but every new primal pair is retained.  Put $k_{-1}\defeq0$.  If level $j$ completes, let $k_j$ denote its accepting iteration; then level $j$ consists of $k\in\{k_{j-1}+1,\ldots,k_j\}$ and
\begin{eqaligned}
\label{eq:ckpt-level}
(\rho^k,L^k,\mu^k)&=(\rho_j,L_j,\mu_j),&
\lambda^k&=\lambda^{[j]},&
T_j&\defeq k_j-k_{j-1},&
K_0&\defeq k_0,
\end{eqaligned}
for every $k$ in that level.  Here $\lambda^{[0]}\defeq\lambda^1=0$.  The finite-arrival result \prettyref{lem:ckpt-arrival} later proves that every level completes, so all of these indices are finite.

Three multiplier roles are distinct.  The \emph{stored} multiplier $\lambda^k=\lambda^{[j]}$ enters every subproblem in level $j$.  Each iteration forms the \emph{virtual} multiplier $\widehat\lambda^{k+1}=\lambda^k+\rho^kr^{k+1}$.  Rejection discards that virtual multiplier, whereas acceptance promotes $\widehat\lambda^{k_j+1}$ to the next stored multiplier $\lambda^{[j+1]}$.  The raw multiplier bounds below cover stored multipliers and accepted virtual multipliers, not rejected virtual multipliers.

At acceptance, the method first records the pre-promotion checkpoint certificate tuple
\begin{eqaligned}
\label{eq:ckpt-certificate-tuple}
(x^{k_j+1},y^{k_j+1},\breve y^{k_j+1},v^{k_j+1},\iota^{k_j+1},\zeta^{k_j+1}).
\end{eqaligned}
It is formed from $(\rho_j,L_j,\mu_j)$ and $\lambda^{[j]}$ by \prettyref{eq:reconstruction} and the definitions preceding \prettyref{eq:crit}, with $\zeta^{k_j+1}=\rho_j\varpi^{k_j+1}$.

\paragraph{Frozen-level merit function.}
The checkpoint proof measures primal progress during level $j$ with
\begin{eqaligned}
\label{eq:ckpt-merit}
\mathcal P_j(x,y)\defeq
\mathcal L(x,y,\lambda^{[j]},L_j,\mu_j,\rho_j)+R_\psi(\mu_j),
\end{eqaligned}
where $\mathcal L$ is defined in~\prettyref{eq:Lagrangian} and $R_\psi$ in~\prettyref{lem:Rpsi-discharge}.  The level-defining state is fixed, so every primal sweep is measured against the same function $\mathcal P_j$; no polynomial drift correction or stored-multiplier update appears within the level.  The primals still change at every iteration.  Only promotion changes the definition of the potential, and \prettyref{lem:ckpt-transition} charges that one-time change at the retained primal pair.

\paragraph{Derived constants.}
The following constants package the support-width, certificate, merit-budget, and work estimates.  They are determined by the explicit choice~\prettyref{eq:checkpoint-parameters} and introduce no additional assumptions:
\begin{eqaligned}
\label{eq:ckpt-derived-constants}
S_\omega&\defeq b_{\rm cp}
+\frac{L_{\max}\tau_{\rm cp}^2}{2(\rho_{\rm cp}^0)^3},\quad
C_{\rm cert}\defeq2Q_{\rm cp}+\eta^{\rm cp}C_\psi+S_\Xi+I_{\rm cp}S_\omega,\\
\Delta_*&\defeq\overline V+\Delta_G,\quad
C_{\rm work}\defeq\frac{\Delta_*}{a_0\tau_{\rm cp}^2}.
\end{eqaligned}

\paragraph{Organization.}
The proof reuses the shared local calculus (Moreau envelopes, walls, model error, block optimality, and reconstruction) and \prettyref{lem:checkpoint-schedule}, but replaces the ordinary mode's drift estimates, Lyapunov functions, and window averaging by frozen-level descent, a two-cap level induction, and promotion charging.  \prettyref{sub:ckpt-descent} and \prettyref{sub:ckpt-feasibility} establish frozen-level descent, the virtual-multiplier identity, wall support, the merit lower bound, and split feasibility under the level hypotheses.  \prettyref{sub:ckpt-invariants} bounds accepted checkpoints, proves finite arrival, and propagates the hypotheses across promotions, and \prettyref{sub:ckpt-certificate} derives the pointwise certificate, the rates, and the scaling results.

\subsection{Frozen-Level Descent and Certificate-Tuple Identities}
\label{sub:ckpt-descent}
Throughout this subsection $k$ is an iteration of level $j$, so that $(\rho^k,L^k,\mu^k)=(\rho_j,L_j,\mu_j)$ and $\lambda^k=\lambda^{[j]}$.  As in \prettyref{sub:lyapunov-descent}, $\Delta x^k=x^{k+1}-x^k$, $d_i^k=y_i^{k+1}-y_i^k$, $r^{k+1}=y^{k+1}-Ax^{k+1}$, and $\upsilon_i^k=m_i^k-\rho^k$, which at level $j$ reads $\upsilon_i^k=m_i^k-\rho_j=(\theta_1^{\rm cp}-1)\rho_j+\theta_1^{\rm cp}\tau_i^k$; $\Upsilon^k=\operatorname{diag}(\upsilon_i^kI_{n_i})_{i\in\IAll}$ is the block-diagonal curvature matrix of~\prettyref{sub:lyapunov-descent}.

\begin{lemma}[Frozen-level one-step descent]
\label{lem:ckpt-descent}
Let $k$ be any iteration of level $j$.  Then, with $a_0$ from~\prettyref{eq:checkpoint-parameters},
\begin{eqaligned}
\label{eq:ckpt-descent}
\mathcal P_j(x^k,y^k)-\mathcal P_j(x^{k+1},y^{k+1})
\geq a_0\rho_j\left(\|d^k\|^2+\|A\Delta x^k\|^2\right).
\end{eqaligned}
\end{lemma}
\begin{proof}
Since $\mu^k=\mu_j$ is frozen, $R_\psi(\mu_j)$ is the same constant at $(x^k,y^k)$ and at $(x^{k+1},y^{k+1})$; it therefore cancels, and it suffices to bound the decrease of $\mathcal L(\cdot,\cdot,\lambda^{[j]},L_j,\mu_j,\rho_j)$.

\emph{Block sweep.}  For fixed $x^k$ the Lagrangian is separable: $\mathcal L(x^k,y,\lambda^{[j]},L_j,\mu_j,\rho_j)=\sum_{i\in\IAll}\Phi_i(y_i)$ with
\[
\Phi_i(u)\defeq\bar f_i(u,L_j,\mu_j)
+\langle\lambda_i^{[j]},u-A_ix^k\rangle
+\frac{\rho_j}{2}\|u-A_ix^k\|^2.
\]
Expanding the quadratic about $y_i^k$ and using $e_i^k=\lambda_i^{[j]}+\rho_j(y_i^k-A_ix^k)$ gives the exact identity
\[
\Phi_i(u)-\Phi_i(y_i^k)
=\bar f_i(u,L_j,\mu_j)-\bar f_i(y_i^k,L_j,\mu_j)
+\langle e_i^k,u-y_i^k\rangle+\frac{\rho_j}{2}\|u-y_i^k\|^2.
\]
Let $\widetilde f_i^k$ be the model-only function of~\prettyref{eq:model-two-sided}; by construction $\widetilde f_i^k(y_i^k)=\bar f_i(y_i^k,L_j,\mu_j)$, and $\widetilde f_i^k(u)+\langle e_i^k,u-y_i^k\rangle+\frac{m_i^k}{2}\|u-y_i^k\|^2$ equals $\mathcal M_i^k(u)$ up to an additive constant independent of $u$.  Applying~\prettyref{eq:model-two-sided} at $u=y_i^{k+1}$,
\[
\Phi_i(y_i^{k+1})-\Phi_i(y_i^k)
\leq\left[\mathcal M_i^k(y_i^{k+1})-\mathcal M_i^k(y_i^k)\right]
-\frac{m_i^k}{2}\|d_i^k\|^2
+\frac{\tau_i^k+\rho_j}{2}\|d_i^k\|^2 .
\]
The model $\mathcal M_i^k$ is $m_i^k$-strongly convex by~\prettyref{ass:obj} on objective blocks and by the convexity of $g(\cdot,L_j)$ in~\prettyref{lem:added_algebraic_3} on constraint blocks.  The block update in~\prettyref{alg:IP-ADMM} minimizes this model at $y_i^{k+1}$.  Hence its bracket is at most $-\frac{m_i^k}{2}\|d_i^k\|^2$.  Since $\tau_i^k+\rho_j=\mathsf B_i^k$ and $m_i^k=\theta_1^{\rm cp}\mathsf B_i^k$,
\[
\Phi_i(y_i^{k+1})-\Phi_i(y_i^k)
\leq-\left(m_i^k-\frac{\mathsf B_i^k}{2}\right)\|d_i^k\|^2
=-\frac{2\theta_1^{\rm cp}-1}{2}\mathsf B_i^k\|d_i^k\|^2
\leq-\frac{2\theta_1^{\rm cp}-1}{2}\rho_j\|d_i^k\|^2,
\]
using $\mathsf B_i^k=\rho_j+\tau_i^k\geq\rho_j$, since $\tau_i^k\geq0$, and $2\theta_1^{\rm cp}-1>0$.  Summing over $i\in\IAll$ shows that the block sweep decreases the Lagrangian by at least $\frac{2\theta_1^{\rm cp}-1}{2}\rho_j\|d^k\|^2$.

\emph{$x$ step.}  Write $\phi(x)\defeq\mathcal L(x,y^{k+1},\lambda^{[j]},L_j,\mu_j,\rho_j)+\frac{\theta_2^{\rm cp}\rho_j}{2}\|A(x-x^k)\|^2$, the exact objective of the $x$ subproblem.  It is a convex quadratic in $x$ with Hessian $(1+\theta_2^{\rm cp})\rho_jA^TA$, so its second-order Taylor expansion about $x^{k+1}$ is an exact identity at every point; evaluated at $x^k$ it reads
\[
\phi(x^k)=\phi(x^{k+1})+\langle\nabla\phi(x^{k+1}),x^k-x^{k+1}\rangle
+\frac{(1+\theta_2^{\rm cp})\rho_j}{2}\|A\Delta x^k\|^2 .
\]
By~\prettyref{eq:opt-x} the gradient $\nabla\phi(x^{k+1})$ lies in the range of $C^T$, while $C\Delta x^k=0$ because $Cx^k=d$ by induction from~\prettyref{ass:feasible} and $Cx^{k+1}=d$ by~\prettyref{eq:opt-x}; hence the inner product vanishes and $\phi(x^k)-\phi(x^{k+1})=\frac{(1+\theta_2^{\rm cp})\rho_j}{2}\|A\Delta x^k\|^2$.  The proximal term of $\phi$ vanishes at $x^k$ and equals $\frac{\theta_2^{\rm cp}\rho_j}{2}\|A\Delta x^k\|^2$ at $x^{k+1}$, so the Lagrangian decrease exceeds the decrease of $\phi$ by exactly that amount; adding it gives
\[
\mathcal L(x^k,y^{k+1},\lambda^{[j]},L_j,\mu_j,\rho_j)
-\mathcal L(x^{k+1},y^{k+1},\lambda^{[j]},L_j,\mu_j,\rho_j)
=\frac{1+2\theta_2^{\rm cp}}{2}\rho_j\|A\Delta x^k\|^2 .
\]
Adding the two sweeps and restoring the cancelled constant $R_\psi(\mu_j)$ gives
\[
\mathcal P_j(x^k,y^k)-\mathcal P_j(x^{k+1},y^{k+1})
\geq\frac{2\theta_1^{\rm cp}-1}{2}\rho_j\|d^k\|^2
+\frac{1+2\theta_2^{\rm cp}}{2}\rho_j\|A\Delta x^k\|^2 .
\]
By~\prettyref{eq:checkpoint-parameters}, $\theta_1^{\rm cp}=101/100$ and $\theta_2^{\rm cp}=1/100$, so both coefficients equal $(2\theta_1^{\rm cp}-1)/2=(1+2\theta_2^{\rm cp})/2=51/100=a_0$, which proves~\prettyref{eq:ckpt-descent}.
\end{proof}

\begin{lemma}[Virtual multiplier and stationarity transfer]
\label{lem:ckpt-virtual-dual}
Let $k$ be any iteration of level $j$.  Then the virtual multiplier obeys the exact identity
\begin{eqaligned}
\label{eq:ckpt-virtual-dual}
\widehat\lambda^{k+1}=\lambda^{[j]}+\rho_jr^{k+1}
=-q^k-\Upsilon^kd^k-\rho_jA\Delta x^k,
\end{eqaligned}
and, with $r_{\rm cp}$ and $c_{\upsilon,{\rm cp}}$ from~\prettyref{eq:checkpoint-parameters},
\begin{eqaligned}
\label{eq:ckpt-upsilon-bound}
0\leq\frac{\tau_i^k}{\rho_j}\leq r_{\rm cp},
\qquad
0<\frac{\upsilon_i^k}{\rho_j}\leq c_{\upsilon,{\rm cp}}
\qquad(i\in\IAll).
\end{eqaligned}
If in addition the step-$k$ test is accepted, that is $E^{k+1}\leq\tau_{\rm cp}/\rho_j$, then the stationarity residual $\Xi^{k+1}$ of~\prettyref{eq:crit} satisfies, with $S_\Xi$ from~\prettyref{eq:checkpoint-parameters},
\begin{eqaligned}
\label{eq:ckpt-stationary-accepted}
\|\Xi^{k+1}\|\leq\frac{S_\Xi}{\rho_j}.
\end{eqaligned}
\end{lemma}
\begin{proof}
\emph{Identity~\prettyref{eq:ckpt-virtual-dual}.}  Substituting $y_i^k-A_ix^k=r_i^{k+1}-d_i^k+A_i\Delta x^k$ into~\prettyref{eq:opt-yi} with $\lambda_i^k=\lambda_i^{[j]}$ and $\rho^k=\rho_j$ gives
\[
0=q_i^k+\lambda_i^{[j]}+\rho_jr_i^{k+1}
+(m_i^k-\rho_j)d_i^k+\rho_jA_i\Delta x^k,
\]
which is exactly $\widehat\lambda_i^{k+1}=-q_i^k-\upsilon_i^kd_i^k-\rho_jA_i\Delta x^k$.  Stacking over $i$ proves~\prettyref{eq:ckpt-virtual-dual}.

\emph{Bound~\prettyref{eq:ckpt-upsilon-bound}.}  For $i\in\IO$, $\tau_i^k=H_i\leq H_{\max}$, and $\rho_j\geq\rho_{\rm cp}^0\geq1$ by~\prettyref{eq:checkpoint-schedule} and~\prettyref{eq:checkpoint-parameters}, so $\tau_i^k/\rho_j\leq H_{\max}\leq r_{\rm cp}$.  For $i\in\IC$, $\tau_i^k=\sqrt{2L_jG^{\rm cp}}\,L_i$, and the second identity of~\prettyref{eq:checkpoint-wall-scaling} gives $\tau_i^k/\rho_j=\sqrt{2G^{\rm cp}\ell_{\rm cp}}\,L_i\leq\sqrt{2G^{\rm cp}\ell_{\rm cp}}\,L_{\max}\leq r_{\rm cp}$.  Nonnegativity is clear.  Consequently
\[
\frac{\upsilon_i^k}{\rho_j}=\frac{\theta_1^{\rm cp}(\rho_j+\tau_i^k)-\rho_j}{\rho_j}
=\theta_1^{\rm cp}-1+\theta_1^{\rm cp}\frac{\tau_i^k}{\rho_j}
\leq\theta_1^{\rm cp}(1+r_{\rm cp})-1=c_{\upsilon,{\rm cp}},
\]
and $\upsilon_i^k\geq(\theta_1^{\rm cp}-1)\rho_j>0$ because $\theta_1^{\rm cp}=101/100$.

\emph{Reconstruction gap.}  Recall the reconstructed-force gaps, defined just before~\prettyref{eq:residual-stationary}, $u_i^k=v_i^{k+1}-q_i^k$ for $i\in\IO$, $u_i^k=\iota_i^{k+1}\nabla f_i(y_i^{k+1})-q_i^k$ for $i\in\IC$, and $u^k=(u_i^k)_{i\in\IAll}$.  For $i\in\IO$, the second line of~\prettyref{eq:recon-error} gives $\|u_i^k\|\leq H_i\|d_i^k\|+\mu_jH_iC_i^\psi=\tau_i^k\|d_i^k\|+\mu_jH_iC_i^\psi$.  For $i\in\IC$, the wall lemma~\prettyref{lem:added_algebraic_3} ($\mathsf m=\mathsf{Ckpt}$, so $G_{\mathsf m}=G^{\rm cp}$) gives $\|u_i^k\|\leq\sqrt{2L_jG^{\rm cp}}\,L_i\|d_i^k\|=\tau_i^k\|d_i^k\|$.  By the previous paragraph $\tau_i^k\leq r_{\rm cp}\rho_j$, so these blockwise estimates and Minkowski's inequality give
\[
\|u^k\|\leq\Big(\sum_{i\in\IAll}(\tau_i^k)^2\|d_i^k\|^2\Big)^{1/2}
+\mu_j\Big(\sum_{i\in\IO}(H_iC_i^\psi)^2\Big)^{1/2}
\leq r_{\rm cp}\rho_j\|d^k\|+\mu_jJ_H,
\]
with $J_H$ from~\prettyref{eq:parameters}.

\emph{Stationarity identity.}  By the definition preceding~\prettyref{eq:crit}, $\Xi^{k+1}=A^T(q^k+u^k)+C^T\zeta^{k+1}$ with $\zeta^{k+1}=\rho_j\varpi^{k+1}$.  Specializing~\prettyref{eq:opt-x} to $\mathsf m=\mathsf{Ckpt}$ and multiplying it by $\rho_j$ and using $y^{k+1}-Ax^k=r^{k+1}+A\Delta x^k$ gives
\[
C^T\zeta^{k+1}
=A^T\left(\lambda^{[j]}+\rho_jr^{k+1}+\rho_jA\Delta x^k\right)
-(1+\theta_2^{\rm cp})\rho_jA^TA\Delta x^k
=A^T\widehat\lambda^{k+1}-\theta_2^{\rm cp}\rho_jA^TA\Delta x^k,
\]
where the last equality uses~\prettyref{eq:ckpt-virtual-dual}.  Also by~\prettyref{eq:ckpt-virtual-dual}, $q^k=-\widehat\lambda^{k+1}-\Upsilon^kd^k-\rho_jA\Delta x^k$.  Adding the two displays, the terms $\pm A^T\widehat\lambda^{k+1}$ cancel and
\[
\Xi^{k+1}=A^T\left[u^k-\Upsilon^kd^k-\rho_jA\Delta x^k\right]
-\theta_2^{\rm cp}\rho_jA^TA\Delta x^k
=A^T\left[u^k-\Upsilon^kd^k\right]-(1+\theta_2^{\rm cp})\rho_jA^TA\Delta x^k .
\]

\emph{Bound~\prettyref{eq:ckpt-stationary-accepted}.}  Because $\Upsilon^k$ is block diagonal, $\|\Upsilon^kd^k\|^2=\sum_i|\upsilon_i^k|^2\|d_i^k\|^2\leq c_{\upsilon,{\rm cp}}^2\rho_j^2\|d^k\|^2$ by~\prettyref{eq:ckpt-upsilon-bound}.  Combining this with the bound on $\|u^k\|$ and $\|A^T\cdot\|\leq\|A\|_2\|\cdot\|$ in the stationarity identity gives
\[
\|\Xi^{k+1}\|\leq\|A\|_2\left[
(r_{\rm cp}+c_{\upsilon,{\rm cp}})\rho_j\|d^k\|
+(1+\theta_2^{\rm cp})\rho_j\|A\Delta x^k\|
+\mu_jJ_H\right].
\]
If the test is accepted, then $E^{k+1}\leq\tau_{\rm cp}/\rho_j$ and, by the Cauchy--Schwarz inequality applied to the two movement channels,
\[
(r_{\rm cp}+c_{\upsilon,{\rm cp}})\rho_j\|d^k\|+(1+\theta_2^{\rm cp})\rho_j\|A\Delta x^k\|
\leq c_{\rm mov}\,\rho_j\sqrt{\|d^k\|^2+\|A\Delta x^k\|^2}
=c_{\rm mov}E^{k+1}
\leq\frac{c_{\rm mov}\tau_{\rm cp}}{\rho_j},
\]
with $c_{\rm mov}$ from~\prettyref{eq:checkpoint-parameters}.  Since $\mu_j=\eta^{\rm cp}/\rho_j$ by~\prettyref{eq:checkpoint-schedule}, substitution in the preceding bound gives $\|\Xi^{k+1}\|\leq\|A\|_2(c_{\rm mov}\tau_{\rm cp}+\eta^{\rm cp}J_H)/\rho_j=S_\Xi/\rho_j$.
\end{proof}

\begin{lemma}[Checkpoint support width]
\label{lem:ckpt-support}
Let $k$ be an iteration of level $j$ whose test is accepted, and suppose $f_i(y_i^{k+1})\geq0$ for every $i\in\IC$.  Then, with $S_\omega$ from~\prettyref{eq:ckpt-derived-constants},
\begin{eqaligned}
\label{eq:ckpt-support}
\iota_i^{k+1}\neq0
\quad\Longrightarrow\quad
0\leq f_i(y_i^{k+1})
\leq\frac{S_\omega}{\rho_j}\leq\delta_{\rm cq}.
\end{eqaligned}
Consequently
\begin{eqaligned}
\label{eq:ckpt-complementarity}
\sum_{i\in\IC}\left|\iota_i^{k+1}f_i(y_i^{k+1})\right|
\leq\frac{S_\omega}{\rho_j}\|\iota^{k+1}\|_1 .
\end{eqaligned}
\end{lemma}
\begin{proof}
By definition, $\iota_i^{k+1}=g'(\ell_i^k(y_i^{k+1}),L_j)$, and the explicit force~\prettyref{eq:band-force} shows that $g'(\cdot,L_j)$ vanishes exactly on $[\beta(L_j),\infty)=[\beta_j,\infty)$, the equality by~\prettyref{eq:checkpoint-wall-scaling}.  Hence $\iota_i^{k+1}\neq0$ forces $\ell_i^k(y_i^{k+1})<\beta_j$.  Since $f_i$ has an $L_i$-Lipschitz gradient by \prettyref{ass:f},
\[
f_i(y_i^{k+1})\leq\ell_i^k(y_i^{k+1})+\frac{L_i}{2}\|d_i^k\|^2
<\beta_j+\frac{L_i}{2}\|d_i^k\|^2 .
\]
The lower bound is the hypothesis.  Acceptance gives $\rho_j\|d_i^k\|\leq\rho_j\|d^k\|\leq E^{k+1}\leq\tau_{\rm cp}/\rho_j$, where the first step is the blockwise bound $\|d_i^k\|\leq\|d^k\|$ and the second discards the nonnegative $\|A\Delta x^k\|$ channel of $E^{k+1}$; hence $\|d_i^k\|^2\leq\tau_{\rm cp}^2/\rho_j^4$, and $L_i\leq L_{\max}$, so
\[
\beta_j+\frac{L_i}{2}\|d_i^k\|^2
\leq\frac{b_{\rm cp}}{\rho_j}+\frac{L_{\max}\tau_{\rm cp}^2}{2\rho_j^4}
=\frac1{\rho_j}\left[b_{\rm cp}+\frac{L_{\max}\tau_{\rm cp}^2}{2\rho_j^3}\right]
\leq\frac1{\rho_j}\left[b_{\rm cp}
+\frac{L_{\max}\tau_{\rm cp}^2}{2(\rho_{\rm cp}^0)^3}\right]
=\frac{S_\omega}{\rho_j},
\]
using $\rho_j\geq\rho_{\rm cp}^0$.  For the last inequality in~\prettyref{eq:ckpt-support}, $\rho_j\geq\rho_{\rm cp}^0$ and~\prettyref{eq:ckpt-derived-constants} give
\[
\frac{S_\omega}{\rho_j}
\leq\frac{b_{\rm cp}}{\rho_{\rm cp}^0}
+\frac{L_{\max}\tau_{\rm cp}^2}{2(\rho_{\rm cp}^0)^4}.
\]
The clause $\rho_{\rm cp}^0\geq2b_{\rm cp}/\delta_{\rm cq}$ of~\prettyref{eq:checkpoint-parameters} bounds the first summand by $\delta_{\rm cq}/2$, and the clause $\rho_{\rm cp}^0\geq(L_{\max}\tau_{\rm cp}^2/\delta_{\rm cq})^{1/4}$ gives $(\rho_{\rm cp}^0)^4\geq L_{\max}\tau_{\rm cp}^2/\delta_{\rm cq}$ and hence bounds the second summand by $\delta_{\rm cq}/2$.  Finally, $\iota_i^{k+1}f_i(y_i^{k+1})$ vanishes for every index outside the support, so~\prettyref{eq:ckpt-complementarity} follows by bounding each surviving factor $f_i(y_i^{k+1})$ by $S_\omega/\rho_j$ and summing $|\iota_i^{k+1}|$.
\end{proof}

\subsection{Merit Bounds under the Stored-Multiplier Cap and Split Feasibility}
\label{sub:ckpt-feasibility}
For $j\geq0$ consider the two \emph{level hypotheses}
\begin{eqaligned}
\label{eq:ckpt-hypotheses}
({\rm I}_j)\quad\|\lambda^{[j]}\|\leq Q_{\rm cp},
\qquad\qquad
({\rm II}_j)\quad
\mathcal P_j\left(x^{k_{j-1}+1},y^{k_{j-1}+1}\right)\leq\overline V .
\end{eqaligned}
\prettyref{lem:ckpt-accepted} and \prettyref{lem:ckpt-arrival} below are proved under $({\rm I}_j)$ and $({\rm II}_j)$ alone, and \prettyref{lem:ckpt-transition} closes the induction that verifies both for every $j$.

\begin{lemma}[Merit lower bound under the stored-multiplier cap and coarse consensus]
\label{lem:ckpt-merit-lower}
Fix $j\geq0$ and suppose $\|\lambda^{[j]}\|\leq Q_{\rm cp}$.  Let $x$ and $y$ be arbitrary.  Then, with $D_{\rm cp}$ from~\prettyref{eq:checkpoint-parameters},
\begin{eqaligned}
\label{eq:ckpt-merit-lower}
\mathcal P_j(x,y)
\geq\sum_{i\in\IC}g(f_i(y_i),L_j)-\frac{D_{\rm cp}}{\rho_j}
\geq-\frac{D_{\rm cp}}{\rho_j}
\geq-\Delta_G .
\end{eqaligned}
If moreover $\mathcal P_j(x,y)\leq\overline V$, then $\|y-Ax\|\leq\delta_{\rm cq}$.
\end{lemma}
\begin{proof}
Write $\nu\defeq\|y-Ax\|$ and $\nu_i\defeq\|y_i-A_ix\|$, so that $\nu^2=\sum_{i\in\IAll}\nu_i^2$.  Split~\prettyref{eq:ckpt-merit} as
\[
\mathcal P_j(x,y)
=\underbrace{\left[\sum_{i\in\IO}\bar f_i(y_i,L_j,\mu_j)+R_\psi(\mu_j)\right]}_{\rm (T1)}
+\underbrace{\sum_{i\in\IC}g(f_i(y_i),L_j)}_{\rm (T2)}
+\underbrace{\langle\lambda^{[j]},y-Ax\rangle+\frac{\rho_j}{2}\nu^2}_{\rm (T3)} .
\]

\emph{Term (T1).}  For $i\in\IO$, $\bar f_i(\cdot,L_j,\mu_j)=h_i+(\psi_i)_{\mu_j}$ is globally $\widehat C_i$-Lipschitz: $h_i$ is $D_i$-Lipschitz by \prettyref{ass:obj}, and $\|\nabla(\psi_i)_{\mu_j}\|\leq C_i^\psi$ by~\prettyref{eq:moreau-bound}, with $\widehat C_i=D_i+C_i^\psi$.  Hence
\[
\sum_{i\in\IO}\bar f_i(y_i,L_j,\mu_j)
\geq\sum_{i\in\IO}\bar f_i(A_ix,L_j,\mu_j)
-\sum_{i\in\IO}\widehat C_i\nu_i
\geq\sum_{i\in\IO}\bar f_i(A_ix,L_j,\mu_j)-\sqrt{\widehat G_0}\,\nu,
\]
the last step by the Cauchy--Schwarz inequality with $\widehat G_0=\sum_{i\in\IO}\widehat C_i^2$, which gives $\sum_{i\in\IO}\widehat C_i\nu_i\leq\sqrt{\widehat G_0}\big(\sum_{i\in\IO}\nu_i^2\big)^{1/2}$, followed by $\sum_{i\in\IO}\nu_i^2\leq\sum_{i\in\IAll}\nu_i^2=\nu^2$ because the omitted squares for $i\in\IC$ are nonnegative.  Since~\prettyref{eq:Rpsi-lower} holds at every $x$, it gives $\sum_{i\in\IO}\bar f_i(A_ix,L_j,\mu_j)+R_\psi(\mu_j)\geq0$.  Therefore ${\rm (T1)}\geq-\sqrt{\widehat G_0}\,\nu$.

\emph{Term (T2).}  The wall sum is retained unchanged and is nonnegative.

\emph{Term (T3).}  By Cauchy--Schwarz and the hypothesis, ${\rm (T3)}\geq-Q_{\rm cp}\nu+\frac{\rho_j}{2}\nu^2$.

\emph{Lower bound.}  Adding the three terms,
\begin{eqaligned}
\label{eq:ckpt-merit-quadratic}
\mathcal P_j(x,y)
\geq\sum_{i\in\IC}g(f_i(y_i),L_j)
-\left(Q_{\rm cp}+\sqrt{\widehat G_0}\right)\nu
+\frac{\rho_j}{2}\nu^2 .
\end{eqaligned}
The minimum of the scalar quadratic $-(Q_{\rm cp}+\sqrt{\widehat G_0})\nu+\frac{\rho_j}{2}\nu^2$ over $\nu\in\R$ is $-(Q_{\rm cp}+\sqrt{\widehat G_0})^2/(2\rho_j)=-D_{\rm cp}/\rho_j$, which proves the first inequality of~\prettyref{eq:ckpt-merit-lower}.  The wall $g$ is nonnegative on both of its branches, so the second inequality follows.  The schedule~\prettyref{eq:checkpoint-schedule} gives $\rho_j\geq\rho_{\rm cp}^0$, while~\prettyref{eq:checkpoint-parameters} gives $\rho_{\rm cp}^0\geq D_{\rm cp}/\Delta_G$, proving the third.

\emph{Coarse consensus.}  Discarding the nonnegative wall sum in~\prettyref{eq:ckpt-merit-quadratic} and using $\mathcal P_j(x,y)\leq\overline V$,
\[
\frac{\rho_j}{2}\nu^2-\left(Q_{\rm cp}+\sqrt{\widehat G_0}\right)\nu-\overline V\leq0 .
\]
The quadratic formula and $\sqrt{a+b}\leq\sqrt a+\sqrt b$ give
\[
\nu\leq\frac{(Q_{\rm cp}+\sqrt{\widehat G_0})
+\sqrt{(Q_{\rm cp}+\sqrt{\widehat G_0})^2+2\rho_j\overline V}}{\rho_j}
\leq\frac{2(Q_{\rm cp}+\sqrt{\widehat G_0})}{\rho_j}
+\sqrt{\frac{2\overline V}{\rho_j}} .
\]
By~\prettyref{eq:checkpoint-schedule} and~\prettyref{eq:checkpoint-parameters}, $\rho_j\geq\rho_{\rm cp}^0\geq4(Q_{\rm cp}+\sqrt{\widehat G_0})/\delta_{\rm cq}$ bounds the first term by $\delta_{\rm cq}/2$, and $\rho_j\geq\rho_{\rm cp}^0\geq8\overline V/\delta_{\rm cq}^2$ gives $2\overline V/\rho_j\leq\delta_{\rm cq}^2/4$ and hence bounds the second term by $\delta_{\rm cq}/2$.  Thus $\|y-Ax\|=\nu\leq\delta_{\rm cq}$.
\end{proof}

\begin{lemma}[Initialization and within-level strict split feasibility]
\label{lem:ckpt-level-feasible}
The following hold.
\begin{enumerate}[label=(\roman*),leftmargin=*]
\item\emph{(Initialization.)}  The starting point satisfies
\begin{eqaligned}
\label{eq:ckpt-init}
\mathcal P_0(x^1,y^1)\leq S_0+\Delta_G=\overline V-\Delta_G .
\end{eqaligned}
\item\emph{(Strict split feasibility from the cap.)}  Fix $j\geq0$ with $\|\lambda^{[j]}\|\leq Q_{\rm cp}$, and let $(x,y)$ be arbitrary with $\mathcal P_j(x,y)\leq\overline V$.  Then $f_i(y_i)>0$ for every $i\in\IC$.
\item\emph{(Propagation along a level.)}  Fix $j\geq0$ with $\|\lambda^{[j]}\|\leq Q_{\rm cp}$ and suppose the first iterate of level $j$ satisfies $\mathcal P_j(x^{k_{j-1}+1},y^{k_{j-1}+1})\leq\overline V$.  Then every generated iterate $(x^k,y^k)$ with $k\geq k_{j-1}+1$ that is produced by an iteration while the level state $(\lambda^{[j]},\rho_j,L_j,\mu_j)$ is in force satisfies $Cx^k=d$, $\mathcal P_j(x^k,y^k)\leq\overline V$, and $f_i(y_i^k)>0$ for every $i\in\IC$.  No termination of the level is presupposed; if the level does terminate at $k_j<\infty$, the range of admissible indices is $k_{j-1}+1\leq k\leq k_j+1$ and includes the checkpoint certificate index $k_j+1$.
\end{enumerate}
\end{lemma}
\begin{proof}
(i)  By construction $y^1=Ax^1$ and $\lambda^{[0]}=\lambda^1=0$, so the linear and quadratic consensus terms of~\prettyref{eq:ckpt-merit} both vanish and
\[
\mathcal P_0(x^1,y^1)
=\sum_{i\in\IO}\left[h_i(A_ix^1)+(\psi_i)_{\mu_0}(A_ix^1)\right]
+\sum_{i\in\IC}g(f_i(A_ix^1),L_0)
+R_\psi(\mu_0).
\]
For $i\in\IO$, the left inequality of~\prettyref{eq:moreau-gap} gives $(\psi_i)_{\mu_0}\leq\psi_i$, whence the first sum is at most $\sum_{i\in\IO}\varphi_i(A_ix^1)=S_0$.  For $i\in\IC$, \prettyref{ass:feasible} gives $f_i(A_ix^1)\geq m>0$, while~\prettyref{eq:checkpoint-parameters} contains the clause $\rho_{\rm cp}^0\geq b_{\rm cp}/m$, so by~\prettyref{eq:checkpoint-wall-scaling} together with $\rho_0=\rho_{\rm cp}^0$ from~\prettyref{eq:checkpoint-schedule}
\[
\beta(L_0)=\beta_0=\frac{b_{\rm cp}}{\rho_{\rm cp}^0}\leq m\leq f_i(A_ix^1).
\]
The feasible branch of $g$ then gives $g(f_i(A_ix^1),L_0)=\frac{L_0}{2}\max\{\beta_0-f_i(A_ix^1),0\}^2=0$, so the whole wall sum vanishes.  Finally $R_\psi(\mu_0)=\eta^{\rm cp}G_\psi/(2\rho_{\rm cp}^0)\leq\Delta_G$ by the clause $\rho_{\rm cp}^0\geq\eta^{\rm cp}G_\psi/(2\Delta_G)$ in~\prettyref{eq:checkpoint-parameters}.  Adding the three estimates and using $\overline V=S_0+2\Delta_G$ proves~\prettyref{eq:ckpt-init}.

(ii)  Suppose, for contradiction, that $f_{i_0}(y_{i_0})\leq0$ for some $i_0\in\IC$.  Both branches of $g$ give $g(s,L)\geq G^{\rm cp}$ for $s\leq0$: for $s<0$ the value is $G^{\rm cp}-\sqrt{2LG^{\rm cp}}\,s>G^{\rm cp}$, and for $s=0$ it is $\frac L2(\sqrt{2G^{\rm cp}/L})^2=G^{\rm cp}$.  Since $g\geq0$ on the remaining constraint blocks, $\sum_{i\in\IC}g(f_i(y_i),L_j)\geq G^{\rm cp}$, and~\prettyref{eq:ckpt-merit-lower} gives
\[
\mathcal P_j(x,y)\geq G^{\rm cp}-\frac{D_{\rm cp}}{\rho_j}\geq G^{\rm cp}-\Delta_G .
\]
By the definition $\Delta_G=(G^{\rm cp}-S_0)/4$ in~\prettyref{eq:checkpoint-parameters}, $G^{\rm cp}=S_0+4\Delta_G$, so $G^{\rm cp}-\Delta_G=S_0+3\Delta_G>S_0+2\Delta_G=\overline V$, where the strict inequality uses $\Delta_G>0$.  This contradicts $\mathcal P_j(x,y)\leq\overline V$.  Thus $f_i(y_i)>0$ for every $i\in\IC$.

(iii)  The equality $Cx^k=d$ holds for $k=1$ by \prettyref{ass:feasible} and for every $k\geq2$ by the $x$-update~\prettyref{eq:opt-x}.  Within level $j$ all four parameters $(\lambda^{[j]},\rho_j,L_j,\mu_j)$ are frozen, so \prettyref{lem:ckpt-descent} applies at every iteration $k\geq k_{j-1}+1$ of the level and yields $\mathcal P_j(x^{k+1},y^{k+1})\leq\mathcal P_j(x^k,y^k)$.  Chaining from the first iterate of the level gives $\mathcal P_j(x^k,y^k)\leq\overline V$ at every index reached while the level state is in force, and, when the level terminates, up to and including $k_j+1$, the index of the pre-promotion checkpoint certificate tuple.  Part (ii) then applies at each such iterate.
\end{proof}

\subsection{Accepted Checkpoints, Global Invariants, and Finite Arrival}
\label{sub:ckpt-invariants}

\begin{lemma}[Accepted-checkpoint multiplier and consensus bounds]
\label{lem:ckpt-accepted}
Suppose \prettyref{ass:MFCQ} holds.  Fix $j\geq0$, assume $({\rm I}_j)$ and $({\rm II}_j)$, and let $k$ be an iteration of level $j$ whose test is accepted, so that $k=k_j$.  Then:
\begin{enumerate}[label=(\alph*),leftmargin=*]
\item\emph{(Coarse consensus and strict feasibility.)}
\begin{eqaligned}
Cx^{k+1}=d,\qquad
f_i(y_i^{k+1})>0\ \ (i\in\IC),\qquad
\|r^{k_j+1}\|\leq\delta_{\rm cq}.
\end{eqaligned}
\item\emph{(Support inside the qualification band.)}
\begin{eqaligned}
\iota_i^{k+1}\neq0
\quad\Longrightarrow\quad
0\leq f_i(y_i^{k+1})\leq\frac{S_\omega}{\rho_j}\leq\delta_{\rm cq}.
\end{eqaligned}
\item\emph{(Uniform Band-MFCQ inequality-multiplier bound.)}  The wall forces satisfy $\iota_i^{k+1}\leq0$, the associated nonnegative inequality multipliers are $-\iota_i^{k+1}\geq0$, and
\begin{eqaligned}
\label{eq:ckpt-acc-iota}
\|\iota^{k+1}\|_1=\sum_{i\in\IC}\left(-\iota_i^{k+1}\right)
\leq\frac{M_F+S_\Xi}{\gamma_{\rm cq}}=I_{\rm cp}.
\end{eqaligned}
\item\emph{(Raw accepted virtual multiplier.)}
\begin{eqaligned}
\label{eq:ckpt-acc-lambda}
\|\widehat\lambda^{k+1}\|
\leq\sqrt{\widehat G_0}+M_{\max}I_{\rm cp}
+(c_{\upsilon,{\rm cp}}+1)\tau_{\rm cp}
=Q_{\rm cp}.
\end{eqaligned}
\item\emph{(Sharp consensus.)}
\begin{eqaligned}
\label{eq:ckpt-acc-sharp}
\|r^{k_j+1}\|\leq\frac{2Q_{\rm cp}}{\rho_j}.
\end{eqaligned}
\end{enumerate}
\end{lemma}
\begin{proof}
\emph{(a).}  By $({\rm I}_j)$, $({\rm II}_j)$, and \prettyref{lem:ckpt-level-feasible}(iii) applied at the checkpoint certificate index $k+1=k_j+1$, we have $Cx^{k+1}=d$, $\mathcal P_j(x^{k+1},y^{k+1})\leq\overline V$, and $f_i(y_i^{k+1})>0$ for every $i\in\IC$.  \prettyref{lem:ckpt-merit-lower} then gives $\|r^{k_j+1}\|=\|y^{k+1}-Ax^{k+1}\|\leq\delta_{\rm cq}$.

\emph{(b).}  Part (a) supplies the hypothesis $f_i(y_i^{k+1})\geq0$ of \prettyref{lem:ckpt-support}, whose conclusion~\prettyref{eq:ckpt-support} is (b).

\emph{(c).}  The explicit force~\prettyref{eq:band-force} gives $\iota_i^{k+1}=g'(\ell_i^k(y_i^{k+1}),L_j)\leq0$ and $-\iota_i^{k+1}\geq0$; consequently $\|\iota^{k+1}\|_1=\sum_{i\in\IC}(-\iota_i^{k+1})$.  By part (a) the pair $(x^{k+1},y^{k+1})$ satisfies every hypothesis of \prettyref{ass:MFCQ}: $Cx^{k+1}=d$, $f_i(y_i^{k+1})\geq0$, and $\max_{i\in\IC}\|A_ix^{k+1}-y_i^{k+1}\|\leq\|r^{k_j+1}\|\leq\delta_{\rm cq}$.  Let $d_{\rm cq}\in\ker C$ be the resulting unit vector of~\prettyref{eq:cq}.

By the definition of $\Xi^{k+1}$ preceding~\prettyref{eq:crit},
\[
\Xi^{k+1}
=\sum_{i\in\IO}A_i^Tv_i^{k+1}
+\sum_{i\in\IC}\iota_i^{k+1}A_i^T\nabla f_i(y_i^{k+1})
+C^T\zeta^{k+1}.
\]
Contract with $d_{\rm cq}$.  Since $Cd_{\rm cq}=0$, the last term drops out and
\[
\langle\Xi^{k+1},d_{\rm cq}\rangle
=\sum_{i\in\IO}\langle v_i^{k+1},A_id_{\rm cq}\rangle
-\sum_{i\in\IC}\left(-\iota_i^{k+1}\right)
\nabla f_i(y_i^{k+1})^TA_id_{\rm cq}.
\]
For $i\in\IO$, the reconstruction~\prettyref{eq:reconstruction} defines $v_i^{k+1}=\nabla h_i(p_i^{k+1})+w_i^{k+1}$, so $\|v_i^{k+1}\|\leq D_i+C_i^\psi=\widehat C_i$, since $\|\nabla h_i\|\leq D_i$ by~\prettyref{eq:objective-derived-bounds} and $\|w_i^{k+1}\|\leq C_i^\psi$ by~\prettyref{eq:moreau-bound}.  Hence, by the Cauchy--Schwarz inequality in the stacked space and $\|d_{\rm cq}\|=1$,
\[
\left|\sum_{i\in\IO}\langle v_i^{k+1},A_id_{\rm cq}\rangle\right|
\leq\Big(\sum_{i\in\IO}\widehat C_i^2\Big)^{1/2}\|Ad_{\rm cq}\|
\leq\|A\|_2\sqrt{\widehat G_0}=M_F,
\]
which is the constant $M_F$ of~\prettyref{eq:parameters}.  For $i\in\IC$ with $\iota_i^{k+1}\neq0$, part (b) puts $f_i(y_i^{k+1})$ in $[0,\delta_{\rm cq}]$, so~\prettyref{eq:cq} gives $\nabla f_i(y_i^{k+1})^TA_id_{\rm cq}\geq\gamma_{\rm cq}$; the indices with $\iota_i^{k+1}=0$ contribute nothing.  Therefore
\[
\gamma_{\rm cq}\|\iota^{k+1}\|_1
\leq\sum_{i\in\IC}\left(-\iota_i^{k+1}\right)
\nabla f_i(y_i^{k+1})^TA_id_{\rm cq}
=\sum_{i\in\IO}\langle v_i^{k+1},A_id_{\rm cq}\rangle
-\langle\Xi^{k+1},d_{\rm cq}\rangle
\leq M_F+\|\Xi^{k+1}\| .
\]
Because the step is accepted, \prettyref{eq:ckpt-stationary-accepted} gives $\|\Xi^{k+1}\|\leq S_\Xi/\rho_j\leq S_\Xi$, using $\rho_j\geq\rho_{\rm cp}^0\geq1$.  Dividing by $\gamma_{\rm cq}>0$ proves~\prettyref{eq:ckpt-acc-iota}.

\emph{(d).}  Bound the three terms of the exact identity~\prettyref{eq:ckpt-virtual-dual}.  For the force, $\|q_i^k\|\leq\|\nabla h_i(y_i^k)\|+\|w_i^{k+1}\|\leq D_i+C_i^\psi=\widehat C_i$ for $i\in\IO$, using $\|\nabla h_i\|\leq D_i$ from~\prettyref{eq:objective-derived-bounds} and $\|w_i^{k+1}\|\leq C_i^\psi$ from~\prettyref{eq:moreau-bound}, while $\|q_i^k\|=|\iota_i^{k+1}|\,\|\nabla f_i(y_i^k)\|\leq M_{\max}|\iota_i^{k+1}|$ for $i\in\IC$ by \prettyref{ass:bg}.  Hence, using $\|\cdot\|_2\leq\|\cdot\|_1$ on the constraint block and part (c),
\[
\|q^k\|\leq\Big(\sum_{i\in\IO}\widehat C_i^2\Big)^{1/2}
+M_{\max}\|\iota^{k+1}\|_2
\leq\sqrt{\widehat G_0}+M_{\max}I_{\rm cp}.
\]
For the movement terms, acceptance gives $\rho_j\|d^k\|\leq E^{k+1}\leq\tau_{\rm cp}/\rho_j\leq\tau_{\rm cp}$ and likewise $\rho_j\|A\Delta x^k\|\leq\tau_{\rm cp}$, so~\prettyref{eq:ckpt-upsilon-bound} gives $\|\Upsilon^kd^k\|\leq c_{\upsilon,{\rm cp}}\rho_j\|d^k\|\leq c_{\upsilon,{\rm cp}}\tau_{\rm cp}$ and $\|\rho_jA\Delta x^k\|\leq\tau_{\rm cp}$.  Adding the three bounds and comparing with the definition of $Q_{\rm cp}$ in~\prettyref{eq:checkpoint-parameters} proves~\prettyref{eq:ckpt-acc-lambda}.

\emph{(e).}  By definition $\rho_jr^{k_j+1}=\widehat\lambda^{k+1}-\lambda^{[j]}$, so $({\rm I}_j)$ and part (d) give $\rho_j\|r^{k_j+1}\|\leq2Q_{\rm cp}$.
\end{proof}

\begin{lemma}[Failure cost and finite arrival]
\label{lem:ckpt-arrival}
Fix $j\geq0$ and assume $({\rm I}_j)$ and $({\rm II}_j)$.  Then:
\begin{enumerate}[label=(\alph*),leftmargin=*]
\item\emph{(Failure cost.)}  If $k$ is an iteration of level $j$ whose test is rejected, then
\begin{eqaligned}
\label{eq:ckpt-failure-cost}
\mathcal P_j(x^k,y^k)-\mathcal P_j(x^{k+1},y^{k+1})
>\frac{a_0\tau_{\rm cp}^2}{\rho_j^3}.
\end{eqaligned}
\item\emph{(Level budget.)}  Along level $j$ the total available merit decrease is at most
\begin{eqaligned}
\overline V+\frac{D_{\rm cp}}{\rho_j}\leq\overline V+\Delta_G=\Delta_* .
\end{eqaligned}
\item\emph{(Finite arrival and per-level work.)}  Level $j$ terminates after finitely many iterations, and, with $C_{\rm work}$ from~\prettyref{eq:ckpt-derived-constants},
\begin{eqaligned}
\label{eq:ckpt-work}
T_j<1+C_{\rm work}\rho_j^3.
\end{eqaligned}
\end{enumerate}
\end{lemma}
\begin{proof}
(a)  Rejection at step $k$ means $E^{k+1}>\tau_{\rm cp}/\rho^k=\tau_{\rm cp}/\rho_j$, that is,
\[
\rho_j\sqrt{\|d^k\|^2+\|A\Delta x^k\|^2}>\frac{\tau_{\rm cp}}{\rho_j},
\qquad\text{hence}\qquad
\|d^k\|^2+\|A\Delta x^k\|^2>\frac{\tau_{\rm cp}^2}{\rho_j^4}.
\]
The state $(\lambda^{[j]},\rho_j,L_j,\mu_j)$ is unchanged by a rejection, so \prettyref{lem:ckpt-descent} applies verbatim at step $k$ and yields a decrease of at least $a_0\rho_j(\|d^k\|^2+\|A\Delta x^k\|^2)>a_0\tau_{\rm cp}^2/\rho_j^3$.

(b)  By $({\rm II}_j)$ the level starts at merit value at most $\overline V$.  The hypothesis $({\rm I}_j)$ supplies the multiplier cap required by~\prettyref{lem:ckpt-merit-lower}, whose conclusion~\prettyref{eq:ckpt-merit-lower} gives $\mathcal P_j(x^k,y^k)\geq-D_{\rm cp}/\rho_j\geq-\Delta_G$ throughout the level.  Since \prettyref{lem:ckpt-descent} makes $\mathcal P_j$ nonincreasing, the accumulated decrease over the whole level cannot exceed $\overline V-(-\Delta_G)=\Delta_*$.

(c)  Suppose first that no test of level $j$ is ever accepted.  Then level $j$ never terminates, so it is the infinite tail $\{k_{j-1}+1,k_{j-1}+2,\ldots\}$ on all of which the level state is frozen and every test is rejected.  Hence, for every $T\geq1$, the first $T$ iterations of the level are all rejections, and summing~\prettyref{eq:ckpt-failure-cost} telescopically over them gives an accumulated decrease exceeding $Ta_0\tau_{\rm cp}^2/\rho_j^3$, which is unbounded in $T$ and contradicts (b).  Hence some test is accepted and $k_j<\infty$.  If $T_j=1$, the claimed bound follows immediately from $\Delta_*>0$.  If $T_j\geq2$, the last iteration of the level, $k=k_j$, is the accepted one and the preceding $T_j-1$ iterations are rejections, so
\[
(T_j-1)\frac{a_0\tau_{\rm cp}^2}{\rho_j^3}<\Delta_*,
\qquad\text{that is,}\qquad
T_j<1+\frac{\Delta_*}{a_0\tau_{\rm cp}^2}\rho_j^3
=1+C_{\rm work}\rho_j^3,
\]
which gives~\prettyref{eq:ckpt-work}.
\end{proof}

\begin{lemma}[Transition charge and global invariants]
\label{lem:ckpt-transition}
Suppose \prettyref{ass:MFCQ} holds.  Then:
\begin{enumerate}[label=(\alph*),leftmargin=*]
\item\emph{(One-step charge.)}  Fix $j\geq0$, assume $({\rm I}_j)$ and $({\rm II}_j)$, and let $k_j<\infty$ be the accepting iteration of level $j$.  Then the promotion at $k_j$ raises the merit value by at most
\begin{eqaligned}
\label{eq:ckpt-charge}
\mathcal P_{j+1}\left(x^{k_j+1},y^{k_j+1}\right)
-\mathcal P_j\left(x^{k_j+1},y^{k_j+1}\right)
\leq\frac{\chi+1}{2}\rho_j\|r^{k_j+1}\|^2
\leq\frac{2(\chi+1)Q_{\rm cp}^2}{\rho_j}
\defeq c_j,
\end{eqaligned}
and the charges are summable with
\begin{eqaligned}
\label{eq:ckpt-charge-sum}
\sum_{j=0}^\infty c_j
=\frac{2(\chi+1)Q_{\rm cp}^2}{\rho_{\rm cp}^0(1-\chi^{-1})}
\leq\Delta_G .
\end{eqaligned}
\item\emph{(Global invariants.)}  $({\rm I}_j)$ and $({\rm II}_j)$ hold, $k_j<\infty$, and $1\leq T_j<\infty$ for every $j\geq0$.  Consequently every generated trial iterate satisfies
\begin{eqaligned}
\label{eq:ckpt-global}
Cx^k=d,
\qquad
f_i(y_i^k)>0\ \ (i\in\IC),
\qquad
\|\lambda^k\|\leq Q_{\rm cp},
\end{eqaligned}
every checkpoint level is reached after finitely many iterations, and the conclusions of \prettyref{lem:ckpt-merit-lower}, \prettyref{lem:ckpt-level-feasible}, \prettyref{lem:ckpt-accepted}, and \prettyref{lem:ckpt-arrival} hold at every level without separately assuming $({\rm I}_j)$ or $({\rm II}_j)$.
\end{enumerate}
\end{lemma}
\begin{proof}
(a)  Write $(x,y)\defeq(x^{k_j+1},y^{k_j+1})$ and $r\defeq r^{k_j+1}=y-Ax$.  Recording precedes promotion, so this primal pair is simultaneously the primal component of the level-$j$ checkpoint output and the first iterate of level $j+1$; only the state changes, from $(\lambda^{[j]},\rho_j,L_j,\mu_j)$ to $(\lambda^{[j+1]},\rho_{j+1},L_{j+1},\mu_{j+1})=(\widehat\lambda^{k_j+1},\chi\rho_j,\chi^2L_j,\mu_j/\chi)$.  Grouping~\prettyref{eq:ckpt-merit} into its objective, wall, and consensus parts, the difference splits exactly as
\[
\mathcal P_{j+1}(x,y)-\mathcal P_j(x,y)
=\underbrace{\Big[\textstyle\sum_{i\in\IO}\bar f_i(y_i,L_{j+1},\mu_{j+1})+R_\psi(\mu_{j+1})
-\sum_{i\in\IO}\bar f_i(y_i,L_j,\mu_j)-R_\psi(\mu_j)\Big]}_{\rm (T1)}
\]
\[
+\underbrace{\sum_{i\in\IC}\left[g(f_i(y_i),L_{j+1})-g(f_i(y_i),L_j)\right]}_{\rm (T2)}
+\underbrace{\langle\lambda^{[j+1]}-\lambda^{[j]},r\rangle
+\frac{\rho_{j+1}-\rho_j}{2}\|r\|^2}_{\rm (T3)} .
\]
Term (T1) is at most zero: $\mu_{j+1}=\mu_j/\chi\leq\mu_j$, so~\prettyref{eq:Rpsi-monotone} of \prettyref{lem:Rpsi-discharge} applies with $(\mu,\mu',L,L')=(\mu_j,\mu_{j+1},L_j,L_{j+1})$.

Term (T2) is at most zero as well.  By $({\rm I}_j)$, $({\rm II}_j)$, and \prettyref{lem:ckpt-level-feasible}(iii), $f_i(y_i)>0$ for every $i\in\IC$, so the final claim of \prettyref{lem:checkpoint-schedule}, namely $g(s,L_{j+1})\leq g(s,L_j)$ for $s\geq0$, applies blockwise.

Term (T3) is the only potentially positive contribution.  Since $\lambda^{[j+1]}-\lambda^{[j]}=\widehat\lambda^{k_j+1}-\lambda^{[j]}=\rho_jr$ and $\rho_{j+1}-\rho_j=(\chi-1)\rho_j$,
\[
{\rm (T3)}=\rho_j\|r\|^2+\frac{(\chi-1)\rho_j}{2}\|r\|^2
=\frac{\chi+1}{2}\rho_j\|r\|^2 .
\]
\prettyref{lem:ckpt-accepted}(e) gives $\|r\|\leq2Q_{\rm cp}/\rho_j$, so ${\rm (T3)}\leq\frac{\chi+1}{2}\rho_j\cdot4Q_{\rm cp}^2/\rho_j^2=c_j$, proving~\prettyref{eq:ckpt-charge}.  Summing $c_j=2(\chi+1)Q_{\rm cp}^2/\rho_j$ with the first identity of~\prettyref{eq:checkpoint-geometric-sums} gives the equality in~\prettyref{eq:ckpt-charge-sum}, and the clause $\rho_{\rm cp}^0\geq2(\chi+1)Q_{\rm cp}^2/[(1-\chi^{-1})\Delta_G]$ of~\prettyref{eq:checkpoint-parameters} gives the inequality.

(b)  We prove by induction on $j$ the joint statement
\[
\mathcal S_j:\qquad
({\rm I}_j)\ \text{holds},
\qquad
\mathcal P_j\left(x^{k_{j-1}+1},y^{k_{j-1}+1}\right)
\leq\overline V-\Delta_G+\sum_{\ell=0}^{j-1}c_\ell,
\qquad
k_{j-1}<\infty .
\]
Because $\sum_{\ell\geq0}c_\ell\leq\Delta_G$ by~\prettyref{eq:ckpt-charge-sum}, $\mathcal S_j$ implies $({\rm II}_j)$.

\emph{Base $j=0$.}  Here $k_{-1}=0<\infty$ and $\lambda^{[0]}=\lambda^1=0$, so $({\rm I}_0)$ holds because $Q_{\rm cp}>0$, its definition in~\prettyref{eq:checkpoint-parameters} being a sum of nonnegative terms that contains the strictly positive term $(c_{\upsilon,{\rm cp}}+1)\tau_{\rm cp}$.  The middle clause with the empty sum is exactly~\prettyref{eq:ckpt-init}.

\emph{Step.}  Assume $\mathcal S_j$; then $({\rm I}_j)$ and $({\rm II}_j)$ hold.  \prettyref{lem:ckpt-arrival}(c) applies and gives $k_j<\infty$, so level $j$ has an accepting iteration.  \prettyref{lem:ckpt-accepted}(d) then gives $\|\lambda^{[j+1]}\|=\|\widehat\lambda^{k_j+1}\|\leq Q_{\rm cp}$, which is $({\rm I}_{j+1})$.  \prettyref{lem:ckpt-descent} makes $\mathcal P_j$ nonincreasing along level $j$, so
\[
\mathcal P_j\left(x^{k_j+1},y^{k_j+1}\right)
\leq\mathcal P_j\left(x^{k_{j-1}+1},y^{k_{j-1}+1}\right)
\leq\overline V-\Delta_G+\sum_{\ell=0}^{j-1}c_\ell,
\]
and part (a) adds at most $c_j$ to obtain the middle clause of $\mathcal S_{j+1}$.  Together with $k_j<\infty$ this proves $\mathcal S_{j+1}$.

Finally, $({\rm I}_j)$ and $({\rm II}_j)$ for all $j$ give~\prettyref{eq:ckpt-global}: the equality $Cx^k=d$ holds for $k=1$ by \prettyref{ass:feasible} and for every $k\geq2$ by~\prettyref{eq:opt-x}, the strict inequalities hold at $k=1$ by~\prettyref{ass:feasible} and $y^1=Ax^1$, and for $k\geq2$ by~\prettyref{lem:ckpt-level-feasible}(iii), and the stored multiplier at any $k$ of level $j$ is $\lambda^k=\lambda^{[j]}$, bounded by $({\rm I}_j)$.  Every level terminates by \prettyref{lem:ckpt-arrival}(c), and its accepting iteration $k_j$ belongs to it, so $1\leq T_j<\infty$ for every $j$; hence $k_j=\sum_{\ell=0}^jT_\ell$ is finite and strictly increasing in $j$, and every checkpoint index $j$ is reached after finitely many iterations.
\end{proof}

\subsection{Checkpoint Certificate, Complexity, and Scaling}
\label{sub:ckpt-certificate}
By \prettyref{lem:ckpt-transition}(b) the level hypotheses~\prettyref{eq:ckpt-hypotheses} hold for every $j\geq0$, so all results of this subsection are unconditional.  Throughout, \prettyref{ass:MFCQ} is in force and the checkpoint certificate is the pre-promotion tuple~\prettyref{eq:ckpt-certificate-tuple} at index $k_j+1$.

\begin{lemma}[Pointwise KKT certificate at a checkpoint]
\label{lem:ckpt-certificate}
Suppose \prettyref{ass:MFCQ} holds.  For every $j\geq0$, the checkpoint certificate tuple~\prettyref{eq:ckpt-certificate-tuple} satisfies every exact side condition of the $\epsilon$-KKT system in~\prettyref{def:criteria}, namely
\begin{eqaligned}
v_i^{k_j+1}\in\partial\varphi_i(\breve y_i^{k_j+1})\ (i\in\IO),\qquad
\breve y_i^{k_j+1}=y_i^{k_j+1}\ (i\in\IC),\\
Cx^{k_j+1}=d,\;
f_i(y_i^{k_j+1})>0\ (i\in\IC),\;
\iota_i^{k_j+1}\leq0\ (i\in\IC),
\end{eqaligned}
together with the aggregate residual bound
\[
\CritKKT^{k_j+1}\leq\frac{C_{\rm cert}}{\rho_j}.
\]
\end{lemma}
\begin{proof}
\emph{Side conditions.}  The subgradient membership is the first line of~\prettyref{eq:recon-error}, while $\breve y_i^{k_j+1}=y_i^{k_j+1}$ for $i\in\IC$ is the definition in~\prettyref{eq:reconstruction}; both hold for the one-call outputs~\prettyref{eq:prox-oracle} at any parameter values and are therefore schedule independent.  The equality $Cx^{k_j+1}=d$ and the strict inequalities $f_i(y_i^{k_j+1})>0$ are~\prettyref{eq:ckpt-global}, and $\iota_i^{k_j+1}\leq0$ was shown in \prettyref{lem:ckpt-accepted}(c).  The checkpoint certificate tuple's inequality multipliers are the raw wall forces, and the corresponding nonnegative multipliers of the original problem are $-\iota_i^{k_j+1}$.

\emph{Consensus.}  Since $\breve y_i^{k_j+1}=y_i^{k_j+1}$ for $i\in\IC$, the first line of~\prettyref{eq:recon-error} gives $\|y^{k_j+1}-\breve y^{k_j+1}\|\leq\mu_jC_\psi=\eta^{\rm cp}C_\psi/\rho_j$.  Combining with the sharp consensus estimate~\prettyref{eq:ckpt-acc-sharp},
\[
\left\|Ax^{k_j+1}-\breve y^{k_j+1}\right\|
\leq\left\|Ax^{k_j+1}-y^{k_j+1}\right\|
+\left\|y^{k_j+1}-\breve y^{k_j+1}\right\|
\leq\frac{2Q_{\rm cp}}{\rho_j}+\frac{\eta^{\rm cp}C_\psi}{\rho_j}.
\]

\emph{Stationarity.}  By~\prettyref{eq:ckpt-stationary-accepted} of \prettyref{lem:ckpt-virtual-dual}, applicable because the step-$k_j$ test is accepted, $\|\Xi^{k_j+1}\|\leq S_\Xi/\rho_j$.

\emph{Complementarity.}  \prettyref{lem:ckpt-accepted}(c) gives $\|\iota^{k_j+1}\|_1\leq I_{\rm cp}$.  Substituting this into the complementarity estimate~\prettyref{eq:ckpt-complementarity} of \prettyref{lem:ckpt-support}, whose hypotheses hold because the step-$k_j$ test is accepted and~\prettyref{eq:ckpt-global} gives $f_i(y_i^{k_j+1})>0$ for every $i\in\IC$, yields $\sum_{i\in\IC}|\iota_i^{k_j+1}f_i(y_i^{k_j+1})|\leq I_{\rm cp}S_\omega/\rho_j$.

\emph{Aggregate.}  By the definition of $\CritKKT$ in~\prettyref{eq:crit}, the three component bounds add to the aggregate bound:
\[
\CritKKT^{k_j+1}
\leq\frac{2Q_{\rm cp}+\eta^{\rm cp}C_\psi}{\rho_j}
+\frac{S_\Xi}{\rho_j}
+\frac{I_{\rm cp}S_\omega}{\rho_j}
=\frac{C_{\rm cert}}{\rho_j},
\]
with $C_{\rm cert}$ exactly as in~\prettyref{eq:ckpt-derived-constants}.
\end{proof}

\begin{corollary}[Completed-checkpoint and anytime complexity]
\label{cor:ckpt-complexity}
Define
\begin{eqaligned}
\label{eq:ckpt-Cany}
C_{\rm any}\defeq\frac{C_{\rm work}+(\rho_{\rm cp}^0)^{-3}}{1-\chi^{-3}} .
\end{eqaligned}
Then, for every $j\geq0$,
\begin{eqaligned}
\label{eq:ckpt-cumulative}
k_j\leq j+1+\frac{C_{\rm work}}{1-\chi^{-3}}\rho_j^3
\leq C_{\rm any}\rho_j^3 .
\end{eqaligned}
Let the integer $K\geq K_0$ be a stopping budget, where the startup budget $K_0$ of~\prettyref{thm:checkpoint-rate} is defined in~\prettyref{eq:ckpt-level}.  Let $J(K)\defeq\max\{j\geq0:k_j\leq K\}$ be the latest completed checkpoint, so that the latest-checkpoint output residual is ${\CritKKT}_{\rm out}(K)=\CritKKT^{k_{J(K)}+1}$.  Then $J(K)$ is well defined and finite, and
\begin{eqaligned}
\label{eq:ckpt-anytime}
{\CritKKT}_{\rm out}(K)
\leq\frac{C_{\rm cert}}{\rho_{J(K)}}
\leq\chi C_{\rm cert}C_{\rm any}^{1/3}K^{-1/3}.
\end{eqaligned}
\end{corollary}
\begin{proof}
\emph{Cumulative bound.}  \prettyref{lem:ckpt-arrival}(c) holds at every level; relaxing its strict bound~\prettyref{eq:ckpt-work} to $T_\ell\leq1+C_{\rm work}\rho_\ell^3$ before summing gives $k_j=\sum_{\ell=0}^jT_\ell\leq(j+1)+C_{\rm work}\sum_{\ell=0}^j\rho_\ell^3$, and the finite geometric sum in~\prettyref{eq:checkpoint-geometric-sums} bounds the last sum by $\rho_j^3/(1-\chi^{-3})$.  For the second inequality of~\prettyref{eq:ckpt-cumulative} it suffices to show
\[
j+1\leq\frac{(\rho_{\rm cp}^0)^{-3}}{1-\chi^{-3}}\rho_j^3
=\frac{\chi^{3j}}{1-\chi^{-3}} .
\]
Put $t\defeq\chi^3>1$.  The claim is $(j+1)(1-t^{-1})\leq t^j$, equivalently $(j+1)(t-1)\leq t^{j+1}$.  This holds because $t^{j+1}-1=(t-1)\sum_{i=0}^{j}t^i\geq(t-1)(j+1)$, each of the $j+1$ terms being at least $1$.

\emph{Anytime bound.}  The set $\{j\geq0:k_j\leq K\}$ contains $0$ by $K\geq K_0$ and is finite because every level executes at least one iteration, $T_\ell\geq1$, so that $k_j=\sum_{\ell=0}^jT_\ell\geq j+1\to\infty$; hence $J\defeq J(K)$ is well defined.  By maximality $K<k_{J+1}$, and~\prettyref{eq:ckpt-cumulative} at level $J+1$ with $\rho_{J+1}=\chi\rho_J$ gives $K<C_{\rm any}\chi^3\rho_J^3$, that is, $1/\rho_J<\chi(C_{\rm any}/K)^{1/3}$.  Combining with \prettyref{lem:ckpt-certificate} at level $J$ proves~\prettyref{eq:ckpt-anytime}.
\end{proof}

\paragraph{Same-$x$ KKT residual in the smooth case.}
The same-$x$ transfer described for the ordinary mode applies verbatim at an accepted checkpoint, with $\|\iota^{k_j+1}\|_1\leq I_{\rm cp}$ from \prettyref{lem:ckpt-accepted}(c) and $\|Ax^{k_j+1}-\breve y^{k_j+1}\|=\mathcal O(\rho_j^{-1})$ from \prettyref{lem:ckpt-certificate}.  Hence, when $\psi_i\equiv0$, the conventional same-$x$ KKT residual is also $\mathcal O(\rho_j^{-1})$ and obeys the latest-checkpoint rate $\mathcal O(K^{-1/3})$.

\begin{lemma}[Block scaling of the checkpoint constants]
\label{lem:ckpt-block-scaling}
Suppose \prettyref{ass:block-scaling} holds, $\chi$, $\ell_{\rm cp}$, and $\tau_{\rm cp}$ are fixed, and $\gamma_{\rm cq}$ and $\delta_{\rm cq}$ are bounded below by positive constants uniformly over the family.  With $D_{\rm cq}$ as in~\prettyref{cor:KKT-MFCQ},
\begin{align}
S_0,\ \widehat G_0,\ G^{\rm cp},\ \Delta_G&=\Theta(N),\qquad
\overline V,\ \Delta_*,\ C_{\rm work},\ G_\psi=\mathcal O(N),
\label{eq:ckpt-scale-a}\\
C_\psi,\ J_H,\ S_\omega,\ b_{\rm cp},\ r_{\rm cp},\ c_{\upsilon,{\rm cp}},\ c_{\rm mov}
&=\mathcal O(\sqrt N),
\label{eq:ckpt-scale-b}\\
M_F,\ S_\Xi=\mathcal O(\|A\|_2\sqrt N),\qquad
I_{\rm cp}&=\mathcal O(D_{\rm cq}),
\label{eq:ckpt-scale-c}\\
Q_{\rm cp}=\mathcal O(\sqrt N+D_{\rm cq}),\qquad
D_{\rm cp}&=\mathcal O(N+D_{\rm cq}^2),\nonumber\\
C_{\rm cert}&=\mathcal O\left(\sqrt N(D_{\rm cq}+\|A\|_2+1)\right).
\label{eq:ckpt-scale-d}
\end{align}
Moreover
\begin{eqaligned}
\label{eq:ckpt-startup}
\rho_{\rm cp}^0=\mathcal O\left(N+\frac{\|A\|_2^2}{\gamma_{\rm cq}^2}\right),\qquad
C_{\rm any}=\mathcal O(N),\qquad
K_0=\mathcal O\left(N\left(N+\frac{\|A\|_2^2}{\gamma_{\rm cq}^2}\right)^3\right),
\end{eqaligned}
and, for $K\geq K_0$,
\begin{eqaligned}
\label{eq:ckpt-block-rate}
{\CritKKT}_{\rm out}(K)
=\mathcal O\left((D_{\rm cq}+\|A\|_2+1)N^{5/6}K^{-1/3}\right)
=\mathcal O\left(\|A\|_2N^{4/3}K^{-1/3}\right).
\end{eqaligned}
\end{lemma}
\begin{proof}
All statements follow by substituting \prettyref{ass:block-scaling} into the exact formulas~\prettyref{eq:checkpoint-parameters} and~\prettyref{eq:ckpt-derived-constants}.

\emph{\prettyref{eq:ckpt-scale-a}.}  $S_0=\sum_{i\in\IO}\varphi_i(A_ix^1)$ is a sum of $N$ terms, each bounded above and below away from zero, so $S_0=\Theta(N)$; likewise $\widehat G_0=\sum_{i\in\IO}\widehat C_i^2=\Theta(N)$ and $G_\psi=\sum_{i\in\IO}(C_i^\psi)^2=\mathcal O(N)$.  Since $m=\Theta(1)$, $G^{\rm cp}=\frac12\max\{m^2,m,3S_0\}=\Theta(N)$.  The objective branch gives $G^{\rm cp}\geq3S_0/2$, hence $S_0\leq2G^{\rm cp}/3$ and
\[
\frac{G^{\rm cp}}{12}\leq\Delta_G=\frac{G^{\rm cp}-S_0}{4}\leq\frac{G^{\rm cp}}4,
\]
so $\Delta_G=\Theta(N)$.  Then $\overline V=S_0+2\Delta_G=\mathcal O(N)$, $\Delta_*=\overline V+\Delta_G=\mathcal O(N)$, and, with $\tau_{\rm cp}$ fixed and $a_0=51/100$, $C_{\rm work}=\Delta_*/(a_0\tau_{\rm cp}^2)=\mathcal O(N)$.

\emph{\prettyref{eq:ckpt-scale-b}.}  $C_\psi=\sqrt{G_\psi}=\mathcal O(\sqrt N)$ and $J_H=(\sum_{i\in\IO}(H_iC_i^\psi)^2)^{1/2}\leq H_{\max}C_\psi=\mathcal O(\sqrt N)$.  Since $\ell_{\rm cp}$ is fixed and $G^{\rm cp}=\Theta(N)$, $b_{\rm cp}=\sqrt{2G^{\rm cp}/\ell_{\rm cp}}=\mathcal O(\sqrt N)$ and $r_{\rm cp}=\max\{H_{\max},\sqrt{2G^{\rm cp}\ell_{\rm cp}}L_{\max}\}=\mathcal O(\sqrt N)$.  Hence $c_{\upsilon,{\rm cp}}=\theta_1^{\rm cp}(1+r_{\rm cp})-1$ and $c_{\rm mov}$ are $\mathcal O(\sqrt N)$.  As $\rho_{\rm cp}^0\geq1$, $S_\omega\leq b_{\rm cp}+L_{\max}\tau_{\rm cp}^2/2=\mathcal O(\sqrt N)$.

\emph{\prettyref{eq:ckpt-scale-c}.}  $M_F=\|A\|_2\sqrt{\widehat G_0}=\mathcal O(\|A\|_2\sqrt N)$, and $S_\Xi=\|A\|_2(c_{\rm mov}\tau_{\rm cp}+\eta^{\rm cp}J_H)=\mathcal O(\|A\|_2\sqrt N)$.  The stacked norm $\|A\|_2$ is not bounded by a family-uniform constant under \prettyref{ass:block-scaling}, so its growth is carried explicitly.  Hence
\[
I_{\rm cp}=\frac{M_F+S_\Xi}{\gamma_{\rm cq}}
=\mathcal O\left(\frac{\|A\|_2\sqrt N}{\gamma_{\rm cq}}\right)
=\mathcal O(D_{\rm cq}),
\]
the last step by $D_{\rm cq}=\Theta(1+\|A\|_2\sqrt N/\gamma_{\rm cq})$.

\emph{\prettyref{eq:ckpt-scale-d}.}  $Q_{\rm cp}=\sqrt{\widehat G_0}+M_{\max}I_{\rm cp}+(c_{\upsilon,{\rm cp}}+1)\tau_{\rm cp}=\mathcal O(\sqrt N+D_{\rm cq})$ since $M_{\max}=\mathcal O(1)$; consequently $D_{\rm cp}=(Q_{\rm cp}+\sqrt{\widehat G_0})^2/2=\mathcal O(N+D_{\rm cq}^2)$.  For the certificate constant,
\[
\begin{aligned}
C_{\rm cert}
&=2Q_{\rm cp}+\eta^{\rm cp}C_\psi+S_\Xi+I_{\rm cp}S_\omega
=\mathcal O\left(\sqrt N+D_{\rm cq}+\|A\|_2\sqrt N+D_{\rm cq}\sqrt N\right)\\
&=\mathcal O\left(\sqrt N(D_{\rm cq}+\|A\|_2+1)\right).
\end{aligned}
\]

\emph{\prettyref{eq:ckpt-startup}.}  Bound the nine branches of $\rho_{\rm cp}^0$ in~\prettyref{eq:checkpoint-parameters} from above.  Write $\Gamma\defeq N+\|A\|_2^2/\gamma_{\rm cq}^2$; by Young's inequality $\|A\|_2\sqrt N/\gamma_{\rm cq}\leq\Gamma/2$, so $D_{\rm cq}=\mathcal O(1+\Gamma)$ and $D_{\rm cq}^2/N=\mathcal O(1+\|A\|_2^2/\gamma_{\rm cq}^2)$.  The floor is $1$; $\eta^{\rm cp}G_\psi/(2\Delta_G)=\mathcal O(1)$ because $\Delta_G=\Theta(N)$; $D_{\rm cp}/\Delta_G=\mathcal O(1+D_{\rm cq}^2/N)=\mathcal O(\Gamma)$; the two $b_{\rm cp}$ branches are $\mathcal O(\sqrt N)$ because $m$ and $\delta_{\rm cq}$ are bounded below; $8\overline V/\delta_{\rm cq}^2=\mathcal O(N)$; the transition branch is $\mathcal O(Q_{\rm cp}^2/\Delta_G)=\mathcal O(1+D_{\rm cq}^2/N)=\mathcal O(\Gamma)$; the coarse-linear branch is $\mathcal O(\sqrt N+D_{\rm cq})=\mathcal O(1+\Gamma)$; and the fourth-root branch is $\mathcal O(1)$.  Hence $\rho_{\rm cp}^0=\mathcal O(\Gamma)$.  Since $\rho_{\rm cp}^0\geq1$, \prettyref{eq:ckpt-Cany} gives $C_{\rm any}\leq(C_{\rm work}+1)/(1-\chi^{-3})=\mathcal O(N)$.  Finally $K_0=T_0\leq1+C_{\rm work}(\rho_{\rm cp}^0)^3$ by~\prettyref{eq:ckpt-work}, which is $\mathcal O(N\Gamma^3)$.

\emph{\prettyref{eq:ckpt-block-rate}.}  By~\prettyref{eq:ckpt-anytime}, the anytime constant is $\chi C_{\rm cert}C_{\rm any}^{1/3}=\mathcal O(\sqrt N(D_{\rm cq}+\|A\|_2+1)N^{1/3})$, which is the first form.  Since $\gamma_{\rm cq}$ is bounded below, $D_{\rm cq}=\mathcal O(1+\|A\|_2\sqrt N)$, and~\prettyref{lem:auto-saturation} gives $\|A\|_2=\Omega(1)$; hence $D_{\rm cq}+\|A\|_2+1=\mathcal O(\|A\|_2\sqrt N)$, which gives the second form.
\end{proof}

\begin{corollary}[Bounded stored and accepted multipliers]
\label{cor:ckpt-multiplier-scaling}
Suppose~\prettyref{ass:MFCQ} holds.  For every $k\geq1$ the stored multiplier satisfies $\|\lambda^k\|\leq Q_{\rm cp}$, and for every $j\geq0$ the accepted virtual multiplier and the accepted wall forces satisfy
\[
\|\widehat\lambda^{k_j+1}\|\leq Q_{\rm cp},\qquad
\|\iota^{k_j+1}\|_1\leq I_{\rm cp}.
\]
In particular $\|\lambda^k\|/\rho^k\leq Q_{\rm cp}/\rho^k\leq Q_{\rm cp}/\rho_{\rm cp}^0$, so no rescaling is needed for boundedness.  No bound is asserted for a rejected virtual multiplier.
\end{corollary}
\begin{proof}
By \prettyref{lem:ckpt-transition}(b), the level hypotheses hold at every level.  The stored bound is then~\prettyref{eq:ckpt-global}, and \prettyref{lem:ckpt-accepted}(c) and (d) give the two accepted bounds~\prettyref{eq:ckpt-acc-iota} and~\prettyref{eq:ckpt-acc-lambda}.  The rescaled bound follows from $\rho^k=\rho_j=\rho_{\rm cp}^0\chi^j\geq\rho_{\rm cp}^0$ on level $j$, by~\prettyref{eq:ckpt-level} and~\prettyref{eq:checkpoint-schedule}.
\end{proof}

\paragraph{Proof of \prettyref{thm:checkpoint-rate}.}
By~\prettyref{lem:ckpt-arrival}(c) and~\prettyref{lem:ckpt-transition}(b), every checkpoint level terminates after finitely many outer iterations and $1\leq T_j<1+C_{\rm work}\rho_j^3$ for every $j\geq0$; in particular the startup budget $K_0$ is finite.  The anytime estimate~\prettyref{eq:ckpt-anytime} of~\prettyref{cor:ckpt-complexity} gives ${\CritKKT}_{\rm out}(K)=\mathcal O(K^{-1/3})$ whenever $K\geq K_0$, and no checkpoint output, hence no certified output, exists for $K<K_0$.  Under the stated block-scaling and family-uniform qualification assumptions, by~\prettyref{lem:ckpt-block-scaling}, \prettyref{eq:ckpt-block-rate} gives ${\CritKKT}_{\rm out}(K)=\mathcal O(\|A\|_2N^{4/3}K^{-1/3})$ for $K\geq K_0$, and~\prettyref{eq:ckpt-startup} gives the stated bound on $K_0$.

\end{document}